%%%%%%%%%%%%%%%%%%%%%%%%%%%%%%%%%%%%%%%%%%%%%%%%%%%%%%%%%%%%%%%%%
%								%
%	Geometry of q-Hypergeometric Functions,			%
%	Quantum Affine Algebras	and Elliptic Quantum Groups	%
%								%
%	V.Tarasov  and  A.Varchenko				%
%								%
%	March 1997, 72 pages					%
%								%
%	Asterisque  246 (1997)					%
%								%
%	amstex.tex (ver. 2.1)  and  amssym.tex  are required	%
%								%
%%%%%%%%%%%%%%%%%%%%%%%%%%%%%%%%%%%%%%%%%%%%%%%%%%%%%%%%%%%%%%%%%

%\mag 1200
\input amstex

\ifx\Asterisque\relax
\mag 1200
\hsize 156truemm
\vsize 231truemm
\hoffset 1.6truemm
\voffset 2.6truemm
\else
\hsize 6.25truein
\vsize 9.63truein
\fi

\expandafter\ifx\csname beta.def\endcsname\relax \else\endinput\fi
\expandafter\edef\csname beta.def\endcsname{%
 \catcode`\noexpand\@=\the\catcode`\@\space}

\let\atbefore @

\catcode`\@=11

\overfullrule\z@

\let\@ft@\expandafter \let\@tb@f@\atbefore

\newif\ifMag
\ifnum\mag>1000 \Magtrue\fi
=\ifMag cmr8\else cmr9\fi

\newdimen\p@@ \p@@\p@
\def\m@ths@r{\ifnum\mathsurround=\z@\z@\else\maths@r\fi}
\def\maths@r{1.6\p@@} \def\mathsurzero{\def\maths@r{\z@}}

\mathsurround\maths@r
\font\Brm=cmr12 \font\Bbf=cmbx12 \font\Bit=cmti12 \font\ssf=cmss10
\font\Bsl=cmsl10 scaled 1200 \font\Bmmi=cmmi10 scaled 1200
\font\BBf=cmbx12 scaled 1200 \font\BMmi=cmmi10 scaled 1440

\def\atletter{\edef\atrestore{\catcode`\noexpand\@=\the\catcode`\@\space}
 \catcode`\@=11}

\newread\@ux \newwrite\@@x \newwrite\@@cd
\let\@np@@\input
\def\@np@t#1{\openin\@ux#1\relax\ifeof\@ux\else\closein\@ux\relax\@np@@ #1\fi}
\def\input#1 {\openin\@ux#1\relax\ifeof\@ux\wrs@x{No file #1}\else
 \closein\@ux\relax\@np@@ #1\fi}
\def\Input#1 {\relax} %% Do not remove the space after #1

\def\wr@@x#1{} \def\wrs@x{\immediate\write\sixt@@n}

\def\readldf{\@np@t{\jobname.ldf}}
\def\writeldf{\def\wr@@x{\immediate\write\@@x}
 \def\cl@selbl{\wr@@x{\string\Snodef{\the\Sno}}\wr@@x{\string\endinput}%
 \immediate\closeout\@@x} \immediate\openout\@@x\jobname.ldf}
\let\cl@selbl\relax

\def\tod@y{\ifcase\month\or
 January\or February\or March\or April\or May\or June\or July\or
 August\or September\or October\or November\or December\fi\space\,
\number\day,\space\,\number\year}

\newcount\c@time
\def\h@@r{hh}\def\m@n@te{mm}
\def\wh@tt@me{\c@time\time\divide\c@time 60\edef\h@@r{\number\c@time}%
 \multiply\c@time -60\advance\c@time\time\edef
 \m@n@te{\ifnum\c@time<10 0\fi\number\c@time}}
\def\t@me{\h@@r\/{\rm:}\m@n@te}  \let\whattime\wh@tt@me
\def\today{\tod@y\wr@@x{\string\todaydef{\tod@y}}}
\def\nowtime{\t@me{\let\/\ic@\wr@@x{\string\nowtimedef{\t@me}}}}
\def\todaydef#1{} \def\nowtimedef#1{}

\def\em#1{{\it #1\/}} \def\emph#1{{\sl #1\/}}

\def\fitem#1{\par\setbox\z@\hbox{#1}\hangindent\wd\z@
 \hglue-2\parindent\kern\wd\z@\indent\llap{#1}\ignore}

\def\itemflat#1{\par\setbox\z@\hbox{\rm #1\enspace}\hang\ifnum\wd\z@>\parindent
 \noindent\unhbox\z@\ignore\else\textindent{\rm#1}\fi}

\newcount\itemlet
\def\newbi{\itemlet 96} \newbi
\def\bitem{\gad\itemlet \par\hangindent1.5\parindent
 \hglue-.5\parindent\textindent{\rm\rlap{\char\the\itemlet}\hp{b})}}
\def\atem{\newbi\bitem}

\newcount\itemrm

\def\iitem{\gad\itemrm \par\hangindent1.5\parindent
 \hglue-.5\parindent\textindent{\rm\hp{v}\llap{\romannumeral\the\itemrm})}}

\def\center{\par\begingroup\leftskip\z@ plus \hsize \rightskip\leftskip
 \parindent\z@\parfillskip\z@skip \def\\{\unskip\break}}
\def\endcenter{\endgraf\endgroup}

\def\Abstract{\begingroup\narrower\nt{\bf Abstract.}\enspace\ignore}
\def\endAbs{\endgraf\endgroup}

\let\b@gr@@\begingroup \let\B@gr@@\begingroup
\def\b@gr@{\b@gr@@\let\b@gr@@\undefined}
\def\B@gr@{\B@gr@@\let\B@gr@@\undefined}

\def\@fn@xt#1#2#3{\let\@ch@r=#1\def\n@xt{\ifx\t@st@\@ch@r
 \def\n@@xt{#2}\else\def\n@@xt{#3}\fi\n@@xt}\futurelet\t@st@\n@xt}

\def\@fwd@@#1#2#3{\setbox\z@\hbox{#1}\ifdim\wd\z@>\z@#2\else#3\fi}
\def\s@twd@#1#2{\setbox\z@\hbox{#2}#1\wd\z@}

\def\r@st@re#1{\let#1\s@v@} \def\s@v@d@f{\let\s@v@}

\def\p@sk@p#1#2{\par\skip@#2\relax\ifdim\lastskip<\skip@\relax\removelastskip
 \ifnum#1=\z@\else\penalty#1\relax\fi\vskip\skip@
 \else\ifnum#1=\z@\else\penalty#1\relax\fi\fi}
\def\sk@@p#1{\par\skip@#1\relax\ifdim\lastskip<\skip@\relax\removelastskip
 \vskip\skip@\fi}

\newbox\p@b@ld
\def\poorbold#1{\setbox\p@b@ld\hbox{#1}\kern-.01em\copy\p@b@ld\kern-\wd\p@b@ld
 \kern.02em\copy\p@b@ld\kern-\wd\p@b@ld\kern-.012em\raise.02em\box\p@b@ld}

\ifx\plainfootnote\undefined \let\plainfootnote\footnote \fi

\let\s@v@\proclaim \let\proclaim\relax
\def\r@R@fs#1{\let#1\s@R@fs} \let\s@R@fs\Refs \let\Refs\relax
\def\r@endd@#1{\let#1\s@endd@} \let\s@endd@\enddocument
\let\bye\relax

\def\myR@fs{\@fn@xt[\m@R@f@\m@R@fs} \def\m@R@fs{\@fn@xt*\m@r@f@@\m@R@f@@}
\def\m@R@f@@{\m@R@f@[References]} \def\m@r@f@@*{\m@R@f@[]}

\def\Twelvepoint{\twelvepoint \let\Bbf\BBf \let\Bmmi\BMmi
\font\Brm=cmr12 scaled 1200 \font\Bit=cmti12 scaled 1200
\font\ssf=cmss10 scaled 1200 \font\Bsl=cmsl10 scaled 1440
\font\BBf=cmbx12 scaled 1440 \font\BMmi=cmmi10 scaled 1728}

\newif\ifamsppt

\newdimen\b@gsize

\def\p@@nt{.\kern.3em} \let\point\p@@nt

\let\proheadfont\bf \let\probodyfont\sl \let\demofont\it

\def\reffont#1{\def\r@ff@nt{#1}} \reffont\rm
\def\keyfont#1{\def\k@yf@nt{#1}} \keyfont\rm
\def\paperfont#1{\def\p@p@rf@nt{#1}} \paperfont\it
\def\bookfont#1{\def\b@@kf@nt{#1}} \bookfont\it
\def\volfont#1{\def\v@lf@nt{#1}} \volfont\bf
\def\issuefont#1{\def\iss@f@nt{#1}} \issuefont{no\p@@nt}

\newdimen\r@f@nd \newbox\r@f@b@x \newbox\adjb@x
\newbox\p@nct@ \newbox\k@yb@x \newcount\rcount
\newbox\b@b@x \newbox\p@p@rb@x \newbox\j@@rb@x \newbox\y@@rb@x
\newbox\v@lb@x \newbox\is@b@x \newbox\p@g@b@x \newif\ifp@g@ \newif\ifp@g@s
\newbox\inb@@kb@x \newbox\b@@kb@x \newbox\p@blb@x \newbox\p@bl@db@x
\newbox\ed@b@x \newif\ifed@ \newif\ifed@s \newif\if@fl@b \newif\if@fn@m
\newbox\p@p@nf@b@x \newbox\inf@b@x \newbox\b@@nf@b@x
\newtoks\@dd@p@n \newtoks\@ddt@ks

\newif\ifp@gen@

\@ft@\ifx\csname amsppt.sty\endcsname\relax

\headline={\hfil}
\footline={\ifp@gen@\ifnum\pageno=\z@\else\hfil\foliorm\folio\fi\else
 \ifnum\pageno=\z@\hfil\foliorm\folio\fi\fi\hfil\global\p@gen@true}
\parindent1pc

\font@\tensmc=cmcsc10
\font@\sevenex=cmex7
\font@\sevenit=cmti7
\font@\eightrm=cmr8
\font@\sixrm=cmr6
\font@\eighti=cmmi8 \skewchar\eighti='177
\font@\sixi=cmmi6 \skewchar\sixi='177
\font@\eightsy=cmsy8 \skewchar\eightsy='60
\font@\sixsy=cmsy6 \skewchar\sixsy='60
\font@\eightex=cmex8
\font@\eightbf=cmbx8
\font@\sixbf=cmbx6
\font@\eightit=cmti8
\font@\eightsl=cmsl8
\font@\eightsmc=cmcsc8
\font@\eighttt=cmtt8
\font@\ninerm=cmr9
\font@\ninei=cmmi9 \skewchar\ninei='177
\font@\ninesy=cmsy9 \skewchar\ninesy='60
\font@\nineex=cmex9
\font@\ninebf=cmbx9
\font@\nineit=cmti9
\font@\ninesl=cmsl9
\font@\ninesmc=cmcsc9
\font@\ninemsa=msam9
\font@\ninemsb=msbm9
\font@\nineeufm=eufm9
\font@\eightmsa=msam8
\font@\eightmsb=msbm8
\font@\eighteufm=eufm8
\font@\sixmsa=msam6
\font@\sixmsb=msbm6
\font@\sixeufm=eufm6

\loadmsam\loadmsbm\loadeufm
\input amssym.tex

\def\footnoterule{\kern-3\p@\hrule width5pc\kern 2.6\p@}
\def\m@k@foot#1{\insert\footins
 {\interlinepenalty\interfootnotelinepenalty
 \eightpoint\splittopskip\ht\strutbox\splitmaxdepth\dp\strutbox
 \floatingpenalty\@MM\leftskip\z@\rightskip\z@
 \spaceskip\z@\xspaceskip\z@
 \leavevmode\footstrut\ignore#1\unskip\lower\dp\strutbox
 \vbox to\dp\strutbox{}}}
\def\ftext#1{\m@k@foot{\vsk-.8>\nt #1}}
\def\pr@cl@@m#1{\p@sk@p{-100}\medskipamount
 \def\endproclaim{\endgroup\p@sk@p{55}\medskipamount}\begingroup
 \nt\ignore\proheadfont#1\unskip.\enspace\probodyfont\ignore}
\outer\def\proclaim{\pr@cl@@m} \s@v@d@f\proclaim \let\proclaim\relax
\def\demo#1{\sk@@p\medskipamount\nt{\ignore\demofont#1\unskip.}\enspace
 \ignore}
\def\enddemo{\sk@@p\medskipamount}

\def\cite#1{{\rm[#1]}} \let\nofrills\relax
 \def\Refs#1#2{\relax}

\def\big@#1#2{{\hbox{$\left#2\vcenter to#1\b@gsize{}%
 \right.\nulldelimiterspace\z@\m@th$}}}
\def\big{\big@\@ne}
\def\Big{\big@{1.5}}
\def\bigg{\big@\tw@}
\def\Bigg{\big@{2.5}}
\normallineskiplimit\p@

\def\tenpoint{\p@@\p@ \normallineskiplimit\p@@
 \mathsurround\m@ths@r \normalbaselineskip12\p@@
 \abovedisplayskip12\p@@ plus3\p@@ minus9\p@@
 \belowdisplayskip\abovedisplayskip
 \abovedisplayshortskip\z@ plus3\p@@
 \belowdisplayshortskip7\p@@ plus3\p@@ minus4\p@@
 \textonlyfont@\rm\tenrm \textonlyfont@\it\tenit
 \textonlyfont@\sl\tensl \textonlyfont@\bf\tenbf
 \textonlyfont@\smc\tensmc \textonlyfont@\tt\tentt
 \ifsyntax@ \def\big##1{{\hbox{$\left##1\right.$}}}%
  \let\Big\big \let\bigg\big \let\Bigg\big
 \else
  \textfont\z@\tenrm \scriptfont\z@\sevenrm \scriptscriptfont\z@\fiverm
  \textfont\@ne\teni \scriptfont\@ne\seveni \scriptscriptfont\@ne\fivei
  \textfont\tw@\tensy \scriptfont\tw@\sevensy \scriptscriptfont\tw@\fivesy
  \textfont\thr@@\tenex \scriptfont\thr@@\sevenex
	\scriptscriptfont\thr@@\sevenex
  \textfont\itfam\tenit \scriptfont\itfam\sevenit
	\scriptscriptfont\itfam\sevenit
  \textfont\bffam\tenbf \scriptfont\bffam\sevenbf
	\scriptscriptfont\bffam\fivebf
  \textfont\msafam\tenmsa \scriptfont\msafam\sevenmsa
	\scriptscriptfont\msafam\fivemsa
  \textfont\msbfam\tenmsb \scriptfont\msbfam\sevenmsb
	\scriptscriptfont\msbfam\fivemsb
  \textfont\eufmfam\teneufm \scriptfont\eufmfam\seveneufm
	\scriptscriptfont\eufmfam\fiveeufm
  \setbox\strutbox\hbox{\vrule height8.5\p@@ depth3.5\p@@ width\z@}%
  \setbox\strutbox@\hbox{\lower.5\normallineskiplimit\vbox{%
	\kern-\normallineskiplimit\copy\strutbox}}%
   \setbox\z@\vbox{\hbox{$($}\kern\z@}\b@gsize1.2\ht\z@
  \fi
  \normalbaselines\rm\dotsspace@1.5mu\ex@.2326ex\jot3\ex@}

\def\eightpoint{\p@@.8\p@ \normallineskiplimit\p@@
 \mathsurround\m@ths@r \normalbaselineskip10\p@
 \abovedisplayskip10\p@ plus2.4\p@ minus7.2\p@
 \belowdisplayskip\abovedisplayskip
 \abovedisplayshortskip\z@ plus3\p@@
 \belowdisplayshortskip7\p@@ plus3\p@@ minus4\p@@
 \textonlyfont@\rm\eightrm \textonlyfont@\it\eightit
 \textonlyfont@\sl\eightsl \textonlyfont@\bf\eightbf
 \textonlyfont@\smc\eightsmc \textonlyfont@\tt\eighttt
 \ifsyntax@\def\big##1{{\hbox{$\left##1\right.$}}}%
  \let\Big\big \let\bigg\big \let\Bigg\big
 \else
  \textfont\z@\eightrm \scriptfont\z@\sixrm \scriptscriptfont\z@\fiverm
  \textfont\@ne\eighti \scriptfont\@ne\sixi \scriptscriptfont\@ne\fivei
  \textfont\tw@\eightsy \scriptfont\tw@\sixsy \scriptscriptfont\tw@\fivesy
  \textfont\thr@@\eightex \scriptfont\thr@@\sevenex
	\scriptscriptfont\thr@@\sevenex
  \textfont\itfam\eightit \scriptfont\itfam\sevenit
	\scriptscriptfont\itfam\sevenit
  \textfont\bffam\eightbf \scriptfont\bffam\sixbf
	\scriptscriptfont\bffam\fivebf
  \textfont\msafam\eightmsa \scriptfont\msafam\sixmsa
	\scriptscriptfont\msafam\fivemsa
  \textfont\msbfam\eightmsb \scriptfont\msbfam\sixmsb
	\scriptscriptfont\msbfam\fivemsb
  \textfont\eufmfam\eighteufm \scriptfont\eufmfam\sixeufm
	\scriptscriptfont\eufmfam\fiveeufm
 \setbox\strutbox\hbox{\vrule height7\p@ depth3\p@ width\z@}%
 \setbox\strutbox@\hbox{\raise.5\normallineskiplimit\vbox{%
   \kern-\normallineskiplimit\copy\strutbox}}%
 \setbox\z@\vbox{\hbox{$($}\kern\z@}\b@gsize1.2\ht\z@
 \fi
 \normalbaselines\eightrm\dotsspace@1.5mu\ex@.2326ex\jot3\ex@}

\def\ninepoint{\p@@.9\p@ \normallineskiplimit\p@@
 \mathsurround\m@ths@r \normalbaselineskip11\p@
 \abovedisplayskip11\p@ plus2.7\p@ minus8.1\p@
 \belowdisplayskip\abovedisplayskip
 \abovedisplayshortskip\z@ plus3\p@@
 \belowdisplayshortskip7\p@@ plus3\p@@ minus4\p@@
 \textonlyfont@\rm\ninerm \textonlyfont@\it\nineit
 \textonlyfont@\sl\ninesl \textonlyfont@\bf\ninebf
 \textonlyfont@\smc\ninesmc \textonlyfont@\tt\ninett
 \ifsyntax@ \def\big##1{{\hbox{$\left##1\right.$}}}%
  \let\Big\big \let\bigg\big \let\Bigg\big
 \else
  \textfont\z@\ninerm \scriptfont\z@\sevenrm \scriptscriptfont\z@\fiverm
  \textfont\@ne\ninei \scriptfont\@ne\seveni \scriptscriptfont\@ne\fivei
  \textfont\tw@\ninesy \scriptfont\tw@\sevensy \scriptscriptfont\tw@\fivesy
  \textfont\thr@@\nineex \scriptfont\thr@@\sevenex
	\scriptscriptfont\thr@@\sevenex
  \textfont\itfam\nineit \scriptfont\itfam\sevenit
	\scriptscriptfont\itfam\sevenit
  \textfont\bffam\ninebf \scriptfont\bffam\sevenbf
	\scriptscriptfont\bffam\fivebf
  \textfont\msafam\ninemsa \scriptfont\msafam\sevenmsa
	\scriptscriptfont\msafam\fivemsa
  \textfont\msbfam\ninemsb \scriptfont\msbfam\sevenmsb
	\scriptscriptfont\msbfam\fivemsb
  \textfont\eufmfam\nineeufm \scriptfont\eufmfam\seveneufm
	\scriptscriptfont\eufmfam\fiveeufm
  \setbox\strutbox\hbox{\vrule height8.5\p@@ depth3.5\p@@ width\z@}%
  \setbox\strutbox@\hbox{\lower.5\normallineskiplimit\vbox{%
	\kern-\normallineskiplimit\copy\strutbox}}%
   \setbox\z@\vbox{\hbox{$($}\kern\z@}\b@gsize1.2\ht\z@
  \fi
  \normalbaselines\rm\dotsspace@1.5mu\ex@.2326ex\jot3\ex@}

\font@\twelverm=cmr10 scaled 1200
\font@\twelveit=cmti10 scaled 1200
\font@\twelvesl=cmsl10 scaled 1200
\font@\twelvebf=cmbx10 scaled 1200
\font@\twelvesmc=cmcsc10 scaled 1200
\font@\twelvett=cmtt10 scaled 1200
\font@\twelvei=cmmi10 scaled 1200 \skewchar\twelvei='177
\font@\twelvesy=cmsy10 scaled 1200 \skewchar\twelvesy='60
\font@\twelveex=cmex10 scaled 1200
\font@\twelvemsa=msam10 scaled 1200
\font@\twelvemsb=msbm10 scaled 1200
\font@\twelveeufm=eufm10 scaled 1200

\def\twelvepoint{\p@@1.2\p@ \normallineskiplimit\p@@
 \mathsurround\m@ths@r \normalbaselineskip12\p@@
 \abovedisplayskip12\p@@ plus3\p@@ minus9\p@@
 \belowdisplayskip\abovedisplayskip
 \abovedisplayshortskip\z@ plus3\p@@
 \belowdisplayshortskip7\p@@ plus3\p@@ minus4\p@@
 \textonlyfont@\rm\twelverm \textonlyfont@\it\twelveit
 \textonlyfont@\sl\twelvesl \textonlyfont@\bf\twelvebf
 \textonlyfont@\smc\twelvesmc \textonlyfont@\tt\twelvett
 \ifsyntax@ \def\big##1{{\hbox{$\left##1\right.$}}}%
  \let\Big\big \let\bigg\big \let\Bigg\big
 \else
  \textfont\z@\twelverm \scriptfont\z@\eightrm \scriptscriptfont\z@\sixrm
  \textfont\@ne\twelvei \scriptfont\@ne\eighti \scriptscriptfont\@ne\sixi
  \textfont\tw@\twelvesy \scriptfont\tw@\eightsy \scriptscriptfont\tw@\sixsy
  \textfont\thr@@\twelveex \scriptfont\thr@@\eightex
	\scriptscriptfont\thr@@\sevenex
  \textfont\itfam\twelveit \scriptfont\itfam\eightit
	\scriptscriptfont\itfam\sevenit
  \textfont\bffam\twelvebf \scriptfont\bffam\eightbf
	\scriptscriptfont\bffam\sixbf
  \textfont\msafam\twelvemsa \scriptfont\msafam\eightmsa
	\scriptscriptfont\msafam\sixmsa
  \textfont\msbfam\twelvemsb \scriptfont\msbfam\eightmsb
	\scriptscriptfont\msbfam\sixmsb
  \textfont\eufmfam\twelveeufm \scriptfont\eufmfam\eighteufm
	\scriptscriptfont\eufmfam\sixeufm
  \setbox\strutbox\hbox{\vrule height8.5\p@@ depth3.5\p@@ width\z@}%
  \setbox\strutbox@\hbox{\lower.5\normallineskiplimit\vbox{%
	\kern-\normallineskiplimit\copy\strutbox}}%
  \setbox\z@\vbox{\hbox{$($}\kern\z@}\b@gsize1.2\ht\z@
  \fi
  \normalbaselines\rm\dotsspace@1.5mu\ex@.2326ex\jot3\ex@}

\font@\twelvetrm=cmr10 at 12truept
\font@\twelvetit=cmti10 at 12truept
\font@\twelvetsl=cmsl10 at 12truept
\font@\twelvetbf=cmbx10 at 12truept
\font@\twelvetsmc=cmcsc10 at 12truept
\font@\twelvettt=cmtt10 at 12truept
\font@\twelveti=cmmi10 at 12truept \skewchar\twelveti='177
\font@\twelvetsy=cmsy10 at 12truept \skewchar\twelvetsy='60
\font@\twelvetex=cmex10 at 12truept
\font@\twelvetmsa=msam10 at 12truept
\font@\twelvetmsb=msbm10 at 12truept
\font@\twelveteufm=eufm10 at 12truept

\def\twelvetruepoint{\p@@1.2truept \normallineskiplimit\p@@
 \mathsurround\m@ths@r \normalbaselineskip12\p@@
 \abovedisplayskip12\p@@ plus3\p@@ minus9\p@@
 \belowdisplayskip\abovedisplayskip
 \abovedisplayshortskip\z@ plus3\p@@
 \belowdisplayshortskip7\p@@ plus3\p@@ minus4\p@@
 \textonlyfont@\rm\twelvetrm \textonlyfont@\it\twelvetit
 \textonlyfont@\sl\twelvetsl \textonlyfont@\bf\twelvetbf
 \textonlyfont@\smc\twelvetsmc \textonlyfont@\tt\twelvettt
 \ifsyntax@ \def\big##1{{\hbox{$\left##1\right.$}}}%
  \let\Big\big \let\bigg\big \let\Bigg\big
 \else
  \textfont\z@\twelvetrm \scriptfont\z@\eightrm \scriptscriptfont\z@\sixrm
  \textfont\@ne\twelveti \scriptfont\@ne\eighti \scriptscriptfont\@ne\sixi
  \textfont\tw@\twelvetsy \scriptfont\tw@\eightsy \scriptscriptfont\tw@\sixsy
  \textfont\thr@@\twelvetex \scriptfont\thr@@\eightex
	\scriptscriptfont\thr@@\sevenex
  \textfont\itfam\twelvetit \scriptfont\itfam\eightit
	\scriptscriptfont\itfam\sevenit
  \textfont\bffam\twelvetbf \scriptfont\bffam\eightbf
	\scriptscriptfont\bffam\sixbf
  \textfont\msafam\twelvetmsa \scriptfont\msafam\eightmsa
	\scriptscriptfont\msafam\sixmsa
  \textfont\msbfam\twelvetmsb \scriptfont\msbfam\eightmsb
	\scriptscriptfont\msbfam\sixmsb
  \textfont\eufmfam\twelveteufm \scriptfont\eufmfam\eighteufm
	\scriptscriptfont\eufmfam\sixeufm
  \setbox\strutbox\hbox{\vrule height8.5\p@@ depth3.5\p@@ width\z@}%
  \setbox\strutbox@\hbox{\lower.5\normallineskiplimit\vbox{%
	\kern-\normallineskiplimit\copy\strutbox}}%
  \setbox\z@\vbox{\hbox{$($}\kern\z@}\b@gsize1.2\ht\z@
  \fi
  \normalbaselines\rm\dotsspace@1.5mu\ex@.2326ex\jot3\ex@}

\font@\elevenrm=cmr10 scaled 1095
\font@\elevenit=cmti10 scaled 1095
\font@\elevensl=cmsl10 scaled 1095
\font@\elevenbf=cmbx10 scaled 1095
\font@\elevensmc=cmcsc10 scaled 1095
\font@\eleventt=cmtt10 scaled 1095
\font@\eleveni=cmmi10 scaled 1095 \skewchar\eleveni='177
\font@\elevensy=cmsy10 scaled 1095 \skewchar\elevensy='60
\font@\elevenex=cmex10 scaled 1095
\font@\elevenmsa=msam10 scaled 1095
\font@\elevenmsb=msbm10 scaled 1095
\font@\eleveneufm=eufm10 scaled 1095

\def\elevenpoint{\p@@1.1\p@ \normallineskiplimit\p@@
 \mathsurround\m@ths@r \normalbaselineskip12\p@@
 \abovedisplayskip12\p@@ plus3\p@@ minus9\p@@
 \belowdisplayskip\abovedisplayskip
 \abovedisplayshortskip\z@ plus3\p@@
 \belowdisplayshortskip7\p@@ plus3\p@@ minus4\p@@
 \textonlyfont@\rm\elevenrm \textonlyfont@\it\elevenit
 \textonlyfont@\sl\elevensl \textonlyfont@\bf\elevenbf
 \textonlyfont@\smc\elevensmc \textonlyfont@\tt\eleventt
 \ifsyntax@ \def\big##1{{\hbox{$\left##1\right.$}}}%
  \let\Big\big \let\bigg\big \let\Bigg\big
 \else
  \textfont\z@\elevenrm \scriptfont\z@\eightrm \scriptscriptfont\z@\sixrm
  \textfont\@ne\eleveni \scriptfont\@ne\eighti \scriptscriptfont\@ne\sixi
  \textfont\tw@\elevensy \scriptfont\tw@\eightsy \scriptscriptfont\tw@\sixsy
  \textfont\thr@@\elevenex \scriptfont\thr@@\eightex
	\scriptscriptfont\thr@@\sevenex
  \textfont\itfam\elevenit \scriptfont\itfam\eightit
	\scriptscriptfont\itfam\sevenit
  \textfont\bffam\elevenbf \scriptfont\bffam\eightbf
	\scriptscriptfont\bffam\sixbf
  \textfont\msafam\elevenmsa \scriptfont\msafam\eightmsa
	\scriptscriptfont\msafam\sixmsa
  \textfont\msbfam\elevenmsb \scriptfont\msbfam\eightmsb
	\scriptscriptfont\msbfam\sixmsb
  \textfont\eufmfam\eleveneufm \scriptfont\eufmfam\eighteufm
	\scriptscriptfont\eufmfam\sixeufm
  \setbox\strutbox\hbox{\vrule height8.5\p@@ depth3.5\p@@ width\z@}%
  \setbox\strutbox@\hbox{\lower.5\normallineskiplimit\vbox{%
	\kern-\normallineskiplimit\copy\strutbox}}%
  \setbox\z@\vbox{\hbox{$($}\kern\z@}\b@gsize1.2\ht\z@
  \fi
  \normalbaselines\rm\dotsspace@1.5mu\ex@.2326ex\jot3\ex@}

\def\m@R@f@[#1]{\mathsurzero{%\let\{\relax
 \s@ct{}{#1}}\wr@@c{\string\Refcd{#1}{\the\pageno}}\B@gr@
 \frenchspacing\rcount\z@\refkey{\k@yf@nt[##1]}\refno{\k@yf@nt[##1]}%
 \widest{AZ}\keyright\let\Key\key\let\refin\relax}
\def\widest#1{\s@twd@\r@f@nd{\r@fk@y{\k@yf@nt#1}\enspace}}
\def\widestno#1{\s@twd@\r@f@nd{\r@fn@{\k@yf@nt#1}\enspace}}
\def\widestlabel#1{\s@twd@\r@f@nd{\k@yf@nt#1\enspace}}
\def\refkey{\def\r@fk@y##1} \def\refno{\def\r@fn@##1}
\def\keyright{\def\r@fit@m{\hang\textindent}}
\def\keyflat{\def\r@fit@m##1{\setbox\z@\hbox{##1\enspace}\hang\noindent
 \ifnum\wd\z@<\parindent\indent\hglue-\wd\z@\fi\unhbox\z@}}

\def\R@fb@x{\global\setbox\r@f@b@x} \def\K@yb@x{\global\setbox\k@yb@x}
\def\ref{\par\b@gr@\r@ff@nt\R@fb@x\box\voidb@x\K@yb@x\box\voidb@x
 \@fn@mfalse\@fl@bfalse\b@g@nr@f}
\def\c@nc@t#1{\setbox\z@\lastbox
 \setbox\adjb@x\hbox{\unhbox\adjb@x\unhbox\z@\unskip\unskip\unpenalty#1}}
\def\adjust#1{\relax\ifmmode\penalty-\@M\null\hfil$\clubpenalty\z@
 \widowpenalty\z@\interlinepenalty\z@\offinterlineskip\endgraf
 \setbox\z@\lastbox\unskip\unpenalty\c@nc@t{#1}\nt$\hfil\penalty-\@M
 \else\endgraf\c@nc@t{#1}\nt\fi}
\def\adjustnext#1{\P@nct\hbox{#1}\ignore}
\def\adjustend#1{\def\@djp@{#1}\ignore}
\def\addtoks#1{\global\@ddt@ks{#1}\ignore}
\def\addnext#1{\global\@dd@p@n{#1}\ignore}

\def\cl@s@{\adjust{\@djp@}\endgraf\setbox\z@\lastbox
 \global\setbox\@ne\hbox{\unhbox\adjb@x\ifvoid\z@\else\unhbox\z@\unskip\unskip
 \unpenalty\fi}\egroup\ifnum\c@rr@nt=\k@yb@x\global\fi
 \setbox\c@rr@nt\hbox{\unhbox\@ne\box\p@nct@}\P@nct\null
 \the\@ddt@ks\global\@ddt@ks{}}
\def\@p@n#1{\def\c@rr@nt{#1}\setbox\c@rr@nt\vbox\bgroup\let\@djp@\relax
 \hsize\maxdimen\nt\the\@dd@p@n\global\@dd@p@n{}}
\def\b@g@nr@f{\bgroup\@p@n\z@}
\def\key{\cl@s@\ifvoid\k@yb@x\@p@n\k@yb@x\k@yf@nt\else\@p@n\z@\fi}
\def\label{\cl@s@\ifvoid\k@yb@x\global\@fl@btrue\@p@n\k@yb@x\k@yf@nt\else
 \@p@n\z@\fi}
\def\no{\cl@s@\ifvoid\k@yb@x\gad\rcount\global\@fn@mtrue
 \K@yb@x\hbox{\k@yf@nt\the\rcount}\fi\@p@n\z@}
\def\labelno{\cl@s@\ifvoid\k@yb@x\gad\rcount\@fl@btrue
 \@p@n\k@yb@x\k@yf@nt\the\rcount\else\@p@n\z@\fi}
\def\by{\cl@s@\@p@n\b@b@x} \def\paper{\cl@s@\@p@n\p@p@rb@x\p@p@rf@nt\ignore}
\def\jour{\cl@s@\@p@n\j@@rb@x} \def\yr{\cl@s@\@p@n\y@@rb@x}
\def\vol{\cl@s@\@p@n\v@lb@x\v@lf@nt\ignore}
\def\issue{\cl@s@\@p@n\is@b@x\iss@f@nt\ignore}
\def\page{\cl@s@\ifp@g@s\@p@n\z@\else\p@g@true\@p@n\p@g@b@x\fi}
\def\pages{\cl@s@\ifp@g@\@p@n\z@\else\p@g@strue\@p@n\p@g@b@x\fi}
\def\inbook{\cl@s@\@p@n\inb@@kb@x}
\def\book{\cl@s@\@p@n\b@@kb@x\b@@kf@nt\ignore}
\def\publ{\cl@s@\@p@n\p@blb@x} \def\publaddr{\cl@s@\@p@n\p@bl@db@x}
\def\ed{\cl@s@\ifed@s\@p@n\z@\else\ed@true\@p@n\ed@b@x\fi}
\def\eds{\cl@s@\ifed@\@p@n\z@\else\ed@strue\@p@n\ed@b@x\fi}
\def\info{\cl@s@\@p@n\inf@b@x} \def\paperinfo{\cl@s@\@p@n\p@p@nf@b@x}
\def\bookinfo{\cl@s@\@p@n\b@@nf@b@x} \let\finalinfo\info
\def\P@nct{\global\setbox\p@nct@} \def\nopunct{\P@nct\box\voidb@x}
\def\p@@@t#1#2{\ifvoid\p@nct@\else#1\unhbox\p@nct@#2\fi}
\def\sp@@{\penalty-50 \space\hskip\z@ plus.1em}
\def\c@mm@{\p@@@t,\sp@@} \def\sp@c@{\p@@@t\empty\sp@@}
\def\p@tb@x#1#2{\ifvoid#1\else#2\@nb@x#1\fi}
\def\@nb@x#1{\unhbox#1\P@nct\lastbox}
\def\endr@f@{\cl@s@\nopunct
 \R@fb@x\hbox{\unhbox\r@f@b@x \p@tb@x\b@b@x\empty
 \ifvoid\j@@rb@x\ifvoid\inb@@kb@x\ifvoid\p@p@rb@x\ifvoid\b@@kb@x
  \ifvoid\p@p@nf@b@x\ifvoid\b@@nf@b@x
  \p@tb@x\v@lb@x\c@mm@ \ifvoid\y@@rb@x\else\sp@c@(\@nb@x\y@@rb@x)\fi
  \p@tb@x\is@b@x\c@mm@ \p@tb@x\p@g@b@x\c@mm@ \p@tb@x\inf@b@x\c@mm@
  \else\p@tb@x \b@@nf@b@x\c@mm@ \p@tb@x\v@lb@x\c@mm@ \p@tb@x\is@b@x\sp@c@
  \ifvoid\ed@b@x\else\sp@c@(\@nb@x\ed@b@x,\space\ifed@ ed.\else eds.\fi)\fi
  \p@tb@x\p@blb@x\c@mm@ \p@tb@x\p@bl@db@x\c@mm@ \p@tb@x\y@@rb@x\c@mm@
  \p@tb@x\p@g@b@x{\c@mm@\ifp@g@ p\p@@nt\else pp\p@@nt\fi}%
  \p@tb@x\inf@b@x\c@mm@\fi
  \else \p@tb@x\p@p@nf@b@x\c@mm@ \p@tb@x\v@lb@x\c@mm@
  \ifvoid\y@@rb@x\else\sp@c@(\@nb@x\y@@rb@x)\fi
  \p@tb@x\is@b@x\c@mm@ \p@tb@x\p@g@b@x\c@mm@ \p@tb@x\inf@b@x\c@mm@\fi
  \else \p@tb@x\b@@kb@x\c@mm@
  \p@tb@x\b@@nf@b@x\c@mm@ \p@tb@x\p@blb@x\c@mm@
  \p@tb@x\p@bl@db@x\c@mm@ \p@tb@x\y@@rb@x\c@mm@
  \ifvoid\p@g@b@x\else\c@mm@\@nb@x\p@g@b@x p\fi \p@tb@x\inf@b@x\c@mm@ \fi
  \else \c@mm@\@nb@x\p@p@rb@x\ic@\p@tb@x\p@p@nf@b@x\c@mm@
  \p@tb@x\v@lb@x\sp@c@ \ifvoid\y@@rb@x\else\sp@c@(\@nb@x\y@@rb@x)\fi
  \p@tb@x\is@b@x\c@mm@ \p@tb@x\p@g@b@x\c@mm@\p@tb@x\inf@b@x\c@mm@\fi
  \else \p@tb@x\p@p@rb@x\c@mm@\ic@\p@tb@x\p@p@nf@b@x\c@mm@
  \c@mm@\@nb@x\inb@@kb@x \p@tb@x\b@@nf@b@x\c@mm@ \p@tb@x\v@lb@x\sp@c@
  \p@tb@x\is@b@x\sp@c@
  \ifvoid\ed@b@x\else\sp@c@(\@nb@x\ed@b@x,\space\ifed@ ed.\else eds.\fi)\fi
  \p@tb@x\p@blb@x\c@mm@ \p@tb@x\p@bl@db@x\c@mm@ \p@tb@x\y@@rb@x\c@mm@
  \p@tb@x\p@g@b@x{\c@mm@\ifp@g@ p\p@@nt\else pp\p@@nt\fi}%
  \p@tb@x\inf@b@x\c@mm@\fi
  \else\p@tb@x\p@p@rb@x\c@mm@\ic@\p@tb@x\p@p@nf@b@x\c@mm@\p@tb@x\j@@rb@x\c@mm@
  \p@tb@x\v@lb@x\sp@c@ \ifvoid\y@@rb@x\else\sp@c@(\@nb@x\y@@rb@x)\fi
  \p@tb@x\is@b@x\c@mm@ \p@tb@x\p@g@b@x\c@mm@ \p@tb@x\inf@b@x\c@mm@ \fi}}
\def\m@r@f#1#2{\endr@f@\ifvoid\p@nct@\else\R@fb@x\hbox{\unhbox\r@f@b@x
 #1\unhbox\p@nct@\penalty-200\enskip#2}\fi\egroup\b@g@nr@f}
\def\endref{\endr@f@\ifvoid\p@nct@\else\R@fb@x\hbox{\unhbox\r@f@b@x.}\fi
 \parindent\r@f@nd
 \r@fit@m{\ifvoid\k@yb@x\else\if@fn@m\r@fn@{\unhbox\k@yb@x}\else
 \if@fl@b\unhbox\k@yb@x\else\r@fk@y{\unhbox\k@yb@x}\fi\fi\fi}\unhbox\r@f@b@x
 \endgraf\egroup\endgroup}
\def\moreref{\m@r@f;\empty}
\def\transl{\m@r@f;{\unskip\space
 {\sl English translation\ic@}:\penalty-66 \space}}
\def\endRefs{\endgraf\goodbreak\endgroup}

\hyphenation{acad-e-my acad-e-mies af-ter-thought anom-aly anom-alies
an-ti-deriv-a-tive an-tin-o-my an-tin-o-mies apoth-e-o-ses
apoth-e-o-sis ap-pen-dix ar-che-typ-al as-sign-a-ble as-sist-ant-ship
as-ymp-tot-ic asyn-chro-nous at-trib-uted at-trib-ut-able bank-rupt
bank-rupt-cy bi-dif-fer-en-tial blue-print busier busiest
cat-a-stroph-ic cat-a-stroph-i-cally con-gress cross-hatched data-base
de-fin-i-tive de-riv-a-tive dis-trib-ute dri-ver dri-vers eco-nom-ics
econ-o-mist elit-ist equi-vari-ant ex-quis-ite ex-tra-or-di-nary
flow-chart for-mi-da-ble forth-right friv-o-lous ge-o-des-ic
ge-o-det-ic geo-met-ric griev-ance griev-ous griev-ous-ly
hexa-dec-i-mal ho-lo-no-my ho-mo-thetic ideals idio-syn-crasy
in-fin-ite-ly in-fin-i-tes-i-mal ir-rev-o-ca-ble key-stroke
lam-en-ta-ble light-weight mal-a-prop-ism man-u-script mar-gin-al
meta-bol-ic me-tab-o-lism meta-lan-guage me-trop-o-lis
met-ro-pol-i-tan mi-nut-est mol-e-cule mono-chrome mono-pole
mo-nop-oly mono-spline mo-not-o-nous mul-ti-fac-eted mul-ti-plic-able
non-euclid-ean non-iso-mor-phic non-smooth par-a-digm par-a-bol-ic
pa-rab-o-loid pa-ram-e-trize para-mount pen-ta-gon phe-nom-e-non
post-script pre-am-ble pro-ce-dur-al pro-hib-i-tive pro-hib-i-tive-ly
pseu-do-dif-fer-en-tial pseu-do-fi-nite pseu-do-nym qua-drat-ic
quad-ra-ture qua-si-smooth qua-si-sta-tion-ary qua-si-tri-an-gu-lar
quin-tes-sence quin-tes-sen-tial re-arrange-ment rec-tan-gle
ret-ri-bu-tion retro-fit retro-fit-ted right-eous right-eous-ness
ro-bot ro-bot-ics sched-ul-ing se-mes-ter semi-def-i-nite
semi-ho-mo-thet-ic set-up se-vere-ly side-step sov-er-eign spe-cious
spher-oid spher-oid-al star-tling star-tling-ly sta-tis-tics
sto-chas-tic straight-est strange-ness strat-a-gem strong-hold
sum-ma-ble symp-to-matic syn-chro-nous topo-graph-i-cal tra-vers-a-ble
tra-ver-sal tra-ver-sals treach-ery turn-around un-at-tached
un-err-ing-ly white-space wide-spread wing-spread wretch-ed
wretch-ed-ly Brown-ian Eng-lish Euler-ian Feb-ru-ary Gauss-ian
Grothen-dieck Hamil-ton-ian Her-mit-ian Jan-u-ary Japan-ese Kor-te-weg
Le-gendre Lip-schitz Lip-schitz-ian Mar-kov-ian Noe-ther-ian
No-vem-ber Rie-mann-ian Schwarz-schild Sep-tem-ber}

\def\leftheadtext#1{} \def\rightheadtext#1{}

\let\nopagenumber\p@gen@false \let\putpagenumber\p@gen@true
\let\pagefirst\nopagenumber \let\pagenext\putpagenumber

\else

\amsppttrue

\let\twelvepoint\relax \let\Twelvepoint\relax \let\putpagenumber\relax
\let\logo@\relax \let\pagefirst\firstpage@true \let\pagenext\firstpage@false
\def\nopagenumber{\let\f@li@ld\folio\def\folio{\global\let\folio\f@li@ld}}

\def\ftext#1{\footnotetext""{\vsk-.8>\nt #1}}

\def\m@R@f@[#1]{\Refs\nofrills{}\m@th\tenpoint
 {%\let\{\relax
 \s@ct{}{#1}}\wr@@c{\string\Refcd{#1}{\the\pageno}}
 \def\k@yf@##1{\hss[##1]\enspace} \let\keyformat\k@yf@
 \def\widest##1{\s@twd@\refindentwd{\tenpoint\k@yf@{##1}}}
 \let\Key\key \def\refin{\kern\refindentwd}}
\let\info\finalinfo \r@R@fs\Refs
\def\adjust#1{#1} \let\adjustend\relax
\let\adjustnext\adjust 

\fi

\outer\def\myRefs{\myR@fs} \r@st@re\proclaim
\def\bye{\par\vfill\supereject\cl@selbl\cl@secd\b@e} \r@endd@\b@e
\let\Cite\cite \let\Key\key \def\endpro{\par\endproclaim}
\let\d@c@\document \def\document{\d@c@\tenpoint}
\hyphenation{ortho-gon-al}

\newtoks\@@tp@t \@@tp@t\output
\output=\@ft@{\let\{\noexpand\the\@@tp@t}
\let\{\relax

\newif\ifVersion

\def\p@n@l#1{\ifnum#1=\z@\else\penalty#1\relax\fi}

\def\s@ct#1#2{\ifVersion
 \skip@\lastskip\ifdim\skip@<1.5\bls\vskip-\skip@\p@n@l{-200}\vsk.5>%
 \p@n@l{-200}\vsk.5>\p@n@l{-200}\vsk.5>\p@n@l{-200}\vsk-1.5>\else
 \p@n@l{-200}\fi\ifdim\skip@<.9\bls\vsk.9>\else
 \ifdim\skip@<1.5\bls\vskip\skip@\fi\fi
 \vtop{\twelvepoint\raggedright\s@cf@nt\vp1\vsk->\vskip.16ex
 \s@twd@\parindent{#1}%
 \ifdim\parindent>\z@\adv\parindent.5em\fi\hang\textindent{#1}#2\strut}
 \else
 \p@sk@p{-200}{.8\bls}\vtop{\s@cf@nt\s@twd@\parindent{#1}%
 \ifdim\parindent>\z@\adv\parindent.5em\fi\hang\textindent{#1}#2\strut}\fi
 \nointerlineskip\nobreak\vtop{\strut}\nobreak\vskip-.6\bls\nobreak}

\def\s@bs@ct#1#2{\ifVersion
 \skip@\lastskip\ifdim\skip@<1.5\bls\vskip-\skip@\p@n@l{-200}\vsk.5>%
 \p@n@l{-200}\vsk.5>\p@n@l{-200}\vsk.5>\p@n@l{-200}\vsk-1.5>\else
 \p@n@l{-200}\fi\ifdim\skip@<.9\bls\vsk.9>\else
 \ifdim\skip@<1.5\bls\vskip\skip@\fi\fi
 \vtop{\elevenpoint\raggedright\s@bf@nt\vp1\vsk->\vskip.16ex%
 \s@twd@\parindent{#1}\ifdim\parindent>\z@\adv\parindent.5em\fi
 \hang\textindent{#1}#2\strut}
 \else
 \p@sk@p{-200}{.6\bls}\vtop{\s@bf@nt\s@twd@\parindent{#1}%
 \ifdim\parindent>\z@\adv\parindent.5em\fi\hang\textindent{#1}#2\strut}\fi
 \nointerlineskip\nobreak\vtop{\strut}\nobreak\vskip-.8\bls\nobreak}

\def\gadv{\global\adv} \def\gad#1{\gadv#1\@ne} \def\gadneg#1{\gadv#1-\@ne}

\newcount\t@@n \t@@n=10 \newbox\testbox

\newcount\Sno \newcount\Lno \newcount\Fno

\def\pr@cl#1{\r@st@re\pr@c@\pr@c@{#1}\global\let\pr@c@\relax}

\def\tagg#1{\tag"\rlap{\rm(#1)}\kern.01\p@"}
\def\l@L#1{\l@bel{#1}L} \def\l@F#1{\l@bel{#1}F} \def\<#1>{\l@b@l{#1}F}
\def\Tag#1{\tag{\l@F{#1}}} \def\Tagg#1{\tagg{\l@F{#1}}}
\def\Rem{\demo{\sl Remark}} \def\Ex{\demo{\bf Example}}
\def\Pf#1.{\demo{Proof #1}} \def\epf{\qed\enddemo}
\def\Ap@x{Appendix}
\def\Appendix{\Sno=64 \t@@n\@ne \wr@@c{\string\Appencd}
 \def\sf@rm{\char\the\Sno} \def\sf@rm@{\Ap@x\space\sf@rm} \def\sf@rm@@{\Ap@x}
 \def\s@ct@n##1##2{\s@ct\empty{\setbox\z@\hbox{##1}\ifdim\wd\z@=\z@
 \if##2*\sf@rm@@\else\if##2.\sf@rm@@.\else##2\fi\fi\else
 \if##2*\sf@rm@\else\if##2.\sf@rm@.\else\sf@rm@.\enspace##2\fi\fi\fi}}}
\def\Appcd#1#2#3{\def\Ap@@{\hglue-\l@ftcd\Ap@x}\ifx\@ppl@ne\empty
 \def\l@@b{\@fwd@@{#1}{\space#1}{}}\if*#2\entcd{}{\Ap@@\l@@b}{#3}\else
 \if.#2\entcd{}{\Ap@@\l@@b.}{#3}\else\entcd{}{\Ap@@\l@@b.\enspace#2}{#3}\fi\fi
 \else\def\l@@b{\@fwd@@{#1}{\c@l@b{#1}}{}}\if*#2\entcd{\l@@b}{\Ap@x}{#3}\else
 \if.#2\entcd{\l@@b}{\Ap@x.}{#3}\else\entcd{\l@@b}{#2}{#3}\fi\fi\fi}

\let\s@ct@n\s@ct
\def\s@ct@@[#1]#2{\@ft@\xdef\csname @#1@S@\endcsname{\sf@rm}\wr@@x{}%
 \wr@@x{\string\labeldef{S}\space{\?#1@S?}\space{#1}}%
 {%\let\{\relax
 \s@ct@n{\sf@rm@}{#2}}\wr@@c{\string\Entcd{\?#1@S?}{#2}{\the\pageno}}}
\def\s@ct@#1{\wr@@x{}{%\let\{\relax
 \s@ct@n{\sf@rm@}{#1}}\wr@@c{\string\Entcd{\sf@rm}{#1}{\the\pageno}}}
\def\s@ct@e[#1]#2{\@ft@\xdef\csname @#1@S@\endcsname{\sf@rm}\wr@@x{}%
 \wr@@x{\string\labeldef{S}\space{\?#1@S?}\space{#1}}%
 {%\let\{\relax
 \s@ct@n\empty{#2}}\wr@@c{\string\Entcd{}{#2}{\the\pageno}}}
\def\s@cte#1{\wr@@x{}{%\let\{\relax
 \s@ct@n\empty{#1}}\wr@@c{\string\Entcd{}{#1}{\the\pageno}}}
\def\theSno#1#2{\dff\?#1@S?{#2}%
 \wr@@x{\string\labeldef{S}\space{#2}\space{#1}}\fi}

\newif\ifd@bn@\d@bn@true
\def\Section{\gad\Sno\ifd@bn@\Fno\z@\Lno\z@\fi\@fn@xt[\s@ct@@\s@ct@}
\def\section{\gad\Sno\ifd@bn@\Fno\z@\Lno\z@\fi\@fn@xt[\s@ct@e\s@cte}
\let\Sect\Section 
\def\subsection{\@fn@xt*\subs@ct@\subs@ct}
\def\subs@ct#1{{%\let\{\relax
 \s@bs@ct\empty{#1}}\wr@@c{\string\subcd{#1}{\the\pageno}}}
\def\subs@ct@*#1{\vsk->\vsk>{%\let\{\relax
 \s@bs@ct\empty{#1}}\wr@@c{\string\subcd{#1}{\the\pageno}}}
\let\subsect\subsection \def\Snodef#1{\Sno #1}

\def\l@b@l#1#2{\def\n@@{\csname #2no\endcsname}%
 \if*#1\gad\n@@ \@ft@\xdef\csname @#1@#2@\endcsname{\l@f@rm}\else\def\t@st{#1}%
 \ifx\t@st\empty\gad\n@@ \@ft@\xdef\csname @#1@#2@\endcsname{\l@f@rm}%
 \else\@ft@\ifx\csname @#1@#2@mark\endcsname\relax\gad\n@@
 \@ft@\xdef\csname @#1@#2@\endcsname{\l@f@rm}%
 \@ft@\gdef\csname @#1@#2@mark\endcsname{}%
 \wr@@x{\string\labeldef{#2}\space{\?#1@#2?}\space\ifnum\n@@<10 \space\fi{#1}}%
 \fi\fi\fi}
\def\labeldef#1#2#3{\dff\?#3@#1?{#2}}
\def\Labeldef#1#2#3{\dff\?#3@#1?{#2}\@ft@\gdef\csname @#3@#1@mark\endcsname{}}

\def\l@bel#1#2{\l@b@l{#1}{#2}\?#1@#2?}

\newcount\c@cite
\def\?#1?{\csname @#1@\endcsname}
\def\[{\@fn@xt:\c@t@sect\c@t@}
\def\c@t@#1]{{\c@cite\z@\@fwd@@{\?#1@L?}{\adv\c@cite1}{}%
 \@fwd@@{\?#1@F?}{\adv\c@cite1}{}\@fwd@@{\?#1?}{\adv\c@cite1}{}%
 \relax\ifnum\c@cite=\z@{\bf ???}\wrs@x{No label [#1]}\else
 \ifnum\c@cite=1\let\@@PS\relax\let\@@@\relax\else\let\@@PS\underbar
 \def\@@@{{\rm<}}\fi\@@PS{\?#1?\@@@\?#1@L?\@@@\?#1@F?}\fi}}
\def\(#1){{\rm(\c@t@#1])}}
\def\c@t@s@ct#1{\@fwd@@{\?#1@S?}{\?#1@S?\relax}%
 {{\bf ???}\wrs@x{No section label {#1}}}}
\def\c@t@sect:#1]{\c@t@s@ct{#1}} \let\SNo\c@t@s@ct

\newdimen\l@ftcd \newdimen\r@ghtcd \let\nlc\relax

\def\d@tt@d{\leaders\hbox to 1em{\kern.1em.\hfil}\hfill}
\def\entcd#1#2#3{\item{\l@bcdf@nt#1}{\entcdf@nt#2}\alb\kern.9em\hbox{}%
 \kern-.9em\d@tt@d\kern-.36em{\p@g@cdf@nt#3}\kern-\r@ghtcd\hbox{}\par}
\def\Entcd#1#2#3{\def\l@@b{\@fwd@@{#1}{\c@l@b{#1}}{}}\vsk.2>%
 \entcd{\l@@b}{#2}{#3}}
\def\subcd#1#2{{\adv\leftskip.333em\entcd{}{\s@bcdf@nt#1}{#2}}}
\def\Refcd#1#2{\def\t@@st{#1}\ifx\t@@st\empty\ifx\r@fl@ne\empty\relax\else
 \R@fcd{\r@fl@ne}{#2}\fi\else\R@fcd{#1}{#2}\fi}
\def\R@fcd#1#2{\sk@@p{.6\bls}\entcd{}{\hglue-\l@ftcd\R@fcdf@nt #1}{#2}}
\def\Refline{\def\r@fl@ne} \def\Refempty{\let\r@fl@ne\empty}
\def\Appencd{\par\adv\leftskip-\l@ftcd\adv\rightskip-\r@ghtcd\@ppl@ne
 \adv\leftskip\l@ftcd\adv\rightskip\r@ghtcd\let\Entcd\Appcd}
\def\appline{\def\@ppl@ne} \def\Appempty{\let\@ppl@ne\empty}
\def\Appline#1{\def\@ppl@ne{\s@bs@ct{}{#1}}}
\def\leftcd#1{\adv\leftskip-\l@ftcd\s@twd@\l@ftcd{\c@l@b{#1}\enspace}
 \adv\leftskip\l@ftcd}
\def\rightcd#1{\adv\rightskip-\r@ghtcd\s@twd@\r@ghtcd{#1\enspace}
 \adv\rightskip\r@ghtcd}
\def\C@nt{Contents} \def\Ap@s{Appendices} \def\R@fcs{References}
\def\contents{\@fn@xt*\cont@@\cont@}
\def\cont@{\@fn@xt[\cnt@{\cnt@[\C@nt]}}
\def\cont@@*{\@fn@xt[\cnt@@{\cnt@@[\C@nt]}}
\def\cnt@[#1]{\c@nt@{M}{#1}{44}{\s@bs@ct{}{\@ppl@f@nt\Ap@s}}}
\def\cnt@@[#1]{\c@nt@{M}{#1}{44}{}}
\def\endco{\par\penalty-500\vsk>\vskip\z@\endgroup}
\def\readcd{\@np@t{\jobname.cd}}
\def\Cde{\@fn@xt*\Cde@@\Cde@}
\def\Cde@{\@fn@xt[\Cd@{\Cd@[\C@nt]}}
\def\Cde@@*{\@fn@xt[\Cd@@{\Cd@@[\C@nt]}}
\def\Cd@[#1]{\cnt@[#1]\readcd\endco}
\def\Cd@@[#1]{\cnt@@[#1]\readcd\endco}
\def\contlabeldef{\def\c@l@b}

\long\def\c@nt@#1#2#3#4{\s@twd@\l@ftcd{\c@l@b{#1}\enspace}
 \s@twd@\r@ghtcd{#3\enspace}\adv\r@ghtcd1.333em
 \def\@ppl@ne{#4}\def\r@fl@ne{\R@fcs}\s@ct{}{#2}\B@gr@\parindent\z@\let\nlc\nl
 \let\nl\relax\parskip.2\bls\adv\leftskip\l@ftcd\adv\rightskip\r@ghtcd}

\def\writecd{\immediate\openout\@@cd\jobname.cd \def\wr@@c{\write\@@cd}
 \def\cl@secd{\immediate\write\@@cd{\string\endinput}\immediate\closeout\@@cd}
 \def\closecd{\cl@secd\global\let\cl@secd\relax}}
\let\cl@secd\relax \def\wr@@c#1{} \let\closecd\relax

\def\dff{\@ft@\d@f} \def\d@f{\@ft@\def}
\def\edff{\@ft@\ed@f} \def\ed@f{\@ft@\edef}
\def\defi#1#2{\def#1{#2}\wr@@x{\string\def\string#1{#2}}}

\def\qed{\hbox{}\nobreak\hfill\nobreak{\m@th$\,\square$}}
\def\back#1 {\strut\kern-.33em #1\enspace\ignore} %% !!! a space after #1 !!!
\def\Text#1{\crcr\noalign{\alb\vsk>\normalbaselines\vsk->\vbox{\nt #1\strut}%
 \nobreak\nointerlineskip\vbox{\strut}\nobreak\vsk->\nobreak}}

\def\hcor#1{\advance\hoffset by #1}
\def\vcor#1{\advance\voffset by #1}
\let\bls\baselineskip \let\ignore\ignorespaces
\ifx\ic@\undefined \let\ic@\/\fi
\def\vsk#1>{\vskip#1\bls} \let\adv\advance
\def\vv#1>{\vadjust{\vsk#1>}\ignore}
\def\vvn#1>{\vadjust{\nobreak\vsk#1>\nobreak}\ignore}
\def\vvv#1>{\vskip\z@\vsk#1>\nt\ignore}
\def\vvgood{\vadjust{\penalty-500}} 
\def\Par{\vsk.5>} \def\setparindent{\edef\Parindent{\the\parindent}}
\def\Type{\vsk.5>\bgroup\parindent\z@\tt\rightskip\z@ plus1em minus1em%
 \spaceskip.3333em \xspaceskip.5em\relax}
\def\endType{\vsk.5>\egroup\nt} 

\let\Hat\widehat \let\Tilde\widetilde \let\dollar\$ \let\ampersand\&
\let\sss\scriptscriptstyle  
\let\vp\vphantom \let\hp\hphantom \let\nt\noindent
\let\cline\centerline \let\lline\leftline \let\rline\rightline
\def\nn#1>{\noalign{\vskip#1\p@@}} \def\NN#1>{\openup#1\p@@}
\def\cnn#1>{\noalign{\vsk#1>}}
 
\let\Lim\lim \def\lim{\Lim\limits} \let\Sum\sum \def\sum{\Sum\limits}
\def\Plus{\bigoplus\limits} 
\let\Prod\prod \def\prod{\Prod\limits} \let\Int\int \def\int{\Int\limits}

\def\tsum{\mathop{\tsize\Sum}\limits} 
\def\tprod{\mathop{\tsize\Prod}\limits}
\def\&{.\kern.1em} \def\>{{\!\;}} \def\]{{\!\!\;}} \def\){\>\]} \def\}{\]\]}
\def\nl{\leavevmode\hfill\break} \def\~{\leavevmode\@fn@xt~\m@n@s\@md@sh}
\def\m@n@s~{\raise.15ex\mbox{-}} \def\@md@sh{\raise.13ex\hbox{--}}
\let\procent\% \def\%#1{\ifmmode\mathop{#1}\limits\else\procent#1\fi}
\let\@ml@t\" \def\"#1{\ifmmode ^{(#1)}\else\@ml@t#1\fi}
\let\@c@t@\' \def\'#1{\ifmmode _{(#1)}\else\@c@t@#1\fi}
\let\colon\: \def\:{^{\vp|}}

\let\texspace\ \def\ {\ifmmode\alb\fi\texspace}
%%% Do not remove a space after \ !!!

\let\n@wp@ge\newpage \def\newpage{\endgraf\n@wp@ge}
\let\=\m@th \def\mbox#1{\hbox{\m@th$#1$}}
\def\mtext#1{\text{\m@th$#1$}} \def\^#1{\text{\m@th#1}}
\def\Line#1{\kern-.5\hsize\line{\m@th$\dsize#1$}\kern-.5\hsize}
\def\Lline#1{\kern-.5\hsize\lline{\m@th$\dsize#1$}\kern-.5\hsize}
\def\Cline#1{\kern-.5\hsize\cline{\m@th$\dsize#1$}\kern-.5\hsize}
\def\Rline#1{\kern-.5\hsize\rline{\m@th$\dsize#1$}\kern-.5\hsize}

\def\Ll@p#1{\llap{\m@th$#1$}} \def\Rl@p#1{\rlap{\m@th$#1$}}
 \def\Cl@p#1{\llap{\m@th$#1$\hss}}
\def\Llap#1{\mathchoice{\Ll@p{\dsize#1}}{\Ll@p{\tsize#1}}{\Ll@p{\ssize#1}}%
 {\Ll@p{\sss#1}}}
\def\Clap#1{\mathchoice{\Cl@p{\dsize#1}}{\Cl@p{\tsize#1}}{\Cl@p{\ssize#1}}%
 {\Cl@p{\sss#1}}}
\def\Rlap#1{\mathchoice{\Rl@p{\dsize#1}}{\Rl@p{\tsize#1}}{\Rl@p{\ssize#1}}%
 {\Rl@p{\sss#1}}}
 
\def\LRtph#1#2{\setbox\z@\hbox{#1}\dimen\z@\wd\z@\hbox{\hbox to\dimen\z@{#2}}}
\def\LRph#1#2{\LRtph{\m@th$#1$}{\m@th$#2$}}

\def\CCph#1#2{\LRph{#1}{\hss#2\hss}}
 \def\RRph#1#2{\LRph{#1}{#2\hss}}
\def\Cph#1#2{\mathchoice{\CCph{\dsize#1}{\dsize#2}}{\CCph{\tsize#1}{\tsize#2}}
 {\CCph{\ssize#1}{\ssize#2}}{\CCph{\sss#1}{\sss#2}}}

\def\Rph#1#2{\mathchoice{\RRph{\dsize#1}{\dsize#2}}{\RRph{\tsize#1}{\tsize#2}}
 {\RRph{\ssize#1}{\ssize#2}}{\RRph{\sss#1}{\sss#2}}}
\def\Lto#1{\setbox\z@\mbox{\tsize{#1}}%
 \mathrel{\mathop{\hbox to\wd\z@{\rightarrowfill}}\limits#1}}
\def\Lgets#1{\setbox\z@\mbox{\tsize{#1}}%
 \mathrel{\mathop{\hbox to\wd\z@{\leftarrowfill}}\limits#1}}
\def\vpb#1{{\vp{\big(}}^{\]#1}} \def\vpp#1{{\vp{\big]}}_{#1}}
\def\lbc{\mathopen{[\![}} \def\rbc{\mathclose{]\!]}}
\def\Lbc{\mathopen{\big[\!\}\big[}} \def\Rbc{\mathclose{\big]\!\}\big]}}
\def\LBc{\mathopen{\Big[\!\}\Big[}} \def\RBc{\mathclose{\Big]\!\}\Big]}}

\let\alb\allowbreak 
\def\ald{\noalign{\alb}} \let\alds\allowdisplaybreaks

\let\o\circ \let\x\times \let\ox\otimes 
\let\sub\subset  \let\tabs\+
\let\le\leqslant \let\ge\geqslant
 \let\8\infty \let\*\star
\let\bra\langle \let\ket\rangle
\let\lto\longrightarrow 
\let\map\mapsto  \let\hto\hookrightarrow
 \let\rto\rightrightarrows
 \def\vert{\ |\ } 
\let\Empty\varnothing 
 
\let\lb\lbrace \let\rb\rbrace
\let\lf\lfloor \let\rf\rfloor
\let\trileft\triangleleft \let\tright\triangleright

\def\lsym#1{#1\alb\ldots\relax#1\alb}
\def\lc{\lsym,}  \def\lx{\lsym\x} \def\lox{\lsym\ox}
\def\llc{\,,\alb\ {\ldots\ ,}\alb\ }
 \def\Im{\mathop{\roman{Im}\>}}
\def\End{\mathop{\roman{End}\>}} \def\Hom{\mathop{\roman{Hom}\>}}
 
 \def\Ker{\mathop{\roman{Ker}\>}}

\def\Res{\mathop{\roman{Res}\>}\limits}

\def\sgn{\mathop{\roman{sgn}\)}\limits}
 
\def\id{\roman{id}}  \def\for{\text{for }\,}
\def\1{^{-1}} \def\_#1{_{\Rlap{#1}}}
\def\vst#1{{\lower1.9\p@@\mbox{\bigr|_{\raise.5\p@@\mbox{\ssize#1}}}}}
\def\vrp#1:#2>{{\vrule height#1 depth#2 width\z@}}
\def\vru#1>{\vrp#1:\z@>} \def\vrd#1>{\vrp\z@:#1>}
\def\qqq{\qquad\quad} 
\def\sscr#1{\raise.3ex\mbox{\sss#1}} \def\@@PS{\bold{OOPS!!!}}
\def\ono{\bigl(1+o(1)\bigr)} 

\def\intcl{\mathop
 {\Rlap{\raise.3ex\mbox{\kern.12em\curvearrowleft}}\int}\limits}
\def\intcr{\mathop
 {\Rlap{\raise.3ex\mbox{\kern.24em\curvearrowright}}\int}\limits}

\def\pms{\raise.25ex\mbox{\ssize\pm}\>}
\def\mps{\raise.25ex\mbox{\ssize\mp}\>}
\def\pss{{\sscr+}} \def\mss{{\sscr-}}

\let\al\alpha
\let\bt\beta
\let\gm\gamma \let\Gm\Gamma \let\Gms\varGamma
\let\dl\delta \let\Dl\Delta 
\let\epe\epsilon \let\eps\varepsilon \let\epsilon\eps
\let\Xj\varXi
\let\zt\zeta
\let\tht\theta \let\Tht\Theta
\let\thi\vartheta \let\Tho\varTheta

\let\ka\kappa
\let\la\lambda \let\La\Lambda

\let\si\sigma 
 \let\Sig\varSigma
\let\pho\phi \let\phi\varphi \let\Pho\varPhi

\let\om\omega \let\Om\Omega \let\Ome\varOmega

\def\C{\Bbb C}

\def\Z{\Bbb Z}

\def\AA{\Bbb A}
\def\BB{\Bbb B}

\def\FF{\Bbb F}

\def\SS{\Bbb S}
\def\TT{\Bbb T}
\def\UU{\Bbb U}

\def\Zp{\Z_{\ge 0}} 
\def\Zn{\Z_{\le 0}} \def\Znn{\Z_{<0}}

\def\difl/{differential} \def\dif/{difference}
\def\cf.{cf.\ \ignore} \def\Cf.{Cf.\ \ignore}
\def\egv/{eigenvector} \def\eva/{eigenvalue} \def\eq/{equation}
\def\lhs/{the left hand side} \def\rhs/{the right hand side}
\def\Lhs/{The left hand side} \def\Rhs/{The right hand side}
\def\gby/{generated by} \def\wrt/{with respect to} \def\st/{such that}
\def\resp/{respectively} \def\off/{offdiagonal} \def\wt/{weight}
\def\pol/{polynomial} \def\rat/{rational} \def\tri/{trigonometric}
\def\fn/{function} \def\var/{variable} \def\raf/{\rat/ \fn/}
\def\inv/{invariant} \def\hol/{holomorphic} \def\hof/{\hol/ \fn/}
\def\mer/{meromorphic} \def\mef/{\mer/ \fn/} \def\mult/{multiplicity}
\def\sym/{symmetric} \def\perm/{permutation} \def\fd/{finite-dimensional}
\def\rep/{representation} \def\irr/{irreducible} \def\irrep/{\irr/ \rep/}
\def\hom/{homomorphism} \def\aut/{automorphism} \def\iso/{isomorphism}
\def\lex/{lexicographical} \def\as/{asymptotic} \def\asex/{\as/ expansion}
\def\ndeg/{nondegenerate} \def\neib/{neighbourhood} \def\deq/{\dif/ \eq/}
\def\hw/{highest \wt/} \def\gv/{generating vector} \def\eqv/{equivalent}
\def\msd/{method of steepest descend} \def\pd/{pairwise distinct}
\def\wlg/{without loss of generality} \def\Wlg/{Without loss of generality}
\def\onedim/{one-dimensional} \def\qcl/{quasiclassical}
\def\hgeom/{hyper\-geometric} \def\hint/{\hgeom/ integral}
\def\hwm/{\hw/ module} \def\emod/{evaluation module} \def\Vmod/{Verma module}
\def\symg/{\sym/ group} \def\sol/{solution} \def\eval/{evaluation}
\def\anf/{analytic \fn/} \def\anco/{analytic continuation}
\def\qg/{quantum group} \def\qaff/{quantum affine algebra}

\def\Rm/{\^{$R$-}matrix} \def\Rms/{\^{$R$-}matrices} \def\YB/{Yang-Baxter \eq/}
\def\Ba/{Bethe ansatz} \def\Bv/{Bethe vector} \def\Bae/{\Ba/ \eq/}
\def\KZv/{Knizh\-nik-Zamo\-lod\-chi\-kov} \def\KZvB/{\KZv/-Bernard}
\def\KZ/{{\sl KZ\/}} \def\qKZ/{{\sl qKZ\/}}
\def\KZB/{{\sl KZB\/}} \def\qKZB/{{\sl qKZB\/}}
\def\qKZo/{\qKZ/ operator} \def\qKZc/{\qKZ/ connection}
\def\KZe/{\KZ/ \eq/} \def\qKZe/{\qKZ/ \eq/} \def\qKZBe/{\qKZB/ \eq/}

\def\h@ph{\discretionary{}{}{-}} \def\$#1$-{\,\^{$#1$}\h@ph}

\def\TFT/{Research Insitute for Theoretical Physics}
\def\HY/{University of Helsinki} \def\AoF/{the Academy of Finland}
\def\CNRS/{Supported in part by MAE\~MICECO\~CNRS Fellowship}
\def\LPT/{Laboratoire de Physique Th\'eorique ENSLAPP}
\def\ENSLyon/{\'Ecole Normale Sup\'erieure de Lyon}
\def\LPTaddr/{46, All\'ee d'Italie, 69364 Lyon Cedex 07, France}
\def\enslapp/{URA 14\~36 du CNRS, associ\'ee \`a l'E.N.S.\ de Lyon,
au LAPP d'Annecy et \`a l'Universit\`e de Savoie}
\def\ensemail/{vtarasov\@ enslapp.ens-lyon.fr}
\def\DMS/{Department of Mathematics, Faculty of Science}
\def\DMO/{\DMS/, Osaka University}
\def\DMOaddr/{Toyonaka, Osaka 560, Japan}
\def\dmoemail/{vt\@ math.sci.osaka-u.ac.jp}
\def\SPb/{St\&Peters\-burg}
\def\home/{\SPb/ Branch of Steklov Mathematical Institute}
\def\homeaddr/{Fontanka 27, \SPb/ \,191011, Russia}
\def\homemail/{vt\@ pdmi.ras.ru}
\def\absence/{On leave of absence from \home/}
\def\UNC/{Department of Mathematics, University of North Carolina}
\def\ChH/{Chapel Hill}
\def\UNCaddr/{\ChH/, NC 27599, USA} \def\avemail/{av\@ math.unc.edu}
\def\grant/{NSF grant DMS\~9501290}	%%% Felder's grant no. 9400841
\def\Grant/{Supported in part by \grant/}

\def\Aomoto/{K\&Aomoto}
\def\Dri/{V\]\&G\&Drin\-feld}
\def\Fadd/{L\&D\&Fad\-deev}
\def\Feld/{G\&Felder}
\def\Fre/{I\&B\&Fren\-kel}
\def\Gustaf/{R\&A\&Gustafson}
\def\Kazh/{D\&Kazhdan} \def\Kir/{A\&N\&Kiril\-lov}
\def\Kor/{V\]\&E\&Kore\-pin}
\def\Lusz/{G\&Lusztig}
\def\MN/{M\&Naza\-rov}
\def\Resh/{N\&Reshe\-ti\-khin} \def\Reshy/{N\&\]Yu\&Reshe\-ti\-khin}
\def\SchV/{V\]\&\]V\]\&Schecht\-man} \def\Sch/{V\]\&Schecht\-man}
\def\Skl/{E\&K\&Sklya\-nin}
\def\Smirn/{F\]\&Smirnov} \def\Smirnov/{F\]\&A\&Smirnov}
\def\Takh/{L\&A\&Takh\-tajan}
\def\VT/{V\]\&Ta\-ra\-sov} \def\VoT/{V\]\&O\&Ta\-ra\-sov}
\def\Varch/{A\&\]Var\-chenko} \def\Varn/{A\&N\&\]Var\-chenko}

\def\AMS/{Amer.\ Math.\ Society}
\def\CMP/{Comm.\ Math.\ Phys.{}}
\def\DMJ/{Duke.\ Math.\ J.{}}
\def\Inv/{Invent.\ Math.{}} %% Inventiones Mathematicae
\def\IMRN/{Int.\ Math.\ Res.\ Notices}
\def\JPA/{J.\ Phys.\ A{}}
\def\JSM/{J.\ Soviet\ Math.{}}
\def\LMP/{Lett.\ Math.\ Phys.{}}
\def\LMJ/{Leningrad Math.\ J.{}}
\def\LpMJ/{\SPb/ Math.\ J.{}}
\def\SIAM/{SIAM J.\ Math.\ Anal.{}}
\def\SMNS/{Selecta Math., New Series}
\def\TMP/{Theor.\ Math.\ Phys.{}}
\def\ZNS/{Zap.\ nauch.\ semin. LOMI}

\def\ASMP/{Advanced Series in Math.\ Phys.{}}

\def\AMSa/{AMS \publaddr Providence}
\def\Birk/{Birkh\"auser}
\def\CUP/{Cambridge University Press} \def\CUPa/{\CUP/ \publaddr Cambridge}
\def\Spri/{Springer-Verlag} \def\Spria/{\Spri/ \publaddr Berlin}
\def\WS/{World Scientific} \def\WSa/{\WS/ \publaddr Singapore}

\newbox\lefthbox \newbox\righthbox

\let\sectsep. \let\labelsep. \let\contsep. \let\labelspace\relax
\let\sectpre\relax \let\contpre\relax
\def\sf@rm{\the\Sno} \def\sf@rm@{\sectpre\sf@rm\sectsep}
\def\c@l@b#1{\contpre#1\contsep}
\def\l@f@rm{\ifd@bn@\sf@rm\labelsep\fi\labelspace\the\n@@}

\def\sectformdef{\def\sf@rm}

\let\DoubleNum\d@bn@true \let\SingleNum\d@bn@false

\def\NoNewNum{\let\writeldf\relax\def\l@b@l##1##2{\if*##1%
 \@ft@\xdef\csname @##1@##2@\endcsname{\mbox{*{*}*}}\fi}}
\def\NoNewTime{\def\todaydef##1{\def\today{##1}}
 \def\nowtimedef##1{\def\nowtime{##1}}}
\def\NoInput{\let\Input\input\let\writeldf\relax}
\def\Fixed{\NoNewTime\NoInput}

\def\sectfont#1{\def\s@cf@nt{#1}} \sectfont\bf
\def\subsectfont#1{\def\s@bf@nt{#1}} \subsectfont\it
\def\Entcdfont#1{\def\entcdf@nt{#1}} \Entcdfont\relax
\def\labelcdfont#1{\def\l@bcdf@nt{#1}} \labelcdfont\relax
\def\pagecdfont#1{\def\p@g@cdf@nt{#1}} \pagecdfont\relax
\def\subcdfont#1{\def\s@bcdf@nt{#1}} \subcdfont\it
\def\applefont#1{\def\@ppl@f@nt{#1}} \applefont\bf
\def\Refcdfont#1{\def\R@fcdf@nt{#1}} \Refcdfont\bf

\tenpoint

\csname beta.def\endcsname

\Fixed

\expandafter\ifx\csname trig.def\endcsname\relax \else\endinput\fi
\expandafter\edef\csname trig.def\endcsname{%
 \catcode`\noexpand\@=\the\catcode`\@\space}
\catcode`\@=11

\font\Brm=cmr12 scaled 1200
\font\Bit=cmti12 scaled 1200

\ifx\twelvemsa\undefined\font\twelvemsa=msam10 scaled 1200\fi
\ifMag\else\fi

\newif\ifPcd \let\noPcd\Pcdfalse \let\makePcd\Pcdtrue
\def\NoPcde{\let\noPcd\relax \let\makePcd\relax}

\def\Th#1{\pr@cl{(\l@F{#1}) Theorem}\ignore}
\def\Lm#1{\pr@cl{(\l@F{#1}) Lemma}\ignore}
\def\Cr#1{\pr@cl{(\l@F{#1}) Corollary}\ignore}
\def\Df#1{\pr@cl{(\l@F{#1}) Definition}\ignore}
\def\Cj#1{\pr@cl{(\l@F{#1}) Conjecture}\ignore}
\def\Prop#1{\pr@cl{(\l@F{#1}) Proposition}\ignore}
\def\Pf#1.{\demo{\let\{\relax Proof #1}\ifPcd\def\t@st@{#1}\ifx\t@st@\empty
 \else\wr@@c{\string\subcd{Proof #1}{\the\pageno}}\fi\fi\ignore}

\def\q{{q\1}} \let\one\id \def\pii{\pi i}
\def\neps{\ne{\)p^s}} \def\nepsr{\ne{\)p^s\eta^r}} \def\netr{\>\ne\)\eta^r}

\def\CC{\bold C}

\def\S{\bold S}
\def\sb{\bold s}

\def\Ec{\Cal E}
\def\Fo{\Cal F}
\def\H{\Cal H}
\def\Qc{\Cal Q}
\def\Rc{\Cal R}
\def\Zc{\Cal Z}

\def\eg{\frak e}
\def\hg{\frak h}
\def\lg{\frak l}
\def\mg{\frak m}
\def\ng{\frak n}

\def\A{\roman A}

\def\Ff{\roman F}
\def\Mu{\roman M}
\def\Nu{\roman N}

\def\Fun{\slanted{Fun}{\)}}

\let\Ks\ka %% \def\Ks{{\ssize K}}

\let\Is I
\def\Ib{\bar I}
\def\Phb{{\,\)\overline{\!\]\Phi\]\!}\)\,}}

\def\Fwh{\Rlap{\,\Hat{\!\phantom\Fo}}\Fo}

\def\Sh{\hat S}

\def\At{\tilde A}
\def\Fqt{\Tilde\Fo_{\!\!\]\sss e\]l\}l}}
\def\Fti{\Tilde F}
\def\Gti{\Tilde G}

\def\Pht{\Tilde\Phi}

\def\Rti{{\>\Tilde{\}R\>}\}}}

\def\Tti{\Tilde T}
\def\TTt{\Tilde\TT}

\def\Wti{\Rlap{\>\Tilde{\}\phantom W\]}}W}
\def\wti{\tilde w}

\def\Xti{\Tilde X}

\def\Cn{\C^{\>n}} \def\Cl{\C^{\,\ell}} \def\Cln{\C^{\,\ell+n}}
 
\def\Cx{\C^\x} \def\Cxn{\C^{\x n}} \def\Cxl{\C^{\x\ell}}
\def\Cxln{\C^{\x(\ell+n)}}

 \def\BE{B\!E}
 \def\Fq{\Fo_{\!\!\]\sss e\]l\}l}}
\def\Fsi{\Fwh_{\vp|\!\}\sss\Sig}} \def\Fz{F_{\]z}}

\def\Mq{M^{\sscr{e\]l\}l}}} \def\Ue{U^{\sscr{e\]l\}l}}}
\def\Sq{S_{\]\sss e\]l\}l}} \def\Omq{\Om_{\]\sss e\]l\}l}}
\def\che{\chi{\vp j}_{\]{\sss e\]l\}l}}}

\def\gge{\mathrel{\vp\gg\smash{\%{\gg}_{\raise1ex\mbox{\sss=\!=}}}}}
\def\lle{\mathrel{\vp\ll\smash{\%{\ll}_{\raise1ex\mbox{\sss=\!=}}}}}

\def\Sl{\S^\ell}  \def\Zpn{\Zp^n}
\def\Zln{\Zc_\ell^n} \def\Zlm{\Zc_\ell^{n-1}} \def\Zll{\Zc_{\ell-1}^n}

\def\HA{H_{\}\sss A}}
\def\Ii{I^{\sscr\o}}
\def\Omb{\Om^{\sscr\bullet}}
\def\qH{q^{-H}}
\def\ve[#1]{v^{\)[#1]}}

\def\Go{\mathchoice{{\%{\vru1.23ex>\smash{G}}^{\>\smash{\sss\o}\}}}}
 {{\%{\vru1.23ex>\smash{G}}^{\>\smash{\sss\o}\}}}}
 {{\%{\vru.92ex>\smash{G}}^{\>\smash{\sss\o}}}}{\@@PS}}
\def\Mo{\mathchoice{{\%{\vru1.19ex>\smash{M}\}}^{\,\smash{\sss\o}\!}}}
 {{\%{\vru1.19ex>\smash{M}\}}^{\,\smash{\sss\o}\!}}}
 {{\%{\vru.893ex>\smash{M}\}}^{\,\smash{\sss\o}\!}}}{\@@PS}}
\def\UUa{\mathchoice{{\%{\vru1.2ex>\smash{\UU}}^{\)\smash{\sss*}\]}}}
 {{\%{\vru1.2ex>\smash{\UU}}^{\)\smash{\sss*}\]}}}
 {{\%{\vru.9ex>\smash{\UU}}^{\)\smash{\sss*}\]}}}{\@@PS}}
\def\UUo{\mathchoice{{\%{\vru1.2ex>\smash{\UU}}^{\)\smash{\sss\o}\]}}}
 {{\%{\vru1.2ex>\smash{\UU}}^{\)\smash{\sss\o}\]}}}
 {{\%{\vru.9ex>\smash{\UU}}^{\)\smash{\sss\o}\]}}}{\@@PS}}
\def\Wo{\mathchoice{{\%{\vru1.2ex>\smash{W}\}}^{\>\smash{\sss\o}\}}}}
 {{\%{\vru1.2ex>\smash{W}\}}^{\>\smash{\sss\o}\}}}}
 {{\%{\vru.9ex>\smash{W}\}}^{\>\smash{\sss\o}\}}}}{\@@PS}}

\def\DWwi{\det\bigl[I(W_\lg,w_\mg)\bigr]}
\def\DWwp{\det\bigl[I(W_\lg,w'_\mg)\bigr]}
\def\Det#1{\mathop{\vp e\smash\det}\limits^{\,\sss\"{#1}\!}}
\def\Dit{\mathop{\italic{Det}\)}}

\def\Zps{Z^\pss} \def\Zms{Z^\mss}
\def\Resp{\Res\nolimits^{\!+}} \def\Resm{\Res\nolimits^{\!-}}
\def\tho{\tht'(1)} \def\Fabc{F(\tell{\)};\alb a,b,c)}

\def\6{{\>\overline{\}\8\!}\,}}
\def\l@inf{\lower.21ex\mbox{\ssize\8}} \def\9{_{\kern-.02em\l@inf}}
\def\l@infb{\lower.21ex\mbox{\ssize\6}} \def\0{_{\kern-.02em\l@infb}}
\def\7#1{^{[#1]}}

\def\smm{\sum_{m=1}} \def\tsmm{\!\tsum_{m=1}}
\def\sun{\smm^n} \def\tsun{\tsmm^n\!}
\def\smk{\sum_{k=0}} 
\def\susi{\sum_{\>\si\in\Sl\}}}
\def\smiZ{\sum_{\,\mg\in\Zln\!}^{\vp.}}
\def\sliZ{\sum_{\,\lg\in\Zln\!}^{\vp.}}
\def\sleZ{\sum_{\,\lg\in\Rph{\Zln}{\Zll}\!}^{\vp.}}

\def\prm{\prod_{m=1}} \def\tprm{\tprod_{m=1}}
\def\pron{\prm^n} \def\tpron{\tprm^n}
\def\prmn{\prm^{n-1}} 
\def\plmn{\prod_{1\le l<m\le n}}
 
\def\prlm{\prod_{1\le l<m}} \def\prllm{\prod_{1\le l\le m}}
\def\tprlm{\!\!\tprod_{1\le l<m}\!} \def\tprllm{\!\!\tprod_{1\le l\le m}\!}
\def\prml{\prod_{m<l\le n}} \def\prmll{\prod_{m\le l\le n}}
\def\tprml{\!\!\tprod_{m<l\le n}\!} \def\tprmll{\!\!\tprod_{m\le l\le n}\!}

\def\pros{\prod_{s=0}^{\ell-1}} \def\prosl{\prod_{s=1}^\ell}
\def\prsl{\prod_{s=1-\ell}^{\ell-1}}
\def\pral{\prod_{a=1}^\ell} \def\tpral{\tprod_{a=1}^\ell}
\def\prkl{\prod_{k=1}^\ell}
\def\prab{\prod_{1\le a<b\le\ell}\!\!} \def\prjk{\prod_{1\le j<k\le\ell}\!\!}
\def\pjkl{\prod_{j=1}^\ell\,\prod_{\tsize{k=1\atop k\ne j}}^\ell}
\def\prmlsi{\prod_{\si_m\}<\si_l\le n}} \def\prlmsi{\prod_{1\le\si_l<\si_m}}
\def\prabsi{\prod_{\tsize{1\le a<b\le\ell\atop\si_a>\si_b}}\!\!}

 \def\kon{k=0\lc n} \def\lcn{l=1\lc n} \def\mn{m=1\lc n}
\def\Lan{\La_1\lc\La_n} \def\phin{\phi_1\lc\phi_n}
\def\xn{x_1\lc x_n} \def\yn{y_1\lc y_n} \def\zn{z_1\lc z_n}
\def\qLan{q^{\La_1}\lc q^{\La_n}} \def\xin{\xi_1\lc\xi_n}
\def\mun{\mu_1\lc\mu_n} \def\zmn{z_1\lc p\)z_m\lc z_n}
\def\zmmn{z_1\lc z_{m+1},z_m\lc z_n}
\def\Gmn{\Gms_1\lc\Gms_n}
\def\Dtel{Dt_1\lsym\wedge Dt_\ell}
 \def\kzy{\ka\1\zy}
\def\yz@{\!\pron\!\xi_m}
\def\yz{\mathchoice{{\tsize\yz@}}{\yz@}{\yz@}{\yz@}}
\def\zy@{\!\pron\!\xi_m\1} 
\def\zy{\mathchoice{{\tsize\zy@}}{\zy@}{\zy@}{\zy@}}
\def\kone{{\eta^{1-\ell}\yz}} \def\koNe{{\eta^{-\ell}\yz}}

 \def\ktwp{{p\1\eta^{\ell-1}\zy\!}}
\def\kOne{\eta^{1-\ell}\!\tpron\!\eta^{\)\La_m}}
\def\kTwp{p\1\eta^{\ell-1}\!\tpron\!\eta^{-\La_m}}
\def\Kone{q^{\,2\!\!\sun\!\!\La_m-2\ell\)+\)2}{\!}}
\def\Ktwp{p\1 q^{-2\!\!\sun\!\!\La_m+\)2\ell-2}{\!}}
\def\aell{a=1\lc\ell} \def\kell{k=1\lc\ell} \def\koll{k=0\lc\ell}
\def\rell{r=1\lc\ell} \def\rtll{r=2\lc\ell} \def\roll{r=0\lc\ell-1}
\def\rll{r=1-\ell\lc\ell-1} \def\rpll{r=\ell-1\lc 2\ell-2}
\def\rolt{r=0\lc 2\ell-2}
\def\tell{t_1\lc t_\ell} \def\tpell{t_1\lc p\)t_a\lc t_\ell}
\def\tall{t_1\lc t_{a+1},t_a\lc t_\ell}
\def\telm{t_1\lc t_{\ell-1}} \def\twll{t_2\lc t_\ell}
\def\tsill{t_{\si_1}\lc t_{\si_\ell}} \def\tsik{t_{\si_1}\ldots t_{\si_k}}
\def\uell{u_1\lc u_\ell}
\def\uxt{u\mathbin{\smash\tright\]}}
\def\uyt{u\mathbin{\]\smash\trileft}}
\def\xt{x\mathbin{\smash\tright\]}}
\def\yt{y\mathbin{\]\smash\trileft}}
\def\taum{\tau\cdot(m,m+1)} \def\taumm{\tau'\cdot(m,m+1)}

\def\VLM{{V^{\!\La}V^{\}\Mu}}} \def\VML{{V^{\}\Mu}V^{\!\La}}}
\def\VLN{{V^{\!\La}V^{\}\Nu}}} \def\VMN{{V^{\}\Mu}V^{\}\Nu}}}
\def\Rq{R^{\vp q\sscr{ell}}} \def\Dlq{\Dl^{\!\sscr{ell}}}
\def\PV{P_{V_1V_2}\:} \def\RV{R_{V_1V_2}}
\def\RVx{\RV\:(x)} \def\RVo{\RV\:(\8)}
\def\VV{V_1\ox V_2} \def\EVV{\End(\VV)}
\def\Vn{V_1\lc V_n}  \def\Vqn{V^e_1\lc V^e_{n{\vp1}}}
\def\Vqnz{V^e_1(z_1)\lc V^e_{n{\vp1}}(z_n)}
\def\Vox{V_1\lox V_n} \def\EV{\End(\Vox)} \def\Vax{(V_1\lox V_n)^*}
\def\Vqx{V^e_1\lox V^e_{n{\vp1}}}

\def\Voxt{V_{\tau_1}\lox V_{\tau_n}} \def\Voxtt{V_{\tau'_1}\lox V_{\tau'_n}}
\def\Vaxt{(V_{\tau_1}\lox V_{\tau_n})^*}
\def\Vqxt{V^e_{\tau_1}\lox V^e_{\tau_n}}
\def\Voxtm{V_{\tau_1}\lox V_{\tau_{m+1}}\ox V_{\tau_m}\lox V_{\tau_n}}
\def\Vaxtm{(V_{\tau_1}\lox V_{\tau_{m+1}}\ox V_{\tau_m}\lox V_{\tau_n})^*}
\def\Vqxtm{V^e_{\tau_1}\lox V^e_{\tau_{m+1}}\ox V^e_{\tau_m}\lox V^e_{\tau_n}}
\def\Voxzt{V_{\tau_1}(z_{\tau_1})\lox V_{\tau_n}(z_{\tau_n})}

\def\Vqxzt{V^e_{\tau_1}(z_{\tau_1})\lox V^e_{\tau_n}(z_{\tau_n})}
\def\Vqxztt{V^e_{\tau'_1}(z_{\tau'_1})\lox V^e_{\tau'_n}(z_{\tau'_n})}
\def\elli{_\ell\:}
\def\Vl{(\Vox)\elli} \def\Vql{(\Vqx)\elli} \def\Val{\Vax_\ell}
\def\Vlt{(\Voxt)\elli} \def\Vqlt{(\Vqxt)\elli} \def\Valt{\Vaxt_\ell}
\def\Vqlto{\ve[0]\ox(V^e_{\tau_2}\lox V^e_{\tau_n})\elli}
\def\Vqltn{(V^e_{\tau_1}\lox V^e_{\tau_{n-1}})\elli\ox\ve[0]}
\def\Vaxtt{(V_{\tau'_1}\lox V_{\tau'_n})^*} \def\Valtt{\Vaxtt_\ell}
\def\Vqxtt{V^e_{\tau'_1}\lox V^e_{\tau'_n}} \def\Vqltt{(\Vqxtt)\elli}
\def\Vqltto{\ve[0]\ox(V^e_{\tau'_2}\lox V^e_{\tau'_n})\elli}
\def\Vqlttn{(V^e_{\tau'_1}\lox V^e_{\tau'_{n-1}})\elli\ox\ve[0]}
 \def\FVal{F\Vax_{\ell-1}}
\def\FzVal{\Fz\Vax_{\ell-1}}
\def\FValt{F\Vaxt_{\ell-1}} \def\FzValt{\Fz\Vaxt_{\ell-1}}
\def\vn{v_1\lc v_n} 
\def\vox{v_1\lox v_n} \def\vqx{\ve[0]\lox\ve[0]}
\def\voxt{v_{\tau_1}\lox v_{\tau_n}}
\def\Fv{F^{\lg_1}v_1\lox F^{\lg_n}v_n} \def\Fvv{F^{\mg_1}v_1\lox F^{\mg_n}v_n}
\def\Fva{(\Fv)^*}
\def\Fvq{\ve[\lg_1]\lox\ve[\lg_n]}
\def\Fvt{F^{\lg_{\tau_1}}v_{\tau_1}\lox F^{\lg_{\tau_n}}v_{\tau_n}}
\def\Fvvt{F^{\mg_{\tau_1}}v_{\tau_1}\lox F^{\mg_{\tau_n}}v_{\tau_n}}
\def\Fvat{(\Fvt)^*}
\def\Fvqt{\ve[\lg_{\tau_1}]\lox\ve[\lg_{\tau_n}]}
\def\Fvqtt{\ve[\lg_{\tau'_1}]\lox\ve[\lg_{\tau'_n}]}

\def\Vls{(\Vox)_\ell\sing} \def\Vlsz{(\Vox)_{\ell,\)z}\sing}
\def\Vlst{(\Voxt)_\ell\sing} \def\Vlszt{(\Voxt)_{\ell,\)z}\sing}
 
\def\VFV{\Val\big/\FVal} \def\VFzV{\Val\big/\FzVal}
\def\VFVt{\Valt\big/\FValt} \def\VFzVt{\Valt\big/\FzValt}

\def\sing{^{\sscr{\italic{sing}}}}
\def\fun#1{{\vp{#1}}^{\sss\italic{Fun}\!}#1}

\def\Tres#1#2#3{\left\lb\matrix\format\l\\#1\\#2\\#3\endmatrix\right\rb}

\def\dt{\,d^\ell t} \def\dtt{(dt/t)^\ell}
\def\ts{t^{\sss\star}} 

 \def\gsl{\frak{sl}_2} \def\Eqg{E_{\rho,\gm}(\gsl)}
\def\Uu{U_q(\gsl)} \def\Uuh{U_q(\Tilde{\gsl})}
\def\Ugh{U_q(\Tilde{\gl})} \def\Ugg{U_q'(\Tilde{\gl})}

\def\azo/{\as/ zone} \def\asol/{\as/ \sol/}
\def\conn/{connection coefficient} \def\sconn/{system of \conn/s}
\def\loc/{local system} \def\dloc/{discrete local system}
\def\prim/{primitive factor} \def\hform/{\hgeom/ form}
\def\hpair/{\hgeom/ pairing} \def\hgf/{\hgeom/ space}
\def\hmap/{\hgeom/ map} \def\hsol/{\hgeom/ \sol/}
\def\thgf/{\tri/ \hgf/} \def\ehgf/{elliptic \hgf/}
\def\hcg/{\hgeom/ cohomology group} \def\hhg/{\hgeom/ homology group}
\def\trib/{trivial bundle} \def\tvb/{the trivial vector bundle}
\def\wtd/{\wt/ decomposition} \def\diag/{diagonalizable}
\def\cosub/{coboundary subspace}\def\bosub/{boundary subspace}
\def\qbeta/{\^{$q$-}beta} \def\qSelberg/{\^{$q$-}Selberg}
\def\dyne/{dynamical elliptic} \def\ecurve/{elliptic curve}
\def\Spair/{Shapo\-va\-lov pairing}

\def\GM/{Gauss-Manin connection} \def\MB/{Mellin-Barnes}

\def\stype/{\^{$\gsl$}\h@ph type}
\def\Umod/{\^{$\Uu$}\h@ph module} \def\Uhmod/{\^{$\Ugg$}\h@ph module}
\def\Emod/{\^{$\Eqg$}\h@ph module} \def\hmod/{\^{$\hg$}\h@ph module}
\def\dhmod/{\diag/ \^{$\hg$}\h@ph module}
\def\eVmod/{evaluation \Vmod/} \def\epoint/{evaluation point}

\def\tenco/{tensor coordinates} \def\traf/{transition \fn/} 
\def\qlo/{quantum loop algebra} \def\phf/{phase \fn/} \def\wtf/{\wt/ \fn/}
\def\eqg/{elliptic \qg/} \def\oalg/{operator algebra}
\def\twf/{\tri/ \wtf/} \def\ewf/{elliptic \wtf/} \def\Atype/{{\sl A\/}-type}
\def\adjf/{adjusting factor} \def\adjm/{adjusting map}

\def\Vval/{\^{$V\!\!\;$-}valued} \def\SL#1{\^{$\,\Sl$}#1}
\def\p#1{\^{$p\>$#1}}

\let\goodbm\relax  \let\mmgood\relax 
\let\goodast\relax \let\agood\relax \let\aagood\relax 

\ifMag\let\goodbm\goodbreak  \let\mmgood\vvgood
 \let\goodbreak\relax  \let\vvgood\relax \fi

\setparindent

\let\fixedpage\newpage

\whattime\readldf\writeldf

\def\lefthead#1{} \def\righthead#1{}
\def\leftsecthead{\@fn@xt*\l@fts@cth@\l@fts@cth}
\def\l@fts@cth#1{} \def\l@fts@cth@*#1{}

\csname trig.def\endcsname

\setbox\lefthbox\hbox{\eightpoint\sl \VT/ \ and \ \Varch/}
\setbox\righthbox\hbox{\=\eightpoint\sl
$q$-Hypergeometric Functions, Quantum Affine Algebras and
Elliptic Quantum Groups}
\leftheadtext{\copy\lefthbox}
\rightheadtext{\copy\righthbox}

\newif\ifAst

\def\Chap/{Section} \def\chap/{section} \def\appx/{section}

\ifx\Asterisque\relax\Asttrue\fi

\expandafter\ifx\csname asterisq.def\endcsname\relax \else\endinput\fi
\expandafter\edef\csname asterisq.def\endcsname{%
 \catcode`\noexpand\@=\the\catcode`\@\space}
\catcode`\@=11

\newbox\l@fthb@x \newbox\r@ghthb@x
\newif\ifn@s@ct \n@s@cttrue
\newcount\m@rkn@

\ifAst

\font\BMmb=cmmib10 scaled 1440
\font\BSlb=cmbxsl10 scaled 1440
\def\BBfcorr{\let\sl\BSlb \BBf \bls 18\p@
 \def\Rms/{{\BMmb R}-matrices} \def\qSelberg/{{\BMmb q}-Selberg}}

\def\headfont{\eightpoint\rm}
\def\foliofont{\eightrm}

\headline={\t@stmarks{\global\setbox\lefthbox\copy\l@fthb@x%
 \global\setbox\righthbox\copy\r@ghthb@x}\ifn@s@ct\foliofont\hfil\ifodd\pageno
 \copy\righthbox\else\copy\lefthbox\fi\hfil\else\hfil\fi\n@s@cttrue}

\footline={\foliofont\hfil\ifp@gen@\folio\hfil\fi\global\p@gen@true}

\def\s@@ct#1#2{\ifodd\pageno\else\global\nopagenumber\strut\par\newpage\fi
 \global\n@s@ctfalse\vp1\vsk->\vtop{\BBfcorr\s@twd@\parindent{#1}%
 \ifdim\parindent>\z@\adv\parindent.5em\fi\hang\textindent{#1}#2\strut}%
 \nointerlineskip\vtop{\strut}\vsk->\vsk>}

\def\t@stmarks#1{\edef\t@stb{\botmark}\edef\t@stt{\topmark}%
 \ifx\t@stb\t@stt\else#1\fi}

\def\s@@@bs@ct#1#2{\p@sk@p{-200}{\bls}%
 \vtop{\twelvepoint\it\s@twd@\parindent{#1}%
 \ifdim\parindent>\z@\adv\parindent.5em\fi\hang\textindent{#1}#2\strut}%
 \nointerlineskip\nobreak\vtop{\strut}\nobreak\vskip-.4\bls\nobreak}

\def\s@@ct@n#1#2{\newpage\s@@ct{#1}{#2}{\def\@@@n{#1}\headfont\pr@@t@
 \edef\s@c@nt@{\uppercase\@ft@{\ifx\@@@n\empty\else\@@@n\enspace\fi#2}}
 \global\setbox\l@fthb@x\hbox{\s@c@nt@}\edef\s@c@nt@{\uppercase{#2}}
 \global\setbox\r@ghthb@x\hbox{\s@c@nt@}\gad\m@rkn@\mark{\the\m@rkn@}}}

\def\s@@bs@ct#1#2{\s@@@bs@ct{#1}{#2}{\def\@@@n{#1}\headfont\pr@@t@
 \edef\s@c@nt@{\uppercase\@ft@{#2}}
 \global\setbox\r@ghthb@x\hbox{\s@c@nt@}\gad\m@rkn@\mark{\the\m@rkn@}}}

\def\lefthead#1{{\headfont\pr@@t@\edef\s@c@nt@{\uppercase\@ft@{#1}}%
 \global\setbox\l@fthb@x\hbox{\s@c@nt@}\gad\m@rkn@\mark{\the\m@rkn@}}}
\def\righthead#1{{\headfont\pr@@t@\edef\s@c@nt@{\uppercase\@ft@{#1}}%
 \global\setbox\r@ghthb@x\hbox{\s@c@nt@}\gad\m@rkn@\mark{\the\m@rkn@}}}

\def\leftsecthead{\@fn@xt*\l@fts@cth@\l@fts@cth}
\def\l@fts@cth#1{\lefthead{\sf@rm@\enspace#1}\righthead{#1}}
\def\l@fts@cth@*#1{\lefthead{#1}\righthead{#1}}

\let\s@ct\s@@ct \let\s@ct@n\s@@ct@n \let\s@bs@ct\s@@bs@ct

\def\Appendix{\Sno=64 \t@@n\@ne \wr@@c{\string\Appencd} \let\s@bs@ct\s@@@bs@ct
 \def\sf@rm{\char\the\Sno} \def\sf@rm@@{\Ap@x}
 \def\s@ct@n##1##2{\newpage\s@@ct{##1}{##2}%
 {\headfont\global\setbox\l@fthb@x\hbox{APPENDIX\enspace\sf@rm}
 \pr@@t@\edef\s@c@nt@{\noexpand\uppercase{##2}}
 \global\setbox\r@ghthb@x\hbox{\s@c@nt@}\gad\m@rkn@\mark{\the\m@rkn@}}}}

\def\endAppendix{\let\s@ct\s@@ct \let\s@ct@n\s@@ct@n \let\s@bs@ct\s@@bs@ct}

\def\myR@fs{\let\s@ct\s@@ct@n\@fn@xt[\m@R@f@\m@R@fs}
\def\R@fcd#1#2{\sk@@p{.8\bls}\entcd{}{\hglue-\l@ftcd\twelvebf #1}{#2}}

\let\c@@nt@@\c@nt@
\long\def\c@nt@#1#2#3#4{\let\s@ct\s@@ct@n \let\s@bs@ct\s@@@bs@ct
 \c@@nt@@{#1}{#2}{#3}{#4}\global\let\s@ct\s@@ct \global\let\s@bs@ct\s@@bs@ct
 \vsk-.56>\vsk0>}

\def\Entcd#1#2#3{\def\l@@b{\@fwd@@{#1}{\c@l@b{#1}}{}}\vsk.36>%
 \entcd{\l@@b}{#2}{#3}}

\def\Eqg{E_{\rho,\gm}(\frak{sl}_2)}
\def\Uu{U_q(\frak{sl}_2)} \def\Uuh{U_q(\Tilde{\frak{sl}_2})}
\def\Ugh{U_q(\Tilde{\frak{gl}_2})} \def\Ugg{U_q'(\Tilde{\frak{gl}_2})}
\def\qSelberg/{\^{$\q@@$-}Selberg} \let\q@@ q

\def\pr@@t@{\let\Uu@\Uu \def\Uu{\noexpand\Uu@}
 \let\Eqg@\Eqg \def\Eqg{\noexpand\Eqg@}
 \let\Ugg@\Ugg \def\Ugg{\noexpand\Ugg@}
 \let\gsl@\gsl \def\gsl{\noexpand\gsl@}
 \let\q@@@\q@@ \def\q@@@{\noexpand\q@@}}

\ifMag
\let\goodast\goodbm \let\agood\mgood \let\aagood\mmgood
\let\goodbm\relax  \let\mmgood\relax
\else
\let\goodast\goodbreak \let\agood\vgood \let\aagood\vvgood
\let\goodbreak\relax  \let\vvgood\relax
\fi

\setbox\lefthbox\hbox{}
\setbox\righthbox\hbox{}

\def\Chap/{Chapter} \def\chap/{chapter} \def\appx/{appendix}

\Magtrue

\fi

\csname asterisq.def\endcsname

\labeldef{S} {1} {1}

\labeldef{S} {2} {2}
\labeldef{F} {2\labelsep \labelspace 1}  {flat}
\labeldef{F} {2\labelsep \labelspace 2}  {period}
\labeldef{F} {2\labelsep \labelspace 3}  {conns}
\labeldef{F} {2\labelsep \labelspace 4}  {Phb}
\labeldef{F} {2\labelsep \labelspace 5}  {Phi}
\labeldef{F} {2\labelsep \labelspace 6}  {tzp}
\labeldef{F} {2\labelsep \labelspace 7}  {Phisym}
\labeldef{F} {2\labelsep \labelspace 8}  {list1}
\labeldef{F} {2\labelsep \labelspace 9}  {act}
\labeldef{F} {2\labelsep \labelspace 10} {thgf}
\labeldef{F} {2\labelsep \labelspace 11} {dis}
\labeldef{F} {2\labelsep \labelspace 12} {subbundle}
\labeldef{F} {2\labelsep \labelspace 13} {npZ}
\labeldef{F} {2\labelsep \labelspace 14} {Lass}
\labeldef{F} {2\labelsep \labelspace 15} {assum}
\labeldef{F} {2\labelsep \labelspace 16} {kanon}
\labeldef{F} {2\labelsep \labelspace 17} {HFo}
\labeldef{F} {2\labelsep \labelspace 18} {kaon}
\labeldef{F} {2\labelsep \labelspace 19} {Zln}
\labeldef{F} {2\labelsep \labelspace 20} {wlga}
\labeldef{F} {2\labelsep \labelspace 21} {wlgp}
\labeldef{F} {2\labelsep \labelspace 22} {wbasis}
\labeldef{F} {2\labelsep \labelspace 23} {wDw}
\labeldef{F} {2\labelsep \labelspace 24} {wDw'}
\labeldef{F} {2\labelsep \labelspace 25} {HFR}
\labeldef{F} {2\labelsep \labelspace 26} {HFR'}
\labeldef{F} {2\labelsep \labelspace 27} {wtau}
\labeldef{F} {2\labelsep \labelspace 28} {1-01}
\labeldef{F} {2\labelsep \labelspace 29} {ect}
\labeldef{F} {2\labelsep \labelspace 30} {Wlga}
\labeldef{F} {2\labelsep \labelspace 31} {Wlgp}
\labeldef{F} {2\labelsep \labelspace 32} {adjust}
\labeldef{F} {2\labelsep \labelspace 33} {Wbasis}
\labeldef{F} {2\labelsep \labelspace 34} {Rests}
\labeldef{F} {2\labelsep \labelspace 35} {xtyt}
\labeldef{F} {2\labelsep \labelspace 36} {ResQ}
\labeldef{F} {2\labelsep \labelspace 37} {ResQ'}
\labeldef{F} {2\labelsep \labelspace 38} {dimFQ}
\labeldef{F} {2\labelsep \labelspace 39} {dimFQ'}
\labeldef{F} {2\labelsep \labelspace 40} {Wtau}

\labeldef{S} {3} {2a}
\labeldef{F} {3\labelsep \labelspace 1}  {VVF}
\labeldef{F} {3\labelsep \labelspace 2}  {VVFz}
\labeldef{F} {3\labelsep \labelspace 3}  {Rdef}
\labeldef{F} {3\labelsep \labelspace 4}  {Rvv}
\labeldef{F} {3\labelsep \labelspace 5}  {Rmore}
\labeldef{F} {3\labelsep \labelspace 6}  {inv}
\labeldef{F} {3\labelsep \labelspace 7}  {Rspec}
\labeldef{F} {3\labelsep \labelspace 8}  {YBq}
\labeldef{F} {3\labelsep \labelspace 9}  {Ughr}
\labeldef{F} {3\labelsep \labelspace 10} {Uintw}
\labeldef{F} {3\labelsep \labelspace 11} {Rij}
\labeldef{F} {3\labelsep \labelspace 12} {Kmz}
\labeldef{F} {3\labelsep \labelspace 13} {FR}
\labeldef{F} {3\labelsep \labelspace 15} {1-02}
\labeldef{F} {3\labelsep \labelspace 16} {Eintw}
\labeldef{F} {3\labelsep \labelspace 18} {YBe}

\labeldef{S} {4} {2b}
\labeldef{F} {4\labelsep \labelspace 1}  {Btau}
\labeldef{F} {4\labelsep \labelspace 2}  {tcfs}
\labeldef{F} {4\labelsep \labelspace 3}  {T12}
\labeldef{F} {4\labelsep \labelspace 4}  {local}
\labeldef{F} {4\labelsep \labelspace 5}  {tcnon}
\labeldef{F} {4\labelsep \labelspace 6}  {tcon}
\labeldef{F} {4\labelsep \labelspace 7}  {qKZ-GM}
\labeldef{F} {4\labelsep \labelspace 9}  {clg}
\labeldef{F} {4\labelsep \labelspace 10} {ecfs}
\labeldef{F} {4\labelsep \labelspace 11} {RespC}
\labeldef{F} {4\labelsep \labelspace 12} {ResmC}
\labeldef{F} {4\labelsep \labelspace 13} {ecfso}
\labeldef{F} {4\labelsep \labelspace 14} {ecfsp}
\labeldef{F} {4\labelsep \labelspace 15} {L21}
\labeldef{F} {4\labelsep \labelspace 16} {localq}
\labeldef{F} {4\labelsep \labelspace 17} {chi}
\labeldef{F} {4\labelsep \labelspace 18} {FFo}
\labeldef{F} {4\labelsep \labelspace 19} {FFq}
\labeldef{F} {4\labelsep \labelspace 20} {Fq1Fqk}

\labeldef{S} {5} {2c}
\labeldef{F} {5\labelsep \labelspace 1}  {IWwt}
\labeldef{F} {5\labelsep \labelspace 2}  {Phis}
\labeldef{F} {5\labelsep \labelspace 3}  {IWw}
\labeldef{F} {5\labelsep \labelspace 4}  {list3}
\labeldef{F} {5\labelsep \labelspace 5}  {list2}
\labeldef{F} {5\labelsep \labelspace 6}  {Ianco}
\labeldef{F} {5\labelsep \labelspace 7}  {IWw=0}
\labeldef{F} {5\labelsep \labelspace 8}  {IWw=0'}
\labeldef{F} {5\labelsep \labelspace 9}  {mu<>0}
\labeldef{F} {5\labelsep \labelspace 10} {mu=0}
\labeldef{F} {5\labelsep \labelspace 11} {mu=02}
\labeldef{F} {5\labelsep \labelspace 12} {AR}
\labeldef{F} {5\labelsep \labelspace 13} {qbeta}
\labeldef{F} {5\labelsep \labelspace 14} {ARl}
\labeldef{F} {5\labelsep \labelspace 15} {sWz}
\labeldef{F} {5\labelsep \labelspace 16} {Vsing}
\labeldef{F} {5\labelsep \labelspace 17} {Vsing'}
\labeldef{F} {5\labelsep \labelspace 18} {PsiW}
\labeldef{F} {5\labelsep \labelspace 19} {qKZsol}
\labeldef{F} {5\labelsep \labelspace 20} {Shat}
\labeldef{F} {5\labelsep \labelspace 21} {Shwd}
\labeldef{F} {5\labelsep \labelspace 22} {qKZfun}
\labeldef{F} {5\labelsep \labelspace 23} {mono}
\labeldef{F} {5\labelsep \labelspace 24} {hmdef}
\labeldef{F} {5\labelsep \labelspace 25} {hmap}
\labeldef{F} {5\labelsep \labelspace 26} {hmapo}
\labeldef{F} {5\labelsep \labelspace 27} {hmps}
\labeldef{F} {5\labelsep \labelspace 28} {hmapp}
\labeldef{F} {5\labelsep \labelspace 29} {hmpsp}
\labeldef{F} {5\labelsep \labelspace 30} {Imm+1}
\labeldef{F} {5\labelsep \labelspace 31} {PsivY}
\labeldef{F} {5\labelsep \labelspace 32} {CCtau}
\labeldef{F} {5\labelsep \labelspace 33} {CCtt}
\labeldef{F} {5\labelsep \labelspace 34} {CCC}
\labeldef{F} {5\labelsep \labelspace 35} {monoR}

\labeldef{S} {6} {2d}
\labeldef{F} {6\labelsep \labelspace 1}  {azone}
\labeldef{F} {6\labelsep \labelspace 2}  {asol}
\labeldef{F} {6\labelsep \labelspace 3}  {Ytau}
\labeldef{F} {6\labelsep \labelspace 4}  {RpR}
\labeldef{F} {6\labelsep \labelspace 5}  {Asol}
\labeldef{F} {6\labelsep \labelspace 6}  {AsIWw}

\labeldef{S} {7} {4}
\labeldef{F} {7\labelsep \labelspace 1}  {Xxy}
\labeldef{F} {7\labelsep \labelspace 2}  {DlL}
\labeldef{F} {7\labelsep \labelspace 3}  {1233}
\labeldef{F} {7\labelsep \labelspace 4}  {3321}
\labeldef{F} {7\labelsep \labelspace 5}  {Xjlg}
\labeldef{F} {7\labelsep \labelspace 6}  {Xixy}
\labeldef{F} {7\labelsep \labelspace 7}  {IWwD}
\labeldef{F} {7\labelsep \labelspace 8}  {IWw!}
\labeldef{F} {7\labelsep \labelspace 9}  {Iz}
\labeldef{F} {7\labelsep \labelspace 10} {wWprim}
\labeldef{F} {7\labelsep \labelspace 11} {wAA}
\labeldef{F} {7\labelsep \labelspace 12} {AsIWwt}
\labeldef{F} {7\labelsep \labelspace 13} {Dit}
\labeldef{F} {7\labelsep \labelspace 14} {w''}

\labeldef{S} {\char 65} {A}
\labeldef{F} {\char 65\labelsep \labelspace 1}  {winFo}
\labeldef{F} {\char 65\labelsep \labelspace 2}  {n=1}
\labeldef{F} {\char 65\labelsep \labelspace 3}  {n=10}
\labeldef{F} {\char 65\labelsep \labelspace 4}  {n=11}
\labeldef{F} {\char 65\labelsep \labelspace 5}  {wDw0}
\labeldef{F} {\char 65\labelsep \labelspace 6}  {Qlg}
\labeldef{F} {\char 65\labelsep \labelspace 7}  {DetM}
\labeldef{F} {\char 65\labelsep \labelspace 8}  {F/R}
\labeldef{F} {\char 65\labelsep \labelspace 9}  {Plg}
\labeldef{F} {\char 65\labelsep \labelspace 10} {DetN}
\labeldef{F} {\char 65\labelsep \labelspace 11} {xtry}
\labeldef{F} {\char 65\labelsep \labelspace 12} {Pxy}
\labeldef{F} {\char 65\labelsep \labelspace 13} {detPxy}
\labeldef{F} {\char 65\labelsep \labelspace 14} {detQxy}
\labeldef{F} {\char 65\labelsep \labelspace 15} {DetNC}

\labeldef{S} {\char 66} {B}
\labeldef{F} {\char 66\labelsep \labelspace 2}  {Tho}
\labeldef{F} {\char 66\labelsep \labelspace 3}  {Glg}
\labeldef{F} {\char 66\labelsep \labelspace 4}  {winFq}
\labeldef{F} {\char 66\labelsep \labelspace 5}  {n=1q}
\labeldef{F} {\char 66\labelsep \labelspace 6}  {n=11q}
\labeldef{F} {\char 66\labelsep \labelspace 7}  {Mqdef}
\labeldef{F} {\char 66\labelsep \labelspace 8}  {DetMq}
\labeldef{F} {\char 66\labelsep \labelspace 9}  {Pilg}
\labeldef{F} {\char 66\labelsep \labelspace 10} {Jxy}
\labeldef{F} {\char 66\labelsep \labelspace 11} {detTxy}
\labeldef{F} {\char 66\labelsep \labelspace 12} {thtomn}

\labeldef{S} {\char 67} {S}
\labeldef{F} {\char 67\labelsep \labelspace 1}  {Omq}
\labeldef{F} {\char 67\labelsep \labelspace 2}  {SWWx}
\labeldef{F} {\char 67\labelsep \labelspace 3}  {SWWy}
\labeldef{F} {\char 67\labelsep \labelspace 4}  {SWW}
\labeldef{F} {\char 67\labelsep \labelspace 5}  {c=1/N}
\labeldef{F} {\char 67\labelsep \labelspace 6}  {ne1}
\labeldef{F} {\char 67\labelsep \labelspace 7}  {Omtri}
\labeldef{F} {\char 67\labelsep \labelspace 8}  {Swwy}
\labeldef{F} {\char 67\labelsep \labelspace 9}  {Sww}
\labeldef{F} {\char 67\labelsep \labelspace 10} {Fform}
\labeldef{F} {\char 67\labelsep \labelspace 11} {Resi}

\labeldef{S} {\char 68} {C}
\labeldef{F} {\char 68\labelsep \labelspace 1}  {Fabc}
\labeldef{F} {\char 68\labelsep \labelspace 2}  {c-eq}
\labeldef{F} {\char 68\labelsep \labelspace 3}  {X1Xl}
\labeldef{F} {\char 68\labelsep \labelspace 4}  {Xtia}
\labeldef{F} {\char 68\labelsep \labelspace 5}  {Xtib}
\labeldef{F} {\char 68\labelsep \labelspace 6}  {Sdeq}
\labeldef{F} {\char 68\labelsep \labelspace 7}  {xident}
\labeldef{F} {\char 68\labelsep \labelspace 8}  {xfact}
\labeldef{F} {\char 68\labelsep \labelspace 9}  {qbetaR}
\labeldef{F} {\char 68\labelsep \labelspace 10} {qbetaS}
\labeldef{F} {\char 68\labelsep \labelspace 11} {FFti}

\labeldef{S} {\char 69} {D}
\labeldef{F} {\char 69\labelsep \labelspace 1}  {IW0}
\labeldef{F} {\char 69\labelsep \labelspace 2}  {IWW0}
\labeldef{F} {\char 69\labelsep \labelspace 3}  {Wal}
\labeldef{F} {\char 69\labelsep \labelspace 4}  {abcx}
\labeldef{F} {\char 69\labelsep \labelspace 5}  {ARl'}
\labeldef{F} {\char 69\labelsep \labelspace 6}  {ARl''}
\labeldef{F} {\char 69\labelsep \labelspace 7}  {points}
\labeldef{F} {\char 69\labelsep \labelspace 8}  {Ascj}
\labeldef{F} {\char 69\labelsep \labelspace 9}  {resfv}
\labeldef{F} {\char 69\labelsep \labelspace 10} {Ascj'}
\labeldef{F} {\char 69\labelsep \labelspace 11} {ysum2}

\labeldef{S} {\char 70} {E}
\labeldef{F} {\char 70\labelsep \labelspace 1}  {xsum}
\labeldef{F} {\char 70\labelsep \labelspace 2}  {ysum}
\labeldef{F} {\char 70\labelsep \labelspace 3}  {PWPW}
\labeldef{F} {\char 70\labelsep \labelspace 4}  {xsumf}

\labeldef{S} {\char 71} {X}
\labeldef{F} {\char 71\labelsep \labelspace 1}  {combi}

\def\ContA{
\Entcd{1}{Introduction}{1}
\Entcd{2}{Discrete flat connections and local systems}{5}
\subcd{Discrete flat connections}{5}
\subcd{Discrete Gauss-Manin connection}{7}
\subcd{Connection coefficients of local systems}{8}
\subcd{The functional space of a trigonometric\/ \stype / local system}{9}
\subcd{Bases in the trigonometric hyper\-geometric space of a fiber}{12}
\subcd{The elliptic hyper\-geometric space}{15}
\Entcd{3}{\Rms / and the \qKZc /}{21}
\subcd{Highest weight \Umod /s}{21}
\subcd{The trigonometric \Rm /}{22}
\subcd{The quantum loop algebra $\Ugg $}{23}
\subcd{The trigonometric \qKZc / associated with $\gsl $}{25}
\subcd{Modules over the elliptic quantum group $\Eqg $ and the elliptic \nl
\nlc \Rms /}{27}
\Entcd{4}{Tensor coordinates on the hyper\-geometric spaces}{31}
\subcd{Tensor coordinates on the trigonometric hyper\-geometric spaces
of fibers}{31}
\subcd{Tensor coordinates on the elliptic hyper\-geometric spaces
of fibers}{34}
\subcd{Tensor products of the hyper\-geometric spaces}{37}
\Entcd{5}{The hyper\-geometric pairing and \nl the hyper\-geometric solutions
of \nlc the \qKZe /}{39}
\subcd{The hyper\-geometric integral}{39}
\subcd{Determinant formulae for the hyper\-geometric pairing}{41}
\subcd{The hyper\-geometric solutions of the \qKZe /}{45}
\subcd{The hyper\-geometric map}{47}
\Entcd{6}{Asymptotic solutions of the \qKZe /}{51}
\Entcd{7}{Proofs}{57}
\subcd{Proof of Lemmas~\[wbasis], \[Wbasis]}{57}
\subcd{Proof of Lemmas~\[wDw], \[wDw']}{57}
\subcd{Proof of Lemma~\[T12]}{57}
\subcd{Proof of Lemma~\[L21]}{59}
\subcd{Proof of Lemmas~\[FFo], \[FFq]}{60}
\subcd{Proof of Lemma~\[IWw=0]}{61}
\subcd{Proof of Lemma~\[IWw=0']}{61}
\subcd{Proof of Theorem~\[sWz]}{61}
\subcd{Proof of Lemmas~\[Vsing], \[Vsing']}{62}
\subcd{Proof of Theorem~\[AsIWw]}{63}
\subcd{Proof of Theorem~\[asol]}{65}
\subcd{Proof of Theorem~\[mu<>0]}{65}
\subcd{Proof of Theorem~\[mu=0]}{66}
\subcd{Proof of Theorem~\[mu=02]}{69}
\subcd{Proof of Theorem~\[hmapo]}{69}
\subcd{Proof of Theorem~\[hmapp]}{70}
\subcd{Proof of Theorem~\[PsivY]}{70}
\Appencd
\Entcd{\char 65}{Basic facts about the trigonometric \nl hyper\-geometric
space}{71}
\subcd{Proof of Lemma~\[wlgp]}{72}
\Entcd{\char 66}{Basic facts about the elliptic hyper\-geometric space}{77}
\subcd{Proof of Lemma~\[Wlgp]}{79}
\subcd{Proof of Lemma~\[dimFQ]}{81}
\subcd{Proof of Lemma~\[dimFQ']}{82}
\subcd{Proof of Lemma~\[ResQ]}{82}
\subcd{Proof of Lemma~\[ResQ']}{82}
\Entcd{\char 67}{The Shapo\-va\-lov pairings of the hyper\-geometric spaces of
fibers}{83}
\Entcd{\char 68}{The \qSelberg / integral}{91}
\subcd{Proof of formula \(qbeta)}{91}
\Entcd{\char 69}{The multidimensional Askey-Roy formula and \nl Askey's
conjecture}{97}
\subcd{Proof of formula \(ARl)}{97}
\Entcd{\char 70}{The Jackson integrals via the hyper\-geometric \nl
integrals}{103}
\Entcd{\char 71}{One useful identity}{107}
\Refcd{References}{109}}

\def\ContMP{
\Entcd{1}{Introduction}{1}
\Entcd{2}{Discrete flat connections and local systems}{5}
\subcd{Discrete flat connections}{5}
\subcd{Discrete Gauss-Manin connection}{7}
\subcd{Connection coefficients of local systems}{7}
\subcd{The functional space of a trigonometric\/ \stype / local system}{9}
\subcd{Bases in the trigonometric hyper\-geometric space of a fiber}{11}
\subcd{The elliptic hyper\-geometric space}{14}
\Entcd{3}{\Rms / and the \qKZc /}{19}
\subcd{Highest weight \Umod /s}{19}
\subcd{The trigonometric \Rm /}{20}
\subcd{The quantum loop algebra $\Ugg $}{21}
\subcd{The trigonometric \qKZc / associated with $\gsl $}{23}
\subcd{Modules over the elliptic quantum group $\Eqg $ and the elliptic \nlc
\Rms /}{24}
\Entcd{4}{Tensor coordinates on the hyper\-geometric spaces}{27}
\subcd{Tensor coordinates on the trigonometric hyper\-geometric spaces of
fibers}{27}
\subcd{Tensor coordinates on the elliptic hyper\-geometric spaces of
fibers}{29}
\subcd{Tensor products of the hyper\-geometric spaces}{32}
\Entcd{5}{The hyper\-geometric pairing and the hyper\-geometric solutions of
\nlc the \qKZe /}{34}
\subcd{The hyper\-geometric integral}{34}
\subcd{Determinant formulae for the hyper\-geometric pairing}{36}
\subcd{The hyper\-geometric solutions of the \qKZe /}{40}
\subcd{The hyper\-geometric map}{41}
\Entcd{6}{Asymptotic solutions of the \qKZe /}{44}
\Entcd{7}{Proofs}{49}
\subcd{Proof of Lemmas~\[wbasis], \[Wbasis]}{49}
\subcd{Proof of Lemmas~\[wDw], \[wDw']}{50}
\subcd{Proof of Lemma~\[T12]}{50}
\subcd{Proof of Lemma~\[L21]}{51}
\subcd{Proof of Lemmas~\[FFo], \[FFq]}{53}
\subcd{Proof of Lemma~\[IWw=0]}{53}
\subcd{Proof of Lemma~\[IWw=0']}{54}
\subcd{Proof of Theorem~\[sWz]}{54}
\subcd{Proof of Lemmas~\[Vsing], \[Vsing']}{55}
\subcd{Proof of Theorem~\[AsIWw]}{55}
\subcd{Proof of Theorem~\[asol]}{57}
\subcd{Proof of Theorem~\[mu<>0]}{57}
\subcd{Proof of Theorem~\[mu=0]}{58}
\subcd{Proof of Theorem~\[mu=02]}{61}
\subcd{Proof of Theorem~\[hmapo]}{61}
\subcd{Proof of Theorem~\[hmapp]}{61}
\subcd{Proof of Theorem~\[PsivY]}{62}
\Appencd
\Entcd{\char 65}{Basic facts about the trigonometric hyper\-geometric
space}{62}
\subcd{Proof of Lemma~\[wlgp]}{63}
\Entcd{\char 66}{Basic facts about the elliptic hyper\-geometric space}{66}
\subcd{Proof of Lemma~\[Wlgp]}{68}
\subcd{Proof of Lemma~\[dimFQ]}{70}
\subcd{Proof of Lemma~\[dimFQ']}{71}
\subcd{Proof of Lemma~\[ResQ]}{71}
\subcd{Proof of Lemma~\[ResQ']}{71}
\Entcd{\char 67}{The Shapo\-va\-lov pairings of the hyper\-geometric spaces of
fibers}{72}
\Entcd{\char 68}{The \qSelberg / integral}{77}
\subcd{Proof of formula \(qbeta)}{77}
\Entcd{\char 69}{The multidimensional Askey-Roy formula and Askey's \nl
conjecture}{82}
\subcd{Proof of formula \(ARl)}{82}
\Entcd{\char 70}{The Jackson integrals via the hyper\-geometric integrals}{86}
\Entcd{\char 71}{One useful identity}{88}
\Refcd{References}{88}}

\def\ContSP{
\Entcd{1}{Introduction}{1}
\Entcd{2}{Discrete flat connections and local systems}{3}
\subcd{Discrete flat connections}{3}
\subcd{Discrete Gauss-Manin connection}{5}
\subcd{Connection coefficients of local systems}{6}
\subcd{The functional space of a trigonometric\/ \stype / local system}{7}
\subcd{Bases in the trigonometric hyper\-geometric space of a fiber}{9}
\subcd{The elliptic hyper\-geometric space}{11}
\Entcd{3}{\Rms / and the \qKZc /}{15}
\subcd{Highest weight \Umod /s}{15}
\subcd{The trigonometric \Rm /}{16}
\subcd{The quantum loop algebra $\Ugg $}{16}
\subcd{The trigonometric \qKZc / associated with $\gsl $}{18}
\subcd{Modules over the elliptic quantum group $\Eqg $ and the elliptic
\Rms /}{19}
\Entcd{4}{Tensor coordinates on the hyper\-geometric spaces}{21}
\subcd{Tensor coordinates on the trigonometric hyper\-geometric spaces of
fibers}{21}
\subcd{Tensor coordinates on the elliptic hyper\-geometric spaces of
fibers}{23}
\subcd{Tensor products of the hyper\-geometric spaces}{25}
\Entcd{5}{The hyper\-geometric pairing and the hyper\-geometric solutions of
the \qKZe /}{27}
\subcd{The hyper\-geometric integral}{27}
\subcd{Determinant formulae for the hyper\-geometric pairing}{28}
\subcd{The hyper\-geometric solutions of the \qKZe /}{31}
\subcd{The hyper\-geometric map}{32}
\Entcd{6}{Asymptotic solutions of the \qKZe /}{35}
\Entcd{7}{Proofs}{38}
\subcd{Proof of Lemmas~\[wbasis], \[Wbasis]}{38}
\subcd{Proof of Lemmas~\[wDw], \[wDw']}{38}
\subcd{Proof of Lemma~\[T12]}{38}
\subcd{Proof of Lemma~\[L21]}{40}
\subcd{Proof of Lemmas~\[FFo], \[FFq]}{41}
\subcd{Proof of Lemma~\[IWw=0]}{41}
\subcd{Proof of Lemma~\[IWw=0']}{42}
\subcd{Proof of Theorem~\[sWz]}{42}
\subcd{Proof of Lemmas~\[Vsing], \[Vsing']}{42}
\subcd{Proof of Theorem~\[AsIWw]}{43}
\subcd{Proof of Theorem~\[asol]}{45}
\subcd{Proof of Theorem~\[mu<>0]}{45}
\subcd{Proof of Theorem~\[mu=0]}{45}
\subcd{Proof of Theorem~\[mu=02]}{47}
\subcd{Proof of Theorem~\[hmapo]}{47}
\subcd{Proof of Theorem~\[hmapp]}{47}
\subcd{Proof of Theorem~\[PsivY]}{48}
\Appencd
\Entcd{\char 65}{Basic facts about the trigonometric hyper\-geometric
space}{48}
\subcd{Proof of Lemma~\[wlgp]}{49}
\Entcd{\char 66}{Basic facts about the elliptic hyper\-geometric space}{52}
\subcd{Proof of Lemma~\[Wlgp]}{53}
\subcd{Proof of Lemma~\[dimFQ]}{55}
\subcd{Proof of Lemma~\[dimFQ']}{55}
\subcd{Proof of Lemma~\[ResQ]}{55}
\subcd{Proof of Lemma~\[ResQ']}{56}
\Entcd{\char 67}{The Shapo\-va\-lov pairings of the hyper\-geometric spaces of
fibers}{56}
\Entcd{\char 68}{The \qSelberg / integral}{60}
\subcd{Proof of formula \(qbeta)}{60}
\Entcd{\char 69}{The multidimensional Askey-Roy formula and Askey's
conjecture}{64}
\subcd{Proof of formula \(ARl)}{64}
\Entcd{\char 70}{The Jackson integrals via the hyper\-geometric integrals}{67}
\Entcd{\char 71}{One useful identity}{69}
\Refcd{References}{69}}

\def\ContM{
\Entcd{1}{Introduction}{1}
\Entcd{2}{Discrete flat connections and local systems}{5}
\subcd{Discrete flat connections}{5}
\subcd{Discrete Gauss-Manin connection}{7}
\subcd{Connection coefficients of local systems}{7}
\subcd{The functional space of a trigonometric\/ \stype / local system}{9}
\subcd{Bases in the trigonometric hyper\-geometric space of a fiber}{11}
\subcd{The elliptic hyper\-geometric space}{14}
\Entcd{3}{\Rms / and the \qKZc /}{19}
\subcd{Highest weight \Umod /s}{19}
\subcd{The trigonometric \Rm /}{20}
\subcd{The quantum loop algebra $\Ugg $}{21}
\subcd{The trigonometric \qKZc / associated with $\gsl $}{23}
\subcd{Modules over the elliptic quantum group $\Eqg $ and the elliptic \nlc
\Rms /}{24}
\Entcd{4}{Tensor coordinates on the hyper\-geometric spaces}{27}
\subcd{Tensor coordinates on the trigonometric hyper\-geometric spaces of
fibers}{27}
\subcd{Tensor coordinates on the elliptic hyper\-geometric spaces of
fibers}{29}
\subcd{Tensor products of the hyper\-geometric spaces}{32}
\Entcd{5}{The hyper\-geometric pairing and the hyper\-geometric solutions of \nlc the \qKZe /}{34}
\subcd{The hyper\-geometric integral}{34}
\subcd{Determinant formulae for the hyper\-geometric pairing}{36}
\subcd{The hyper\-geometric solutions of the \qKZe /}{40}
\subcd{The hyper\-geometric map}{41}
\Entcd{6}{Asymptotic solutions of the \qKZe /}{44}
\Entcd{7}{Proofs}{49}
\Appencd
\Entcd{\char 65}{Basic facts about the trigonometric hyper\-geometric
space}{62}
\Entcd{\char 66}{Basic facts about the elliptic hyper\-geometric space}{66}
\Entcd{\char 67}{The Shapo\-va\-lov pairings of the hyper\-geometric spaces of
fibers}{72}
\Entcd{\char 68}{The \qSelberg / integral}{77}
\Entcd{\char 69}{The multidimensional Askey-Roy formula and Askey's \nl
conjecture}{82}
\Entcd{\char 70}{The Jackson integrals via the hyper\-geometric integrals}{86}
\Entcd{\char 71}{One useful identity}{88}
\Refcd{References}{88}}

\def\ContS{
\Entcd{1}{Introduction}{1}
\Entcd{2}{Discrete flat connections and local systems}{3}
\subcd{Discrete flat connections}{3}
\subcd{Discrete Gauss-Manin connection}{5}
\subcd{Connection coefficients of local systems}{6}
\subcd{The functional space of a trigonometric\/ \stype / local system}{7}
\subcd{Bases in the trigonometric hyper\-geometric space of a fiber}{9}
\subcd{The elliptic hyper\-geometric space}{11}
\Entcd{3}{\Rms / and the \qKZc /}{15}
\subcd{Highest weight \Umod /s}{15}
\subcd{The trigonometric \Rm /}{16}
\subcd{The quantum loop algebra $\Ugg $}{16}
\subcd{The trigonometric \qKZc / associated with $\gsl $}{18}
\subcd{Modules over the elliptic quantum group $\Eqg $ and the elliptic
\Rms /}{19}
\Entcd{4}{Tensor coordinates on the hyper\-geometric spaces}{21}
\subcd{Tensor coordinates on the trigonometric hyper\-geometric spaces of
fibers}{21}
\subcd{Tensor coordinates on the elliptic hyper\-geometric spaces of
fibers}{23}
\subcd{Tensor products of the hyper\-geometric spaces}{25}
\Entcd{5}{The hyper\-geometric pairing and the hyper\-geometric solutions of
the \qKZe /}{27}
\subcd{The hyper\-geometric integral}{27}
\subcd{Determinant formulae for the hyper\-geometric pairing}{28}
\subcd{The hyper\-geometric solutions of the \qKZe /}{31}
\subcd{The hyper\-geometric map}{32}
\Entcd{6}{Asymptotic solutions of the \qKZe /}{35}
\Entcd{7}{Proofs}{38}
\Appencd
\Entcd{\char 65}{Basic facts about the trigonometric hyper\-geometric
space}{48}
\Entcd{\char 66}{Basic facts about the elliptic hyper\-geometric space}{52}
\Entcd{\char 67}{The Shapo\-va\-lov pairings of the hyper\-geometric spaces of
fibers}{56}
\Entcd{\char 68}{The \qSelberg / integral}{60}
\Entcd{\char 69}{The multidimensional Askey-Roy formula and Askey's
conjecture}{64}
\Entcd{\char 70}{The Jackson integrals via the hyper\-geometric integrals}{67}
\Entcd{\char 71}{One useful identity}{69}
\Refcd{References}{69}}

\document

\ifAst
\hfuzz 10pt
\else
\ifMag
\hfuzz 60pt
\else
\hfuzz 8pt
\fi\fi

\def\abstext{\Abstract
The \tri/ quantized \KZv/ (\qKZ/) \eq/ associated with the \qg/ $\Uu$ is
a system of linear \deq/s with values in a tensor product of $\Uu\]$ \Vmod/s.
We solve the equation in terms of multidimensional \^{$q$-}\hgeom/ \fn/s and
define a natural \iso/ of the space of \sol/s and the tensor product of
the corresponding \eVmod/s over the \eqg/ $\Eqg$ where parameters $\rho$
and $\gm$ are related to the parameter $q$ of the \qg/ $\Uu$ and the step $p$
of the \qKZe/ via $p=e^{2\pii\rho}$ and $q=e^{-2\pii\gm}$.
\par
We construct \asol/s associated with suitable \azo/s and compute the \traf/s
between the \asol/s in terms of the dynamical elliptic \Rms/. This description
of the \traf/s gives a connection between \rep/ theories of the \qlo/ $\Ugh$
and the \eqg/ $\Eqg$ and is analogous to the Kohno-Drinfeld theorem on
the monodromy group of the \difl/ \KZe/.
\par
In order to establish these results we construct a discrete \GM/,
in particular, a suitable discrete \loc/, discrete homology and cohomology
groups with coefficients in this \loc/, and identify an associated \deq/
with the \qKZe/.
\endAbs}

\Sno 0

\ifAst
\center
\=\twelvepoint
{\BBfcorr \bls 20pt
Geometry of {\Bmmi q\)}-Hypergeometric Functions,
\\
Quantum Affine Algebras and
\\
Elliptic Quantum Groups
\par}
\vsk1.5>
\VT/$^{\,\star}$ \ and \ \Varch/$^{\,*}$
\vsk1.5>
{\tenpoint\it
$^\star$\DMO/\\ \DMOaddr/
\vsk.35>
$^*$\UNC/\\ \UNCaddr/
\vsk1.8>
\sl March \,19, 1997}
\endcenter
\ftext{\=\bls11pt $\]^\star\>$\absence/\vv-.05>\nl
{\tenpoint\hp{$^\star$}\sl E-mail\/{\rm:} \dmoemail/\,, \;\homemail/}\nl
$\]^*\>$\Grant/\nl
{\tenpoint\hp{$^*$}\sl E-mail\/{\rm:} \avemail/}}
\vsk2.4>
\abstext

\newpage

\pageno -1
{\raggedbottom
\contents
\ContA
\endco
\newpage}

\else
\vp1
\vsk-3>
\rline{\sl Ast\'erisque \,{\bf 246} {\rm(}1997\/\rm)}
\vsk>

\center
\=
{\bls 16pt
\Bbf
Geometry of {\Bmmi q\)}-Hypergeometric Functions,
\\
Quantum Affine Algebras and Elliptic Quantum Groups
\par}
\vsk1.5>
\VT/$^{\,\star}$ \ and \ \Varch/$^{\,*}$
\vsk1.5>
{\it
\ifMag
$^\star$\DMO/\\ \DMOaddr/
\else
$^\star$\DMO/, \DMOaddr/
\fi
\vsk.35>
\ifMag
$^*$\UNC/\\ \UNCaddr/
\else
$^*$\UNC/, \UNCaddr/
\fi}
\vsk1.8>
{\sl March \,19, 1997}
\endcenter
\ftext{\=\bls11pt $\]^\star\>$\absence/\vv-.05>\nl
{\tenpoint\hp{$^\star$}\sl E-mail\/{\rm:} \dmoemail/\,, \;\homemail/}\nl
$\]^*\>$\Grant/\nl
{\tenpoint\hp{$^*$}\sl E-mail\/{\rm:} \avemail/}}
\fi

%% \vsk1.8>\vsk0>\abstext\vsk1.5>\vsk0>

\ifAst\newpage\pageno1\fi
\Sect[1]{Introduction}
In this paper we solve the \tri/ quantized \KZv/ (\qKZ/) \eq/ associated with
the \qg/ $\Uu$. The \tri/ \qKZe/ associated with $\Uu$ is a system of \deq/s
for a \fn/ $\Psi(\zn)$ with values in a tensor product $\Vox$ of \Umod/s.
The system of \eq/s has the form
$$
\align
\Psi(\zmn)\, &{}=\,R_{m,m-1}(p\)z_m/z_{m-1})\ldots R_{m,1}(p\)z_m/z_1)\>
\ka^{-H_m}\;\x
\\
\nn3>
&{}\>\x\,R_{m,n}(z_m/z_n)\ldots R_{m,m+1}(z_m/z_{m+1})\>\Psi(\zn)\,,
\endalign
$$
$\mn$, where $p$ and $\ka$ are parameters of the \qKZe/, $H$ is a generator of
the Cartan subalgebra of $\Uu$, $H_m$ is $H$ acting in the \^{$m$-th} factor,
$R_{l,m}(x)$ is the \tri/ \Rm/ $R_{V_lV_m}(x)\in\End(V_l\ox V_m)$ acting in
the \^{$l$-th} and \^{$m$-th} factors of the tensor product of \Umod/s.
In this paper we consider only the steps $p$ with absolute value less than 1.
\par
The \qKZe/ is an important system of \deq/s. The \qKZe/s had been introduced in
\Cite{FR} as \eq/s for matrix elements of vertex operators of the \qaff/.
An important special case of the \qKZe/ had been introduced earlier in \Cite{S}
as \eq/s for form factors in massive integrable models of quantum field theory;
relevant \sol/s of these \eq/s had been given therein. Later the \qKZe/s were
derived as \eq/s for correlation \fn/s in lattice integrable models,
\cf. \Cite{JM} and references therein.
\par
In the \qcl/ limit the \qKZe/ turns into the \difl/ \KZv/ \eq/ for conformal
blocks of the Wess-Zumino-Witten model of conformal field theory on the sphere.
\par
Asymptotic \sol/s of the \qKZe/ as the step $p$ tends to $1$ are closely
related to diagonalization of the transfer-matrix of the corresponding lattice
integrable model by the algebraic \Ba/ method \Cite{TV2}.
\Par
We describe the space of \sol/s of the \qKZe/ in terms of \rep/ theory.
Namely, we consider the \eqg/ $\Eqg$ with parameters $\rho$ and $\gm$ defined
by ${p=e^{2\pii\rho}}$, ${q=e^{-2\pii\gm}}$ and \)the \Emod/s $\Vqnz$ where
$V^e_m(z_m)$ is the \eVmod/ over $\Eqg$ which corresponds to the \Umod/ $V_m$.
Notice that as a vector space the \eVmod/ $V^e_m(z_m)$ does not depend on
$z_m$. For every \perm/ $\tau\in\S^n$ we consider the tensor product $\Vqxt$
and establish a natural \iso/ of the space $\SS$ of \sol/s of the \qKZe/ with
values in $\Vox$ and the space $\Vqxt\ox\FF$, where $\FF$ is the space of
\fn/s of $\zn$ which are \p-periodic \wrt/ each of the \var/s,
$$
\CC_\tau:\Vqxt\ox\FF\,\to\,\SS\,,
$$
\cf. \(CCtau). Notice that if $\Psi(z)$ is a \sol/ of the \qKZe/ and $F(z)$
is a \p-periodic \fn/, then also $F(z)\Psi(z)$ is a \sol/ of the \qKZe/.
\par
We call the \iso/s $\CC_\tau$ the \tenco/ on the space of \sol/s.
The compositions of the \iso/s are linear maps
$$
\CC_{\tau,\tau'}(\zn):\Vqxtt\,\to\,\Vqxt
$$
depending on $\zn$ and \p-periodic \wrt/ all \var/s. We call these compositions
the \traf/s. It turns out that the \traf/s are defined in terms of the elliptic
\Rms/ ${\Rq_{V^e_lV^e_m}(x,\la)\in\End(V^e_l\ox V^e_m)}$ acting in tensor
products of \Emod/s. Namely, for any \perm/ $\tau$ and for any transposition
$(m,m+1)$ the \traf/
$$
\CC_{\tau,\tau\cdot(m,m+1)}(\zn):\Vqxtm\,\to\,\Vqxt
$$
equals the operator
\ifMag
$$
P_{V^e_{\tau_{m+1}}\}V^e_{\tau_m}}\:
\Rq_{V^e_{\tau_{m+1}}\}V^e_{\tau_m}}\!\bigl(z_{\tau_{m+1}}/z_{\tau_m},
(\eta^{\)H}\}\lox\eta^{\)H}\}\ox\%{\eta^{-H}}_{\^{$m$-th}}\}\lox\eta^{-H})
\>\eta\1\ka\bigr)
$$
\else
${P_{V^e_{\tau_{m+1}}\}V^e_{\tau_m}}\:
\Rq_{V^e_{\tau_{m+1}}\}V^e_{\tau_m}}\!\bigl(z_{\tau_{m+1}}/z_{\tau_m},
(\eta^{\)H}\}\lox\eta^{\)H}\}\ox\%{\eta^{-H}}_{\^{$m$-th}}\}\lox\eta^{-H})
\>\eta\1\ka\bigr)}$
\fi
twisted by certain \adjm/s, \cf. \(CCC) and Theorem~\[localq].
Here $P_{V^e_lV^e_m}\:$ is the transposition of the tensor factors.
\Par
We consider \azo/s $|z_{\tau_m}/z_{\tau_{m+1}}|\ll 1$, $\mn-1$, labelled by
\perm/s ${\tau\in\S^n}$. For every \azo/ we define a basis of \asol/s of the
\qKZe/. We show that for every \perm/ $\tau$ the basis of the corresponding
\asol/s is the image of the standard monomial basis in $\Vqxt$ under the map
$$
\Vqxt\,\to\,\Vqxt\ox 1\,\hookrightarrow\,\Vqxt\ox\FF
\,\ {\overset{C_\tau}\to\lto}\,\ \SS\,,
$$
\cf. Theorem~\[asol]. The last two statements express the \traf/s between
the \asol/s via the elliptic \Rms/.
\Par
The \tri/ \Rm/ $R_{V_lV_m}\:(x)\in\End(V_l\ox V_m)$ is defined in terms of the
action of the \qlo/ $\Ugg$ in the tensor product of \Umod/s. The \qlo/ $\Ugg$
is a Hopf algebra which contains the \qg/ $\Uu$ as a Hopf subalgebra and has
a family of \hom/s $\Ugg\to\Uu$ depending on a parameter. Therefore, each
\Umod/ $V_m$ carries a \Uhmod/ structure $V_m(x)$ depending on a parameter.
For \Vmod/s $V_l,\,V_m$ over $\Uu$ the \qlo/ modules $V_l(x)\ox V_m(y)$ and
$V_m(y)\ox V_l(x)$ are isomorphic for generic $x,y$. Moreover, for \irr/
\Umod/s $V_l,\,V_m$ the \qlo/ modules $V_l(x)\ox V_m(y)$ and $V_m(y)\ox V_l(x)$
are \irr/ and isomorphic for generic $x,y$. The map
$$
P_{V_lV_m}\:R_{V_lV_m}\:(x/y):V_l(x)\ox V_m(y)\,\to\,V_m(y)\ox V_l(x)
$$
is the unique suitably normalized intertwiner \Cite{T}, \Cite{CP}.
\Par
The elliptic \Rm/ $\Rq_{V^e_lV^e_m}(x,\la)\in\End(V^e_l\ox V^e_m)$ is defined
in terms of the action of the \eqg/ $\Eqg$ in the tensor product of \eVmod/s.
For \eVmod/s $V^e_l(x),\,V^e_m(y)$ over $\Eqg$, the \Emod/s
$V^e_l(x)\ox V^e_m(y)$ and $V^e_m(y)\ox V^e_l(x)$ are isomorphic for generic
$x,y$. The map
$$
P_{V^e_lV^e_m}\:\Rq_{V^e_lV^e_m}(x/y,\la):
V^e_l(x)\ox V^e_m(y)\,\to\,V^e_m(y)\ox V^e_l(x)
$$
is the unique suitably normalized intertwiner \Cite{F}, \Cite{FV}.
\Par
Our result on the \traf/s between \asol/s together with the indicated
construction of \Rms/ shows that the \qKZe/ establishes a connection between
\rep/ theories of the \qlo/ $\Ugg$ and the \eqg/ $\Eqg$. Our result is
analogous to the Kohno-Drinfeld theorem on the monodromy group of the \difl/
\KZv/ \eq/ \Cite{K}, \Cite{D2}. The Kohno-Drinfeld theorem establishes
a connection between \rep/ theories of a Lie algebra and the corresponding
\qg/, see \Cite{D2}. Using the ideas of the Kohno-Drinfeld result it was proved
in \Cite{KL} that the category of \rep/s of a \qg/ is \eqv/ to a suitably
defined fusion category of \rep/s of the corresponding affine Lie algebra.
Similarly to the Kazhdan-Lusztig theorem one could expect that our result for
the \dif/ \qKZe/ could be a base for a Kazhdan-Lusztig type result connecting
certain categories of \rep/s of \qaff/s and \eqg/s, \cf. \Cite{KS}.
\Par
In this paper we consider the \tri/ \qKZe/. There are other types of the
\qKZe/: the \rat/ \qKZe/ \Cite{FR} and the elliptic \qKZBe/ \Cite{F}.
Here \KZB/ stands for \KZv/-Bernard, and the \dif/ \qKZBe/ is a discretization
of the \difl/ \KZB/ \eq/ for conformal blocks on the torus.
\par
The \rat/ \qKZe/ was considered in \Cite{TV3}. It is a system of \deq/s
analogous to the \tri/ \qKZe/ in which the role of the \tri/ \Rm/ is played by
the \rat/ \Rm/ defined in terms of the Yangian \rep/s. In \Cite{TV3} we solved
the \rat/ \qKZe/ in terms of multidimensional \hint/s of \MB/ type, introduced
\azo/s and described \traf/s between \asol/s in terms of \tri/ \Rms/,
thus showing that the \rat/ \qKZe/ gives a connection between \rep/s of
Yangians and \qaff/s.
\par
The elliptic \qKZBe/ is an analogue of the \tri/ \qKZe/ in which the role of
the \tri/ \Rm/ is played by the \dyne/ \Rm/ \Cite{F}. The \dyne/ \Rm/ is
a matrix acting in the tensor product of two \emod/s over the \eqg/ $\Eqg$,
and the \eqg/ is an elliptic analogue of the \qaff/ $\Uuh$ associated with
$\gsl$ \Cite{F}, \Cite{FV}. The \eqg/ depends on two parameters: an \ecurve/
${\C/(\Z+\rho\)\Z)}$ and Planck's constant $\gm$. The elliptic \qKZBe/ is
considered in \Cite{FTV1}, \Cite{FTV2} where we solve it in terms of
multidimensional elliptic \^{$q$-}\hint/s. We formulate and solve in
\Cite{FTV2} a ``monodromy'' problem for the \qKZBe/ analogous to the problem
of description of the \traf/s for the \rat/ and \tri/ \qKZe/s. We describe in
\Cite{FTV2} the ``monodromy'' matrices of the elliptic \qKZBe/ associated with
the \eqg/ $\Eqg$ in terms of the elliptic \Rm/ associated with the elliptic
quantum group $E_{\al,\gm}(\gsl)$ where $\al$ is the step of the initial \qKZB/
\deq/, thus showing the \sym/ role of the \ecurve/s ${\C/(\Z+\rho\)\Z)}$ and
${\C/(\Z+\al\)\Z)}$ in the story.
\Par
In this paper, in order to establish a connection between \rep/ theories of the
\qlo/ $\Ugg$ and the \eqg/ $\Eqg$ we define a discrete analogue of a locally
\trib/ and a \loc/ on the space of bundle. We define a discrete analogue of the
\GM/ for the discrete locally \trib/ with a discrete \loc/ and consider the
corresponding \deq/. We identify that \dif/ Gauss-Manin \eq/ with the \dif/
\qKZe/. To realize this idea we introduce a suitable discrete de~Rham complex
and its cohomology group in the spirit of \Cite{A}, then we define the homology
group as the dual space to the cohomology group and construct a family of
discrete cycles, elements of the discrete homology group, using ideas of
\Cite{S}. We construct the space of discrete cycles as a certain space of
\fn/s. Having a representative of a discrete cohomology class (a \fn/) and
a discrete cycle (a \fn/ again) we define the pairing (the \hpair/) between the
cohomology class and the cycle as an integral of their product with a certain
fixed ``\hgeom/ \phf/'' over a certain fixed contour of the middle dimension.
We show that there are enough discrete cycles and they form the space dual to
the quotient space of the space of our discrete closed forms modulo discrete
coboundaries. To prove this we compute the determinant of the period matrix
and get an explicit formula \(mu<>0) for the determinant analogous to the
determinant formulae for the \hgeom/ \fn/s of \MB/ type \Cite{TV3} and for the
``continuous'' \hgeom/ \fn/s \Cite{V1}, \cf. Loeser's determinant formula for
the Frobenius transformation \Cite{L}.
\par
The form of our discrete cycles suggests a natural identification of the space
of our discrete cycles with a tensor product of \Emod/s and this identification
allows us to prove the result on \traf/s between \asol/s.
\par
The paper is organized as follows. \Chap/s~{\SNo{2}\,--\,\SNo{2d}} contain
constructions and statements and \Chap/~\SNo{4} contains proofs. We give
necessary prelimimnaries, proofs of technical results and some applications
in Appendices.
\Par
Parts of this work had been written when the authors visited the University of
Tokyo, the Kyoto University, the Osaka University, the University Paris VI,
the IHES at Bures~sur~Yvette, \'Ecole Normale Sup\'erieure de Lyon, the MSRI at
Berkeley, the ETH at Z\"urich. The authors thank those institutions for
hospitality. The authors thank \Feld/, P\]\&Etingof and E\&Mukhin for valuable
discussions. The authors also thank the referee for important comments.

\ifAst\newpage\fi
\Sect[2]{Discrete flat connections and \loc/s}
In this \chap/ we recall basic notions introduced in \Cite{TV3}.

\subsect{Discrete flat connections}
Let ${\Cx\!=\C\setminus\lb 0\rb}$. Consider a complex vector space $\Cn$ called
the \em{base space}. Let $\Cxn$ be the complement of the coordinate hyperplanes
in the base space. Fix a complex number $p$, \st/ $p\ne 0,1$, which is called
the \em{step}.
\ifMag\else\newpage\fi
The lattice $\Z^n$ acts on the base space by dilations:
$$
l:(\zn)\,\map\,(p^{l_1}z_1\lc p^{l_n}z_n)\,,\qqq l\in\Z^n\,.
$$
Let $\BB\sub\Cxn$ be an \inv/ subset of the base space. Say that there is a
\em{bundle with a discrete connection} over $\BB$ if for any $z\in\BB$ there
are a vector space $V(z)$ and linear \iso/s
$$
A_m(\zn):V(\zmn)\,\to\,V(\zn)\,,\qqq\mn\,.
$$
The connection is called \em{flat} (or \em{integrable}) if the \iso/s
$A_1\lc A_n$ commute:
\ifMag
$$
\align
A_l(\zn)\> & A_m(z_1\lc p\)z_l\lc z_n)\;=
\Tag{flat}
\\
\nn4>
=\;{} & A_m(\zn)\>A_l(\zmn)\,.
\endalign
$$
\else
$$
A_l(\zn)\>A_m(z_1\lc p\)z_l\lc z_n)\,=\,A_m(\zn)\>A_l(\zmn)\,.
\Tag{flat}
$$
\fi
Say that a \em{discrete subbundle} in $\BB$ is given if a subspace in every
fiber is distinguished and the family of subspaces is \inv/ \wrt/ the
connection.
\par
A section $s:z\map s(z)$ is called \em{periodic} (or \em{horizontal}) if
its values are \inv/ \wrt/ the connection:
$$
A_m(\zn)\>s(\zmn)\,=\,s(\zn)\,,\qqq\mn\,.
\Tag{period}
$$
A \fn/ $f(\zn)$ on the base space is called a \em{quasiconstant} if
$$
f(\zmn)\,=\,f(\zn)\,,\qqq\mn\,.
$$
Periodic sections form a module over the ring of quasiconstants.
\par
The \em{dual bundle} with the \em{dual connection} has fibers
$V^*(z)$ and \iso/s
$$
A_m^*(\zn):V^*(\zn)\,\to\,V^*(\zmn)\,.
$$
Let $s_1\lc s_N$ be a basis of sections of the initial bundle. Then the
\iso/s $A_m$ of the connection are given by matrices $\A\"m$:
$$
A_m(\zn)\>s_k(\zmn)\,=\,\sum_{l=1}^N \A\"m_{kl}(\zn)\>s_l(\zn)\,.
$$
\par
For any section ${\psi:z\map\psi(z)}$ of the dual bundle denote by
${\Psi:z\map\Psi(z)}$ its coordinate vector,
$\Psi_k(z)=\bra\psi(z),s_k(z)\ket$. The section $\psi$ is periodic
if and only if its coordinate vector satisfies the system of \deq/s
$$
\Psi(\zmn)\,=\,\A\"m(\zn)\>\Psi(\zn)\,,\qqq\mn\,.
$$
This system of \deq/s is called the \em{periodic section \eq/}.
\Par
Say that \fn/s $\phin$ in \var/s $\zn$ form a \em{\sconn/} if
$$
\phi_l(\zmn)\>\phi_m(\zn)\,=\,\phi_m(z_1\lc p\)z_l\lc z_n)\>\phi_l(\zn)
$$
for all $l,\mn$. These \fn/s define a connection on the trivial complex
\onedim/ vector bundle.
\par
The \sconn/ is called \em{decomposable} if it has the form
\ifAst
$$
\Line{\phi_m(\zn)\,=\,
\ka_m\>\Bigl[\,\prod_{1\le l<m}\!\pho_{lm}(p\1z_l/z_m)\>\Bigr]\1\!
\prml\pho_{ml}(z_m/z_l)\,,\hfil\mn\,,}
$$
\else
$$
\phi_m(\zn)\,=\,
\ka_m\>\Bigl[\,\prod_{1\le l<m}\!\pho_{lm}(p\1z_l/z_m)\>\Bigr]\1\!
\prml\pho_{ml}(z_m/z_l)\,,\qqq\mn\,,
$$
\fi
for certain \fn/s $\pho_{lm}$, $l<m$, in one \var/ and nonzero complex numbers
$\ka_m$. The \fn/s $\pho_{lm}$ are called \em{\prim/s} and $\ka_m$ are called
\em{scaling parameters}.
\Par
A \fn/ $\Phb(\zn)$ is called a \em{\phf/} of a \sconn/ if
$$
\Phb(\zmn)\,=\,\phi_m(\zn)\>\Phb(\zn)\,,\qqq\mn\,.
$$
Similarly, a \fn/ $\Pho(x)$ is called a \em{\phf/} of a \fn/ $\pho(x)$ in one
\var/ if ${\Pho(p\)x)={}}\alb\pho(x)\>\Pho(x)$. Notice that the \phf/s are not
unique.
\Par
For any \fn/ $f(\zn)$ define new \fn/s $Q_1f\lc Q_nf$ and $D_1f\lc D_nf$ by
the rule:
$$
(Q_mf)(\zn)\,=\,\phi_m(\zn)\>f(\zmn)\,,
$$
and $D_mf\,=\,Q_mf-f$, $\mn$. The \fn/s $D_1f\lc D_nf$ are the \em{discrete
partial derivatives} of the \fn/ $f$. We have that $D_lD_mf=D_mD_lf$ for any
$l,\mn$.
\par
Let $\Ff$ be a vector space of \fn/s of $\zn$ \st/ the operators $Q_1\lc Q_n$
induce linear \iso/s of $\Ff$:
$$
Q_m:\Ff\,\to\,\Ff\,.
$$
Say that the space $\Ff$ and the \conn/s $\phin$ form a one-dimen\-sional
%% \onedim/
\em{\dloc/} on $\Cxn\!$. $\Ff$ is called the \em{\fn/al space} of the \loc/.
\Par
Define the \em{de~Rham complex} $\bigl(\Omb(\Ff),D\bigr)$ of the \loc/ in a
standard way. Namely, set
$$
\Om^a\,=\,\bigl\lb\>\om=\tsum_{k_1\lc k_a}f_{k_1\lc k_a}\>
Dz_{k_1}\lsym\wedge Dz_{k_a}\>\bigr\rb
$$
where $Dz_1\lc Dz_n$ are formal symbols and the coefficients $f_{k_1\lc k_a}$
belong to $\Ff$. Define the \difl/ of a \fn/ by $Df=\sun D_mf\>Dz_m$,
and the \difl/ of a form by
$$
D\)\om\;=
\tsum_{k_1\lc k_a}Df_{k_1\lc k_a}\wedge Dz_{k_1}\lsym\wedge Dz_{k_a}\,.
$$
The cohomology groups $H^1\lc H^n$ of this complex are called the
\em{cohomology groups of $\Cxn\!$ with coefficients in the \dloc/}.
In particular, the top cohomology group is $H^n=\Ff/D\Ff$ where
$D\Ff=\sun D_m\Ff$. The dual spaces $H_a=(H^a)^*$ are called
the \em{homology groups}.

\subsect{Discrete \GM/}
There is a geometric construction of bundles with discrete flat connections,
a discrete version of the \GM/ construction.
\par
Let $\pi:\Cln\to\Cn$ be an affine projection onto the base with fiber $\Cl\!$.
$\Cln$ will be called the \em{total space}. Let $\zn$ be coordinates on the
base, $\tell$ coordinates on the fiber, so that $\tell$, $\zn$ are coordinates
on the total space. When it is convenient, we denote the coordinates $\zn$ by
$t_{\ell+1}\lc t_{\ell+n}$.
\par
Let $\Ff$, $\phi_1\lc\phi_{\ell+n}$ be a \loc/ on $\Cxln\!$. For a point
$z\in\Cxn\!$ define a \loc/ $\Ff(z)$, $\phi_a(\cdot;z)$, $\aell$, on the fiber
over $z$. Set
$$
\Ff(z)\,=\,\lb\>f\vst{\pi\1(z)}\vert f\in\Ff\>\rb\qquad
\text{and}\qquad
\phi_a(\cdot;z)\,=\,\phi_a\vst{\pi\1(z)}\,.
$$
The de~Rham complex, cohomology and homology groups of the fiber are denoted by
$\bigl(\Omb(z),D(z)\bigr)$, $H^a(z)$ and $H_a(z)$, \resp/.
\par
There is a natural \hom/ of the de~Rham complexes
$$
(\Omb(\Cxln\!,\Ff),D)\,\to\,(\Omb(z),D(z))\,,\qqq \om\map\om|_{\pi\1(z)}\,,
$$
where the restriction of a form is defined in a standard way: all symbols
$Dz_1\lc Dz_n$ are replaced by zero and all coefficients of the remaining
monomials $Dt_{k_1}\lsym\wedge Dt_{k_a}$ are restricted to the fiber.
\par
For a fixed $a$ the vector spaces $H^a(z)$ form a bundle with a discrete flat
connection. The linear maps
$$
A_m(\zn):H^a(\zmn)\,\to\,H^a(\zn)
$$
are defined as follows. Define $Q_m:\Om^a(\Cxln,\Ff)\to\Om^a(\Cxln,\Ff)$ by
$$
\om\;\map\tsum_{k_1\lc k_a}Q_mf_{k_1\lc k_a}\>
Dz_{k_1}\lsym\wedge Dz_{k_a}\,.
\vvgood
$$
Then $Q_m$ induces a homomorphism of the de~Rham complexes
$$
\bigl(\Omb(\zmn),D(\zmn)\bigr)\,\to\,\bigl(\Omb(\zn),D(\zn)\bigr)\,.
$$
We set $A_m(\zn)$ to be equal to the induced map of the cohomology spaces.
This connection is called the \em{discrete \GM/}.
\par
The \GM/ on the cohomological bundle induces the dual flat connection on the
homological bundle:
$$
A^*_m(\zn):H_a(\zn)\,\to\,H_a(\zmn)\,.
$$

\subsect{Connection coefficients of \loc/s}
In this paper we study the \GM/ for a class of \loc/s with decomposable
\conn/s, namely for \tri/ \stype/ \loc/s \Cite{TV3}, for the \rat/ case see
\Cite{TV3} and for the elliptic case see \Cite{FTV1}, \Cite{FTV2}.
\par
Primitive factors and scaling parameters of a \tri/ \stype/ \loc/ have the
following form:
\vvn->
$$
\NN3>
\alignat2
\pho_{ab}(x)\,&{}=\;{x-\eta\over\eta\)x-1}\qqq &&\for\quad a<b\le\ell\,,
\\
\pho_{ab}(x)\,&{}=\;{\xi_{b-\ell}x-1\over x-\xi_{b-\ell}}\qqq &&
\for\quad a\le\ell<b\,,
\\
\pho_{ab}(x)\,&{}=\;1\qqq &&\for\quad \ell<a<b\,.
\\
\ka_a\, &{}=\,\ka\qqq &&\for\quad a\le\ell\,,
\\
\ka_a\, &{}=\,1\qqq &&\for\quad \ell<a\,.
\endalignat
$$
Such a \sconn/ depends on $n+2$ nonzero complex numbers $\xin,\,\eta,\ka$.
The \conn/s of a \tri/ \stype/ \loc/ have the form
\ifMag
$$
\NN4>
\gather
\phi_a(t,z)\,=\,\ka\pron{\xi_mt_a-z_m\over t_a-\xi_mz_m}\,
\prod_{a<b\le\ell}{t_a-\eta\)t_b\over\eta\)t_a-t_b}\;
\prod_{1\le b<a}{p\)t_a-\eta\)t_b\over p\)\eta\)t_a-t_b}\;,
\Tag{conns}
\\
\nn-4>
\Rline{\aell\,,}\\
\\
\nn-24>
\kern1em
\phi_{\ell+m}(t,z)\,=\,
\pral\,{t_a-p\>\xi_mz_m\over\xi_mt_a-p\)z_m}\;,\qqq\mn\,.
\kern-1em
\endgather
$$
\else
$$
\NN4>
\gather
\kern1em
\phi_a(t,z)\,=\,\ka\pron{\xi_mt_a-z_m\over t_a-\xi_mz_m}\,
\prod_{a<b\le\ell}{t_a-\eta\)t_b\over\eta\)t_a-t_b}\;
\prod_{1\le b<a}{p\)t_a-\eta\)t_b\over p\)\eta\)t_a-t_b}\;,\qqq\aell\,,
\kern-1em
\Tag{conns}
\\
\kern1em
\phi_{\ell+m}(t,z)\,=\,
\pral\,{t_a-p\>\xi_mz_m\over\xi_mt_a-p\)z_m}\;,\qqq\mn\,.
\kern-1em
\endgather
$$
\fi
Let $\Pho(x;\al)$ be a \phf/ of the \fn/ $(x\al-1)/(x\al\1-1)$.
Then a \phf/ of the \sconn/ is given by
$$
\gather
\Phb(\tell,\zn)\,=\,f(\tell,\zn)\>\Phi(\tell,\zn)
\ifMag\kern-2em\fi
\Tag{Phb}
\\
\nn2>
\Text{where}
\nn-6>
\Phi(\tell,\zn)\,=\,\pron\,\pral\>\Pho(t_a/z_m;\xi_m)\prab\Pho(t_a/t_b;\eta\1)
\ifMag\kern-2em\fi
\Tag{Phi}
\endgather
$$
and $f(\tell,\zn)$ is an arbitrary \fn/ \st/
$$
\gather
\ifMag\qquad\fi
f(\tpell,\zn)\,=\,\ka\>\eta^{\ell-2a+1}\zy\,f(\tell,\zn)\,,
\Tagg{tzp}
\\
\nn6>
\ifMag\qquad\fi
f(\tell,\zmn)\,=\,\xi_m^\ell\>f(\tell,\zn)\,.
\endgather
$$
Later in this \chap/ we will describe a convenient space of such \fn/s called
the \ehgf/.
\Par
The \fn/ $\Phi(\tell,\zn)$ will be called a \em{short \phf/}.
\Ex
Let ${0<|p|<1}$. Let $(u)\9=\prod_{k=0}^\8(1-p^ku)$. We can take
$$
\Pho(x;\al)={(x\al\1)\9\!\over(x\al)\9}\;.
$$
Then the \phf/ $\Pho(x;\al)$ has a symmetry property
$$
\Pho(-x;\al)\,=\,\Pho(x;\al)\;{(x\al-1)\over(x-\al)}\,
{\tht(x\al)\over\al\>\tht(x/\al)}
$$
where ${\tht(u)=(u)\9(p\)u\1)\9(p)\9}$ is the Jacobi theta-\fn/.
%% with the transformation property: $\tht({p\)u})=-u\1\tht(u)$.
The symmetry property for $\Pho(x;\al)$ leads to a symmetry property
\ifMag
$$
\align
& \Phi(\tall,\zn)\;=
\Tag{Phisym}
\\
\nn6>
&\qqq{}=\,\Phi(\tell,\zn)\;{(t_a-\eta\)t_{a+1})\over(\eta\)t_a-t_{a+1})}\,
{\eta\>\tht(\eta\1t_a/t_{a+1})\over\tht(\eta\)t_a/t_{a+1})}
\endalign
$$
\else
$$
\Phi(\tall,\zn)\,=\,
\Phi(\tell,\zn)\;{(t_a-\eta\)t_{a+1})\over(\eta\)t_a-t_{a+1})}\,
{\eta\>\tht(\eta\1t_a/t_{a+1})\over\tht(\eta\)t_a/t_{a+1})}
\kern-.25em
\Tag{Phisym}
$$
\fi
of the short \phf/ of the \sconn/. This property motivates definitions \(act)
and \(ect) of certain actions of the \symg/.
\enddemo

\subsect{The \fn/al space of a \tri/\{\/ \{\stype/ \loc/}
Define the \fn/al space $\Fwh$ of a \tri/ \stype/ \loc/ as the space of \raf/s
on the total space $\Cln$ with at most simple poles at the following
hyperplanes:
$$
t_a=p^{s+1}\xi_m\1z_m\,,\qquad t_a=p^{-s}\xi_mz_m\,,\qquad
t_a=p^{s+1}\eta\)t_b\,,\qquad t_b=p^s\eta\)t_a\,,
\ifMag\kern-1em\else\kern-2em\fi
\Tag{list1}
$$
$1\le a<b\le\ell$, $\mn$, $s\in\Zp$, and any poles at the coordinate
\aagood
hyperplanes ${t_a=0}$, $\aell$, and ${z_m=0}$, $\mn$. It is easy to check
that the \fn/al space
\vv.1>
is \inv/ \wrt/ all the shift operators $Q_1^{\pm1}\lc Q_n^{\pm1}\!$.
\par
Define an action of the \symg/ $\Sl$ on the \fn/al space:
$$
\gather
\si:\Fwh\,\to\,\Fwh,\qqq f\map[\)f\)]\vpp\si\,,\qqq\si\in\Sl,
\Tag{act}
\\
\nn4>
[\)f\)]\vpp\si(\tell,\zn)\,=\,f(\tsill,\zn)
\prabsi\;{t_{\si_b}-\eta\)t_{\si_a}\over\eta\)t_{\si_b}-t_{\si_a}}\;.
\endgather
$$
The operators $Q_1\lc Q_{\ell+n}$ and $D_1\lc D_{\ell+n}$ commute with
the action of the \symg/.
\par
We extend the \SL-action to the de~Rham complex assuming that it respects
the exterior product and
$$
\si:Dt_a\,\map\,Dt_{\si_a}\,,\qqq\si:Dz_m\,\map\,Dz_m\,,\qqq\si\in\Sl.
$$
The same formulae define an action of the \symg/ on the de~Rham complex of a
fiber. The \hom/ of the restriction of the de~Rham complex of the total space
to the de~Rham complex of a fiber commutes with the action of the \symg/.
The action of the \symg/ induces an action of the \symg/ on the homology and
cohomology groups. The \GM/ commutes with this action.
\par
If a \symg/ acts on a vector space $V\!$, we denote by $V_{\]\sss\Sig}$
the subspace of \inv/ vectors and by $V_{\}\sss A}$ the subspace of skew-\inv/
vectors.
\par
In this paper we are interested in the skew-\inv/ part $\HA^\ell(z)$ of
the top cohomology group of a fiber. This subspace is \gby/ forms
${f\)\Dtel}$ where $f$ runs through the space $\Fsi(z)$ of \inv/ \fn/s.
\par
Introduce an important \em{\thgf/} $\Fo\sub\Fsi$ as the subspace of
\fn/s of the form
$$
P(\tell,\zn)\,\tpral\}t_a\pron\,\pral\,{1\over t_a-\xi_mz_m}\;
\prab\ \,{t_a-t_b\over\eta\)t_a-t_b}
\Tag{thgf}
$$
where $P$ is a \pol/ with complex coefficients which is \sym/ in \var/s $\tell$
and has degree less than $n$ in each of the \var/s $\tell$. The restriction of
the \thgf/ to a fiber defines the \em{\thgf/} $\Fo(z)\sub\Fsi(z)$ of the fiber
which is a complex \fd/ vector space. If required, we will write down
explicitly dependence of the \thgf/s on $\zn$, $\xin$ and $\ell$, that is
$$
\Fo=\Fo[\zn;\xin;\ell\,]\qquad\text{and}\qquad
\Fo(z)=\Fo[\zn;\xin;\ell\,](z)\,. \vp{\Big|}
$$
\par
A form ${f\)\Dtel}$ with the coefficient \ifAst\^{$f\)$}\else $f$\fi{}
in the \thgf/ of a fiber is called a \em{\hform/}. The subspace
$\H(z)\sub\HA^\ell(z)$ of the top cohomology group of a fiber \gby/
the \hform/s is called the \em{\hcg/}.
\goodast
\Par
The union of the hyperplanes
$$
\xi_l\>\xi_mz_l/z_m\,=\,p^s\eta^r\,,\qqq\roll\,,\qquad s\in\Z\,,
\Tag{dis}
$$
$l,\mn$, $l\ne m$, in the base space $\Cn$ is called the \em{discriminant}.
The complement to the discriminant in $\Cxn\!$ will be denoted by $\BB$.
\Th{subbundle}
\back\Cite{V3},\;\Cite{TV1}
The family of subspaces ${\lb\>\H(z)\>\rb}_{z\in\BB}$ is \inv/ \wrt/ the
\GM/ and, therefore, defines a discrete subbundle.
\endpro
\nt
This subbundle is called the \em{\hgeom/ subbundle}.
\Par
Later on we make the following assumptions. We always assume that the step $p$
is \st/ $0<|p|<1$, and the parameters $\eta$, $\ka$, $\xin$, $\zn$ are nonzero.
We often assume that the parameter $\eta$ is \st/
$$
\eta^r\!\ne\)p^s\,,\qqq\rell\,,\qquad s\in\Z\,,
\Tag{npZ}
$$
the parameters $\xin$ are \st/
$$
\xi_m^2\nepsr\,,\qqq \mn\,,\qquad\rll\,,\qquad s\in\Z\,,
\Tag{Lass}
$$
and the coordinates $\zn$ obey the condition
$$
\xi_l^{\pm1}\xi_m^{\pm1}z_l/z_m\,\ne\,p^s\eta^r\,,\qqq l,\mn\,,\quad l\ne m\,,
\qquad s\in\Z\,,
\Tag{assum}
$$
for any $\rll$ and for an arbitrary combination of signs.
\Th{kanon}
Let ${\ka\neps\eta^{-r}\!\yz}$ and ${\)\ka\ne p^{-s-1}\eta^r\zy\!}$, \,
\vv.1>
$\roll$, $s\in\Zp$. Let $0<|p|<1$. Let \(npZ)\;--\;\(assum) hold. Then
$$
\dim\H(z)\,=\,\dim\Fo(z)\,=\,{n+\ell-1\choose n-1}\,.
$$
\vvv->
\endpro
\nt
This means that
\ifMag\vv-.2>\else\vv-.5>\fi
$$
\H(z)\,\simeq\,\Fo(z)\,.
\Tag{HFo}
$$
In what follows we will consider in detail only two of the special values of
the scaling parameter $\ka$ mentioned in Theorem~\[kanon], namely $\ka=\kone$
and $\ka=\ktwp$,
\vvn.1>
which correspond to the values $r=\ell-1$ and $s=0$. In principle, all other
special values of the scaling parameter $\ka$ can be considered similarly.
\Th{kaon}
Let either $\ka=\kone$ or $\ka=\ktwp$. Let $0<|p|<1$. Let
\(npZ)\;--\;\(assum) hold. If ${\pron\!\xi_m^2\nepsr}$ for all $\rpll$
and $s\in\Znn$, then $\dsize\dim\H(z)\,=\,{n+\ell-2\choose n-2}$.
\endpro
\nt
Theorems~\[kanon] and \[kaon] follow from Theorems~\[mu<>0] and \[mu=0],
\[mu=02], \resp/, and Lemma~\[IWwD].
\Par
Theorem~\[kanon] means that if the value of the scaling parameter $\ka$ is not
special, then every nonzero \hform/ defines a nonzero cohomology class.
\vv.15>
If either ${\ka=\kone}$ or ${\ka=\ktwp}$, then Theorem~\[kaon] says that
there exist exact \hgeom/ forms. We describe them in Lemma \[wDw].

\subsect{Bases in the \thgf/ of a fiber}
The \fd/ \thgf/ $\Fo(z)$ of a fiber has $n!$ remarkable bases. These bases
will allow us to identify the geometry of an \stype/ \loc/ with \rep/ theory.
The bases are labelled by elements of the \symg/ $\S^n\!$. First we define
the basis corresponding to the unit element of the \symg/.
\Par
Let
\vv-.5>
$$
\Zln\,=\,\lb\>\lg\in\Zpn\vert\tsun\lg_m=\ell\>\rb\,.
\Tag{Zln}
$$
\vvv-.5>
Set $\lg^m=\sum_{k=1}^m\lg_k$. In particular, $\lg^0=0$, $\lg^n=\ell$.
For any $\lg\in\Zln\!$ define a \raf/ $w_\lg\in\Fo$ as follows:
\ifMag
$$
\align
\kern1.5em
w_\lg( & \tell,\zn)\;=
\Tagg{wlga}
\\
\nn2>
&{}=\>\prod_{k=1}^n\,\prod_{s=1}^{\lg_k}\>{1-\eta\over 1-\eta^s}\,
\susi\,\Bigl[\,\pron\>\prod_{a\in\Gm_{\Rph lm}}
\Bigl(\,{t_a\over t_a-\xi_mz_m}
\ \prlm\,{\xi_l\)t_a-z_l\over t_a-\xi_lz_l}\,\Bigr)\>\Bigr]_\si
\kern-1.5em
\endalign
$$
\else
$$
\kern1.5em
w_\lg(\tell,\zn)\,=\>\prod_{k=1}^n\,\prod_{s=1}^{\lg_k}\>
{1-\eta\over1-\eta^s}\,\susi\,\Bigl[\,\pron\>\prod_{a\in\Gm_{\Rph lm}}
\Bigl(\,{t_a\over t_a-\xi_mz_m}
\ \prlm\,{\xi_l\)t_a-z_l\over t_a-\xi_lz_l}\,\Bigr)\>\Bigr]_\si
\kern-1.5em
\Tagg{wlga}
$$
\fi
where $\Gm_m=\lb 1+\lg^{m-1}\,\lc\lg^m\rb$, $\mn$.
The \fn/s $w_\lg$ are called the \em{\twf/s}.
\goodbm
\Lm{wlgp}
Let $\lg\in\Zln\!$. Then
$$
\align
w_\lg(\tell,\zn)\,=\!{}& \prab\ \,{t_a-t_b\over\eta\)t_a-t_b}\ \x
\\
\nn4>
\x\sum_{\Gmn}\biggl\lb\,\pron\>\prod_{a\in\Gms_{\Rph lm}}
& \>\Bigl(\,{t_a\over t_a-\xi_mz_m}
\ \prlm\,{\xi_l\)t_a-z_l\over t_a-\xi_lz_l}\,\Bigr)\!
\prod_{\tsize{1\le l<m\le n\atop\vru1.4ex>a\in\Gms_l,\>b\in\Gms_{\Rph lm}}}
\!{\eta\)t_a -t_b\over t_a-t_b}\,\biggr\rb
\endalign
$$
where the summation is over all \^{$\,n$-}tuples $\Gmn$ of disjoint subsets of
${\)\lb 1\lc\ell\)\rb}$ \st/ $\Gms_m$ has $\lg_m$ elements.
\endpro
\nt
The lemma is proved in Appendix~\SNo{A}.
\Ex
For $\ell=1$ the \twf/s have the form
$$
w_{\eg(m)}(t_1,\zn)\,=\,
{t_1\over t_1-\xi_mz_m}\ \prlm\,{\xi_l\)t_1-z_l\over t_1-\xi_lz_l}
$$
where $\eg(m)=(0\lc\%1_{\sss\^{$m$-th}}\lc 0)$, $\,\mn$.
\enddemo
\Ex
For $n=1$ the \fn/ $w\'\ell$ has the form
$$
w\'\ell(\tell,z_1)\,=\,
\pral\,{t_a\over t_a-\xi_1z_1}\;\prab\ {t_a-t_b\over\eta\)t_a-t_b}\;.
$$
\enddemo
\Ex
For $\ell=2$ and $n=2$ the \fn/s $w_\lg$ have the form
$$
\NN6>
\align
w_{(2,0)}(t_1,t_2,z_1,z_2)\, &{}=\;
{t_1t_2\over(t_1-\xi_1z_1)\>(t_2-\xi_1z_1)}\,{t_1-t_2\over\eta\)t_1-t_2}\;,
\\
w_{(1,1)}(t_1,t_2,z_1,z_2)\, &{}=\;
{t_1t_2\over(t_1-\xi_1z_1)\>(t_2-\xi_2z_2)}\,
{\xi_1t_2-z_1\over t_2-\xi_1z_1}\;+{}
\\
\nn-2>
&{}\)+\;{t_1t_2\over(t_2-\xi_1z_1)\>(t_1-\xi_2z_2)}\,
{\xi_1t_1-z_1\over t_1-\xi_1z_1}\,{t_1-\eta\)t_2\over\eta\)t_1-t_2}\;,
\\
w_{(0,2)}(t_1,t_2,z_1,z_2)\, &{}=\;
{t_1t_2\over(t_1-\xi_2z_2)\>(t_2-\xi_2z_2)}\,
{(\xi_1t_1-z_1)\>(\xi_1t_2-z_1)\over(t_1-\xi_1z_1)\>(t_2-\xi_1z_1)}\,
{t_1-t_2\over\eta\)t_1-t_2}\;.
\endalign
$$
\enddemo
\Lm{wbasis}
The \fn/s $w_\lg$, $\lg\in\Zln\!$, restricted to a fiber over $z$ form
a basis in the \thgf/ $\Fo(z)$ of the fiber provided that
\ ${\xi_l\>\xi_m z_m/z_l\netr}$, \ $1\le l<m\le n$ \ for any $\roll$.
\endpro
\nt
Lemma~\[wbasis] is proved in \Chap/~\SNo{4}.
\Par\nt
Let $\,\botsmash{\eg(m)=(0\lc\%1_{\sss\^{$m$-th}}\lc 0)}$, $\,\mn$.
\Lm{wDw}
Let $\ka=\kone$. Then for any $\lg\in\Zll$ the following relation holds:
\ifAst
$$
\align
\\
\nn-24>
\sun w_{\lg+\eg(m)}\>(1-\eta^{\lg_m+1})\> & (\xi_m-\eta^{\lg_m}\xi_m\1)\>
\tprllm\!\eta^{-\lg_l}\xi_l\,={}
\\
\nn1>
& {}\!=\,(1-\eta)\>
\sum_{a=1}^\ell\,D_a\bigl[\)w_\lg(\twll)\)\bigr]_{(1,a)}\,,
\endalign
$$
\else
$$
\sun w_{\lg+\eg(m)}\>(1-\eta^{\lg_m+1})\>(\xi_m-\eta^{\lg_m}\xi_m\1)\>
\tprllm\!\eta^{-\lg_l}\xi_l\,=\,(1-\eta)\>
\sum_{a=1}^\ell\,D_a\bigl[\)w_\lg(\twll)\)\bigr]_{(1,a)}\,,
$$
\fi
where ${(1,a)\in\Sl}$ are transpositions. Moreover, if $\Rc(z)$ is the
\mmgood
\aagood
subspace in $\Fo(z)$ \gby/ the elements in \lhs/ of the relations, then
$$
\dim\Fo(z)/\Rc(z)\,=\,{n+\ell-2\choose n-2}
$$
provided that \ ${\xi_l\>\xi_m z_m/z_l\netr}$, \ $1\le l\le m\le n$,
\ for any $\roll$.
\endpro
\Lm{wDw'}
Let $\ka=\ktwp$. Then for any $\lg\in\Zll$ the following relation holds:
\ifAst
$$
\align
\sun w_{\lg+\eg(m)}\>(1-\eta^{\lg_m+1})\> & (\xi_m-\eta^{\lg_m}\xi_m\1)\,
z_m\tprlm\!\eta^{\lg_l}\xi_l\1\>={}
\\
& {}\!=\,
(1-\eta)\>\sum_{a=1}^\ell\,D_a\bigl[\)t_1\)w_\lg(\twll)\)\bigr]_{(1,a)}\,,
\endalign
$$
\else
$$
\sun w_{\lg+\eg(m)}\>(1-\eta^{\lg_m+1})\>(\xi_m-\eta^{\lg_m}\xi_m\1)\,
z_m\tprlm\!\eta^{\lg_l}\xi_l\1\>=\,(1-\eta)\>
\sum_{a=1}^\ell\,D_a\bigl[\)t_1\)w_\lg(\twll)\)\bigr]_{(1,a)}\,,
$$
\fi
where ${(1,a)\in\Sl}$ are transpositions. Moreover, if $\Rc'(z)$ is the
subspace in $\Fo(z)$ \gby/ the elements in \lhs/ of the relations, then
$$
\dim\Fo(z)/\Rc'(z)\,=\,{n+\ell-2\choose n-2}
$$
provided that \ ${\xi_l\>\xi_m z_m/z_l\netr}$, \ $1\le l\le m\le n$,
\ for any $\roll$.
\endpro
\nt
Lemmas~\[wDw] and \[wDw'] are proved in \Chap/~\SNo{4}.
\Par
The subspaces $\Rc(z),\Rc'(z)\sub\Fo(z)$ are called the \em{\cosub/s}.
\Par
For $\ka=\kone$ relations \(wDw) induce the relations
\ifAst
$$
\Line{\sun\lf w_{\lg+\eg(m)}\)\Dtel\rf\>
(1-\eta^{\lg_m+1})\>(\xi_m-\eta^{\lg_m}\xi_m\1)\>
\tprllm\!\eta^{-\lg_l}\xi_l\,=\,0\,,\hfil\lg\in\Zll\,,}
$$
\else
$$
\sun\lf w_{\lg+\eg(m)}\)\Dtel\rf\>
(1-\eta^{\lg_m+1})\>(\xi_m-\eta^{\lg_m}\xi_m\1)\>
\tprllm\!\eta^{-\lg_l}\xi_l\,=\,0\,,\qqq\lg\in\Zll\,,
$$
\fi
in the cohomology group $H^\ell(z)$, where $\lf\al\rf$ denotes the
cohomological class of a form $\al$. Under assumptions of Theorem~\[kaon]
we have
$$
\H(z)\,\simeq\,\Fo(z)/\Rc(z)\,.
\Tag{HFR}
$$
Similarly, for $\ka=\ktwp$ relations \(wDw') induce the relations
\ifAst
$$
\align
\sun\lf w_{\lg+\eg(m)}\)\Dtel\rf\>
(1-\eta^{\lg_m+1})\>(\xi_m-\eta^{\lg_m}\xi_m\1)\,
z_m\tprlm\!\eta^{\lg_l}\xi_l\1\>={} & \,0\,,
\\
\lg\, & \!{}\in\Zll\,,
\endalign
$$
\else
$$
\sun\lf w_{\lg+\eg(m)}\)\Dtel\rf\>
(1-\eta^{\lg_m+1})\>(\xi_m-\eta^{\lg_m}\xi_m\1)\,
z_m\tprlm\!\eta^{\lg_l}\xi_l\1\>=\,0\,,\qqq\lg\in\Zll\,,
$$
\fi
in the cohomology group $H^\ell(z)$, and under assumptions of Theorem~\[kaon]
we have
$$
\H(z)\,\simeq\,\Fo(z)/\Rc'(z)\,.
\Tag{HFR'}
$$
\par
For any \perm/ ${\tau\in\S^n}$ we define a basis
${\lb\)w^\tau_\lg\)\rb_{\lg\)\in\Zln}\:}$ in the \thgf/ of a fiber
by formulae similar to \(wlga). Namely, set
$$
\ifMag\kern2em\else\kern1em\fi
w^\tau_\lg(\tell,\zn;\xin)\,=\,w_{{\vp(}^\tau\!\lg}^{\vp1}
(\tell,z_{\tau_1}\lc z_{\tau_n};\xi_{\tau_1}\lc\xi_{\tau_n})
\ifMag\kern-2em\else\kern-1em\fi
\Tag{wtau}
$$
where $^\tau\}\lg=(\lg_{\tau_1}\lc\lg_{\tau_n})$.
\goodbreak
\Ex
For ${\ell=1}$ and a \perm/ ${\tau=(n,n-1\lc 1)}$ the \twf/s have the form
$$
w^\tau_{\eg(m)}(t_1,\zn)\,=\,
{t_1\over t_1-\xi_mz_m}\ \prml\,{\xi_l\)t_1-z_l\over t_1-\xi_lz_l}\;.
$$
\enddemo

\subsect{The \ehgf/}
In our study of the \GM/ for the discrete \loc/ \(conns) an important role
is played by the following \em{\ehgf/}. The \ehgf/ is an elliptic counterpart
of the \thgf/ introduced above.
\Par
All over this section we assume that $0<|p|<1$. Let
\ifMag
$$
\tht(u)\,=\,(u;p)\9(p\)u\1;p)\9(p;p)\9
$$
\else
\vv.1>
${\tht(u)=(u;p)\9(p\)u\1;p)\9(p;p)\9}$
\fi
be the Jacobi theta-\fn/.
\Par
The \ehgf/ $\Fq$ is the space of \fn/s of \var/s $\tell,\zn$ spanned over $\C$
by \fn/s which have the form
$$
Y(\zn)\>\Tho(\tell,\zn)\,\pron\,\pral\,{1\over\tht(\xi_m\1t_a/z_m)}\;
\prab\ {\tht(t_a/t_b)\over\tht(\eta\)t_a/t_b)}\;.
$$
Here $Y$ is a \mef/ on ${\Cxn\!}$ and $\Tho$ is a \hof/ on ${\Cxln\!}$,
\sym/ in the \var/s $\tell$; we assume that the functions $Y$ and $\Tho$ have
the properties
\vv-.5>
$$
\alignat2
\Tho &(\tpell,\zn)\,=\,(-t_a)^{-n}\,\ka\}\tpron\!z_m\;
\Tho(\tell,\zn)\,, &&
\\
\nn8>
Y &(\zmn)\>\Tho(\tell,\zmn)\;= &&
\\
\nn3>
&& \Llap{=\;(-p/z_m)^n\)\tpral t_a\;Y(\zn)\>\Tho(\tell,\zn)\,.} &
\endalignat
$$
The restriction of the \ehgf/ to a fiber defines the \em{\ehgf/} $\Fq(z)$ of
the fiber. The \ehgf/ $\Fq(z)$ is a complex \fd/ vector space of the same
dimension as the \thgf/ of the fiber, see Appendix~\SNo{B}.
\par
All elements $f(t,z)$ of the \ehgf/ satisfy the periodicity conditions
$$
\gather
f(\tpell,\zn)\,=\,\ka\>\eta^{\ell-2a+1}\zy\,f(\tell,\zn)\,,
\ifMag\kern-2.5em\fi
\Tag{1-01}
\\
\nn6>
f(\tell,\zmn)\,=\,\xi_m^\ell\>f(\tell,\zn)\,.
\ifMag\kern-2em\fi
\endgather
$$
(\cf. \(tzp)\)). Therefore, if $\Phi(t,z)$ is a short \phf/ given by \(Phi)
and $f(t,z)$ is any element of the \ehgf/, then $\Phi(t,z)f(t,z)$ is a \phf/
of our discrete \loc/.
\par
Transformation properties \(1-01) also mean that the \ehgf/s of fibers over
$z$ and $z'$ are naturally identified if the points $z$ and $z'$ lie in
the same orbit of the \^{$\Z^n\]$-}action on the base space.
\par
We give basic facts about the \ehgf/ in Appendix~\SNo{B}.
\Par
If required, we will write down explicitly the dependence of the \ehgf/s on
$\ka$, $\zn$, $\xin$ and $\ell$, that is
\ifAst
$$
\Line{\Fq=\Fq[\)\ka;\zn;\xin;\ell\,]\hfil\text{and}\hfil
\Fq(z)=\Fq[\)\ka;\zn;\xin;\ell\,](z)\,.}
$$
\else
$$
\Fq=\Fq[\)\ka;\zn;\xin;\ell\,]\qquad\text{and}\qquad
\Fq(z)=\Fq[\)\ka;\zn;\xin;\ell\,](z)\,.
$$
\fi
\par
Introduce a new action of the \symg/ $\Sl$ on \fn/s,
$$
\gather
f\,\map\,\lbc\)f\)\rbc\vpp\si\,,\qqq\si\in\Sl,
\Tag{ect}
\\
\nn4>
\lbc\)f\)\rbc\vpp\si(\tell,\zn)\,=\,f(\tsill,\zn)
\prabsi{\eta\)\tht(\eta\1t_{\si_b}/t_{\si_a})\]\over
\tht(\eta\)t_{\si_b}/t_{\si_a})}\;.
\endgather
$$
The \ehgf/ is \inv/ \wrt/ this action. The action commutes with the restriction
of \fn/s to a fiber.
\Par
The \ehgf/ of a fiber has $n!$ remarkable bases. The bases are labelled by
elements of the \symg/ $\S^n\!$. First we define the basis corresponding to the
unit element of the \symg/.
For any $\lg\in\Zln\!$ define a \fn/ $W_\lg(t,z)$ as follows:
\ifMag
$$
\gather
\Lline{\hp{(2.22)}\quad W_\lg(\tell,\zn)\;=}
\Tagg{Wlga}
\\
\nn4>
\Rline{{}=\>\prod_{k=1}^n\,\prod_{s=1}^{\lg_k}\,{\tht(\eta)\over\tht(\eta^s)}\,
\susi\,\LBc\,\pron\,\prod_{a\in\Gm_{\Rph lm}}
\Bigl(\,{\tht(\eta^{\)2a-\ell-1}\ka_m\1t_a/z_m)\over\tht(\xi_m\1t_a/z_m)}
\prlm\,{\tht(\xi_l\)t_a/z_l)\over\tht(\xi_l\1t_a/z_l)}\,\Bigr)\>\RBc_\si}
\endgather
$$
\else
$$
\align
\kern1.8em
& W_\lg(\tell,\zn)\;=
\Tagg{Wlga}
\\
\nn4>
& \kern1.5em
{}=\>\prod_{k=1}^n\,\prod_{s=1}^{\lg_k}\,{\tht(\eta)\over\tht(\eta^s)}\,
\susi\,\LBc\,\pron\,\prod_{a\in\Gm_{\Rph lm}}
\Bigl(\,{\tht(\eta^{\)2a-\ell-1}\ka_m\1t_a/z_m)\over\tht(\xi_m\1t_a/z_m)}
\prlm\,{\tht(\xi_l\)t_a/z_l)\over\tht(\xi_l\1t_a/z_l)}\,\Bigr)\>\RBc_\si
\kern-1.8em
\endalign
$$
\fi
\goodast\nt
where $\Gm_m=\lb 1+\lg^{m-1}\,\lc\lg^m\rb$ and
$\ka_m=\ka\!\prlm\!\xi_l\:\prml\xi_l\1\!$,\quad$\mn$.
The \fn/s $W_\lg$ are called the \em{\ewf/s}.
\Lm{Wlgp}
Let $\lg\in\Zln\!$. Then
\ifAst
$$
\gather
\Lline{W_\lg(\tell,\zn)\,=\!
\prab\ \,{\tht(t_a/t_b)\over\tht(\eta\)t_a/t_b)}\ \x}
\\
\nn4>
\Rline{\x\sum_{\Gmn}\biggl\lb\,\pron\>\prod_{a\in\Gms_{\Rph lm}}
\>\Bigl(\,{\tht(\ka_{\lg,m}\1\)t_a/z_m)\over\tht(\xi_m\1t_a/z_m)}
\ \prlm\,{\tht(\xi_l\)t_a/z_l)\over\tht(\xi_l\1t_a/z_l)}\,\Bigr)\!
\prod_{\tsize{1\le l<m\le n\atop\vru1.4ex>a\in\Gms_l,\>b\in\Gms_{\Rph lm}}}
\!{\tht(\eta\)t_a/t_b)\over\tht(t_a/t_b)}\,\biggr\rb}
\endgather
$$
\else
$$
\align
W_\lg(\tell,\zn)\,=\!{}& \prab\ \,{\tht(t_a/t_b)\over\tht(\eta\)t_a/t_b)}\ \x
\\
\nn4>
\x\sum_{\Gmn}\biggl\lb\,\pron\>\prod_{a\in\Gms_{\Rph lm}}
& \>\Bigl(\,{\tht(\ka_{\lg,m}\1\)t_a/z_m)\over\tht(\xi_m\1t_a/z_m)}
\ \prlm\,{\tht(\xi_l\)t_a/z_l)\over\tht(\xi_l\1t_a/z_l)}\,\Bigr)\!
\prod_{\tsize{1\le l<m\le n\atop\vru1.4ex>a\in\Gms_l,\>b\in\Gms_{\Rph lm}}}
\!{\tht(\eta\)t_a/t_b)\over\tht(t_a/t_b)}\,\biggr\rb
\endalign
$$
\fi
where ${\,\ka_{\lg,m}=\ka\!\prod_{1\le i<m}\!\}\eta^{-\lg_i}\xi_i\:\!
\prod_{m<i\le n}\!\}\eta^{\lg_i}\xi_i\1}$
\vvn.15>
and the summation is over all \^{$\,n$-}tuples $\Gmn$ of disjoint subsets of
${\)\lb 1\lc\ell\)\rb}$ \st/ $\Gms_m$ has $\lg_m$ elements.
\endpro
\nt
The lemma is proved in Appendix~\SNo{B}.
\Par
Let $Y_\lg(z)$ be any \mef/ \st/
$$
\gather
Y_\lg\)(\zmn)\,=\,\al_{\lg,m}\,Y_\lg\)(\zn)
\Tag{adjust}
\\
\nn8>
\Line{\rlap{where}\hfil\tsize\al_{\lg,m}\,=\,\ka^{\lg_m}\!\!\!
\prlm\!\eta^{-\lg_l\lg_m}\xi_l^{\lg_m}\xi_m^{\lg_l}\!
\prml\!\eta^{\lg_l\lg_m}\xi_l^{-\lg_m}\xi_m^{-\lg_l},\hfil}
\endgather
$$
$\mn$. Then the product $Y_\lg(z)W_\lg(t,z)$ is an element of the \ehgf/.
The function $Y_\lg$ will be called an \em{\adjf/} for the \wtf/ $W_\lg$.
The \adjf/s can be chosen to be \mef/s in parameters $\eta$, $\xin$ and $\ka$.
\Ex
Let $c_1\lc c_n$ be arbitrary nonzero complex numbers.
Let $\al_{\lg,1}\lc \al_{\lg,n}$ be the same as in \(adjust). Then the \fn/
$$
Y_\lg(\zn)\,=\,\pron\,{\tht(c_mz_m/\al_{\lg,m})\over\tht(c_mz_m)}
$$
is an \adjf/ for the \wtf/ $W_\lg$.
\enddemo
\nt
Notice that an \adjf/ is not unique. In what follows we never need to know the
\adjf/s explicitly.
\Ex
For $\ell=1$ the \ewf/s have the form
$$
W_{\eg(m)}(t_1,\zn)\,=\,{\tht(\ka_m\1t_1/z_m)\over\tht(\xi_m\1t_1/z_m)}
\,\prlm\,{\tht(\xi_l\)t_1/z_l)\over\tht(\xi_l\1t_1/z_l)}
$$
where $\ka_m=\ka\!\prlm\!\xi_l\:\prml\xi_l\1\!$,\quad $\mn$.
\enddemo
\Ex
For $n=1$ the \fn/ $W\'\ell$ has the form
$$
W\'\ell(\tell,z_1)\,=\,\pral\,
{\tht(\ka\1t_a/z_1)\over\tht(\xi_1\1t_a/z_1)}\;
\prab\ {\tht(t_a/t_b)\over\tht(\eta\)t_a/t_b)}\;.
$$
\enddemo
\Ex
For $\ell=2$ and $n=2$ the \fn/s $W_\lg$ have the form
\ifAst
$$
\NN8>
\gather
\Lline{W_{(2,0)}(t_1,t_2,z_1,z_2)\,=\;
{\tht(\ka\1\xi_2t_1/z_1)\>\tht(\ka\1\xi_2t_2/z_1)\over
\tht(\xi_1\1t_1/z_1)\>\tht(\xi_1\1t_2/z_1)}\,
{\tht(t_1/t_2)\over\tht(\eta\)t_1/t_2)}\;,}
\\
\ald
\Lline{\aligned\!
W_{(1,1)}(t_1,t_2,z_1,z_2)\, &{}=\;
{\tht(\eta\1\ka\1\xi_2t_1/z_1)\>\tht(\eta\)\ka\1\xi_1\1t_2/z_2)\>
\tht(\xi_1t_2/z_1)\over
\tht(\xi_1\1t_1/z_1)\>\tht(\xi_2\1t_2/z_2)\>\tht(\xi_1\1t_2/z_1)}\ +{}
\\
\nn-2>
&{}\>+\;{\tht(\eta\)\ka\1\xi_1\1t_1/z_2)\>\tht(\eta\1\ka\1\xi_2t_2/z_1)\>
\tht(\xi_1t_1/z_1)\over
\tht(\xi_1\1t_1/z_1)\>\tht(\xi_2\1t_1/z_2)\>\tht(\xi_1\1t_2/z_1)}\,
{\eta\)\tht(\eta\1t_1/t_2)\over\tht(\eta\)t_1/t_2)}\;,
\endaligned}
\\
\ald
\Lline{W_{(0,2)}(t_1,t_2,z_1,z_2)\,=\;
{\tht(\ka\1\xi_1\1t_1/z_2)\>\tht(\ka\1\xi_1\1t_2/z_2)\>
\tht(\xi_1t_1/z_1)\>\tht(\xi_1t_2/z_1)\over
\tht(\xi_2\1t_1/z_2)\>\tht(\xi_2\1t_2/z_2)\>
\tht(\xi_1\1t_1/z_1)\>\tht(\xi_1\1t_2/z_1)}\,
{\tht(t_1/t_2)\over\tht(\eta\)t_1/t_2)}\;.}
\endgather
$$
\else
$$
\NN8>
\align
W_{(2,0)}(t_1,t_2,z_1,z_2)\, &{}=\;
{\tht(\ka\1\xi_2t_1/z_1)\>\tht(\ka\1\xi_2t_2/z_1)\over
\tht(\xi_1\1t_1/z_1)\>\tht(\xi_1\1t_2/z_1)}\,
{\tht(t_1/t_2)\over\tht(\eta\)t_1/t_2)}\;,
\\
\ald
W_{(1,1)}(t_1,t_2,z_1,z_2)\, &{}=\;
{\tht(\eta\1\ka\1\xi_2t_1/z_1)\>\tht(\eta\)\ka\1\xi_1\1t_2/z_2)\>
\tht(\xi_1t_2/z_1)\over
\tht(\xi_1\1t_1/z_1)\>\tht(\xi_2\1t_2/z_2)\>\tht(\xi_1\1t_2/z_1)}\ +{}
\\
\nn-2>
&{}\>+\;{\tht(\eta\)\ka\1\xi_1\1t_1/z_2)\>\tht(\eta\1\ka\1\xi_2t_2/z_1)\>
\tht(\xi_1t_1/z_1)\over
\tht(\xi_1\1t_1/z_1)\>\tht(\xi_2\1t_1/z_2)\>\tht(\xi_1\1t_2/z_1)}\,
{\eta\)\tht(\eta\1t_1/t_2)\over\tht(\eta\)t_1/t_2)}\;,
\\
\ald
W_{(0,2)}(t_1,t_2,z_1,z_2)\, &{}=\;
{\tht(\ka\1\xi_1\1t_1/z_2)\>\tht(\ka\1\xi_1\1t_2/z_2)\>
\tht(\xi_1t_1/z_1)\>\tht(\xi_1t_2/z_1)\over
\tht(\xi_2\1t_1/z_2)\>\tht(\xi_2\1t_2/z_2)\>
\tht(\xi_1\1t_1/z_1)\>\tht(\xi_1\1t_2/z_1)}\,
{\tht(t_1/t_2)\over\tht(\eta\)t_1/t_2)}\;.
\endalign
$$
\fi
\enddemo
\Lm{Wbasis}
The \fn/s $W_\lg$, $\lg\in\Zln\!$, restricted to a fiber over $z$ form
a basis in the \ehgf/ $\Fq(z)$ of the fiber, provided that
\ $\xi_l\>\xi_m z_m/z_l\>\nepsr$, \ $1\le l<m\le n$, \ for any $\roll$,
$s\in\Z$, and
\ $\,\ka\tprllm\xi_l\:\!\prml\!\xi_l\1\!\nepsr$ \ for any $\mn-1$
and $\rll$, $s\in\Z$.
\endpro
\nt
Lemma~\[Wbasis] is proved in \Chap/~\SNo{4}.
\Par
Let $\Qc(z)$ be the space of \fn/s of the form
$$
\susi\,\Lbc\>W(\twll)\>\Rbc_\si\,,
$$
where ${W\in\Fq{[\)\eta\1\ka;\zn;\xin;\ell-1\)]}}(z)$.
Let $\Qc'(z)$ be the space of \fn/s of the form
$$
\susi\,\LBc\>W(t_1\lc t_{\ell-1})\;t_\ell\1\!
\pron\,{\tht(\xi_mt_\ell/z_m)\over\tht(\xi_m\1t_\ell/z_m)}\;\RBc_\si\,,
$$
where ${W\in\Fq[\)\eta\)\ka;\zn;\xin;\ell-1\)]}(z)$. In general, the spaces
$\Qc(z)$ and $\Qc'(z)$ are not subspaces of the \ehgf/ of the fiber $\Fq(z)$,
because their elements do not have the required quasiperiodicity properties.
However, if $\ka=\kone$, then the space $\Qc(z)$ is a subspace of $\Fq(z)$,
and if $\ka=\ktwp$, then the space $\Qc'(z)$ is a subspace of $\Fq(z)$.
The spaces $\Qc(z), \Qc'(z)$ are called the \em{\bosub/s}.
\Par
For any \fn/ $f(\tell)$ and a point ${\ts\}=(\ts_1\lc\ts_\ell)}$ we define
the multiple residue $\Res f(t)\vst{t=\ts}$ by the formula
$$
\Res f(t)\vst{t=\ts}\)=\,\Res\bigl(\,\ldots\,
\Res f(\tell)\vst{t_\ell=\ts_\ell}\ldots\,\bigr)\vst{t_1=\ts_1}\,.
\Tag{Rests}
$$
For any $\lg\in\Zln\!$ define the points $\xt\lg\,,\,\yt\lg\in\Cxl\!$
as follows:
\ifMag
$$
\NN2>
\align
\xt\lg\,=\,
(\eta^{1-\lg_1}\xi_1z_1,\>\eta^{\)2-\lg_1}\xi_1z_1\lc\xi_1z_1,
\>\eta^{1-\lg_2}\xi_2z_2\lc\xi_2z_2,\,\ldots\,,{} &
\Tagg{xtyt}
\\
\eta^{1-\lg_n}\xi_nz_n\lc{} & \xi_nz_n)\,,
\\
\nn4>
\yt\lg\,=\,(\eta^{\lg_1-1}\xi_1\1z_1,\>\eta^{\lg_1-2}\xi_1\1z_1\lc\xi_1\1\]z_1,
\>\eta^{\lg_2-1}\xi_2\1z_2\lc\xi_2\1\]z_2,\,\ldots\,,{} &
\\
\eta^{\lg_n-1}\xi_n\1z_n\lc{} &\xi_n\1z_n)\,.
\endalign
$$
\else
$$
\gather
\enspace
\xt\lg\,=\,
(\eta^{1-\lg_1}\xi_1z_1,\>\eta^{\)2-\lg_1}\xi_1z_1\lc\xi_1z_1,
\>\eta^{1-\lg_2}\xi_2z_2\lc\xi_2z_2,\,\ldots\,,
\)\eta^{1-\lg_n}\xi_nz_n\lc\xi_nz_n)\,,
\Tag{xtyt}
\\
\nn4>
\yt\lg\,=\,(\eta^{\lg_1-1}\xi_1\1z_1,\>\eta^{\lg_1-2}\xi_1\1z_1\lc\xi_1\1\]z_1,
\>\eta^{\lg_2-1}\xi_2\1z_2\lc\xi_2\1\]z_2,\,\ldots\,,
\)\eta^{\lg_n-1}\xi_n\1z_n\lc\xi_n\1z_n)\,.
\endgather
$$
\fi
\Ex
Let ${\ell=1}$ and ${\lg=(0\lc\%1_{\sss\^{$m$-th}}\lc 0)}$.
Then $\xt\lg=\xi_mz_m$ and $\yt\lg=\xi_m\1z_m$.
\enddemo
\Lm{ResQ}
Let ${\eta^r\!\neps}$ for any $\rell-1$, $s\in\Z$. Assume that
\vv.1>
\ $\xi_l\1\xi_mz_l\1z_m\nepsr$, \ $l,\mn$, $l\ne m$, \ for any $\roll$,
$s\in\Z$. Then ${\Res f(t)\vst{t=\xt\mg}\!=0}$ for any $f\in\Qc(z)$ and
$\mg\in\Zln\!$. Moreover, if $\,\ka=\kone\,$ and $\,\eta^\ell\!\neps\!$,
$s\in\Z$, then
$$
\Qc(z)\,=\,\lb\)f\in\Fq(z)\vert\Res f(t)\vst{t=\xt\mg}\!=0\quad\text{for any}
\ \ \mg\in\Zln\)\rb\,.
$$
\endpro
\Lm{ResQ'}
Let ${\eta^r\!\neps}$ for any $\rell-1$, $s\in\Z$. Assume that
\vv.1>
\ ${\xi_l\)\xi_m}\1z_l\1z_m\nepsr$, \ $l,\mn$, $l\ne m$, \ for any $\roll$,
$s\in\Z$. Then $f(\yt\mg)=0$ for any $f\in\Qc'(z)$ and $\mg\in\Zln\!$.
Moreover, if $\ka=\ktwp$ and $\,\eta^\ell\!\neps\!$, $s\in\Z$, then
$$
\Qc'(z)\,=\,\lb\)f\in\Fq(z)\vert f(\yt\mg)=0\quad\text{for any}
\ \ \mg\in\Zln\)\rb\,.
$$
\endpro
\nt
Lemmas~\[ResQ], \[ResQ'] are proved in Appendix~\SNo{B}.
\Ex
Let ${\ell=1}$. Then the spaces $\Qc(z)$ and $\Qc'(z)$ are \onedim/.
The space $\Qc(z)$ is the space of constant \fn/s in one \var/,
and the space $\Qc'(z)$ is spanned by the \fn/
${W(t_1)=\dsize t_1\1\pron\,{\tht(\xi_m t_1/z_m)\over\tht(\xi_m\1t_1/z_m)\}}}$.
\par
Let $\ka=\yz$. Then \fn/s of the \ehgf/ of a fiber $\Fq(z)$ are \p-periodic.
The space $\Fq(z)$ has a \onedim/ subspace of constant \fn/s which is the
\bosub/ $\Qc(z)$. The constant \fn/s are the only \fn/s in $\Fq(z)$ which are
regular in $\Cx\!$. The regular \ewf/ is $W\'{1,0\lc 0}=1$.
\vsk.3>
Let $\ka=p\1\!\zy$. Then the \ewf/ $W\'{0\lc 0,1}$ equals ${-\xi_n\1z_n\)W}$
and generates the \bosub/ ${\Qc'(z)\sub\Fq(z)}$. The \fn/ $W$ is the only \fn/
in $\Fq(z)$ which vanishes at all the points $\xi_1\1z_1\lc\xi_n\1z_n$.
\enddemo
\Lm{dimFQ}
Let $\ka=\kone$. Let ${\eta^r\!\neps}$ for any $\rtll$, $s\in\Z$. Then
$$
\dim\Fq(z)/\Qc(z)\,=\,{n+\ell-2\choose n-2}\,.
$$
Moreover, the equivalence classes of \fn/s $W_\lg$, ${\lg_1\]=0}$,
$\lg\in\Zln\!$, restricted to a fiber over $z$ form a basis in the space
$\Fq(z)/\Qc(z)$,
\vv.1>
provided that \ ${\xi_l\>\xi_m}z_m/z_l\>\nepsr$,
\ $1\le l<m\le n$, \ for any $\roll$, $s\in\Z$,
and $\prllm\!\xi_l^2\nepsr$, \ $\mn-1$, fro any $\rolt$, $s\in\Z$.
\endpro
\Lm{dimFQ'}
Let $\ka=\ktwp$. Let ${\eta^r\!\neps}$ for any $\rtll$, $s\in\Z$. Then
$$
\dim\Fq(z)/\Qc'(z)\,=\,{n+\ell-2\choose n-2}\,.
$$
Moreover, the equivalence classes of \fn/s $W_\lg$, ${\lg_n\]=0}$,
$\lg\in\Zln\!$, restricted to a fiber over $z$ form a basis in the space
$\Fq(z)/\Qc'(z)$,
\vv.1>
provided that \ ${\xi_l\>\xi_m}z_m/z_l\>\nepsr$,
\ $1\le l<m\le n$ for any $\roll$, $s\in\Z$, and
$\prml\!\xi_l^2\nepsr$, \ $\mn-1$, fro any $\rolt$, $s\in\Z$.
\endpro
\nt
Lemmas~\[dimFQ], \[dimFQ'] are proved in Appendix~\SNo{B}.
\Par
For any \perm/ ${\tau\in\S^n}$ we define a basis
${\lb\)W^\tau_\lg\)\rb_{\lg\)\in\Zln}\:}$ in the \ehgf/ of a fiber
by formulae similar to \(Wlga). Namely, set
$$
\ifMag\kern2em\else\kern1em\fi
W^\tau_\lg(\tell,\zn;\xin)\,=\,W_{{\vp(}^\tau\!\lg}^{\vp1}
(\tell,z_{\tau_1}\lc z_{\tau_n};\xi_{\tau_1}\lc\xi_{\tau_n})
\ifMag\kern-2em\else\kern-1em\fi
\Tagg{Wtau}
$$
where $^\tau\}\lg=(\lg_{\tau_1}\lc\lg_{\tau_n})$.
\Ex
For ${\ell=1}$ and a \perm/ ${\tau=(n,n-1\lc 1)}$ the \ewf/s have the form
$$
W^\tau_{\eg(m)}(t_1,\zn)\,=\,
{\tht(\tilde\ka_m\1t_1/z_m)\over\tht(\xi_m\1t_1/z_m)}
\,\prml\,{\tht(\xi_l\)t_1/z_l)\over\tht(\xi_l\1t_1/z_l)}
$$
where $\tilde\ka_m=\ka\!\prlm\!\xi_l\1\!\!\prml\xi_l\:$,\quad $\mn$.
\enddemo

\ifAst\newpage\fi
\Sect[2a]{\{\Rms/ and the \{\qKZc/}
\subsect*{Highest \wt/ \{\Umod/s}
Let $q$ be a nonzero complex number which is not a root of unity.
Consider the \qg/ $\Uu$ with generators $E,F,q^{\pm H}$ and relations:
\vv-1>
$$
\NN3>
\gather
q^H\qH=\>\qH q^H=\>1\,,
\\
q^H\]E=q\)E\>q^H\,,\qqq q^H\]F=\q F\>q^H\,,
\\
[E,F]\;=\,{q^{2H}-q^{-2H}\over q-\q}\;.
\endgather
$$
Let the coproduct $\Dl:\Uu\,\to\,\Uu\ox\Uu$ be given by
$$
\NN2>
\align
\Dl(q^H\])=q^H\]\ox q^H, &\qqq \Dl(\qH\])=\qH\]\ox\qH,
\\
\Dl(E)=E\ox\qH\]+ q^H\]\ox E\,, &\qqq
\Dl(F)=F\ox\qH\]+ q^H\]\ox F\,.
\endalign
$$
The coproduct defines a \Umod/ structure on a tensor product of \Umod/s.
\Par
For a \Umod/ $V$ let $V=\Plus_\la V_\la$ be its \wtd/.
Let $V^*=\Plus_\la V^*_\la$ be its restricted dual. Define a structure of
a \Umod/ on $V^*\!$ by
$$
\bra E\phi,x\ket=\bra\phi,Fx\ket\,,\qqq
\bra F\phi,x\ket=\bra\phi,Ex\ket\,,\qqq
\bra q^{\pm H}\]\phi,x\ket=\bra\phi,q^{\pm H}\]x\ket\,.
$$
This \Umod/ structure on $V^*\!$ will be called the \em{dual} module
structure. For any \Umod/s $V_1,V_2$, the tautological map
$V_1^*\ox V_2^*\to (V_1\ox V_2)^*$ is an \iso/ of \Umod/s.
\Par
Let $V$ be the \Umod/ with \hw/ $q^\La\!$. Let
\ifMag $\topsmash{V=\Plus_{l=0}^\8 V_{\La-l}}$ \else
${V=\Plus_{l=0}^\8 V_{\La-l}}$ \fi
\vv.1>
be its \wtd/. For any nonzero complex number $u$ define an operator
${u^{\La-H}\}\in{}}\End(V)$ by $u^{\La-H}v=u^l v$ for any $v\in V_{\La-l}$.
\ifMag\Par\else\par\fi
Let $\Vn$ be \Umod/s with \hw/s $\qLan\!$, \resp/. We have the \wtd/s
\ifMag\vv-.2>\else\vv-.5>\fi
$$
\Vox\,=\,\Plus_{\ell=0}^\8\,\Vl\qquad\text{and}\qquad
\Vax\,=\,\Plus_{\ell=0}^\8\,\Val
$$
\ifMag\vvv-.5>\else\vvv->\fi
where $()\elli$ denotes the eigenspace of $q^H\!$ with \eva/
$q^{\,\sun\!\!\La_m-\ell}\!$.
\Par
Let $\FVal\sub\Val$ be the image of the operator $F$. Let $\Vls\sub\Vl$ be
the kernel of the operator $E$. There is a natural pairing
$$
\Vls\ox\VFV\,\to\,\C\,.
\Tag{VVF}
$$
\par
Let $\zn$ be nonzero complex numbers. Set
\ifMag
$$
\NN-5>
\align
E_z\, &{}=\>\sun\,\qH\]\lox\)\%{z_m E}_{\^{$m$-th}}\>\lox\)q^H
\\
\Text{and}
\Fz\, &{}=\>\sun\,\qH\]\lox\)\%{z_m F}_{\^{$m$-th}}\>\lox\)q^H.
\endalign
$$
\else
\vv-.5>
$$
E_z\,=\>\sun\,\qH\]\lox\)\%{z_m E}_{\^{$m$-th}}\)\lox q^H\qquad\text{and}
\qquad\Fz\,=\>\sun\,\qH\]\lox\)\%{z_m F}_{\^{$m$-th}}\)\lox q^H.
$$
\vvv-.5>
\fi
Let $\FzVal\sub\Val$ be the image of the operator $\Fz$.
\vv.2>
Let $\Vlsz\alb\sub\Vl$ be the kernel of the operator $E_z$.
There is a natural pairing
$$
\Vlsz\ox\VFzV\,\to\,\C\,.
\Tag{VVFz}
$$
Let $\Vn$ be \Vmod/s, then pairings \(VVF) and \(VVFz) are nondegenerate
provided
$$
\pron\,\pros\,(1-q^{4\La_m-2s})\ne 0\,.
$$
\vv->
\vv->

\subsect{The \tri/ \{\Rm/}
Let $V_1,V_2$ be \Vmod/s for $\Uu$ with \hw/s $q^{\La_1}\!,\,q^{\La_2}$ and
\gv/s $v_1,v_2$, \resp/. Consider an \^{$\EVV$-}valued \mef/ $\RVx$ with the
following properties:
$$
\align
\kern1em
\RVx\>(F\ox\qH\]+\>q^H\]\ox F)\, {}& =\,
(F\ox q^H\]+\>\qH\]\ox F)\>\RVx\,,
\kern-1em
\Tagg{Rdef}
\\
\nn4>
\kern1em
\RVx\>(F\ox q^H\]+\>x\)\qH\]\ox F)\, {}& =
\,(F\ox\qH\]+\>x\)q^H\]\ox F)\>\RVx
\kern-1em
\endalign
$$
in $\EVV$ and
$$
\RVx\>v_1\ox v_2=v_1\ox v_2\,.
\Tag{Rvv}
$$
Such a \fn/ $\RVx$ exists and is uniquely determined. $\RVx$ is called
the \em{$\gsl\!$ \tri/ \Rm/} for the tensor product $\VV$.
\Par
The \tri/ \Rm/ $\RVx$ also satifies the following relations:
$$
\NN4>
\align
\kern1em
\RVx\>(E\ox\qH\]+\>q^H\]\ox E)\, {}& =\,
(E\ox q^H\]+\>\qH\]\ox E)\>\RVx\,,
\kern-1em
\Tagg{Rmore}
\\
\kern1em
\RVx\>(x\)E\ox q^H\]+\>\qH\]\ox E)\, {}& =
\,(x\)E\ox\qH\]+\>q^H\]\ox E)\>\RVx\,,
\kern-1em
\\
\RVx\,q^H\]\ox q^H {}& =\,q^H\]\ox q^H\RVx\,.
\endalign
$$
In particular, $\RVx$ respects the \wtd/ of $\VV$.
\par
$\RVx$ satisfies the inversion relation
$$
\PV\>\RVx\,=\,\bigl(R_{V_2V_1}\:(x\1)\bigr)\1\>\PV
\Tag{inv}
$$
where $\PV:\VV\to V_2\ox V_1$ is the \perm/ map.
\Par
Let $\VV=\Plus_{l=0}^\8 V\"l$ be the decomposition of the \Umod/ $\VV$ into the
direct sum of \irr/s, where the \irr/ module $V\"l\!$ is \gby/ a singular
vector of \wt/ $q^{\La_1+\La_2-l}\!$. Let $\Pi\"l$ be the projector onto
$V\"l\!$ along the other summands. Then we have
$$
\gather
\RVx,=\,\RVo\,\sum_{l=0}^\8\,\Pi\"l\cdot\prod_{s=0}^{l-1}\,
{x-q^{2s-2\La_1-2\La_2}\over x-q^{2\La_1+2\La_2-2s}}
\Tag{Rspec}
\\
\Text{where}
\nn-4>
\RVo\,=\,q^{2\La_1\La_2-2H\ox H}\smk^\8(q^2-1)^{2k}
\tprod_{s=1}^k(1-q^{2s})\1\>(\qH\]E\ox q^H\]F)^k\,.
\endgather
$$
\par
Let $V_1,\,V_2,\,V_3$ be \Vmod/s. The corresponding \Rms/ satisfy the \YB/:
$$
\RV\:(x/y)\>R_{V_1V_3}\:(x)\>R_{V_2V_3}\:(y)\,=\,
R_{V_2V_3}\:(y)\>R_{V_1V_3}\:(x)\>\RV\:(x/y)\,.
\Tag{YBq}
$$
All of the properties of $\RVx$ given above are well known (\cf. \Cite{T},
\Cite{D1}, \Cite{J}, \Cite{CP}\)).

\subsect{The \qlo/ $\{\Ugg$}
The \tri/ \Rm/ is connected with an action of the \qlo/ $\Ugg$ in a tensor
product of \Umod/s. The \qlo/ $\Ugg$ is a Hopf algebra which contains $\Uu$ as
a Hopf subalgebra. We give the necessary facts about $\Ugg$ in this section.
\Par
Let $q$ be a complex number, $q\ne\pm1$. The \em{\qlo/} $\Ugh$ is a unital
associative algebra with generators $L_{ij}\"{+0},\ L_{ji}\"{-0}\!$,
$1\le j\le i\le 2$, and $L_{ij}\"s\!$, $i,j=1,2$, $s=\pm 1,\pm 2,\ldots$,
subject to relations \(Ughr)\,\Cite{RS}, \Cite{DF}.
\par
Let $e_{ij}$, $i,j=1,2$, be the ${2{\x}2}$ matrix with the only nonzero entry
$1$ at the intersection of the \^{$i$-}th row and \^{$j$-}th column. Set
$$
\align
R(x)\,=\,{} &(xq-\q)\>(e_{11}\ox e_{11}+e_{22}\ox e_{22})\;+
\\
\nn3>
{} +\,{} &(x-1)\>(e_{11}\ox e_{22}+e_{22}\ox e_{11})\,+\,
x(q-\q)\>e_{12}\ox e_{21}\,+\,(q-\q)\>e_{21}\ox e_{12}\,.
\endalign
$$
Introduce the generating series
$L^{\pm}_{ij}(u)=L_{ij}\"{\pm0}+\sum_{s=1}^\8 L_{ij}\"{\pm s}u^{\pm s}\!$.
\aagood
The relations in $\Ugh$ have the form
\vvn-.7>
$$
\NN4>
\gather
L_{ii}\"{+0}L_{ii}\"{-0}\,=\,1\,,\qquad L_{ii}\"{-0}L_{ii}\"{+0}\,=\,1\,,
\Rlap{\qqq i=1,2\,,}
\Tag{Ughr}
\\
{\align
& R(x/y)\>L^+\'1(x)\>L^+\'2(y)\,=\,
L^+\'2(y)\>L^+\'1(x)\>R(x/y)\,,
\\
& R(x/y)\>L^+\'1(x)\>L^-\'2(y)\,=\,
L^-\'2(y)\>L^+\'1(x)\>R(x/y)\,,
\\
& R(x/y)\>L^-\'1(x)\>L^-\'2(y)\,=\,
L^-\'2(y)\>L^-\'1(x)\>R(x/y)\,,
\endalign}
\endgather
$$
where $L^\nu\'1(u)=\sum_{ij} e_{ij}\ox 1\ox L_{ij}^\nu(u)$ and
$L^\nu\'2(u)=\sum_{ij}1\ox e_{ij}\ox L_{ij}^\nu(u)$,\quad $\nu=\pm$.
\par
Elements $\,L_{11}\"{+0}L_{22}\"{+0}\!$, $L_{22}\"{+0}L_{11}\"{+0}\!$,
$L_{11}\"{-0}L_{22}\"{-0}\!$, $L_{22}\"{-0}L_{11}\"{-0}$ are central in $\Ugh$.
Impose the following relations:
$$
L_{11}\"{+0}L_{22}\"{+0}\>=\>1\,,\qquad
L_{22}\"{+0}L_{11}\"{+0}\>=\>1\,,\qquad
L_{11}\"{-0}L_{22}\"{-0}\>=\>1\,,\qquad
L_{22}\"{-0}L_{11}\"{-0}\>=\>1\,,
$$
in addition to relations \(Ughr). Denote the corresponding quotient algebra
by $\Ugg$.
\Par
\ifMag\nt\fi
The \qlo/ $\Ugg$ is a Hopf algebra with a coproduct $\Dl:\Ugg\to\Ugg\ox\Ugg$:
$$
\Dl:L^\nu_{ij}(u)\,\map\,\tsum_k L^\nu_{ik}(u)\ox L^\nu_{kj}(u)\,,\qqq
\nu=\pm\,.
$$
There is an important one-parametric family of \aut/s
${\rho_x:\Ugg\to{}}\alb\Ugg$:
$$
\gather
\rho_x:L_{ij}^\nu(u)\,\map\,L_{ij}^\nu(u/x)\,,\qqq \nu=\pm\,,
\\
\nn-2>
\Text{that is}
\nn-2>
\rho_x:L_{ij}\"{\pm0}\,\map\,L_{ij}\"{\pm0}\qquad\text{and}\qquad
\rho_x:L_{ij}\"{s}\map x^{-s}L_{ij}\"{s}\,,\qquad s\in\Z_{\ne0}\,.
\endgather
$$
The \qlo/ $\Ugg$ contains $\Uu$ as a Hopf subalgebra; the embedding is given
by
$$
E\>\map\>-L_{21}\"{+0}\big/(q-\q)\,,\qqq
F\>\map\>L_{12}\"{-0}\big/(q-\q)\,,\qqq q^H\>\map\>L_{11}\"{-0}\,.
$$
There is also an \em{evaluation \hom/} $\epe:\Ugg\to\Uu$:
$$
\gather
\NN4>
\\
\ifMag\nn-18>\else\nn-20>\fi
\alignedat2
& \epe\):\)L^+_{11}(u)\,\map\,\qH-q^Hu\,, &&
\epe\):\)L^+_{12}(u)\,\map\,-F\)(q-\q)\>u\,,
\\
& \epe\):\)L^+_{21}(u)\,\map\,-E\)(q-\q)\,, &&
\epe\):\)L^+_{22}(u)\,\map\,q^H-\qH u\,,
\\
\ald
\nn4>
& \epe\):\)L^-_{11}(u)\,\map\,q^H-\qH u\1\,, &&
\epe\):\)L^-_{12}(u)\,\map\,F\)(q-\q),
\\
& \epe\):\)L^-_{21}(u)\,\map\,E\)(q-\q)\>u\1\,,\qqq
\ifMag\quad\else\qqq\fi &&
\epe\):\)L^-_{22}(u)\,\map\,\qH-q^Hu\1,
\endalignedat
\\
\ifMag\nn2>\else\nn1>\fi
\Text{that is}
\ifMag\nn3>\else\nn2>\fi
\alignedat4
& \epe\):\)L_{11}\"{+0}\,\map\,\qH,\qquad &&
\epe\):\)L_{11}\"1\,\map\,-q^H, &&
\qqq\Rlap{\epe\):\)L_{12}\"1\,\map\,-F\)(q-\q)\,,} &&
\\
&\qqq \Rlap{\epe\):\)L_{21}\"{+0}\,\map\,-E\)(q-\q)\,,} &&&&
\epe\):\)L_{22}\"{+0}\,\map\,q^H, &&
\epe\):\)L_{22}\"1\,\map\,-\qH,
\\
\ald
\nn4>
& \epe\):\)L_{11}\"{-0}\,\map\,q^H, &&
\epe\):\)L_{11}\"{-1}\,\map\,-\qH, \ifMag\ \else\qqq\fi &&
\qqq\Rlap{\epe\):\)L_{12}\"{-0}\,\map\,F\)(q-\q)\,,} &&
\\
&\qqq \Rlap{\epe\):\)L_{21}\"{-1}\,\map\,E\)(q-\q)\,,} &&&&
\epe\):\)L_{22}\"{-0}\,\map\,\qH,\qquad &&
\epe\):\)L_{22}\"{+1}\,\map\,-q^H
\endalignedat
\endgather
$$
and ${\epe\):\)L_{ij}\"s\,\map\,0}$ for all other generators $L_{ij}\"s$.
\par
Both the \aut/s $\rho_x$ and $\epe$ restricted to the subalgebra $\Uu$
are the identity maps.
\par
For any \Umod/ $V$ denote by $V(x)$ the \Uhmod/ which is obtained from the
module $V$ via the \hom/ ${\epe\o\rho_x}$. The module $V(x)$ is called
the \em{\emod/}.
\par
Let $V_1,V_2$ be \Vmod/s for $\Uu$ with \gv/s $v_1,v_2$, \resp/. For generic
complex numbers $x,y$ the \Uhmod/s $V_1(x)\ox V_2(y)$ and
$V_2(y)\ox V_1(x)$ are isomorphic and the \tri/ \Rm/ $\PV\RV\:(x/y)$
intertwines them \Cite{T}, \Cite{CP}. The vectors $v_1\ox v_2$ and
$v_2\ox v_1$ are respective \gv/s of the \Uhmod/s $V_1(x)\ox V_2(y)$ and
$V_2(y)\ox V_1(x)$. The \tri/ \Rm/ $\RV\:(x/y)$ can be defined as
the unique element of $\EVV$ with property \(Rvv) and \st/
$$
\PV\RV\:(x/y)\>:\>V_1(x)\ox V_2(y)\,\to\,V_2(y)\ox V_1(x)
\Tag{Uintw}
$$
is an \iso/ of the \Uhmod/s.

\subsect{The \tri/ \{\qKZc/ associated with $\{\gsl$}
Let $\Vn$ be \Umod/s. The \qKZc/ is a discrete connection on the \trib/ over
$\Cxn\!$ with fiber $\Vox$. We define it below.
\Par
Let $\Vn$ be \Vmod/s with \hw/s $\qLan\!$, \resp/. Let $R_{V_iV_j}\:(x)$ be
the \tri/ \Rms/. Let $R_{ij}(x)\in\EV$ be defined in a standard way:
$$
R_{ij}(x)\,=\>\sum\,\one\lox\%{r(x)}_{\^{$i$-th}\,}\lox\%{r'(x)}_{\^{$j$-th}\,}
\lox\one
\Tag{Rij}
$$
provided that $R_{V_iV_j}\:(x)=\sum r(x)\ox r'(x)\in\End(V_i\ox V_j)$.
For any $X\in\Uu$ set
$$
X_m\,=\,\one\lox\%{X}_{\sss\^{$m$-th}}\lox\one\,.
$$
\goodast\nt
Let $p,\Ks$ be complex numbers. For any $\mn$ set
$$
\align
\kern1em
K_m(\zn)\,=\, R_{m,m-1}(p\)z_m/z_{m-1})\ldots R_{m,1}(p\)z_m/z_1)\>
\Ks^{\La_m-H_m}_{\vp1}\;\x &
\Tagg{Kmz}
\\
\nn4>
\x\; R_{m,n}(z_m/z_n)\ldots R_{m,m+1}(z_m/z_{m+1}) & \,.
\kern-1em
\endalign
$$
\Th{FR}
\back\Cite{FR}
The linear maps $K_1(z)\lc K_n(z)$ obey the flatness conditions
$$
K_l(\zmn)\>K_m(\zn)\,=\,K_m(z_1\lc p\)z_l\lc z_n)\>K_l(\zn)\,,
$$
$l,\mn$.
\endpro
The maps $K_1(z)\lc K_n(z)$ define a flat connection on the \trib/
over $\Cxn\!$ with fiber $\Vox$. This connection is called the \em{\qKZc/}.
\Par
By \(Rmore) the operators $K_1(z)\lc K_n(z)$ commute with the action of
$q^H$ in $\Vox$:
$$
[K_m(\zn)\>,\)q^H\)]=0\,,\qqq\mn\,,
$$
and, therefore, preserve the \wtd/ of $\Vox$. Hence, the \qKZc/ induces
the dual flat connection on the \trib/ over $\Cxn$ with fiber $\Vax$.
This connection will be called the \em{dual \qKZc/}.
\Par
\Lm*
For any $z\in\Cxn\!$ \st/ $q^{2\La_l+2\La_m-2r}z_l/z_m\neps$ for all
\vvn.1>
$\roll$, and $l,\mn$, $l\ne m$, $s\in\Z$, the linear maps
$K^*_1(z)\lc K^*_n(z)$ define \iso/s of $\Val$.
\endpro
\nt
This statement follows from formulae \(Rspec) and \(Kmz).
\par
If ${\Ks=\Kone\!}$,
\vv.1>
then the dual \qKZc/ admits a trivial discrete subbundle with fiber $\FVal$
and,
\vvn.1>
therefore, it induces a flat connection on the \trib/ with fiber
$\Val\big/\FVal$.
\par
If ${\Ks=\Ktwp\!}$,
\vv.1>
then the dual \qKZc/ admits a discrete subbundle with fiber $\FzVal$ and,
\vvn.1>
therefore, it induces a flat connection on the discrete vector bundle with
fiber $\Val\big/\FzVal$.
\Par
Let $\Vn$ be \Umod/s. The \em{\qKZe/} for a \^{$\Vox$-}valued \fn/
$\Psi(\zn)$ is the following system of \eq/s:
$$
\Psi(\zmn)\,=\,K_m(\zn)\>\Psi(\zn)\,, \qqq \mn\,.
$$
The \qKZe/ is a remarkable \deq/, see \Cite{S}, \Cite{FR}, \Cite{JM},
\Cite{Lu}.

\goodast
\subsect{Modules over the \eqg/ $\{\Eqg$ and the elliptic
\ifAst\{\nl\fi \ifMag\{\nlc\fi \{\Rms/}
\righthead{Modules over $\{\Eqg$ and the elliptic \{\Rms/}
In this section we recall the definitions concerning the \eqg/ $\Eqg$, \Emod/s
and the \Rms/ associated with tensor products of \Emod/s. For a more detailed
exposition on the subject and proofs see \Cite{F}, \Cite{FV}.
\Par
Fix two complex numbers $\rho, \gm$ \st/ $\Im\rho>0$.
Set $p=e^{2\pii\rho}$ and $\eta=e^{-4\pii\gm}$.
\ifMag\else\nl\fi
Let $\,\tht(u)=(u;p)\9\>(p\)u\1;p)\9\>(p;p)\9$
be the Jacobi theta-\fn/ and
$$
\al(x,\la)\,=\;{\eta\>\tht(x)\>\tht(\la/\eta)\over\tht(\eta\)x)\>\tht(\la)}\;,
\qqq \bt(x,\la)\,=\;{\tht(\eta)\>\tht(x\la)\over\tht(\eta\)x)\tht(\la)}\;.
$$
Let $e_{ij}$, $i,j=1,2$, be the ${2{\x}2}$ matrix with the only nonzero entry
$1$ at the intersection of the \^{$i$-}th row and \^{$j$-}th column. Set
$$
\align
R(x,\la)\,&{}=\,e_{11}\ox e_{11}+\>e_{22}\ox e_{22}\)+\)
\al(x,\la)\>e_{11}\ox e_{22}\,\,+
\\
\nn4>
&{}+\,\al(x,\la\1)\>e_{22}\ox e_{11}+\)
\bt(x,\la)\>e_{12}\ox e_{21}+\)\bt(x,\la\1)\>e_{21}\ox e_{12}\,.
\endalign
$$
\Rem
In \Cite{FV} the \eqg/ $\Eqg$ is described in terms of the additive theta-\fn/
$$
\Tht(u)\,=\,-\!\sum_{m=-\8}^\8\!
\exp\bigl(\pii(m+{1/2})^2\rho+2\pii(m+1/2)(u+1/2)\bigr)\,,
$$
which is related to the multiplicative theta-\fn/ $\tht(x)$ by the equality
$$
\Tht(u)\,=\,i\)\exp(\pii\rho/4-\pii u)\>\tht(e^{2\pii u})\,.
$$
\enddemo
Let $\hg$ be the \onedim/ Lie algebra with the generator $H$. Let $V$ be
an \hmod/. Say that $V$ is \em {\diag/} if $V$ is a direct sum of \fd/
eigenspaces of $H$:
$$
V\,=\,\Plus_\mu\;V_\mu\,,\qquad\qquad Hv=\mu\)v\qquad\for v\in V_\mu\,.
\mmgood
$$
For a \fn/ $X(\mu)$ taking values in $\End(V)$ we set $X(H)v=X(\mu)v$ for
any $v\in V_\mu$.
\par
Let $\Vn$ be \dhmod/s. We have the decomposition
$$
\Vox\,=\Plus_{\mun}(V_1)_{\mu_1}\lox(V_n)_{\mu_n}\,.
$$
Set $H_m=\one\lox\%{H}_{\sss\^{$m$-th}}\lox\one$.
\vv.1>
For a \fn/ $X(\mun)$ taking values in $\End(\Vox)$ we set
$X(H_1\lc H_n)v=X(\mun)v$ for any $v\in(V_1)_{\mu_1}\lox(V_n)_{\mu_n}$.
\Par
By definition \Cite{FV}, a \em{module over the \eqg/ $\Eqg$} is a \dhmod/ $V$
together with four \^{$\End(V)$-}valued \fn/s $T_{ij}(u,\la)$, $i,j=1,2$, which
are \mer/ in $u,\la\in\Cx$ and obey the following relations
\ifAst
$$
\gather
[\)T_{ij}(u,\la)\),H\)]\,=\,(j-i)\>H\,,\qqq i,j=1,2\,, \kern-1em
\\
\nn10>
\aligned
R(x/y,\eta^{\)1\ox 1\ox 2H\]}\la)\> & T\'1(x,\la)\>
T\'2(y,\eta^{\)2H\ox 1\ox 1\]}\la)\,={}
\\
\nn4>
{}={} &\,T\'2(y,\la)\>T\'1(x,\eta^{\)1\ox 2H\ox 1\]}\la)\>R(x/y,\la)\,.
\endaligned
\endgather
$$
\else
$$
\gather
[\)T_{ij}(u,\la)\),H\)]\,=\,(j-i)\>H\,,\qqq i,j=1,2\,, \kern-1em
\\
\nn10>
R(x/y,\eta^{\)1\ox 1\ox 2H\]}\la)\>T\'1(x,\la)\>
T\'2(y,\eta^{\)2H\ox 1\ox 1\]}\la)\,=\,
T\'2(y,\la)\>T\'1(x,\eta^{\)1\ox 2H\ox 1\]}\la)\>R(x/y,\la)\,.
\endgather
$$
\fi
Here \ $T\'1(u,\la)=\sum_{ij}^{\vp1}e_{ij}\ox\one\ox T_{ij}(u,\la)$,
\ $T\'2(u,\la)=\sum_{ij}\one\ox e_{ij}\ox T_{ij}(u,\la)$,
\ifMag\else\vvn-.1>\fi
\ and $H$ acts in $\C^2\!$ as $(e_{11}-e_{22})/2$.
\Ex
Fix a complex number $\La$. Consider a \dhmod/
$V^\La={\!\Plus_{\smash{\,\,k\in\Zp\!\!}}\!\C\)\ve[k]}$, \st/
$$
H\ve[k]\,=\,(\La-k)\>\ve[k]\,.
$$
Let $x$ be a nonzero complex number. Set
$$
\NN8>
\align
T_{11}(u,\la)\>\ve[k]\, &{}=\;{\tht(\eta^{\)\La-k}u/x)\>\tht(\eta^{-k}\la)
\over\tht(\eta^{\)\La}u/x)\>\tht(\la)}\;\eta^k \ve[k]\,,
\\
T_{12}(u,\la)\>\ve[k]\, &{}=\;{\tht(\eta^{\)\La-k-1}\la\)u/x)\>\tht(\eta)
\over\tht(\eta^{\)\La}u/x)\>\tht(\la)}\;\ve[k+1]\,,
\\
T_{21}(u,\la)\>\ve[k]\, &{}=\;
{u\>\tht(\eta^{\)\La-k+1}\la\)x/u)\>\tht(\eta^{\)2\La-k+1})\>\tht(\eta^k)\over
x\>\tht(\eta^{\)\La}u/x)\>\tht(\la)\>\tht(\eta)}\;\eta^{k-1-\La}\ve[k-1]\,,
\\
T_{22}(u,\la)\>\ve[k]\, &{}=\;
{\tht(\eta^{k-\La}\)u/x)\>\tht(\eta^{\)2\La-k}\la)\over
\tht(\eta^{\)\La}u/x)\>\tht(\la)}\;\ve[k]\,.
\endalign
$$
These formulae make $V^\La$ into an \Emod/ $V^\La(x)$, which is called the
\em{\eVmod/} with \em{\hw/} $\La$ and \em{\epoint/} $x$ \Cite{FV}.
\enddemo
For any complex vector space $V$ denote by $\Fun(V)$ the space of \Vval/ \mef/s
on $\Cx\!$. The space $V$ is naturally embedded in $\Fun(V)$ as the subspace of
constant \fn/s.
\par
Let $V_1,V_2$ be complex vector spaces. Any \fn/
$\phi\in\Fun\bigl(\Hom(V_1,V_2)\bigr)$ induces a linear map
$$
\fun\phi\):\)\Fun(V_1)\,\to\,\Fun(V_2)\,,\qqq
\fun\phi\):\)f(\la)\,\map\,\phi(\la)\)f(\la)\,.
$$
\par
For any \Emod/ $V$ we define the associated \em{\oalg/} acting on the space
$\Fun(V)$. The \oalg/ is \gby/ \mef/s in $\la,\eta^H\!$ acting pointwise,
$$
\phi(\la,\eta^H)\):\)f(\la)\,\map\,\phi(\la,\eta^H)f(\la)\,,
\Tag{1-02}
$$
and by values and residues \wrt/ $x$ of the operator-valued \mef/s
$\Tti_{ij}(x)$, $i,j=1,2$, defined below:
$$
\Tti_{i1}(x)\):\)f(\la)\,\map\,T_{i1}(x,\la)\)f(\eta\la)\,,\qqq
\Tti_{i2}(x)\):\)f(\la)\,\map\,T_{i2}(x,\la)\)f(\eta\1\la)\,.
$$
The relations obeyed by the generators of the \oalg/ are described in detail
in \Cite{FV}.
\Par
Let $V_1,V_2$ be \Emod/s. An element $\phi\in\Fun\bigl(\Hom(V_1,V_2)\bigr)$
\st/ the induced map $\fun\phi$ intertwines the actions of the respective
\oalg/s is called a \em{morphism} of \Emod/s $V_1,V_2$. A morphism $\phi$ is
called an \em{\iso/} if the linear map $\phi(\la)$ is \ndeg/ for generic $\la$.
\Ex
Evaluation \Vmod/s $V^\La(x)$ and $V^\Mu(x)$ are isomorphic if
$\eta^{\)\La}=\eta^{\)\Mu}\!$ with the tautological \iso/.
\enddemo
The \eqg/ $\Eqg$ has the coproduct $\Dlq$:
\ifMag
$$
\gather
\Dlq\!:H\,\map\,H\ox 1 + 1\ox H\,,
\\
\nn6>
\Dlq\!:T_{ij}(u,\la)\,\map\,
\tsum_k \bigl(1\ox T_{ik}(u,\eta^{\)2H\ox 1}\la)\bigr)\)
\bigl(T_{kj}(u,\la)\ox 1\bigr)\,.
\endgather
$$
\else
$$
\Dlq\!:H\,\map\,H\ox 1 + 1\ox H\,,\qquad\qquad
\Dlq\!:T_{ij}(u,\la)\,\map\,
\tsum_k \bigl(1\ox T_{ik}(u,\eta^{\)2H\ox 1}\la)\bigr)\)
\bigl(T_{kj}(u,\la)\ox 1\bigr)\,.
$$
\fi
The precise meaning of the coproduct is that it defines an \Emod/ structure on
the tensor product $\VV$ of \Emod/s $V_1,V_2$. If $V_1,V_2,V_3$ are \Emod/s,
then the modules $(V_1\ox V_2)\ox V_3$ and $V_1\ox(V_2\ox V_3)$ are naturally
isomorphic \Cite{F}.
\Rem
Notice that we take the coproduct $\Dlq\!$ which is opposite to the
coproduct used in \Cite{F}, \Cite{FV}, \cite{FTV1}. The coproduct $\Dlq$ is
in a sence opposite to the coproduct $\Dl$ taken for the \qlo/ $\Ugg$.
\enddemo
\Th{Eintw}
\back\Cite{FV},\;\Cite{FTV1}
Let $V^\La(x), V^\Mu(y)$ be \eVmod/s. Then for any $\La,\Mu$ and generic $x,y$
there is a unique \iso/ $\Rti_{(x,y)}$ of \,
\ifMag $\Eqg$ modules \else\Emod/s \fi
$V^\La(x)\ox V^\Mu(y)$, $V^\Mu(y)\ox V^\La(x)$ \st/
$\Rti_{(x,y)}(\la)\ve[0]\ox\ve[0]=\ve[0]\ox\ve[0]$. Moreover, $\Rti_{(x,y)}$
has the form
$$
\Rti_{(x,y)}(\la)\,=\,P_\VLM\:\>\Rq_\VLM(x/y,\la)
$$
where $P_\VLM\::V^\La\ox V^\Mu\to V^\Mu\ox V^\La$ is the \perm/ map and
$\Rq_\VLM(u,\la)$ is a \mef/ of $u,\la\in\Cx\!$ with values in
$\End(V^\La\ox V^\Mu)$.
\endpro
\Cr*
The \fn/ $\Rq_\VLM(x,\la)$ satisfies the inversion relation
$$
P_\VLM\:\>\Rq_\VLM(x,\la)\,=\,\bigl(\Rq_\VML(x\1,\la)\bigr)\1\>P_\VLM\:\,.
$$
\endpro
The \fn/ $\Rq_\VLM(x,\la)$ is called the $\gsl\!$ \em{dynamical elliptic \Rm/}
for the tensor product $V^\La\ox V^\Mu\!$.
\Th{YBe}
\back\Cite{FV},\;\Cite{FTV1}
For any complex numbers $\La,\Mu,\Nu$ the corresponding elliptic \Rms/ satisfy
the dynamical \YB/ in the space\ifAst\break\fi{}
$\End(V^\La\ox V^\Mu \ox V^\Nu)${\rm:}
$$
\align
\Rq_\VLM(x/y &,\eta^{\)1\ox 1\ox 2H\]}\la)\>
\Rq_\VLN(x,\la)\>\Rq_\VMN(y,\eta^{\)2H\ox 1\ox 1\]}\la)\,={}
\\
\nn6>
&{}\>=\,\Rq_\VMN(y,\la)\>\Rq_\VLN(x,\eta^{\)1\ox 2H\ox 1\]}\la)\>
\Rq_\VLM(x/y,\la)\,.
\endalign
$$
\endpro
One can associate a discrete flat connection with the dynamical elliptic \Rm/.
This connection is studied in \Cite{FTV1}, \Cite{FTV2}.

\Sect[2b]{Tensor coordinates on the \hgeom/ spaces}
In this \chap/ we identify the \GM/ and the \qKZc/.

\subsect{Tensor coordinates on the \thgf/s of fibers}
\righthead{Tensor coordinates on the \thgf/s}
Let $\Vn$ be $\Uu$ Verma modules with \hw/s $\qLan\!$ and \gv/s $\vn$, \resp/.
Consider the \wt/ subspace $\Vl$ with a basis given by monomials $\Fv$,
$\lg\in\Zln$. The dual space $\Val$ has the dual basis denoted by $\Fva$,
$\lg\in\Zln$.
\par
Consider the \tri/ \stype/ \loc/ with \conn/s \(conns) where the parameters
$\xin$ and $\eta$ are related to the parameter $q$ and the \hw/s $\qLan\!$
as follows:
$$
\eta=q^2\,,\qquad\qquad\xi_m=q^{2\La_m},\qqq\mn\,.
$$
Let ${\Fo=\Fo[\zn;\xin;\ell\,]}$ be the corresponding \tri/ hypergeomet\-ric
%% \thgf/
space. For any $z\in\Cxn\!$ and for any $\tau\in\S^n$ denote by $B_\tau(z)$
the following linear map:
\ifMag\vv-.5>\fi
$$
\gather
B_\tau(z):\Valt\,\to\,\Fo(z)\,,
\Tag{Btau}
\\
\nn4>
B_\tau(z):\Fvat\,\map\,b_\lg\>w^\tau_\lg(t,z)\,,
\endgather
$$
where $\Fo(z)$ is the \thgf/ of the fiber and
\ifMag
$$
b_\lg\,=\tpron\) q^{\lg_m(\lg_m-1)/2\,+\)\lg_m\La_m},
$$
\else
${b_\lg\)=\!\}\pron q^{\lg_m(\lg_m-1)/2\,+\)\lg_m\La_m}\!}$,
\vvn.1>
\fi
(\cf. \(wlga), \(wtau)\)). The linear maps $B_\tau(z)$ are called
the \em{\tenco/} on the \thgf/ of a fiber. The composition maps
$$
B_{\tau,\tau'}(z):\Valtt\,\to\,\Valt\,,\ifMag\else\quad\fi\qquad
B_{\tau,\tau'}(z)\,=\,B\1_\tau(z)\o B_{\tau'}(z)\,,
$$
are called the \em{\traf/s}, \cf. \Cite{V3}, \Cite{TV3}.
\Lm{tcfs}
Let ${q^{2\La_l+2\La_m-2r}z_l/z_m\ne 1}$ for any $\roll$, and $l,\mn$,
$l\ne m$. Then for any \perm/ $\tau$ the linear map
$B_\tau(z):\Valt\to\Fo(z)$ is \ndeg/.
\endpro
\nt
The statement follows from Lemma~\[wbasis] since
${\xi_l\>\xi_m z_l/z_m\]\netr}$ for any $\roll$, and $l,\mn$, $l\ne m$.
\goodbreak
\Par
Consider a tensor product $\Voxzt$ of \emod/s over $\Ugg$ coinciding
with $\Voxt$ as a \Umod/.
\Lm{T12}
For any $\phi\in\Valt$ we have
$$
\NN4>
\align
\bra\phi,L^+_{12}(t_1)\ldots & L^+_{12}(t_\ell)\>\voxt\ket\,{}=\,
\bigl(B_\tau(z)\)\phi\bigr)(\tell)\ \x
\\
&\!\]\x\;(q-\q)^\ell\>q^{\ell(1-\ell)/2-\ell\!\!\sun\!\!\La_m}
\pral\,\pron\,(\xi_m-t_a/z_m)\prab\,{\eta\)t_a-t_b\over t_a-t_b}\;,
\\
\nn6>
\bra\phi,L^-_{12}(t_1)\ldots & L^-_{12}(t_\ell)\>\voxt\ket\,{}=\,
\bigl(B_\tau(z)\)\phi\bigr)(\tell)\ \x
\\
&\!\]\x\;(q-\q)^\ell\>q^{\ell(1-\ell)/2-\ell\!\!\sun\!\!\La_m}
\pral\,\pron\,(1-\xi_m z_m/t_a)\prab\,{\eta\)t_a-t_b\over t_a-t_b}\;.
\endalign
$$
\endpro
\nt
It is easy to see that \rhs/s of the formulae above are \pol/s in $\tell$ and
$t_1\1\lc t_\ell\1\!$, \resp/, \cf. \(thgf) and \(Btau). So both the formulae
make sense without additional prescriptions.
\Par\nt
Lemma~\[T12] is proved in \Chap/~\SNo{4}.
\Th{local}
\back\Cite{V3}
For any ${\tau\in\S^n}$ and any transposition $(m,m+1)$, $\mn-1$, the \traf/
$$
B_{\tau,\taum}(z):\Vaxtm_\ell\,\to\,\Valt
$$
equals the operator
${\bigl(P_{V_{\tau_m}\}V_{\tau_{m+1}}}\:
R_{V_{\tau_m}\}V_{\tau_{m+1}}}\:\!(z_{\tau_m}/z_{\tau_{m+1}})\bigr)^*}$ acting
in the \^{$m$-th} and
\ifMag\nl\fi
\^{$(m+1)$-th} factors.
\endpro
\nt
The theorem follows from Lemma~\[T12] and formula \(Uintw).
\Par
Each $B_\tau(z)$ induces a linear map $\Valt\,\to\,\H(z)$ which also will be
denoted by $B_\tau(z)$.
\par
\Th{tcnon}
Let ${\ka\neps\eta^{-r}\!\yz}$ and ${\ka\ne p^{-s-1}\eta^r\zy\!}$
\vv.2>
for any $\roll$, $s\in\Zp$. Let $0<|p|<1$. Let \(npZ)\;--\;\(assum) hold.
Then for any $\tau\in\S^n$ the map $B_\tau(z):\Valt\to\H(z)$ is an \iso/.
\endpro
\nt
This statement follows from Theorem~\[kanon] and Lemma~\[tcfs]. The assumption
of the theorem means that ${\ka\neps q^{\,2\!\!\sun\!\!\La_m-2r}\!}$ and
${\ka\ne p^{-s-1}q^{-2\!\!\sun\!\!\La_m+\)2r}\!}$
\vv.2>
for any $\roll$, $s\in\Zp$.
\Par
It is easy to see that for any $\tau\in\S^n$ the images of $\FValt$ and
\vv.1>
$\FzValt$ under the map $B_\tau(z)$ coincide \resp/ with the \cosub/s
$\Rc(z)$ and $\Rc'(z)$ in the \thgf/ of the fiber $\Fo(z)$.
\Th{tcon}
Let $0<|p|<1$. Let \(npZ)\;--\;\(assum) hold. Let ${\pron\!\xi_m^2\nepsr}$
\ifMag\else\vv-.2>\fi
for any $\rpll$, $s\in\Znn$. If $\ka=\kone$, that is $\ka=\Kone$,
\ifMag\else\vv-.1>\fi
then for any $\tau\in\S^n$ the map $B_\tau(z)$ induces an \iso/
\ifMag\else\vv-.3>\fi
$$
\VFVt\)\to\,\H(z)\,.
$$
\ifMag\else\vvv-.1>\fi
Similarly, if $\ka=\ktwp$, that is $\ka=\Ktwp$, then for any $\tau\in\S^n$
the map $B_\tau(z)$ induces an \iso/
\ifMag\else\vv-.2>\fi
$$
\VFzVt\)\to\,\H(z)\,.
$$
\endpro
\nt
The statements follow from Theorem~\[kaon] and Lemmas~\[wbasis]\;--\;\[wDw'].
The assumption of the theorem means that
$q^{\,4\!\!\sun\!\!\La_m-2r}\!\neps$ for any $\rpll$, $s\in\Znn$.
\Par
Taking into account formulae \(VVF) and \(VVFz) we get an \iso/
$$
\gather
\bigl(\Vlst\bigr)^*\,\to\,\H(z)
\\
\nn3>
\Text{for $\ka=\kone$ and an \iso/}
\nn3>
\bigl(\Vlszt\bigr)^*\,\to\,\H(z)
\endgather
$$
\ifMag\vvv-.7>\else\vvv-.85>\fi
for $\ka=\ktwp$.
\Th{qKZ-GM}
\back\Cite{V3},\;\Cite{TV1}
For any $\mn$, the following diagram is commutative:
$$
\kern1.4em
\CD
\Valt @>{\tsize\ K^*_m(z_{\tau_1}\lc z_{\tau_n})\ }>> \Valt
\\
@V{\tsize{\vrp24pt:18pt>}B_\tau(\zmn)\ }VV
@VV{\tsize\ B_\tau(\zn{\vrp24pt:18pt>})}V
\\
\H(\zmn) @>>{\tsize\ \Cph{K^*_m(z_{\tau_1}\lc z_{\tau_n})}{A_m(\zn)}\ }>
\H(\zn)
\endCD
\kern-1.4em
$$
Here $A_m(z)$ are the operators of the \GM/, $K^*_m(z)$ are the operators dual
to $K_m(z)$, and $K_m(z)$ are the operators of the \qKZc/ in $\Vlt$ defined by
\(Kmz).
\endpro
\Cr*
The construction above identifies the \qKZc/ and the \GM/ restricted to the
\hgeom/ subbundle.
\endpro

\subsect{Tensor coordinates on the \ehgf/s of fibers}
Let $\Vqnz$ be \eVmod/s over $\Eqg$ with \hw/s $\Lan$ and \epoint/s $\zn$,
\resp/. Let $\Vqn$ be the corresponding \hmod/s. The \wt/ subspace $\Vql$ has
a basis given by the monomials $\Fvq$, $\lg\in\Zln$.
\par
Let ${\Fq=\Fq[\ka;\zn;\xin;\ell\,]}$ be the \ehgf/ where the parameters $\xin$
and the \hw/s $\Lan$ are related as follows:
$$
\xi_m=\eta^{\)\La_m},\qqq\mn\,.
$$
For any $z\in\BB$ and for any $\tau\in\S^n$ denote by $C_\tau(z)$ the following
linear map:
$$
\align
C_\tau(z):\Vqlt\, &{}\to\,\Fq(z)\,,
\\
\nn3>
C_\tau(z):\Fvqt\, &{}\map\,c^\tau_\lg\>W^\tau_\lg(t,z)\,.
\endalign
$$
Here $\Fq(z)$ is the \ehgf/ of the fiber and
$c^\tau_\lg=c_{{\vp(}^\tau\!\lg}\:(\xi_{\tau_1}\lc\xi_{\tau_n})$ where
$^\tau\}\lg=(\lg_{\tau_1}\lc\lg_{\tau_n})$ and
\ifMag\vv-.6>\else\vv->\fi
$$
\alignat2
c_\lg\:(\xin)\, &{}=\,\pron\,\prod_{s=1}^{\lg_m}\,
{\tht(\eta^s)\>\tht(\eta^{1-s}\xi_m^2)\over\tht(\eta)}
\plmn\!\}\eta^{-\lg_l\lg_m}\xi_m^{2\lg_l}\,\x{} &
\Tag{clg}
\\
\nn4>
&{}\>\x\,\Bigl(\,\prod_{s=1}^{\lg_1}\>
\tht\bigl(\eta^{s-\ell}\ka\1\yz\bigr)\,
\prod_{s=1}^{\lg_n}\>\tht\bigl(\eta^{\ell-s}\kzy\bigr)
\,\x{}&
\\
\nn2>
&& \Llap{\x\;\prmn\>\prod_{\tsize{s=-\lg\_m\atop s\ne 0}}^{\lg_{m+1}}
\tht\bigl(\eta^s\ka\1\!\!\!\tprllm\!\!\eta^{\)\lg_l}\xi_l\1\!\!\!
\tprml\!\!\eta^{-\lg_l}\xi_l\bigr)\Bigr)\1}\! &.
\endalignat
$$
The linear maps $C_\tau(z)$ are called the \em{\tenco/} on the \ehgf/ of a
fiber. The composition maps
$$
C_{\tau,\tau'}(z):\Vqltt\,\to\,\Vqlt\,,\ifMag\else\quad\fi\qquad
C_{\tau,\tau'}(z)\,=\,C\1_\tau(z)\o C_{\tau'}(z)\,,
$$
are called the \em{\traf/s}, \cf. \Cite{V3}, \Cite{TV3}.
\Ex
For $\ell=1$ the coefficients $c_\lg$ have the form
$$
\tsize c_{\eg(m)}(\xin)\,=\,\tht(\xi_m^2)\,
\Bigl(\tht\bigl(\ka\1\!\!\!\prllm\!\!\xi_l\1\!\!\!\prml\!\!\xi_l\bigr)\>
\tht\bigl(\ka\1\!\!\!\prlm\!\!\xi_l\1\!\!\!\prmll\!\!\xi_l\bigr)\Bigr)\1
\!\!\!\prml\!\!\xi_l^2
$$
where $\eg(m)=(0\lc\%1_{\sss\^{$m$-th}}\lc 0)$, $\,\mn$.
\enddemo
\Ex
For $\ell=2$ and $n=2$ the coefficients $c_\lg$ have the form
\ifMag
$$
\NN3>
\align
c\'{2,0}(\xi_1,\xi_2)\, &{}=\,\xi_2^4\;{\tht(\eta^2)\over\tht(\eta)}\;
\tht(\xi_1^2)\>\tht(\eta\1\xi_1^2)\,\x{}
\\
& {}\>\x\,\bigl(\tht(\ka\1\xi_1\xi_2)\>\tht(\eta\1\ka\1\xi_1\xi_2)\>
\tht(\ka\1\xi_1\1\xi_2)\>\tht(\eta\>\ka\1\xi_1\1\xi_2)\bigr)\1,
\\
\nn3>
c\'{1,1}(\xi_1,\xi_2)\, &{}=\,\eta\1\xi_2^2\;\tht(\xi_1^2)\>\tht(\xi_2^2)\,\x{}
\\
& {}\>\x\,\bigl(\tht(\eta\1\ka\1\xi_1\xi_2)\>\tht(\eta\>\ka\1\xi_1\1\xi_2)\>
\tht(\eta\1\ka\1\xi_1\1\xi_2)\>\tht(\eta\>\ka\1\xi_1\1\xi_2\1)\bigr)\1,
\\
\nn3>
c\'{0,2}(\xi_1,\xi_2)\, &{}=\;{\tht(\eta^2)\over\tht(\eta)}\;
\tht(\xi_2^2)\>\tht(\eta\1\xi_2^2)\,\x{}
\\
& {}\>\x\,\bigl(\tht(\ka\1\xi_1\1\xi_2)\>\tht(\eta\1\ka\1\xi_1\1\xi_2)\>
\tht(\ka\1\xi_1\1\xi_2\1)\>\tht(\eta\>\ka\1\xi_1\1\xi_2\1)\bigr)\1.
\endalign
$$
\else
$$
\NN4>
\align
& c\'{2,0}(\xi_1,\xi_2)\,=\,\xi_2^4\;{\tht(\eta^2)\over\tht(\eta)}\;
\tht(\xi_1^2)\>\tht(\eta\1\xi_1^2)\>
\bigl(\tht(\ka\1\xi_1\xi_2)\>\tht(\eta\1\ka\1\xi_1\xi_2)\>
\tht(\ka\1\xi_1\1\xi_2)\>\tht(\eta\>\ka\1\xi_1\1\xi_2)\bigr)\1,
\\
& c\'{1,1}(\xi_1,\xi_2)\,=\,\eta\1\xi_2^2\;\tht(\xi_1^2)\>\tht(\xi_2^2)\>
\bigl(\tht(\eta\1\ka\1\xi_1\xi_2)\>\tht(\eta\>\ka\1\xi_1\1\xi_2)\>
\tht(\eta\1\ka\1\xi_1\1\xi_2)\>\tht(\eta\>\ka\1\xi_1\1\xi_2\1)\bigr)\1,
\\
& c\'{0,2}(\xi_1,\xi_2)\,=\;{\tht(\eta^2)\over\tht(\eta)}\;
\tht(\xi_2^2)\>\tht(\eta\1\xi_2^2)\>
\bigl(\tht(\ka\1\xi_1\1\xi_2)\>\tht(\eta\1\ka\1\xi_1\1\xi_2)\>
\tht(\ka\1\xi_1\1\xi_2\1)\>\tht(\eta\>\ka\1\xi_1\1\xi_2\1)\bigr)\1.
\endalign
$$
\fi
\enddemo
\Lm{ecfs}
Let ${\eta^r\!\neps}$ for any $\rtll$, $s\in\Z$. Let
${\eta^{\)\La_l+\La_m}z_l/z_m\nepsr}$ for any $l,\mn$, and $\roll$, $s\in\Z$.
Let ${\ka^{\)\pm1}\!\pron\!\eta^{\)\La_m}\!\nepsr}$
\,for any $\roll$, $s\in\Z$. Assume that for a \perm/ $\tau$ we have
\ifMag\vv.15>\fi
${\ka\!\!\prllm^{}\!\!\eta^{\)\La_{\tau_l}}\!\!\!\prml\!\!\eta^{-\La_{\tau_l}}
\!\nepsr}$ \,for any $\mn-1$, and $\rll$, $s\in\Z$.
Then the linear map $C_\tau(z):\Vqlt\to\Fq(z)$ is an \iso/.
\goodbm
\endpro
\nt
The statement follows from Lemma~\[Wbasis].
\Rem
The map $C_\tau(z)$ considered as a \fn/ of $\ka$ has a simple pole at
$\ka=\kone$
\vv-.3>
\,because of the factor $\,\tht\bigl(\eta^{1-\ell}\ka\1\yz\bigr)$
in formula \(clg) for the coefficients $c_\lg\:$, and
\ifMag
$$
\gather
\Ker\bigl(\Resp C_\tau(z)\bigr)\,=\,\Vqlto\,,
\Tag{RespC}
\\
\nn4>
\Im\bigl(\Resp C_\tau(z)\bigr)\,=\,\Qc(z)\,,
\endgather
$$
\else
$$
\Ker\bigl(\Resp C_\tau(z)\bigr)\,=\,\Vqlto\,,\qqq
\Im\bigl(\Resp C_\tau(z)\bigr)\,=\,\Qc(z)\,,
\Tag{RespC}
$$
\fi
where $\Resp C_\tau(z)$ is the residue of $C_\tau(z)$ at $\ka=\kone$ and
$\Qc(z)$ is the \bosub/.
\par
Similarly, the map $C_\tau(z)$ has a simple pole at $\,\ka=\ktwp$
\vv-.3>
\;because of the factor
\ifMag\else\nl\fi
$\,\tht\bigl(\eta^{\ell-1}\kzy\bigr)$ in formula \(clg) for
the coefficients $c_\lg\:$, and
\ifMag
$$
\gather
\Ker\bigl(\Resm C_\tau(z)\bigr)\,=\,\Vqltn\,,
\Tag{ResmC}
\\
\nn4>
\Im\bigl(\Resm C_\tau(z)\bigr)\,=\,\Qc'(z)\,,
\endgather
$$
\else
$$
\Ker\bigl(\Resm C_\tau(z)\bigr)\,=\,\Vqltn\,,\qqq
\Im\bigl(\Resm C_\tau(z)\bigr)\,=\,\Qc'(z)\,,
\Tag{ResmC}
$$
\fi
where $\Resm C_\tau(z)$ is the residue of $C_\tau(z)$ at $\ka=\ktwp$ and
$\Qc'(z)$ is the \bosub/.
\enddemo
\Lm{ecfso}
Let ${\ka=\kOne}\!$. Let ${\eta^r\!\neps}$ for any $\rtll$, $s\in\Z$.
Let $\eta^{\)\La_l+\La_m}z_l/z_m\nepsr$ for any $l,m=2\lc n$, and $\roll$,
$s\in\Z$. Let ${\eta^{\)2\La_1}\nepsr}$ for any $\roll$, $s\in\Z$.
Let ${\pron\!\eta^{\)2\La_m}\!\nepsr}$ \,for any $\rpll$, $s\in\Z$.
Assume that for a \perm/ $\tau$ we have
${\!\prllm^{}\!\eta^{\)2\La_{\tau_l}}\!\nepsr}$ \,for any $m=2\lc n-1$, and
$\rolt$, $s\in\Z$. Then the linear map $C_\tau(z):\Vqlto\to\Fq(z)/\Qc(z)$
is an \iso/.
\endpro
\nt
The statement follows from Lemma~\[dimFQ].
\goodbreak
\Lm{ecfsp}
Let $\ka=\kTwp\!$. Let ${\eta^r\!\ne p^s}$ for any $\rtll$, $s\in\Z$. Let
$\eta^{\)\La_l+\La_m}z_l/z_m\alb\nepsr$ for any $l,\mn-1$, and $\roll$,
$s\in\Z$. Let ${\eta^{\)2\La_n}\ne p^s\eta^r}$ for any $\roll$, $s\in\Z$.
Let ${\pron\!\eta^{\)2\La_m}\!\ne p^s\eta^r}$ \,for any $\rpll$, $s\in\Z$.
Assume that for a \perm/ $\tau$ we have
${\!\prml^{}\!\!\eta^{\)2\La_{\tau_l}}\!\nepsr}$ \,for any $\mn-2$, and
$\rolt$, $s\in\Z$. Then the linear map $C_\tau(z):\Vqltn\to\Fq(z)/\Qc'(z)$ is
an \iso/.
\endpro
\nt
The statement follows from Lemma~\[dimFQ'].
\Par
Consider a tensor product $\Vqxzt$ of \eVmod/s over $\Eqg$ coinciding with
$\Vqxt$ as an \hmod/.
\Lm{L21}
Let $\la=\ka\zy$. Then for any $v\in\Vqlt$ we have
$$
\align
& T_{21}(t_1,\la)\ldots T_{21}(t_\ell,\eta^{\ell-1}\la)\>v\,=\,
\bigl(C_\tau(z)\)v\bigr)(\tell)\;\x
\\
\nn4>
&\quad\x\,\pros\>\tht(\eta^s\kzy)\,
\pral\,\pron\,{\tht(\xi_m\1t_a/z_m)\over\tht(\xi_m t_a/z_m)}
\prab\,{\tht(\eta\)t_a/t_b)\over\tht(t_a/t_b)}\ \vqx\;,
\endalign
$$
\endpro
\nt
The lemma is proved in \Chap/~\SNo{4}.
\Th{localq}
For any $\tau\in\S^n$ and any transposition $(m,m+1)$, $\mn-1$, the \traf/
$$
C_{\tau,\taum}(z):\Vqxtm\,\to\,\Vqxt
$$
equals the operator
\ifMag\nl\fi
${P_{V^e_{\tau_{m+1}}\}V^e_{\tau_m}}\:
\Rq_{V^e_{\tau_{m+1}}\}V^e_{\tau_m}}\!\bigl(z_{\tau_{m+1}}/z_{\tau_m},
(\eta^{\)H}\}\lox\eta^{\)H}\}\ox\%{\eta^{-H}}_{\^{$m$-th}}\}\lox\eta^{-H})
\>\eta\1\ka\bigr)}$.
\endpro
\nt
The theorem follows from Lemma~\[L21] and Theorem~\[Eintw].
\goodast
\Rem
The elliptic \Rm/ in Theorem~\[localq] has an operator
$$
K\,=\,(\eta^{\)H}\lox\eta^{\)H}\ox
\%{\eta^{-H}}_{\^{$m$-th}}\lox\eta^{-H})\>\eta\1\ka
$$
at the place of the second argument, and 
\vv.2>
$\Rq_{V^e_{\tau_{m+1}}\}V^e_{\tau_m}}\!(x,\la)$ commutes with $K$ for any
values of $x,\la$. The operator $\Rq_{V^e_{\tau_{m+1}}\}V^e_{\tau_m}}\!(x,K)$
is understood in the standard way:
$$
\Rq_{V^e_{\tau_{m+1}}\}V^e_{\tau_m}}\!(x,K)\>v\,=\,
\Rq_{V^e_{\tau_{m+1}}\}V^e_{\tau_m}}\!(x,\la)\>v
$$
for any $v\in\Vqxtm$ \st/ ${Kv\>=\)\la\)v}$.
\enddemo

\subsect{Tensor products of the \hgf/s}
Let ${\Fo[z_1\lc z_m;\xi_1\lc\xi_m;l\,]}$ and
${\Fq[\)\al;z_1\lc z_m;\xi_1\lc\xi_m;l\,]}$
\ifMag\else\vvn.2>\fi
be \resp/ the \tri/ and the \ehgf/s defined for the projection
${\C^{\,l+m}\!\to\C^{\>m}\!}$. In particular, in our previous notations we have
$$
\Fo=\Fo[\zn;\xin;\ell\,]\qquad\text{and}\qquad\Fq=\Fq[\ka;\zn;\xin;\ell\,]\,.
$$
There are maps
\ifMag\vv-.6>\fi
$$
\NN6>
\gather
{\align
\chi:\)\Fo[z_1\lc z_k;\xi_1\lc\xi_k;j\)] \ox
\Fo[z_{k+1}\lc z_{k+m};\xi_{k+1}\lc\xi_{k+m};l\,]\,\to{} \ifMag\kern-2em\fi &
\Tagg{chi}
\\
{}\to\,\Fo[z_1\lc z_{k+m};\xi_1\lc\xi_{k+m};j+l\,]\> \ifMag\kern-2em\fi &
\endalign}
\\
\nn-10>
\Text{and}
\ifMag\nn-8>\else\nn-10>\fi
\alignedat2
\ifMag\kern2em\fi \che:{} &
\Fq[\)\al\)\eta^{\)l}\!\tprod_{i=1}^m\xi_{i+k}\1;z_1\lc z_k;\xi_1\lc\xi_k;j\)]
\bigl((z_1\lc z_k)\bigr)\;\ox &&
\\
\nn-3>
&\ox\;\Fq[\)\al\)\eta^{-j}\!\tprod_{i=1}^k\xi_i\);
z_{k+1}\lc z_{k+m};\xi_{k+1}\lc\xi_{k+m};l\,]
\bigl((z_{k+1}\lc z_{k+m})\bigr)\,\to{} &&
\\
&& \Llap{{}\to\,\Fq[\)\al;z_1\lc z_{k+m};\xi_1\lc\xi_{k+m};j+l\,]
\bigl((z_1\lc z_{k+m})\bigr)\>} &
\endalignedat
\endgather
$$
which are \resp/ defined by ${\chi:f\ox g\map f\*g}$ and
${\che\]:f\ox g\map f*g}$ where
\ifMag
$$
\align
\kern.7em
(f\*g) &(t_1\lc t_{j+l})\;=
\\
&\!{}=\;{1\over j!\>l\)!}
\sum_{\si\in\S^{j+l}}\,\Bigl[\>
f(t_1\lc t_j)\>g(t_{j+1}\lc t_{j+l})\,\prod_{i=1}^k\,\prod_{a=1}^l\,
{\xi_i\)t_{a+j}-z_i\over t_{a+j}-\xi_i z_i}\,\Bigr]_\si
\\
\Text{and}
\nn4>
(f*g) &(t_1\lc t_{j+l})\;=
\\
&\!{}=\;{1\over j!\>l\)!} \sum_{\si\in\S^{j+l}}\,\LBc\>
f(t_1\lc t_j)\>g(t_{j+1}\lc t_{j+l})\,\prod_{i=1}^k\,\prod_{a=1}^l\,
{\tht(\xi_i\)t_{a+j}/z_i)\over\tht(\xi_i\1t_{a+j}/z_i)}\,\RBc_\si.
\kern-.7em
\endalign
$$
\else
$$
\gather
(f\*g)(t_1\lc t_{j+l})\,=\,{1\over j!\,l\)!}
\sum_{\si\in\S^{j+l}}\,\Bigl[\>
f(t_1\lc t_j)\>g(t_{j+1}\lc t_{j+l})\,\prod_{i=1}^k\,\prod_{a=1}^l
{\xi_i\)t_{a+j}-z_i\over t_{a+j}-\xi_i z_i}\,\Bigr]_\si,
\\
\Text{and}
(f*g)(t_1\lc t_{j+l})\,=\,{1\over j!\,l\)!}
\sum_{\si\in\S^{j+l}}\,\Bigl[\!\}\Bigl[\>
f(t_1\lc t_j)\>g(t_{j+1}\lc t_{j+l})\,\prod_{i=1}^k\,\prod_{a=1}^l
{\tht(\xi_i\)t_{a+j}/z_i)\over\tht(\xi_i\1t_{a+j}/z_i)}\,\Bigr]\!\}\Bigr]_\si,
\endgather
$$
\fi
We have the next lemmas.
\Lm{FFo}
Assume that ${\xi_i\>\xi_{j+k}z_{j+k}/z_i\]\netr}$ for any $i=1\lc k$,
$j=1\lc m$, $r=0\lc l-1$. Then the map
\ifMag\else\vv->\fi
$$
\alignat2
\chi:\)\Plus_{i+j=l} &
\Fo[z_1\lc z_k;\xi_1\lc\xi_k;i\)]\)\bigl((z_1\lc z_k)\bigr)\ox{} &&
\\
\nn-2>
&\!{}\ox\Fo[z_{k+1}\lc z_{k+m};\xi_{k+1}\lc\xi_{k+m};j\)]\)
\bigl((z_{k+1}\lc z_{k+m})\bigr)\,\to{} &&
\\
\nn6>
&&\Llap{{}\to\,\Fo[z_1\lc z_{k+m};\xi_1\lc\xi_{k+m};l\,]\)
\bigl((z_1\lc z_{k+m})\bigr)}\! &
\endalignat
$$
defined by linearity is bijective.
\endpro
\Lm{FFq}
Assume that ${\xi_i\>\xi_{j+k}z_{j+k}/z_i\nepsr}$ and \,
\ifMag\vv.15>\fi
${\al\}}\tprod_{a=1}^k\xi_a\tprod_{a=1}^m\xi_{a+k}\1\neps\eta^{\pm r}$
for any $i=1\lc k$, $j=1\lc m$, $r=0\lc l-1$, $s\in\Z$. Then the map
$$
\alignat2
\che\]:\)\Plus_{i+j=l} &
\Fq[\)\al\)\eta^{\)j}\!\tprod_{a=1}^m\xi_{a+k}\1;
z_1\lc z_k;\xi_1\lc\xi_k;i\)]\)\bigl((z_1\lc z_k)\bigr)\ox{} &&
\\
\nn-2>
&\!{}\ox\Fq[\)\al\)\eta^{-i}\!\tprod_{a=1}^k\xi_a\);
z_{k+1}\lc z_{k+m};\xi_{k+1}\lc\xi_{k+m};j\)]\)
\bigl((z_{k+1}\lc z_{k+m})\bigr)\,\to{} &&
\\
\nn6>
&&\Llap{{}\to\,\Fq[\)\al;z_1\lc z_{k+m};\xi_1\lc\xi_{k+m};l\,]\)
\bigl((z_1\lc z_{k+m})\bigr)}\! &
\endalignat
$$
defined by linearity is bijective
\endpro
\nt
Lemmas~\[FFo] and \[FFq] are proved in \Chap/~\SNo{4}.
\Par
It is clear that for any \fn/s $f,g,h$ we have $(f\*g)\*h\,=\,f\*(g\*h)$ and
for any \fn/s $f,g,h$ we have ${(f*g)*h\,=\,f*(g*h)}$. Lemmas~\[FFo], \[FFq]
can be extended naturally to an arbitrary number of factors.
\Par
The map $\che$ admits the following generalization.
Fix a nonnegative integer $k$. Let $n_0\lc n_k$ be integers \st/
$$
0=n_0<n_1\lsym< n_k\,.
$$
Fix nonnegative integers $l_1\lc l_k$. Let
$$
\Fq^i\,=\,
\Fq[\)\al_i;z_{n_{i-1}+1}\lc z_{n_i};\xi_{n_{i-1}+1}\lc\xi_{n_i};l_i\)]
$$
where ${\,\al_i=\al\)\eta^{\,\sum_{j>i}\]l_j\)-\]\sum_{j<i}\]l_j}\!\}
\prod_{\!\!\!j\le n\_{i-1}\;\;}\xi_j\,\prod_{j>n\_i}\xi_j\1}\!$.
Let $h(z_1\lc z_{n_k})$ be a \mef/ on $\C^{\x n_k}\!$ \st/
$$
h(z_1\lc p\)z_i\lc z_{n_k})\,=\,
\xi_i^{\,\sum_{j<i}\]l_j\)-\]\sum_{j>i}\]l_j}\]h(z_1\lc z_{n_k})
$$
for any $i=1\lc n_k$. Then we have a well defined map
$$
\gather
\Fq^1\lox\Fq^k\,\to\,
\Fq[\)\al;z_1\lc z_{n_k};\xi_1\lc\xi_{n_k};l_1\lsym+ l_k\)]\,,
\Tag{Fq1Fqk}
\\
\nn8>
f_1\lox f_k\,\map\,(f_1\]\lsym* f_k)\>h
\endgather
$$
of the \ehgf/s. We call the \fn/ $h(z_1\lc z_{n_k})$ an \em{\adjf/} for the
tensor product of the \ehgf/s $\Fq^1\lox\Fq^k$.

\ifAst\newpage\fi
\Sect[2c]{The \hpair/ and \ifAst\{\nl\fi the \hsol/s of
\ifMag\{\nlc\fi the \{\qKZe/}
\leftsecthead{The \hpair/ and the \hsol/s of the \{\qKZe/}
In this \chap/ we define the main object of this paper, the \hpair/.
We define a pairing between the \tri/ and the \ehgf/s of a fiber.

\subsect{The \hint/}
For any \fn/s $w\in\Fo(z)$ and $W\in\Fq(z)$ we define the \em{\hint/} by
\ifMag\vv-.5>\fi
$$
I(W,w)\,=\int_{\){\TTt^\ell}\]}\Phi(t)\>w(t)\>W(t)\;\dtt
\Tag{IWwt}
$$
\ifMag\vvv-.8>\fi
where ${\dtt=\!\pral dt_a/t_a}$, $\Phi(t)$ is the short \phf/ defined
in \(Phis) and $\TTt^\ell$ is a suitable deformation of the torus
$$
\TT^\ell=\lb\)t\in\Cl\vert\ |t_1|=1\llc|t_\ell|=1\)\rb\,.
$$
\par
Recall that we always have $0<|p|<1$. Let $(u)\9=\prod_{k=0}^\8(1-p^ku)$.
We take \,$\dsize\Pho(x;\al)={(x\al\1)\9\!\over(x\al)\9}$ \ in \(Phi) so that
the short \phf/ has the form
\ifMag\vv-.3>\fi
$$
\Phi(\tell,\zn)\,=\,\pron\,\pral\,{(\xi_m\1t_a/z_m)\9\!\over(\xi_mt_a/z_m)\9}\,
\prab\;{(\eta\)t_a/t_b)\9\!\over(\eta\1t_a/t_b)\9\!}\;.
\Tag{Phis}
$$
\ifMag\vsk.1>\else\par\fi
We define the \hint/ as follows.
Assume that $|\eta|>1$ and $|z_m|=1$, $|\xi_m|<1$, $\mn$. Set
$$
I(W,w)\,=\int_{\){\TT^\ell}\]}\Phi(t)\>w(t)\>W(t)\;\dtt\,.
\Tag{IWw}
\mmgood
$$
Notice that the integrand has simple poles at the hyperplanes
$$
\alignat2
\NN3>
& t_a/z_m=(p^s\xi_m)^{\pm1}\,, && \aell\,,\qquad\mn\,,
\Tag{list3}
\\
& t_a/\)t_b=(p^s\eta\1)^{\pm1}\,,\qqq && 1\le a<b\le\ell\,,
\endalignat
$$
for $s\in\Zp$ and essential singularities at the coordinate hyperplanes.
The set of hyperplanes \(list3) could be decomposed into subsets corresponding
to couples ${\lb a,m\rb}$ or ${\lb a,b\rb}$. Under the above assumptions for
each subset of the hyperplanes the torus $\TT^\ell$ separates the hyperplanes
corresponding to different choices of the sign.
\par
The \hint/ for generic $\xin$, $\zn$ and arbitrary $\eta$ is defined by \anco/
\wrt/ $\xin$, $\zn$ and $\eta$. This \anco/ makes sence since the integrand
is analytic in $\xin$, $\zn$ and $\eta$, \cf. \(Phi), \(wlga), \(Wlga).
More precisely, first we define the \hint/ for basis \fn/s ${w_\lg\,,\ W_\mg}$
and then extend the definition by linearity to arbitrary \fn/s $w\in\Fo(z)$,
$W\in\Fq(z)$. The result of the \anco/ can be represented as an integral of
the integrand over a suitably deformed torus. Namely, the poles of
the integrand of the \hint/ $I(W_\lg,w_\mg)$ are located at the hyperplanes
$$
t_a=p^s\xi_mz_m\,,\qquad t_a=p^{-s}\xi_m\1z_m\,,\qquad
t_a=p^s\eta\1t_b\,,\qquad t_a=p^{-s}\eta\)t_b\,,
\Tag{list2}
$$
${1\le b<a\le\ell}$, $\mn$, $s\in\Zp$. We deform $\xin$, $\zn$ and $\eta$
in such a way that the topology of the complement in $\Cxl$ to the union of
hyperplanes \(list2) does not change. We deform accordingly the torus
$\TT^\ell$ so that it does not intersect the hyperplanes \(list2) at every
moment of the deformation. The deformed torus is denoted by $\TTt^\ell$.
Then the \anco/ of the integral \(IWw) is given by formula \(IWwt).
\Th{Ianco}
For any $\lg,\mg\in\Zln$ the \hint/ $I(W_\lg,w_\mg)$ can be analytically
continued as a \hol/ univalued \fn/ of complex \var/s $\eta$, $\xin$,
$\zn$ to the region:
$$
\gather
\NN4>
\eta\ne 0\,,\qquad\qquad\xi_m\ne 0\,,\qquad z_m\ne 0\,,\qquad\mn\,,
\\
\eta^{r+1}\neps\,,\qquad\qquad \xi_m^2\>\nepsr\,,\qquad\mn\,,
\\
\xi_l^{\pm1}\xi_m^{\pm1}z_l/z_m\,\ne\,p^s\eta^r\,,\qqq
l,\mn\,, \qquad l\ne m\,,
\endgather
$$
where ${\roll}$, ${s\in\Z}$, and the combination of signs $\pm1$ can be
arbitrary {\rm(}\cf. \(npZ)\;--\;\(assum)\){\rm)}.
\endpro
\nt
The proof of the theorem is the same as the proof of Theorem 5.7 in
\Cite{TV3}.
\Par
Let $\Rc(z),\Rc'(z)\sub\Fo(z)$ be the \cosub/s and
let $\Qc(z),\Qc'(z)\sub\Fq(z)$ be the \bosub/s.
\Lm{IWw=0}
Let $\ka=\kone$. Let \(npZ)\;--\;\(assum) hold. Then
\vsk.2>
\atem For any ${w\in\Rc(z)}$, $\!{W\in\Fq(z)}$, the \hint/ $I(W,w)$ equals
zero.
\bitem For any ${w\in\Fo(z)}$, $\!{W\in\Qc(z)}$, the \hint/ $I(W,w)$ equals
zero.
\endpro
\Ex
Let $\ell=1$ and $\ka=\!\pron\!\xi_m$.
\ifMag\vv.05>\fi
Then the space $\Qc(z)$ is \onedim/ and is spanned by the \fn/ $W(t_1)=1$.
\mmgood
\vv.2>
Assume for simplicity that ${|z_m|=1}$, ${|\xi_m|<1}$ for any $\mn$.
Then the \hint/ $I(W,w)$ is given by
$$
I(W,w)\,=\int_{|t_1|=1}w(t_1)\)\pron\)
{(\xi_m\1t_1/z_m)\9\!\over(\xi_mt_1/z_m)\9}\;{dt_1\over t_1}\;.
$$
Since $w(0)=0$ for any $w\in\Fo(z)$, the integrand is regular in the disk
${|t_1|\le 1}$. Hence, $I(W,w)=0$ for any $w\in\Fo(z)$.
\enddemo
\Lm{IWw=0'}
Let $\ka=\ktwp$. Let \(npZ)\;--\;\(assum) hold. Then
\vsk.1>
\atem For any ${w\in\Rc'(z)}$, $\!{W\in\Fq(z)}$, the \hint/ $I(W,w)$
\hbox{equals zero.}
\bitem For any ${w\in\Fo(z)}$, $\!{W\in\Qc'(z)}$, the \hint/ $I(W,w)$ equals
zero.
\endpro
\Ex
Let $\ell=1$ and $\ka=\!\pron\!\xi_m\1$. Then the space $\Qc'(z)$ is \onedim/
and is spanned by the \fn/ ${\dsize
W(t_1)=\>t_1\1\!\pron\,{\tht(\xi_mt_1/z_m)\over\tht(\xi_m\1t_1/z_m)}}$. Assume
for simplicity that ${|z_m|=1}$, ${|\xi_m|<1}$ for any $\mn$. Then the \hint/
$I(W,w)$ is given by
$$
I(W,w)\,=\int_{|t_1|=1}w(t_1)\)\pron\)
{(p\>\xi_m\1z_m/t_1)\9\!\over(p\>\xi_mz_m/t_1)\9}\;{dt_1\over t_1^2}\;.
$$
Since $w(t_1)=O(1)$ as ${t_1\to\8}$ for any $w\in\Fo(z)$, the integrand is
regular in the domain ${|t_1|\ge 1}$ and behaves as $O(t_1^{-2})$
as ${t_1\to\8}$. Hence, $I(W,w)=0$ for any $w\in\Fo(z)$.
\enddemo
Lemmas~\[IWw=0] and \[IWw=0'] are proved in \Chap/~\SNo{4}.

\subsect{Determinant formulae for the \hpair/}
The \hint/ defines the \em{\hpair/}
$$
I:\Fq(z)\ox\Fo(z)\,\to\,\C
$$
which induces the \hpair/s
$$
\gather
\Ii\!:\Fq(z)/\Qc(z)\ox\Fo(z)/\Rc(z)\,\to\,\C
\\
\nn-4>
\Text{for \ $\ka=\kone$ \ and}
\nn-4>
I':\Fq(z)/\Qc'(z)\ox\Fo(z)/\Rc'(z)\,\to\,\C
\endgather
$$
for \ $\ka=\ktwp$. According to \(HFo), \(HFR) and
\ifMag\vvn.1>\fi
Lemmas~\[dimFQ], \[dimFQ'] this can be \resp/ written as
\vv->
$$
\gather
I:\Fq(z)\ox\H(z)\,\to\,\C\,,
\\
\nn8>
\Ii\!:\Fq(z)/\Qc(z)\ox\H(z)\,\to\,\C \qqq\text{and}\qqq
I':\Fq(z)/\Qc'(z)\ox\H(z)\,\to\,\C\,.
\endgather
$$
\par
Set $\dsize d(n,m,\ell,s)\,=
\sum_{\tsize{\vp( i,j\ge 0\atop\tsize{\vp( i+j<\ell\atop i-j=s\vp(}}}
{m-1+i\choose m-1}\,{n-m-1+j\choose n-m-1}$.
\goodbm
\goodast
\ifamsppt\goodbreak\fi
\Th{mu<>0}
Let ${\ka\neps\eta^{-r}\!\yz}$ and ${\ka\ne p^{-s-1}\eta^r\zy\!}$, \,
\vv.1>
$\roll$, $s\in\Zp$. Let \(npZ)\;--\;\(assum) hold. Then the \hpair/
\ifMag\vv.2>\nl\fi
${I:{}}\Fq(z)\ox\Fo(z)\to\C$ is \ndeg/. Moreover,
\ifMag
$$
\NN8>
\gather
\Lline{\DWwi_{\lg,\mg\in\Zln}\>=\,(2\pii)\vpb{\ell\tsize{n+\ell-1\choose n-1}}
\,\ell\)!\vpb{\tsize{n+\ell-1\choose n-1}}\>
\eta\vpb{-n\tsize{n+\ell-1\choose n+1}}\ \x}
\\
\Lline{\qquad{}\x\>
\pron\xi_{m\vp1}^{(n-m){\tsize{n+\ell-1\choose n}}}
\prsl\,\prmn\>\tht\bigl({\tsize\eta^s
\ka\1\!\!\!\prllm\!\!\xi_l\1\!\!\!\prml\!\!\xi_l\)}\bigr)^{d(n,m,\ell,s)}\ \x}
\\
\Rline{{}\x\,\pros\,\Bigl[\>
{(\eta\1)\9^n\>(\eta^{s+1-\ell}\ka\1\}\prod\xi_m)\9\>
(p\)\eta^{s+1-\ell}\)\ka\prod\xi_m)\9\!\over
(\eta^{-s-1})\9^n\>(p)\9^{2n-1}\>\prod(\eta^{-s}\xi_m^2)\9}\ \x\qquad}
\\
\nn-2>
\Rline{\x\;\plmn{(\eta^s\xi_l\1\xi_m\1z_l/z_m)\9\!\over
(\eta^{-s}\xi_l\)\xi_m z_l/z_m)\9}\,\Bigr]^{\tsize{n+\ell-s-2\choose n-1}}.}
\endgather
$$
\else
$$
\NN8>
\align
& \DWwi_{\lg,\mg\in\Zln}\,=\,(2\pii)\vpb{\ell\tsize{n+\ell-1\choose n-1}}\,
\ell\)!\vpb{\tsize{n+\ell-1\choose n-1}}\>
\eta\vpb{-n\tsize{n+\ell-1\choose n+1}}\ \x
\\
& \qquad\,\) {}\x\>
\pron\xi_{m\vp1}^{(n-m){\tsize{n+\ell-1\choose n}}}
\prsl\,\prmn\>\tht\bigl({\tsize\eta^s
\ka\1\!\!\!\prllm\!\!\xi_l\1\!\!\!\prml\!\!\xi_l\)}\bigr)^{d(n,m,\ell,s)}\ \x
\\
& \qquad\,\) {}\x\;
\pros\,\Bigl[\>{(\eta\1)\9^n\>(\eta^{s+1-\ell}\ka\1\}\prod\xi_m)\9\>
(p\)\eta^{s+1-\ell}\)\ka\prod\xi_m)\9\!\over
(\eta^{-s-1})\9^n\>(p)\9^{2n-1}\>\prod(\eta^{-s}\xi_m^2)\9}
\plmn{(\eta^s\xi_l\1\xi_m\1z_l/z_m)\9\!\over
(\eta^{-s}\xi_l\)\xi_m z_l/z_m)\9}\,\Bigr]^{\tsize{n+\ell-s-2\choose n-1}}.
\endalign
$$
\fi
\endpro
\Th{mu=0}
Let $\ka=\kone$. Let \(npZ)\;--\;\(assum) hold.
\vv.1>
If $\!{\pron\!\xi_m^2\nepsr}$ for all $\rpll$ and $s\in\Znn$, then the \hpair/
\ifMag\vv.2>\nl\fi
${\Ii{\!}:\Fq(z)/\Qc(z)\ox\Fo(z)/\Rc(z)\to\C}$ is \ndeg/. Moreover,
\ifMag
$$
\NN8>
\gather
\Lline{\DWwi_{\%{\Rlap{\lg,\mg\in\Zln}}_{\ssize\Rlap{\lg_1=\mg_1=0}}}
\hp{\ssize\lg_n=\mg_n=0}=\,(2\pii)^{\ell\tsize{n+\ell-2\choose n-2}}\,
\ell\)!^{\tsize{n+\ell-2\choose n-2}}\>
\eta\vpb{\)(1-n)\tsize{n+\ell-2\choose n}}\ \x}
\\
\Rline{
\aligned
&{}\x\>\pron\xi_{m\vp1}^{(n-m)\tsize{n+\ell-2\choose n-1}}
\prsl\,\prod_{m=2}^{n-1}\>\tht\bigl({\tsize\eta^{s+\ell-1}\!\!\!
\prllm\!\!\xi_l^{-2}}\bigr)^{d(n-1,m-1,\ell,s)}\ \x
\\
& {}\x\;\pros\,\Bigl[\>{(\eta\1)\9^{n-1}\>(p\)\eta^{s+2-2\ell}\prod\xi_m^2)\9
\>(\eta^s\xi_1^{-2})\9\!\over
(\eta^{-s-1})\9^{n-1}\>(p)\9^{2n-3}\!\prod_{1<m\le n}\!(\eta^{-s}\xi_m^2)\9}
\plmn{(\eta^s\xi_l\1\xi_m\1z_l/z_m)\9\!\over
(\eta^{-s}\xi_l\)\xi_m z_l/z_m)\9}\,\Bigr]^{\tsize{n+\ell-s-3\choose n-2}}.
\endaligned}
\endgather
$$
\else
$$
\NN8>
\align
& \DWwi_{\%{\Rlap{\lg,\mg\in\Zln}}_{\ssize\Rlap{\lg_1=\mg_1=0}}}
\hp{\ssize\lg_n=\mg_n=0}=\,(2\pii)\vpb{\ell\tsize{n+\ell-2\choose n-2}}\,
\ell\)!\vpb{\tsize{n+\ell-2\choose n-2}}\>
\eta\vpb{\)(1-n)\tsize{n+\ell-2\choose n}}\ \x
\\
& \qquad\,\) {}\x\>
\pron\xi_{m\vp1}^{(n-m)\tsize{n+\ell-2\choose n-1}}
\prsl\,\prod_{m=2}^{n-1}\>\tht\bigl({\tsize\eta^{s+\ell-1}\!\!\!
\prllm\!\!\xi_l^{-2}}\bigr)^{d(n-1,m-1,\ell,s)}\ \x
\\
& \qquad\,\) {}\x\,
\pros\,\Bigl[\>{(\eta\1)\9^{n-1}\>(p\)\eta^{s+2-2\ell}\prod\xi_m^2)\9\>
(\eta^s\xi_1^{-2})\9\!\over
(\eta^{-s-1})\9^{n-1}\>(p)\9^{2n-3}\!\prod_{1<m\le n}\!(\eta^{-s}\xi_m^2)\9}
\plmn{(\eta^s\xi_l\1\xi_m\1z_l/z_m)\9\!\over
(\eta^{-s}\xi_l\)\xi_m z_l/z_m)\9}\,\Bigr]^{\tsize{n+\ell-s-3\choose n-2}}.
\endalign
$$
\fi
\endpro
\Th{mu=02}
Let $\ka=\ktwp$. Let \(npZ)\;--\;\(assum) hold.
\vv.1>
If ${\pron\!\xi_m^2\nepsr}$ for all $\rpll$ and $s\in\Znn$, then the \hpair/
\ifMag\vv.2>\nl\fi
${I'{\}}:\Fq(z)/\Qc'(z)\ox\Fo(z)/\Rc'(z)\to\C}$ is \ndeg/. Moreover,
\ifAst\fixedpage\strut\vsk-2>\nt\fi
\ifMag
$$
\gather
\NN8>
\Lline{\DWwi_{\%{\Rlap{\lg,\mg\in\Zln}}_{\ssize\Rlap{\lg_n=\mg_n=0}}}
\hp{\ssize\lg_n=\mg_n=0}=\,(2\pii)\vpb{\ell\tsize{n+\ell-2\choose n-2}}\,
\ell\)!\vpb{\tsize{n+\ell-2\choose n-2}}\>
\eta\vpb{\)(1-n)\tsize{n+\ell-2\choose n}}\ \x}
\\
\Rline{
\aligned
&{}\x\,
\prmn\xi_{m\vp1}^{(n-m-1)\tsize{n+\ell-2\choose n-1}}
\prsl\,\prod_{m=1}^{n-2}\>\tht\bigl({\tsize p\)\eta^{s+1-\ell}\!\!\!
\prml\!\!\xi_l^2}\bigr)^{d(n-1,m,\ell,s)}\ \x
\\
& {}\x\;\pros\,\Bigl[\>{(\eta\1)\9^{n-1}\>(p\)\eta^{s+2-2\ell}\prod\xi_m^2)\9
\>(\eta^s\xi_n^{-2})\9\!\over
(\eta^{-s-1})\9^{n-1}\>(p)\9^{2n-3}\!\prod_{1\le m<n}\!(\eta^{-s}\xi_m^2)\9}
\plmn{(\eta^s\xi_l\1\xi_m\1z_l/z_m)\9\!\over
(\eta^{-s}\xi_l\)\xi_m z_l/z_m)\9}\,\Bigr]^{\tsize{n+\ell-s-3\choose n-2}}.
\endaligned}
\endgather
$$
\else
$$
\NN8>
\align
& \DWwi_{\%{\Rlap{\lg,\mg\in\Zln}}_{\ssize\Rlap{\lg_n=\mg_n=0}}}
\hp{\ssize\lg_n=\mg_n=0}=\,(2\pii)\vpb{\ell\tsize{n+\ell-2\choose n-2}}\,
\ell\)!\vpb{\tsize{n+\ell-2\choose n-2}}\>
\eta\vpb{\)(1-n)\tsize{n+\ell-2\choose n}}\ \x
\\
& \qquad\,\) {}\x\,
\prmn\xi_{m\vp1}^{(n-m-1)\tsize{n+\ell-2\choose n-1}}
\prsl\,\prod_{m=1}^{n-2}\>\tht\bigl({\tsize p\)\eta^{s+1-\ell}\!\!\!
\prml\!\!\xi_l^2}\bigr)^{d(n-1,m,\ell,s)}\ \x
\\
& \qquad\,\) {}\x\,
\pros\,\Bigl[\>{(\eta\1)\9^{n-1}\>(p\)\eta^{s+2-2\ell}\prod\xi_m^2)\9\>
(\eta^s\xi_n^{-2})\9\!\over
(\eta^{-s-1})\9^{n-1}\>(p)\9^{2n-3}\!\prod_{1\le m<n}\!(\eta^{-s}\xi_m^2)\9}
\plmn{(\eta^s\xi_l\1\xi_m\1z_l/z_m)\9\!\over
(\eta^{-s}\xi_l\)\xi_m z_l/z_m)\9}\,\Bigr]^{\tsize{n+\ell-s-3\choose n-2}}.
\endalign
$$
\fi
\endpro
\nt
In Theorems~\[mu<>0]\;--\;\[mu=02], the product $\prod$ without limits stands
for $\pron$.
\par\nt
Theorems~\[mu<>0]\;--\;\[mu=02] are proved in \Chap/~\SNo{4}.
\Ex
Theorem~\[mu<>0] for $n=1$, $\ell=1$ gives
$$
\int_{\!C\,}\!{\tht(c\)t)\over(at)\9(b/t)\9\!}\,{dt\over t}\;=\,
2\pii\,{(p\)a/c)\9\>(bc)\9\!\over(ab)\9}\;.
$$
Here $C$ is an anticlockwise oriented contour around the origin $t=0$
separating the sets ${\lb\)p^s\!/a\vert s\in\Zn\)\rb}$ and
${\lb\)p^s b\vert s\in\Zp\)\rb}$.
\enddemo
\Ex
Theorems~\[mu=0] and \[mu=02] for $n=2$, $\ell=1$ give particular cases of
the Askey-Roy formula \Cite{GR, (4.11.2)}
$$
\int_{\!C\,}\!
{\tht(p\)t/c)\>\tht(abc\)t)\over(at)\9(b\)t)\9(\al/t)\9(\bt/t)\9\!}\,
{dt\over t}\;=\,2\pii\,{(ab\al\bt)\9\>\tht(ac)\>\tht(bc)\over
(p)\9(a\al)\9(a\bt)\9(b\al)\9(b\bt)\9\!}\;.
\Tag{AR}
$$
Here $C$ is an anticlockwise oriented contour around the origin $t=0$
separating the sets ${\lb\)p^s\!/a,\;p^s\!/\)b}\vert{s\in\Zn\)\rb}$ and
${\lb\)p^s\al,\;p^s\bt}\vert{s\in\Zp\)\rb}$.
\enddemo
\Ex
There are \p-analogues of the gamma-\fn/ and the power \fn/:
$$
\Gm_{\!p}(x)\)=\)(1-p)^{1-x}(p)\9/(p^x)\9\,,\qqq
(1-u)_p^{2x}=\)(p^{-x}u)\9/(p^xu)\9\,.
$$
Introduce new \fn/s $\,\lb\]-u\rb_p^{2x}=\)\tht(p^{-x}u)/\tht(p^xu)\,$ and
$$
\sin_p(\pi x)\,=\;{\pi\over\Gm_{\!p}(x)\>\Gm_{\!p}(1-x)}\;=\;
{\pi\>\tht(p^x)\over(1-p)\>(p)\9^3}\;.
$$
We have $(1-u)_p^{2x}=\)\lb\]-u\rb_p^{2x}\)(1-pu\1)_p^{2x}\!$.
\Par
Theorems~\[mu=0] and \[mu=02] for $\ell=1$ give \resp/ the following formulae:
$$
\gather
\Lline{\det\biggl[\;\)\int_{\!C\,}{p^{\)2\Dl_l}-1\over t-p^{\)\Dl_l}z_l}\;
\tpron(1-t/z_m)_p^{2\Dl_m}\>
\bigl\lb\]-p^{-\!\!\sum_{1\le j\le k}\!\!\Dl_j}t/z_k\bigr\rb
^{2\!\!\sum_{1\le j<k}\!\!\Dl_j}_{p\vp1}\;\x}
\\
\Cline{\x\,\tprod_{1\le j<k}\lb\]-t/z_j\rb^{-2\Dl_j}_{p\vp1}
\prod_{1\le j<l}{t-p^{-\Dl_j}z_j\over t-p^{\)\Dl_j}z_j}\
{dt\over p-1\]}\;\)\biggr]_{k,l=2}^n={}\qqq}
\\
\nn8>
\ifMag
\Rline{
\aligned
=\;\Gm_{\!p}(1+2\tsun\]\Dl_m)\1\pron\Gm_{\!p}(1+2\Dl_m)\!
\plmn\!\!(1-z_l/z_m)_{p\vp1}^{2(\Dl_l+\Dl_m)}\;\x &
\\
\nn6>
\x\;\prmn 2i\)\sin_p\bigl(-2\pi\!\!\}\tsum_{1\le j\le m}\!\!\!\Dl_j\bigr) &
\endaligned}
\else
\Rline{=\;\Gm_{\!p}(1+2\tsun\]\Dl_m)\1\pron\Gm_{\!p}(1+2\Dl_m)\!
\plmn\!\!(1-z_l/z_m)_{p\vp1}^{2(\Dl_l+\Dl_m)}\>
\prmn 2i\)\sin_p\bigl(-2\pi\!\!\}\tsum_{1\le j\le m}\!\!\!\Dl_j\bigr)}
\fi
\\
\Text{and}
\nn8>
\Lline{\det\biggl[\;\)\int_{\!C\,}
{p^{\)\Dl_l}-p^{-\Dl_l}\over t/z_l-p^{\)\Dl_l}}\,
\tpron(1-t/z_m)_p^{2\Dl_m}\>
\bigl\lb\]-p^{\>\sum_{k<j\le n}\!\!\Dl_j}t/z_k\bigr\rb
^{-2\!\!\sum_{k\le j\le n}\!\!\Dl_j}_{p\vp1}\;\x}
\\
\Cline{\x\,\tprod_{1\le j<k}\lb-t/z_j\rb^{-2\Dl_j}_{p\vp1}
\prod_{1\le j<l}{t-p^{-\Dl_j}z_j\over t-p^{\)\Dl_j}z_j}\
{dt\over t\>(p-1)\}}\;\>\biggr]_{k,l=1}^{n-1}={}\qqq}
\\
\nn8>
\ifMag
\Rline{
\aligned
=\;\Gm_{\!p}(1+2\tsun\]\Dl_m)\1\pron\Gm_{\!p}(1+2\Dl_m)\!
\plmn\!\!(1-z_l/z_m)_{p\vp1}^{2(\Dl_l+\Dl_m)}\;\x &
\\
\nn6>
\x\;\prmn 2i\)\sin_p\bigl(2\pi\!\!\}\tsum_{m<j\le n}\!\!\!\Dl_j\bigr) &
\endaligned}
\else
\Rline{=\;\Gm_{\!p}(1+2\tsun\]\Dl_m)\1\pron\Gm_{\!p}(1+2\Dl_m)\!
\plmn\!\!(1-z_l/z_m)_{p\vp1}^{2(\Dl_l+\Dl_m)}\>
\prmn 2i\)\sin_p\bigl(2\pi\!\!\}\tsum_{m<j\le n}\!\!\!\Dl_j\bigr)}
\fi
\endgather
$$
where $C$ is an anticlockwise oriented contour around the origin $t=0$
separating the sets ${\lb\)p^{s-\Dl_m}z_m}\vert\mn\,,\ {s\in\Zn\)\rb}$ and
${\lb\)p^{s+\Dl_m}z_m}\vert\mn\,,\ {s\in\Zp\)\rb}$. These formulae are
analogues of the next formula \Cite{V1}\>:
\ifAst
$$
\align
\det\Bigl[\int_{z_k}^{z_{k+1}}\!{\,\la_l\over t-z_l}\,
& \tpron(t-z_m)^{\la_m}\>dt\>\Bigr]_{k,l=1}^{n-1}={}
\\
&\Cph{{}\tpron{}}{{=}}\}
\Gm(1+\tsun\]\la_m)\1\pron\Gm(1+\la_m)\,\prod_{l\ne m}(z_l-z_m)^{\la_m}.
\endalign
$$
\else
$$
\det\Bigl[\int_{z_k}^{z_{k+1}}\!{\,\la_l\over t-z_l}\,
\tpron(t-z_m)^{\la_m}\>dt\>\Bigr]_{k,l=1}^{n-1}=\,
\Gm(1+\tsun\]\la_m)\1\pron\Gm(1+\la_m)\,\prod_{l\ne m}(z_l-z_m)^{\la_m}.
$$
\fi
\enddemo
\ifAst\fixedpage\fi
\Ex
Theorem~\[mu<>0] for $n=1$ and arbitrary $\ell$ gives the following \qbeta/
integral
\ifMag
$$
\align
\int_{\)\TT^\ell\]}\,
\prkl\,{\tht(c\)t_k)\over t_k\>(at_k)\9\>(b/t_k)\9\!}\;
\pjkl{(t_j/t_k)\9\over(x\)t_j/t_k)\9\!}\;\dt\,=\qqq &
\Tag{qbeta}
\\
\nn4>
=\;(2\pii)^\ell\,\ell\)!\,
\pros{(x)\9\>(x^sbc)\9\>(p\)x^s a/c)\9\!\over(x^{s+1})\9\>(x^s ab)\9} &
\endalign
$$
\else
$$
\int_{\)\TT^\ell\]}\,
\prkl\,{\tht(c\)t_k)\over t_k\>(at_k)\9\>(b/t_k)\9\!}\;
\pjkl{(t_j/t_k)\9\over(x\)t_j/t_k)\9\!}\;\dt\,=\,(2\pii)^\ell\,\ell\)!\,
\pros{(x)\9\>(x^sbc)\9\>(p\)x^s a/c)\9\over(x^{s+1})\9\>(x^s ab)\9}
\Tag{qbeta}
$$
\fi
where $|a|<1$, $|b|<1$, $|x|<1$. In \Chap/~\SNo{4} we use this formula to
prove Theorems~\[mu<>0]\;--\;\[mu=02]. We give a proof of the formula in
Appendix~\SNo{C} and show there that the calculation of the integral by
residues inside the torus $\TT^\ell$ implies the formula for the \qSelberg/
integral proved by Aomoto \Cite{AK, Theorem 3.2}, see formula \(qbetaR).
\enddemo
\Ex
The formulae of Theorems~\[mu=0] and \[mu=02] for $n=2$ are particular cases of
the following formula
$$
\align
\ifMag\kern2.5em\fi
\int_{\)\TT^\ell\]}\,\prkl\,{\tht(p\)t_k/c)\>\tht(x^{\ell-1}abc\)t_k)\over
t_k\>(at_k)\9\>(b\)t_k)\9\>(\al/t_k)\9\>(\bt/t_k)\9\!}\;
\pjkl{(t_j/t_k)\9\over(x\)t_j/t_k)\9\!}\;\dt\,=\qqq &
\Tagg{ARl}
\\
\nn4>
=\;(2\pii)^\ell\,\ell\)!\,\pros\,
{(x)\9\>(x^{\ell+s-1}ab\al\bt)\9\>\tht(x^s ac)\>\tht(x^s bc)\over
(x^{s+1})\9\>(x^s a\al)\9\>(x^s a\bt)\9\>(x^s b\al)\9\>(x^s b\bt)\9\>(p)\9\!}&
\ifMag\kern-2.5em\fi
\endalign
$$
where $|a|<1$, $|b|<1$, $|\al|<1$, $|\bt|<1$, $|x|<1$. This formula is a
multidimensional generalization of the Askey-Roy formula \(AR).
\mmgood
In Appendix~\SNo{D} we give a proof of this formula and show that
the calculation of the integral by residues outside the torus $\TT^\ell$
implies the formula for the most general multidimensional \qbeta/ integral
conjectured by Askey \Cite{As, Conjecture 8}, see formula \(Ascj).
\enddemo
\Rem
It is plausible that the assumptions on $p$, $\xin$, $\zn$ of
Theorems~\[mu<>0]\;--\;\[mu=02] as well as of Theorems and Lemmas~\[kanon],
\[kaon], \[tcnon], \[tcon], \[IWw=0], \[IWw=0'], \[asol], \[Asol], \[AsIWw]
could be replaced by the following weaker assumptions: the step $p$ and
the parameter $\eta$ are \st/ \(npZ) holds, and the parameters $\xin$
and $\zn$ are \st/
$$
\xi_l\)\xi_m z_l/z_m\,\ne\,p^s\eta^r\,,\qqq l,\mn\,,
\qquad s\in\Z\,,
$$
for any $\roll$ and $s\in\Z$.
\enddemo

\subsect{The \hsol/s of the \{\qKZe/}
Let $W$ be any element of the \ehgf/ $\Fq$. The restriction of the \fn/ $W$ to
a fiber defines an element $W|_z\in\Fq(z)$ of the \ehgf/ of the fiber.
The \hpair/ allows us to consider the element $W|_z\in\Fq(z)$ as an element
$s_W\:(z)$ of the space $\H^*(z)$ dual to the \hcg/ $\H(z)$. This construction
defines a section of the bundle over $\Cxn$ with fiber $\H^*(z)$.
\par
There is a simple but important statement.
\Th{sWz}
Let $\xin$ obey \(Lass). Then the section $s_W\:$ is a periodic section \wrt/
the \GM/.
\endpro
\nt
The theorem is proved in \Chap/~\SNo{4}.
\Par
Consider the \hpair/ as a map $\,{\Ib(z):\Fq(z)\to\bigl(\Fo(z)\bigr)^*\!}$
\,so that for any ${W\in\Fq}$ we have ${s_W\:=\Ib(z)\>W|_z}$.
Let $\Vn$ be \Vmod/s over $\Uu$ with \hw/s $\qLan\!$. Recall that
$$
q^2=\eta\,,\qqq\eta^{\)\La_m}=\xi_m\,,\qquad\mn\,.
$$
The map $\Ib(z)$ and the \tenco/ $B_\tau(z)$ induce a map
$$
\bigl(B_\tau(z)\bigr)^*\!\o\Ib(z)\>:\>\Fq(z)\,\to\,\Vlt\,.
$$
\Lm{Vsing}
Let $\ka=\kone$. Let \(npZ)\;--\;\(assum) hold. Then for any $W\in\Fq(z)$ we
have that $\bigl(B_\tau(z)\bigr)^*\!\o\Ib(z)\cdot W\in\Vlst\!$.
\endpro
\Lm{Vsing'}
Let $\ka=\ktwp$. Let \(npZ)\;--\;\(assum) hold. Then for any $W\in\Fq(z)$ we
have that $\bigl(B_\tau(z)\bigr)^*\!\o\Ib(z)\cdot W\in\Vlszt\!$.
\endpro
\nt
Lemmas~\[Vsing], \[Vsing'] follow from Lemmas~\[IWw=0], \[IWw=0'], \resp/.
\Par
A section $s_W\:$ and the \tenco/ $B_\tau$ induce a section
$$
\Psi_W^\tau\):\>z\,\map\,B_\tau^*(z^{\tau\1})\,s_W\:(z^{\tau\1})\in\Vlt
\Tag{PsiW}
$$
of the \trib/ with fiber ${\Vlt}$. Here ${z^{\si}\]=(z_{\si_1}\lc z_{\si_n})}$.
If $\tau$ is the identity \perm/, then we write $\Psi_W$ instead of
$\Psi_W^\id$.
\goodbm
\Par
Theorems~\[sWz], \[qKZ-GM] and Lemmas~\[Vsing], \[Vsing'] imply the
following statement.
\Cr{qKZsol}
The section $\Psi_W^\tau$ is a \sol/ of the \qKZe/ with values in $\Vlt$.
Moreover, if $\ka=\kone$, \,then $\Psi_W^\tau$ takes values in $\Vlst\!$,
and if $\ka=\ktwp$, \,then $\Psi_W^\tau$ takes values in $\Vlszt\!$.
\endpro
\nt
We call \sol/s $\Psi_W^\tau$ the \em{\hsol/s} of the \qKZe/.
\Par
Let $\,{S_{m,\tau}(\zn):\Voxt\to\>\Voxtm}$ equal the operator
\nl
$\,P_{V_{\tau_m}\}V_{\tau_{m+1}}}\:
R_{V_{\tau_m}\}V_{\tau_{m+1}}}\:\!(z_m/z_{m+1})$
acting in the \^{$m$-th} and \^{$(m+1)$-th} factors. Define operators
$\Sh_{\tau,\tau'}$ acting on \fn/s of $\zn$ by the following formulae:
\ifMag
$$
\gather
{\align
& \bigl(\Sh_{\taum,\tau}\>f\)\bigr)\)(\zn)\,={}
\Tagg{Shat}
\\
\nn6>
& \quad {}=\,S_{m,\tau}(\zmmn)\>f(\zmmn)\,,
\endalign}
\\
\nn8>
\Sh_{\tau,\tau'}\>\Sh_{\tau'\!,\tau''}\,=\,\Sh_{\tau,\tau''}
\endgather
$$
\else
$$
\gather
\kern1em
\bigl(\Sh_{\taum,\tau}\>f\)\bigr)\)(\zn)\,=\,S_{m,\tau}(\zmmn)\>f(\zmmn)\,,
\kern-1em
\Tag{Shat}
\\
\nn6>
\Sh_{\tau,\tau'}\>\Sh_{\tau'\!,\tau''}\,=\,\Sh_{\tau,\tau''}
\endgather
$$
\fi
where $(m,m+1)$ is a transposition, $\mn-1$, and $\tau,\tau'\!,\tau''\}\in\S^n$
are arbitrary \perm/s. The operator $\Sh_{\tau,\tau'}$ acts on a \fn/ taking
values in $\Voxtt$ and the result is a \fn/ taking values in $\Voxt$.
\Lm{Shwd}
Formulae \(Shat) define operators $\Sh_{\tau,\tau'}$ selfconsistently.
\endpro
\nt
The statement follows from the inversion relation \(inv) and the \YB/ \(YBq).
\Par
The \qKZe/ has the following important property.
\Th{qKZfun}
The \qKZe/ is functorial. Namely, for any \perm/s $\tau,\tau'\in\S^n$ and any
\sol/ $\Psi$ of the \qKZe/ with values in $\Voxtt$, the \fn/
$\Sh_{\tau,\tau'}\Psi$ is a \sol/ of the \qKZe/ with values in $\Voxt$.
\endpro
\Th{mono}
The \hsol/s $\Psi^\tau_W$ of the \qKZe/ are functorial. Namely, for any \perm/s
$\tau,\tau'\in\S^n$ and any \fn/ $W\in\Fq$ we have that
$\Sh_{\tau,\tau'}\Psi_W^{\tau'}=\Psi_W^\tau$.
\endpro
\nt
The statement follows from Theorem~\[local].

\subsect{The \hmap/}
Let $\Vqnz$ be \eVmod/s over $\Eqg$ with \hw/s $\Lan$ and \epoint/s $\zn$.
Let $\Vqn$ be the corresponding \hmod/s.
The \tenco/ $B_\tau(z)$, $C_{\tau'}(z)$ induce the \em{\hmap/}
$$
\gather
\Is_{\tau,\tau'}(z)\>:\>\Vqltt\,\to\,\Vlt\,,
\Tag{hmdef}
\\
\nn4>
\Is_{\tau,\tau'}(z)\,=\,(B_\tau(z)\bigr)^*\!\o\Ib(z)\o C_{\tau'}(z)\,.
\endgather
$$
\Th{hmap}
Let \(npZ)\;--\;\(assum) hold.
Let ${\ka^{\)\pm1}\!\pron\!\eta^{\)\La_m}\!\nepsr}$ \,for any $\roll$,
$s\in\Z$. Assume that for a \perm/ $\tau'$ we have
${\ka\!\!\prllm^{}\!\!\eta^{\)\La_{\tau'_l}}\!\!\!\prml\!\!
\eta^{-\La_{\tau'_l}}\!}{}\nepsr$ \,for any $\mn-1$, and $\rll$, $s\in\Z$.
Then the \hmap/ $\Is_{\tau,\tau'}(z)$ is well defined and \ndeg/.
\endpro
\nt
The statement follows from Lemmas~\[tcfs],~\[ecfs] and Theorem~\[mu<>0].
\Th{hmapo}
Let \(npZ)\;--\;\(assum) hold. Let ${\pron\!\eta^{\)2\La_m}\!\nepsr}$
\,for any $\rpll$, $s\in\Z$. Assume that for a \perm/ $\tau'$ we have
${\!\prllm^{}\!\eta^{\)2\La_{\tau'_l}}\!\nepsr}$ \,for any $\mn-1$, and
$\rolt$, $s\in\Z$. Then the \hmap/ $\Is_{\tau,\tau'}(z)$ is well defined and
\ndeg/ for any $\ka$ in the punctured \neib/ of $\kOne$. Moreover,
$\Is_{\tau,\tau'}(z)$ considered as a \fn/ of $\ka$ has a finite limit as
$\,\ka\to\kOne$ \,and the limit $\;\lim\>\Is_{\tau,\tau'}(z)\,$ is \ndeg/.
\endpro
Define the \hmap/ $\Is_{\tau,\tau'}(z)$ at $\ka=\kOne$ by the \anco/ \wrt/
$\ka$. Notice that the restriction of the map $\Is_{\tau,\tau'}(z)$ to
the subspace $\Vqltto$ \,is regular at $\ka=\kOne$, since in this case all
the maps involved in definition \(hmdef) are well defined at $\ka=\kOne$.
\Cr{hmps}
Let $\ka=\kOne$ and the assumptions of Theorem~\[hmapo] hold. Then
$$
\Is_{\tau,\tau'}(z)\bigl(\Vqltto\bigr)\,=\,\Vlst.
$$
\endpro
\nt
The statement follows from Lemma~\[Vsing].
\Th{hmapp}
Let \(npZ)\;--\;\(assum) hold. Let ${\pron\!\eta^{\)2\La_m}\!\nepsr}$
\,for any $\rpll$, $s\in\Z$. Assume that for a \perm/ $\tau'$ we have
${\!\prml^{}\!\!\eta^{\)2\La_{\tau'_l}}\!\nepsr}$ \,for any $\mn-1$, and
$\rolt$, $s\in\Z$. Then the \hmap/ $\Is_{\tau,\tau'}(z)$ is well defined and
\ndeg/ for any $\ka$ in the punctured \neib/ of $\kTwp$. Moreover,
$\Is_{\tau,\tau'}(z)$ considered as a \fn/ of $\ka$ has a finite limit as
$\,\ka\to\kTwp$ \,and the limit $\;\lim\>\Is_{\tau,\tau'}(z)\,$ is \ndeg/.
\endpro
Define the \hmap/ $\Is_{\tau,\tau'}(z)$ at $\ka=\kTwp$ by the \anco/ \wrt/
$\ka$. Notice that the restriction of the map $\Is_{\tau,\tau'}(z)$ to
the subspace $\Vqlttn$ \,is regular at $\ka=\kTwp$, since in this case all
the maps involved in definition \(hmdef) are well defined at $\ka=\kTwp$.
\Cr{hmpsp}
Let $\ka=\kTwp$ and the assumptions of Theorem~\[hmapp] hold. Then
$$
\Is_{\tau,\tau'}(z)\bigl(\Vqltn\bigr)\,=\,\Vlszt.
$$
\endpro
\nt
The statement follows from Lemma~\[Vsing'].
\Par\nt
Theorems~\[hmapo] and~\[hmapp] are proved in \Chap/~\SNo{4}.
\Par
Therefore, we constructed the \hmap/s
$$
\Is_{\tau,\tau'}(z):\Vqxztt\,\to\,\Voxzt
$$
from modules over the \eqg/
\aagood
to modules over the \qlo/. The maps have the following properties:
\ifMag
$$
\gather
\Is_{\taum,\tau'}(z)\,=\,P_{V_{\tau_m}\}V_{\tau_{m+1}}}\:
R_{V_{\tau_m}\}V_{\tau_{m+1}}}\:\!(z_{\tau_m}/z_{\tau_{m+1}})\,
\Is_{\tau,\tau'}(z)\,,
\Tagg{Imm+1}
\\
\nn6>
{\align
& \Is_{\tau,\taumm}(z)\,=\,\Is_{\tau,\tau'}(z)\;\x
\\
\nn3>
& \ \,\,{}\x\,P_{V^e_{\tau'_{m+1}}\}V^e_{\tau'_m}}\:
\Rq_{V^e_{\tau'_{m+1}}\}V^e_{\tau'_m}}\!\bigl(z_{\tau'_{m+1}}/z_{\tau'_m},
(\eta^{\)H}\lox\eta^{\)H}\ox\%{\eta^{-H}}_{\^{$m$-th}}\lox\eta^{-H})
\>\eta\1\ka\bigr)
\endalign}
\endgather
$$
\else
$$
\align
& \Is_{\taum,\tau'}(z)\,=\,P_{V_{\tau_m}\}V_{\tau_{m+1}}}\:
R_{V_{\tau_m}\}V_{\tau_{m+1}}}\:\!(z_{\tau_m}/z_{\tau_{m+1}})\,
\Is_{\tau,\tau'}(z)\,,
\Tagg{Imm+1}
\\
\nn6>
& \Is_{\tau,\taumm}(z)\,=\,\Is_{\tau,\tau'}(z)\;\x
\\
\nn3>
& \ \,\,{}\x\,P_{V^e_{\tau'_{m+1}}\}V^e_{\tau'_m}}\:
\Rq_{V^e_{\tau'_{m+1}}\}V^e_{\tau'_m}}\!\bigl(z_{\tau'_{m+1}}/z_{\tau'_m},
(\eta^{\)H}\lox\eta^{\)H}\ox\%{\eta^{-H}}_{\^{$m$-th}}\lox\eta^{-H})
\>\eta\1\ka\bigr)\,,
\endalign
$$
\fi
where $(m,m+1)$ is a transposition.
\Par
{\bls 1.15\bls
For any \ewf/ $W_\lg^\tau(t,z)$ let $Y_\lg^\tau(z)$ be the corresponding
\adjf/. Recall that this means that the product $Y_\lg^\tau(z)W_\lg^\tau(t,z)$
is an element of the \ehgf/. Define a map $Y^\tau(z)\in\End(\Vqxt)$ by
the rule:
\par\ifAst\vsk->\fi}
$$
Y^\tau(z)\>:\>\Fvqt\,\to\,Y_\lg^\tau(z)\>\Fvqt\,.
$$
The map $Y^\tau(z)$ is called an \em{\adjm/} for the tensor product $\Vqxt$.
\Par
If $v\in\Vqxtt$, \,then the \hmap/ $\Is_{\tau,\tau'}(z)$ and the \adjm/
$Y^{\tau'}\!(z)$ define a section
$$
\Psi_{v,Y^{\tau'}_{\vp|}}^\tau\):\>z\,\map\,
\Is_{\tau,\tau'}(z^{\tau\1})\>Y^{\tau'}\!(z^{\tau\1})\!\cdot v\in\Voxt
$$
where ${z^{\si}\]=(z_{\si_1}\lc z_{\si_n})}$.
\Th{PsivY}
For any \adjm/ $Y^{\tau'}\!(z)$ and any $v\in\Vqltt$ the section
$\Psi_{v,Y^{\tau'}_{\vp|}}^\tau\!$ is a \sol/ of the \qKZe/ with values in
$\Vlt$. Under assumptions of each of the Theorems~\[mu<>0]\,--\,\[mu=02] all
\sol/s are constructed in this way.
\endpro
\nt
The theorem is proved in \Chap/~\SNo{4}.
\Rem
Theorem~\[PsivY] can be reformulated as follows. For a given \adjm/ $Y^\tau(z)$
the assignment
${v\)\map\)\Psi_{v,Y^{\tau}_{\vp|}}^\id}\!$ defines an \iso/ of the space
$\SS$ of \sol/s of the \qKZe/ with values in $\Vox$ and the space
$\Vqxt\ox\FF$, where $\FF$ is the space of \fn/s of $\zn$ which are
\p-periodic \wrt/ each of the \var/s,
$$
\CC_\tau:\Vqxt\ox\FF\,\to\,\SS.
\Tag{CCtau}
$$
The compositions of the \iso/s: $\CC_{\tau,\tau'}=\CC\1_\tau\>\CC_{\tau'}$,
define linear maps
$$
\CC_{\tau,\tau'}(z):\Vqxtt\,\to\,\Vqxt
\Tag{CCtt}
$$
depending on $\zn$ and \p-periodic \wrt/ all the \var/s. We call these
compositions the \traf/s. Theorem~\[asol] in the next \chap/ shows that
$\CC_{\tau,\tau'}(z)$ is a \traf/ from the \asol/ of the \qKZe/ in
the \azo/ $\AA_\tau$ to the \asol/ in the \azo/ $\AA_{\tau'}$, \cf. \(azone).
\par
Notice that $\CC_{\tau,\tau'}(z)$ differs from the \traf/ $C_{\tau,\tau'}(z)$
defined in \Chap/~\[:2b], namely
$$
\CC_{\tau,\tau'}(z)\,=\,\bigl(Y^\tau(z)\bigr)\1 C_{\tau,\tau'}(z)\>
Y^{\tau'}\!(z)\,.
\Tag{CCC}
$$
\enddemo
Let $\Sh_{\tau,\tau'}$ be operators defined by formulae \(Shat). We extend
their action to matrix-valued \fn/s in a natural way.
\par
The maps ${\Ib_\tau(z)=I_{\tau,\tau}(z^{\tau\1})\>Y^\tau(z^{\tau\1})}$ satisfy
the \qKZe/s with values in $\Voxt$, \resp/. The following theorem describes
their ``monodromy'' properties \wrt/ \perm/s of the \var/s $\zn$ in terms of
the elliptic \Rms/.
\Th{monoR}
For any \perm/ $\tau\in\S^n$ and any transposition $(m,m+1)$, $\mn-1$,
we have that
$$
\alignat2
\smash{\bigl(\Sh_{\tau,\taum}\)\Ib_{\taum}\bigr)\)(z)\,=\,
\Ib_{\tau}(z)\,\bigl(Y^\tau(z^{\tau\1})\bigr)\1}\>{}\x{}\!\!\!&
\\
\nn7>
{}\x\,P_{V^e_{\tau_{m+1}}\}V^e_{\tau_m}}\:
\Rq_{V^e_{\tau_{m+1}}\}V^e_{\tau_m}}\!\bigl(z_{m+1}/z_m, (\eta^{\)H}\}\lox
\eta^{\)H}\} &&& {}\ox \%{\eta^{-H}}_{\^{$m$-th}}\lox\eta^{-H})
\>\eta\1\ka\bigr)\,\x{}
\\
\nn4>
&&& \hp{\ox}{}\x\,Y^{\taum}(z^{(m,m+1)\cdot\tau\1})\,.
\endalignat
$$
\endpro
\nt
The statement follows from formulae \(Imm+1).

\ifAst\newpage\fi
\Sect[2d]{Asymptotic \sol/s of the \{\qKZe/}
One of the most important characteristics of a \difl/ \eq/ is the monodromy
group of its \sol/s. For the \difl/ \KZe/ with values in a tensor product
of \rep/s of a simple Lie algebra its monodromy group is described in terms of
the corresponding \qg/. This fact establishes a remarkable connection between
\rep/ theories of simple Lie algebras and their quantum groups, see \Cite{K},
\Cite{D2}, \Cite{KL}, \Cite{SV2}, \Cite{V2}, \Cite{V4}.
\par
The analogue of the monodromy group for \deq/s is the set of \traf/s
between \asol/s. For a \deq/ one defines suitable \azo/s in the domain of the
definition of the \eq/ and then an \asol/ for every zone. Thus, for
every pair of \azo/s one gets a \traf/ between the corresponding \asol/s,
\cf. \Cite{TV3}.
\goodbm
\par
In this \chap/ we describe \azo/s, \asol/s, and their \traf/s for the \qKZe/
with values in a tensor product of \Umod/s. A remarkable fact is that the
\traf/s are described in terms of the elliptic \Rms/ acting in the tensor
product of the corresponding \Emod/s. This fact establishes a correspondence
between \rep/ theories of \qlo/s and \eqg/s, since the \qKZe/ is defined in
terms of the \tri/ \Rm/ action in the tensor product of \Umod/s (and,
therefore, in terms of the \qlo/ action), and the elliptic \Rm/ action in
the tensor product of \Emod/s is defined in terms of the action of the \eqg/.
\vsk0>
\Par
Consider the \qKZe/ with values in $\Vl$. For every \perm/ $\tau\in\S^n$ we
consider an \azo/ $\AA_\tau$ in $\Cxn\!$ given by
$$
\AA_\tau\,=\,\lb\)z\in\Cxn\vert\ |\)z_{\tau_m}/z_{\tau_{m+1}}|\ll 1\,,
\quad\mn-1\)\rb\,.
\Tag{azone}
$$
Say that $z$ tends to limit in the \azo/, $z\rto\AA_\tau$, if
${z_{\tau_m}/z_{\tau_{m+1}}\!\to 0}$ for all $\mn-1$.
\par
Say that a basis $\Psi_1\lc \Psi_N$ of \sol/s of \qKZe/ form an \em{\asol/} in
the \azo/ $\AA_\tau$ if
$$
\Psi_j(z)\,=\,h_j(z)\>\bigl(v_j+o(1)\bigr)\,,
$$
where $h_1(z)\lc h_N(z)$ are \mef/s \st/
$$
h_j(\zmn)\,=\,a_{jm}\>h_j(\zn)
$$
for suitable numbers $a_{jm}$, $v_1\lc v_N$ are constant vectors which form
a basis in $\Vl$, and $o(1)$ tends to $0$ as $z\rto\AA_\tau$.
\Par
For every \perm/ $\tau\in\S^n$ we constructed the \fn/s $W_\lg^\tau$,
$\lg\in\Zln$, whose restriction to a fiber gives a basis in the \ehgf/ of
the fiber. Let $Y_\lg^\tau$, $\lg\in\Zln$, be the corresponding \adjf/s so
that the \fn/s $Y_\lg^\tau W_\lg^\tau$ are in the \ehgf/. These \fn/s define
a basis $\Psi_{Y_\lg^\tau W_\lg^\tau}$, $\lg\in\Zln$, of \sol/s of the \qKZe/,
\cf. \(PsiW).
\ifMag\Par\else\par\fi
Recall that ${^\tau\}\lg=(\lg_{\tau_1}\lc\lg_{\tau_n})}$ for any $\lg\in\Zln$.
For $\lg,\mg\in\Zln$ say that ${\lg\ll\mg}$ if $\lg\ne\mg$ and
\ifMag\else\vv-.5>\fi
$\,\sum_{i=1}^m\lg_i\le\sum_{i=1}^m\mg_i\;$ for any $\mn-1$.
\Th{asol}
Let the parameters $\xin$ obey condition \(Lass). Then for any \perm/
$\tau\in\S^n$ the basis $\Psi_{Y_\lg^\tau W_\lg^\tau}$, $\lg\in\Zln$, is an
\asol/ in the \azo/ $\AA_\tau$. Namely,
$$
\Psi_{Y_\lg^\tau W_\lg^\tau}(z)\,=\,Y_\lg^\tau(z)\>
\bigl(\Ome_\lg^\tau+o(1)\bigr)
$$
as $z$ tends to limit in the \azo/, $z\rto\AA_\tau$, so that at any moment
assumption \(assum) holds. Here
$$
\Ome_\lg^\tau\>=\,\Xi_{\lg}^\tau\)
\Fv\,+\!\sum_{\]{}^\tau\}\mg\gg\]{\vp(}^\tau\!\lg\)}N^\tau_{\lg\mg}\)\Fvv
$$
for suitable constant coefficients $N^\tau_{\lg\mg}$ and $\,\Xi_\lg^\tau$,
and the constant $\,\Xi_\lg^\tau$ is given by
$$
\align
\Xi_\lg^\tau\,=\,(2\pii)^\ell\,\ell\)!\,\pron &\;
\Bigl(\)q^{\lg_m(1-\lg_m)/2\,+\)\lg_m\La_m}
\prod_{\tsize{1\le l<m\atop \si_l<\si_m}}\!\xi_l^{\lg_m}\!\!
\prod_{\tsize{1\le l<m\atop \si_l>\si_m}}\!
\eta^{\lg_l\lg_m}\xi_l^{-\lg_m}\ \x
\\
\nn6>
\x\;\, & \!\prod_{s=0}^{\lg_m-1}\,{(\eta\1)\9\>
(\eta^{-s}(\ka_{\lg,m}^\tau)\1\xi_m)\9\>
(p\)\eta^{-s}\ka_{\lg,m}^\tau\)\xi_m)\9\!\over
(\eta^{-s-1})\9\>(\eta^{-s}\xi_m^2)\9\>(p)\9}\,\Bigr)
\endalign
$$
where \,$\si=\tau\1$ and \,$\ka_{\lg,m}^\tau=
\ka\!\prlmsi\!\!\!\eta^{-\lg_l}\xi_l\:\!\prmlsi\!\eta^{\lg_l}\xi_l\1\]$.
\endpro
\nt
The theorem is proved in \Chap/~\SNo{4}.
\Rem
The \tri/ \Rm/ $R(x)$ has finite limits as $x$ tends to zero or infinity,
\cf. \(Rspec).
Thus, the \qKZo/s $K_1(z)\lc K_n(z)$ have finite limits as $z$ tends to limit
in an \azo/:
$$
\vp{K}\smash{K_m(z)\,=\,K_m^\tau\,\ono\,,\qqq z\rto\AA_\tau\,,\qqq\mn\,.}
$$
where $K_m^\tau$ are some operators independent of $z$. The vectors
$\Ome_\lg^\tau$ form an eigenbasis of the operators $K_m^\tau$ with \eva/s
$a_{\lg,m}^\tau$, \cf. \(Ytau).
\enddemo
\Rem
Recall that the \adjf/s $Y_\lg^\tau(z)$ have the following properties:
$$
\gather
Y_\lg^\tau(\zmn)\,=\,a_{\lg,m}^\tau\>Y_\lg^\tau(\zn)\,,
\Tag{Ytau}
\\
\nn6>
\tsize a_{\lg,m}^\tau\,=\,\ka^{\lg_m}\!\!\!
\prlmsi\!\!\!\eta^{-\lg_l\lg_m}\xi_l^{\lg_m}\xi_m^{\lg_l}\!
\prmlsi\!\eta^{\lg_l\lg_m}\xi_l^{-\lg_m}\xi_m^{-\lg_l},
\endgather
$$
for any $\tau\in\Sl$, $\lg\in\Zln$, $\mn$. Here $\si=\tau\1$.
\enddemo
\Rem
If the absolute value of $\ka$ is sufficiently small, then the relation
$Y_\mg^\tau(z)=o\bigl(Y_\lg^\tau(z)\bigr)$ as $z$ tends to limit in the \azo/,
$z\rto\AA_\tau$, implies that ${^\tau\}}\mg\gg{^\tau\}}\lg$, \cf. \(Ytau).
Similarly, if the absolute value of $\ka$ is sufficiently large, then
the relation $Y_\mg^\tau(z)=o\bigl(Y_\lg^\tau(z)\bigr)$ as $z$ tends to limit
in the \azo/, $z\rto\AA_\tau$, implies that ${^\tau\}}\mg\ll{^\tau\}}\lg$.
\par
For example, assume that $|\eta|=1$, $|\xi_m|=1$, $\mn$, $|\ka|<1$, and all
the \adjf/s $Y_\mg^\tau(z)$ are regular at point $(1\lc 1)\in\Cxn\!$. Let
$z=(p^{s_1}\lc p^{s_n})$ where $s_1\lc s_n$ are integers. Then the relation
$Y_\mg^\tau(z)=o\bigl(Y_\lg^\tau(z)\bigr)$ as $z\rto\AA_\tau$ implies that
the sum $\sum_{i=1}^n s_i(\mg_i-\lg_i)$
\ifMag\vv-.4>\else\vv-.2>\fi
is large positive if all the differences ${s_{\tau_m}\!\!-s_{\tau_{m+1}}}$,
$\,\mn-1$ are large positive. Since $\sum_{i=1}^n \lg_i=\sum_{i=1}^n \mg_i\,$
we have that
$$
\vp{\sum}\smash{
\sum_{i=1}^n s_i(\mg_i-\lg_i)\,=\,\sum_{i=1}^{n-1}\,(s_{\tau_i}-s_{\tau_{i+1}})
\tsum_{j=1}^i(\mg_{\tau_j}\}-\lg_{\tau_j})\,.}
$$
Therefore, $\,\sum_{j=1}^i\mg_{\tau_j}\ge\sum_{j=1}^i\lg_{\tau_j}\,$ for any
$i=1\lc n-1$, and $\mg\ne\lg$, \, that is ${\,^\tau\}}\mg\gg{^\tau\}}\lg$.
\enddemo
\Rem
The \qKZe/ depends \mer/ally on parameters $\ka$, $\xin$. Let the \adjf/s
$Y_\lg^\tau$ depend \mer/ally on $\ka$, $\xin$. Then the basis of \sol/s
$\Psi_{Y_\lg^\tau W_\lg^\tau},\ \lg\in\Zln$ also depends \mer/ally on $\ka$,
$\xin$. The \as/s of the basis $\Psi_{Y_\lg^\tau W_\lg^\tau},\ \lg\in\Zln$,
described in Theorem~\[asol] determine the basis uniquely. Namely, if a basis
of \sol/s \mer/ally depends on the parameters $\ka$, $\xin$ and has \as/s in
$\AA_\tau$ described in Theorem~\[asol], then such a basis coincides with the
basis $\Psi_{Y_\lg^\tau W_\lg^\tau}\!$. In fact, elements of any such a basis
are linear combinations of the \fn/s $\Psi_{Y_\lg^\tau W_\lg^\tau}\!$ with
coefficients \mer/ally depending on $\ka$, $\xin$ and \p-periodic in $\zn$.
To preserve the \as/s one can add to an element
$\Psi_{Y_\lg^\tau W_\lg^\tau}\!$ any other \fn/s
$\Psi_{Y_{\mg\vp1}^\tau W_{\mg\vp1}^\tau}\!$ having smaller \as/s.
If the absolute value of $\ka$ is sufficiently small, then one can add only the
\fn/s $\Psi_{Y_{\mg\vp1}^\tau W_{\mg\vp1}^\tau}\!$ with
${^\tau\}}\mg\gg{^\tau\}}\lg$, and if the absolute value of $\ka$ is
sufficiently large, then one can add only the \fn/s
$\Psi_{Y_{\mg\vp1}^\tau W_{\mg\vp1}^\tau}\!$ with
${^\tau\}}\mg\ll{^\tau\}}\lg$, see the previous Remark.
Since the coefficients of added terms are \mer/ they have to be zero.
\enddemo
\Rem
The \asol/ $\Psi_{Y_\lg^\tau W_\lg^\tau}$, $\lg\in\Zln$, in the \azo/
$\AA_\tau$ of the \qKZe/ with values in $\Vl$, \cf. Theorem~\[asol], is
an image of the monomial basis $\Fvqt$, $\lg\in\Zln$, of $\Vqlt$ under
the composition ${\Is_{\id,\tau}(z)\>Y^\tau(z)}$ of the \hmap/ and the \adjm/,
\cf. Theorem~\[PsivY]. The \traf/s between the \asol/s are linear maps
$\CC_{\tau,\tau'}$, see \(CCtt). Formula \(CCC) and Theorem~\[localq] show
that the \traf/s between \asol/s corresponding to neighbouring \azo/s are given
by the dynamical elliptic \Rms/ twisted by the corresponding \adjm/s.
\enddemo
\Ex
Theorem~\[asol], formula \(CCC) and Theorem~\[localq] allow us to write
an elliptic \Rm/ as an infinite product of \tri/ \Rms/. Namely, consider
the \qKZe/ with values in the tensor product of two $\Uu$ \Vmod/s $\VV$
with \hw/s $q^{\La_1},q^{\La_2}\!$, \resp/. Then there are two \azo/s
$|z_1/z_2|\gg 1$ and $|z_1/z_2|\ll 1$. Our result on the \traf/ from
the first \azo/ to the second one is the following statement.
\par
Let $V^e_1, V^e_2$ be the \eVmod/ over $\Eqg$ with \hw/s $\La_1,\La_2$, \resp/.
Then we have
\vv-.5>
$$
\bigl(\Rq_{V^e_1V^e_2}(x)\bigr)\1\,=
\prod_{s=-\8}^\8\!\ka^{-H\ox\one}\)\RV\:(p^s x)\,,
\Tag{RpR}
$$
provided the infinite product in \rhs/ is suitably regularized and the factors
of the product are ordered in such a way that $s$ grows from the right to
the left, see the example below. The restriction of formula \(RpR) to the \wt/
subspace $(\VV)_1\:$ of \wt/ $q^{\La_1+\La_2-1}$ can be transformed into
the infinite product formula for $2{\)\x}2$ matrices, which looks as follows.
\par
Let $a,b,c,d,\al,\dl,p$ be nonzero complex numbers \st/ $|p|<1$ and
$$
\al/\dl\ne p^s,\qqq a/d\ne p^s,\qqq bc\ne(1-p^s)\>(p^{-s}\al d-\dl a)
$$
for any $s\in\Z$. Set
$\dsize A(u)\,=\,\pmatrix a-\al u & b\)u \\ c & d-\dl u \endpmatrix$.
Let $\la,\mu$ be two \sol/s of the quadratic \eq/ $\det A(u)=0$. Then
\ifMag
\mmgood
$$
\gather
\Lline{\lim_{s\to\8}\, \pmatrix a & 0 \\ c & d \endpmatrix ^{\tsize\!-s}
\Bigl(\prod_{r=-s}^s A(p^r\} u)\Bigr)\,
\pmatrix -\al & \>b\> \\ 0 & \>\Llap{-}\dl \endpmatrix ^{\tsize\!-s}
\!u^{-s}\>p^{\)s(s+1)/2}\;=\;
\pmatrix 1&0{\vp{\dsize{1\over1}}}\,\\ \dsize{\vp1c\over a-d}&1\,\endpmatrix
\;\x}
\\
\nn14>
\Lline{\x\;\pmatrix \dsize\,
{a\>\tht(\al u/a)\>(p\)\la\)\dl/a)\9\>(p\)\mu\)\dl/a)\9\!\over
(p\)d/a)\9\>(p\)\dl/\al)\9\>(p)\9} &
\dsize {b\)u\>\tht(u\1a/\dl)\>(p\)\la\)\al/a)\9\>(p\)\mu\)\al/a)\9\!\over
(p\)d/a)\9\>(\al/\dl)\9\>(p)\9}\,
\\
\nn8>
\dsize\, {c\>\tht(\al u/d)\>(p\)\la\)\dl/d)\9\>(p\)\mu\)\dl/d)\9\!\over
(a/d)\9\>(p\)\dl/\al)\9\>(p)\9} &
\dsize {d\>\tht(\dl u/d)\>(\la\)\al/d)\9\>(\mu\)\al/d)\9\!\over
(a/d)\9\>(\al/\dl)\9\>(p)\9}\, \endpmatrix\!
\pmatrix\,1& \dsize{b\over\dl-\al}\\ \,0{\vp{\dsize{1\over1}}}& 1 \endpmatrix}
\endgather
$$
\else
$$
\gather
\lim_{s\to\8}\, \pmatrix a & 0 \\ c & d \endpmatrix ^{\tsize\!-s}
\Bigl(\prod_{r=-s}^s A(p^r\} u)\Bigr)\,
\pmatrix -\al & \>b\> \\ 0 & \>\Llap{-}\dl \endpmatrix ^{\tsize\!-s}
\! u^{-s}\>p^{\)s(s+1)/2}\;=
\\
\ald
\nn14>
=\ \!
\pmatrix 1&0{\vp{\dsize{1\over1}}}\,\\ \dsize{\vp1c\over a-d}&1\,\endpmatrix\!
\pmatrix \dsize\,
{a\>\tht(\al u/a)\>(p\)\la\)\dl/a)\9\>(p\)\mu\)\dl/a)\9\!\over
(p\)d/a)\9\>(p\)\dl/\al)\9\>(p)\9} &
\dsize {b\)u\>\tht(u\1a/\dl)\>(p\)\la\)\al/a)\9\>(p\)\mu\)\al/a)\9\!\over
(p\)d/a)\9\>(\al/\dl)\9\>(p)\9}\,
\\
\nn8>
\dsize\, {c\>\tht(\al u/d)\>(p\)\la\)\dl/d)\9\>(p\)\mu\)\dl/d)\9\!\over
(a/d)\9\>(p\)\dl/\al)\9\>(p)\9} &
\dsize {d\>\tht(\dl u/d)\>(\la\)\al/d)\9\>(\mu\)\al/d)\9\!\over
(a/d)\9\>(\al/\dl)\9\>(p)\9}\, \endpmatrix\!
\pmatrix\,1& \dsize{b\over\dl-\al}\\ \,0{\vp{\dsize{1\over1}}}& 1 \endpmatrix
\endgather
$$
\fi
\vvv>
where the factors of the product are ordered in such a way that $r$ grows from
the right to the left.
\enddemo
Theorem~\[asol] admits the following generalization. Fix a nonnegative integer
$k$ not greater than $n$. Fix nonnegative integers $n_0\lc n_k$ \st/
$$
0=n_0<n_1\lsym< n_k=n
$$
and consider an \azo/ in $\Cxn$ given by
\ifMag
$$
\align
\AA\,=\,\lb\) z\in\Cxn\vert & |\)z_{m_i}/z_{m_{i+1}}|\ll 1
\\
\nn2>
&\,\text{for all $\,m_1\lc m_k\,$ \st/ $\,n_{i-1}<m_i\le n_i\,$,
$\,i=1\lc k\)\rb\,$.}
\endalign
$$
\else
$$
\AA\,=\,\lb\) z\in\Cxn\vert |\)z_{m_i}/z_{m_{i+1}}|\ll 1 \quad
\text{for all $\,m_1\lc m_k\,$ \st/ $\,n_{i-1}<m_i\le n_i\,$,
$\,i=1\lc k\)\rb\,$.}
$$
\fi
We say that $z$ tends to limit in the \azo/, $z\rto\AA$, if
${z_l/z_m\}\to 0}$ for all $l,m$ \st/
\vv.1>
$n_{i-1}<l\le n_i<m\le n_{i+1}$ for some $i=1\lc k-1$, and $|z_l/z_m|$ remains
bounded for all $l,m$ \st/ $n_{i-1}<l,m\le n_i$ for some $i=1\lc k$.
\par
Fix nonnegative integers $\ell_1\lc\ell_k$ \st/ $\sum_{i=1}^k\ell_i=\ell$. Let
\vv-.5>
$$
\gather
\Fo^i[\)\ell_i\)]\,=\,
\Fo[z_{n_{i-1}+1}\lc z_{n_i};\xi_{n_{i-1}+1}\lc\xi_{n_i};\ell_i\)]
\\
\nn-2>
\Text{and}
\nn-4>
\Fq^i[\)\ell_i\)]\,=\,
\Fq[\)\ka_i;z_{n_{i-1}+1}\lc z_{n_i};\xi_{n_{i-1}+1}\lc\xi_{n_i};\ell_i\)]
\endgather
$$
where ${\,\ka_i=\ka\)\eta^{\,\sum_{j>i}\]\ell_j\)-\]\sum_{j<i}\]\ell_j}\!\}
\prod_{l\le n_{\Rph{i}{i-1}}}\>\xi_l\:\;\prod_{l>n\_i}\xi_l\1}\!$.
\vv.1>
Let $Y(\zn)$ be an \adjf/ for the tensor product
${\Fq^1[\)\ell_1\)]\lox\Fq^k[\)\ell_k\)]}$, \cf. \(Fq1Fqk).
Then we have linear maps
$$
\NN4>
\gathered
\Fo^1[\)\ell_1\)]\lox\Fo^k[\)\ell_k\)]\,\hto\,\Fo
\\
f_1\lox f_k\,\map\,f_1\lsym\*f_k
\endgathered
\qqq\text{and}\qqq
\aligned
& \Fq^1[\)\ell_1\)]\lox\Fq^k[\)\ell_k\)]\,\hto\,\Fq
\\
& f_1\lox f_k\,\map\,(f_1\lsym*f_k)Y
\endaligned
$$
\wrt/ the tensor products introduced in \Chap/~\SNo{2b}.
\par
Let $\mg\in\Zln$. Say that ${\mg\gg\}(\ell_1\lc\ell_k)\,}$ if
\vv-.5>
$\,\sum_{l=1}^{n_i}\mg_{\)l}\ge\sum_{j=1}^i\ell_j\,$ for any $i=1\lc k-1$,
and at least one of the inequalities is strict.
\goodbreak
\Par
For any ${W\in\Fq^i[\)\ell_i\)]}$ let $\Psi_W(z_{n_{i-1}+1}\lc z_{n_i})$ be
the \sol/ of the \qKZe/ with values in $(V_{n_{i-1}+1}\lox V_{n_i})_{\ell_i}\:$
corresponding to $W$ (\cf. \(PsiW)\)).
\Th{Asol}
Let the parameters $\xin$ obey condition \(Lass).
\vvn.2>
Let ${W_i\in{}}\alb{\Fq^i[\)\ell_i\)]}$, $i=1\lc k$. Let $W=W_1\lsym*W_k$
and let $Y$ be an \adjf/ for the tensor product
${\Fq^1[\)\ell_1\)]\lox\Fq^k[\)\ell_k\)]}$.
Then the \sol/ $\Psi_{YW}(\zn)$ of the \qKZe/ with values in $\Vl$ has
the following \as/s as $z$ tends to limit in the \azo/, $z\rto\AA$,
so that at any moment assumption \(assum) holds\/{\rm:}
\ifMag\else\vv-.8>\fi
$$
\gather
\Psi_{YW}(z)\,=\;{\ell\>!\over\ell_1!\ldots\ell_k!}\;
\prod_{i=1}^k\,\prod_{1\le\,j\le n\_{i-1}}\xi_j^{\ell_i}\>\>Y(z)\,
\bigl(\Ome_{\>W}(z)+o(1)\bigr)
\\
\Text{where}
\ifMag\nn2>\fi
\aligned
\Ome_{\>W}(\zn)\,& {}=\,
\Psi_{W_1}(z_1\lc z_{n_1})\lox\Psi_{W_k}(z_{n_{k-1}+1}\lc z_n)\;+
\\
\nn6>
&\>+\!\sum_{\mg\gg(\ell_1\lc\ell_k)}\)
\tprod_{i=1}^k N\"i_{W\!,\mg}(z_{n_{i-1}+1}\lc z_{n_i})\,\Fvv
\endaligned
\endgather
$$
for suitable coefficients
$N\"i_{W\!,\mg}(z_{n_{i-1}+1}\lc z_{n_i})$.
\endpro
\nt
Theorem~\[Asol] follows from Theorem~\[AsIWw] below which describes \as/s of
the \hpair/.
\par
Let $\ell'_1\lc\ell'_k$ be nonnegative integers \st/
$\sum_{i=1}^k\ell'_i=\ell$.
\vvn-.4>
Say that ${(\ell'_1\lc\ell'_k)\]}\gg\alb{\](\ell_1\lc\ell_k)}\,$ if
$\sum_{j=1}^i\ell'_j\ge\sum_{j=1}^i\ell_j$ for any $i=1\lc k-1$, and
$\,(\ell'_1\lc\ell'_k)\ne(\ell_1\lc\ell_k)$.
\Th{AsIWw}
{\bls 1.1\bls
Let the parameters $\xin$ obey condition \(Lass).
Let ${w_i\in\Fo^i[\)\ell'_i\)]}$ and ${W_i\in\Fq^i[\)\ell_i\)]}$, $i=1\lc k$.
Let $w=w_1\lsym\*w_k$ and $W=W_1\lsym*W_k$. Then the \hint/ $I(W,w)$ has the
following \as/s as $z$ tends to limit in the \azo/, $z\rto\AA$, so that
at any moment assumption \(assum) holds\/{\rm:}\par\nt}
\ifMag\else\vv-.5>\fi
$$
I(W,w)\,=\;{\ell\>!\over\ell_1!\ldots\ell_k!}\;
\prod_{i=1}^k\,\prod_{1\le\,j\le n\_{i-1}}\xi_j^{\ell_i}\;
\bigl(\>\tprod_{i=1}^k I(W_i,w_i)+o(1)\bigr)
$$
for $(\ell'_1\lc\ell'_k)=(\ell_1\lc\ell_k)$ and
$$
\alignat2
& I(W,w)\,=\,O(1) && \for (\ell'_1\lc\ell'_k)\]\gg\](\ell_1\lc\ell_k)\,,
\\
\nn2>
& I(W,w)\,=\,o(1)\,,\qqq && \text{otherwise}\,.
\endalignat
$$
\endpro
\nt
The theorem is proved in \Chap/~\SNo{4}.

\ifAst\newpage\fi
\makePcd
\Sect[4]{Proofs}
This \chap/ contains proofs of the statements formulated in
\Chap/s~{\SNo{2}\,--\,\SNo{2d}}. Basic facts about the \tri/ and \ehgf/s are
given in Appendices~\hbox{\SNo{A} and \SNo{B}}, \resp/. In particular, we give
there proofs of Lemmas~\[wlgp], \[Wlgp], \[ResQ]\,\)--\,\[dimFQ'].
\par
\Pf of Lemmas\{~\{\[wbasis], \{\[Wbasis].
The statements immediately follows from Theorems~\[DetM], \[DetMq], \resp/.
\epf
\Pf of Lemmas\{~\{\[wDw], \{\[wDw'].
The first claims of the lemmas are \resp/ \eqv/ to the first and second
formulae in \(wDw0). The second claims of the lemmas are the same as
Corollary~\[F/R].
\epf
Lemma~\[T12] follows from formula (A.3) in \Cite{IK} and the definition of the
\emod/s by induction \wrt/ the number of factors in the tensor product.
Nevertheless, we give here an independent proof of Lemma~\[T12] which is
similar to the proof of Lemma~\[L21] in the elliptic case.
\Pf of Lemma\{~\{\[T12].
\Wlg/ we can assume that $\tau$ is the identity \perm/. We give a proof of
the second formula of the lemma. The proof of the first formula is similar.
\par
It suffices to prove the formula for generic values of parameters $\eta$,
$\xin$, $\zn$, since both sides of the formula are \anf/s of the parameters.
\par
Define \fn/s $X_\lg(\tell)$, $\lg\in\Zln$, by the rule:
\ifMag
$$
\align
L^-_{12}(t_1) &{}\ldots L^-_{12}(t_\ell)\>\vox\;=
\\
&\Cph{{}\ldots{}}{{=}}\!\!\sum_{\,\)\lg\in\Zln\!\]}^{\vp1} X_\lg(\tell)\!\!
\plmn\!\!\! q^{\lg_m\La_l-\lg_l\La_m-\lg_l\lg_m}\pral t_a^{\ell(1-n)}\>\Fv\,.
\endalign
$$
\else
$$
L^-_{12}(t_1)\ldots L^-_{12}(t_\ell)\>\vox\,=
\sum_{\,\)\lg\in\Zln\!\]} X_\lg(\tell)\!\!
\plmn\!\!\! q^{\lg_m\La_l-\lg_l\La_m-\lg_l\lg_m}\pral t_a^{\ell(1-n)}\Fv\,.
$$
\fi
The claim of the lemma means that for any $\lg\in\Zln$ the \fn/ $X_\lg$
coincides with the \pol/ $P_\lg$ defined by formula \(Plg).
\par
For $\lg,\mg\in\Zln$ say that ${\lg\lle\mg}$ if
\vv-.4>
$\,\sum_{i=1}^m\lg_i\le\sum_{i=1}^m\mg_i\;$ for any $\mn-1$.
Say that $\lg\ll\mg$ if $\lg\ne\mg$ and $\lg\lle\mg$.
\vsk.1>
By the definition of \fn/s $X_\lg$, they are \sym/ \pol/s in $\ell$ \var/s of
degree less than $n$ in each of the \var/. Hence, they are linear combinations
of \pol/s $P_\lg$:
$$
X_\lg(t)\,=\,(q-\q)\vpb\ell\smiZ U_{\lg\mg}\>P_\mg(t)\,,\qqq\lg\in\Zln.
$$
By Lemmas~\[Xxy], \[Pxy] the matrix $U$ is upper triangular: $U_{\lg\mg}=0$
unless $\lg\lle\mg$,
\vv.1>
as a ratio of upper triangular matrices, and simultaneously lower triangular:
$U_{\lg\mg}=0$ unless $\lg\gge\mg$, as a ratio of lower triangular matrices.
Therefore, the matrix $U$ is diagonal. Moreover, by the same lemmas
$U_{\lg\lg}=1$ for any $\lg\in\Zln$. Hence, $U$ is the unit matrix.
Lemma~\[T12] is proved.
\epf
For any $\lg\in\Zln$ let $\xt\lg,\,\yt\lg\in\Cxl\!$ be the following points:
\ifMag
$$
\NN2>
\align
\xt\lg\,=\,
(\eta^{1-\lg_1}\xi_1z_1,\>\eta^{\)2-\lg_1}\xi_1z_1\lc\xi_1z_1,
\>\eta^{1-\lg_2}\xi_2z_2\lc\xi_2z_2,\,\ldots\,,{} &
\\
\eta^{1-\lg_n}\xi_nz_n\lc{} & \xi_nz_n)\,,
\\
\nn4>
\yt\lg\,=\,(\eta^{\lg_1-1}\xi_1\1z_1,\>\eta^{\lg_1-2}\xi_1\1z_1\lc\xi_1\1\]z_1,
\>\eta^{\lg_2-1}\xi_2\1z_2\lc\xi_2\1\]z_2,\,\ldots\,,{} &
\\
\eta^{\lg_n-1}\xi_n\1z_n\lc{} &\xi_n\1z_n)\,,
\endalign
$$
\else
$$
\gather
\xt\lg\,=\,
(\eta^{1-\lg_1}\xi_1z_1,\>\eta^{\)2-\lg_1}\xi_1z_1\lc\xi_1z_1,
\>\eta^{1-\lg_2}\xi_2z_2\lc\xi_2z_2,\,\ldots\,,
\)\eta^{1-\lg_n}\xi_nz_n\lc\xi_nz_n)\,,
\\
\nn4>
\yt\lg\,=\,(\eta^{\lg_1-1}\xi_1\1z_1,\>\eta^{\lg_1-2}\xi_1\1z_1\lc\xi_1\1\]z_1,
\>\eta^{\lg_2-1}\xi_2\1z_2\lc\xi_2\1\]z_2,\,\ldots\,,
\)\eta^{\lg_n-1}\xi_n\1z_n\lc\xi_n\1z_n)\,,
\endgather
$$
\fi
\cf. \(xtyt), \(xtry).
\Lm{Xxy}
$X_\lg(\xt\mg)=0$ unless $\lg\lle\mg$. $X_\lg(\yt\mg)=0$ unless $\lg\gge\mg$.
Moreover,
\vvn-.5>
$$
X_\lg(\xt\lg)\,=\,(q-\q)\vpb\ell\>\pron\,\prod_{s=0}^{\lg_m-1}
\Bigl(\>\prlm(\eta^{-s}\xi_mz_m-\xi_l\1z_l)
\prml(\eta^{\lg_l-s}\xi_mz_m-\xi_lz_l)\>\Bigr)\,.
$$
\endpro
\Pf.
The proof is given by the straightforward calculation based on the
definition of the coproduct in the \qlo/ $\Ugg$. We illustrate the calculation
by the following example.
\par
Let $n=3$. Let ${\Dl\vpb{\!(3)}\}=(\Dl\ox\id)\o\Dl:\Ugg\to\Ugg^{\ox3}}$ be the
iterated coproduct. We have that
$$
\align
\ifMag\qquad\fi
\Dl\vpb{\!(3)}(L^-_{12}(t))\,=\,L^-_{11}(t)\ox L^-_{11}(t)\ox L^-_{12}(t)+
L^-_{11}(t)\ox L^-_{12}(t)\ox L^-_{22}(t)+{} &
\Tagg{DlL}
\\
\nn4>
{}+L^-_{12}(t)\ox L^-_{22}(t)\ox L^-_{22}(t)+
L^-_{12}(t)\ox L^-_{21}(t)\ox L^-_{12}(t) &\,.
\endalign
$$
Let $\ell=4$, $\mg=(1,1,2)$. Recall that $\eta=q^2$ and $\xi_m=q^{2\La_m}$.
\par
We have to calculate the following expressions
$$
\gather
L^-_{12}(\xi_1z_1)\>L^-_{12}(\xi_2z_2)\>L^-_{12}(\eta\1\xi_3z_3)\>
L^-_{12}(\xi_3z_3)\>v_1\ox v_2\ox v_3\,,
\Tag{1233}
\\
\nn4>
L^-_{12}(\xi_1\1z_1)\>L^-_{12}(\xi_2\1z_2)\>L^-_{12}(\eta\)\xi_3\1z_3)\>
L^-_{12}(\xi_3\1z_3)\>v_1\ox v_2\ox v_3\,.
\Tag{3321}
\endgather
$$
To compute ${L^-_{12}(\eta\1\xi_3z_3)\>L^-_{12}(\xi_3z_3)\>v_1\ox v_2\ox v_3}$
we need only the first term in \rhs/ of formula \(DlL), since all other terms
vanish. Then to compute ${L^-_{12}(\xi_2z_2)\>L^-_{12}(\eta\1\xi_3z_3)\>
L^-_{12}(\xi_3z_3)\>v_1\ox v_2\ox v_3}$ we need only the first and second terms
in \rhs/ of formula \(DlL) for the same reason. Finally, at the last step we
have to use all four terms in \(DlL) and we find that expression \(1233) is
a linear combination of vectors
\ifMag
$$
\gather
Fv_1\ox Fv_2\ox F^2v_3\,,\quad Fv_1\ox v_2\ox F^3v_3\,,
\\
\nn3>
v_1\ox F^2v_2\ox F^2v_3\,,\quad v_1\ox Fv_2\ox F^3v_3\,,\quad
v_1\ox v_2\ox F^4v_3\,.
\endgather
$$
\else
$$
Fv_1\ox Fv_2\ox F^2v_3\,,\quad Fv_1\ox v_2\ox F^3v_3\,,\quad
v_1\ox F^2v_2\ox F^2v_3\,,\quad v_1\ox Fv_2\ox F^3v_3\,,\quad
v_1\ox v_2\ox F^4v_3\,.
$$
\fi
The coefficient of vector ${Fv_1\ox Fv_2\ox F^2v_3}$ can be easily
calculated and it has the prescribed form.
\par
To calculate expression \(3321) we first use the commutativity of the factors
in the product and transform the expression as follows:
$$
L^-_{12}(\eta\)\xi_3\1z_3)\>L^-_{12}(\xi_3\1z_3)\>L^-_{12}(\xi_2\1z_2)\>
L^-_{12}(\xi_1\1z_1)\>v_1\ox v_2\ox v_3\,.
$$
Then to compute ${L^-_{12}(\xi_1\1z_1)\>v_1\ox v_2\ox v_3}$ we need only the
third term in \rhs/ of formula \(DlL), and to compute
${L^-_{12}(\xi_2\1z_2)\>L^-_{12}(\xi_1\1z_1)\>v_1\ox v_2\ox v_3}$ we need only
the second and third terms. The rest of the calculation is clear.
\par
The general case can be considered similarly.
\epf
\Pf of Lemma\{~\{\[L21].
The proof is similar to the proof of Lemma~\[T12]. So we give only the main
points of the proof.
\par
\Wlg/ we can assume that parameters $\eta$, $\xin$, $\zn$ are generic, and
$\tau$ is the identity \perm/. Define \fn/s $\Xj_\lg(\tell)$, $\lg\in\Zln$,
by the rule:
\ifMag
$$
\align
\qquad& T_{21}(t_1,\la)\ldots T_{21}(t_\ell,\eta^{\ell-1}\la)\>\Fvq\;=
\Tagg{Xjlg}
\\
\nn3>
& {}=\,\Xj_{\lg}(\tell)\,\pros\>\tht(\eta^s\kzy)\,
\pral\,\pron\>\tht(\xi_m\1t_a/z_m)\,\vqx\,.
\endalign
$$
\else
$$
\align
T_{21}(t_1,\la)\ldots T_{21}(t_\ell,\eta^{\ell-1}\la) &\>\Fv\;=
\Tag{Xjlg}
\\
\nn3>
{}=\,\Xj_{\lg}(\tell) &\,\pros\>\tht(\eta^s\kzy)\,
\pral\,\pron\>\tht(\xi_m\1t_a/z_m)\,\vox\,.
\endalign
$$
\vvv-.5>
\fi
Here $\la=\ka\zy$. We have to show that $\Xj_\lg=c_\lg\:\>J_\lg$ for any
$\lg\in\Zln$, where the constant $c_\lg\:$ and the \fn/ $J_\lg$ are given
by formulae \(clg) and \(Pilg), \resp/.
\par
By the definition of \fn/s $\Xj_\lg$, they are \hof/ on $\Cxl\!$ having
the property
$$
\Xj_\lg(\tpell)\,=\,\ka\tpron z_m\,(-t_a)^{-n}\>\Xj_\lg(\tell)\,.
$$
Therefore, they are linear combinations of the \fn/s $J_\lg$:
$\dsize\Xj_\lg(t)\,=\sum_{\mg\in\Zln}\Ue_{\lg\mg}\>J_\mg(t)$,
$\ \lg\in\Zln$.
Lemmas~\[Xixy] and~\[Jxy] imply that the matrix $\Ue$ is diagonal and
$\Ue_{\lg\lg}=c_\lg\:$, $\,\lg\in\Zln$. Lemma~\[L21] is proved.
\epf
\Lm{Xixy}
$\Xj_\lg(\xt\mg)=0$ unless $\lg\lle\mg$.
$\Xj_\lg(\yt\mg)=0$ unless $\lg\gge\mg$. Moreover,
$$
\NN6>
\align
\Xj_\lg(\xt\lg)\,=\plmn\!\}\eta^{-\lg_l\lg_m}\xi_m^{2\lg_l}\,
\pron\,\prod_{s=0}^{\lg_m-1}\Bigl(\>
{\tht(\eta^s)\>\tht(\eta^{1-s}\xi_m^2)\over\tht(\eta)\>
\tht\bigl(\eta^{s+1}\ka\1\!\!\!\prllm\!\!\eta^{\lg_l}\xi_l\1\!
\prml\!\!\eta^{-\lg_l}\xi_l\)\bigr)}\ \x{} &
\\
\x\;\prlm\!\}\tht(\eta^{-s}\xi_l\xi_mz_m/z_l)\!
\prml\!\}\tht(\eta^{\lg_l-s}\xi_l\1\xi_mz_m/z_l)\>\Bigr) & \,.
\endalign
$$
\ifMag\else\vvv-.6>\fi
\endpro
\Pf.
The proof is similar to the proof of Lemma~\[Xxy] and is given by the
straightforward calculation based on the definition of the coproduct in the
\eqg/ $\Eqg$. Two remarks on the calculation is to be done.
\par
The calculation becomes more transparent if it is done in the dual picture,
that is if we replace formula \(Xjlg) by the dual one
\ifMag
$$
\gather
\Lline{\bigl(T_{21}(t_1,\la)\ldots T_{21}(t_\ell,\eta^{\ell-1}\la)\bigr)^*\>
\bigl(\vqx\bigr)^*\,=}
\\
\nn3>
\Rline{=\sum_{\,\)\lg\in\Zln\!\]}\Xj_{\lg}(\tell)\,\pros\>\tht(\eta^s\kzy)\,
\pral\,\pron\>\tht(\xi_m\1t_a/z_m)\,\bigl(\Fvq\bigr)^*.}
\endgather
$$
\else
$$
\align
\bigl(T_{21}(t_1,\la)\ldots T_{21}(t_\ell,\eta^{\ell-1}\la) &\bigr)^*\>
(\vox)^*\,=
\\
\nn3>
=\sum_{\,\)\lg\in\Zln\!\]}\Xj_{\lg}(\tell) &\,\pros\>\tht(\eta^s\kzy)\,
\pral\,\pron\>\tht(\xi_m\1t_a/z_m)\,\Fva.
\endalign
$$
\fi
\par
The factors in the product in \lhs/ of this formula should be put in the
suitable order, which can be done using commutation relations in the \eqg/
\vvgood
$\Eqg$. For instance, if $\ell=4$ and $\mg=(1,1,2)$, then for the point
$\xt\mg$ the suitable form of the corresponding product is
$$
\align
T_{21}(\xi_3z_3,\la)\> &T_{21}(\eta\1\xi_3z_3,\eta\)\la)\>
T_{21}(\xi_2z_2,\eta^2\la)\>T_{21}(\xi_1z_1,\eta^3\la)\;=
\\
\nn3>
=\; &T_{21}(\xi_1z_1,\la)\>T_{21}(\xi_2z_2,\eta\)\la)\>
T_{21}(\eta\1\xi_3z_3,\eta^2\la)\>T_{21}(\xi_3z_3,\eta^3\la)
\endalign
$$
and for the point $\yt\mg$ the suitable form of the corresponding product is
$$
\align
T_{21}(\xi_1\1z_1,\la)\> &T_{21}(\xi_2\1z_2,\eta\)\la)\>
T_{21}(\xi_3\1z_3,\eta^2\la)\>T_{21}(\eta\)\xi_3\1z_3,\eta^3\la)\;=
\\
\nn3>
=\; &T_{21}(\xi_1\1z_1,\la)\>T_{21}(\xi_2\1z_2,\eta\)\la)\>
T_{21}(\eta\)\xi_3\1z_3,\eta^2\la)\>T_{21}(\xi_3\1z_3,\eta^3\la)\,.
\endalign
$$
The necessary transformation of the product in the general case is similar.
\epf
\Pf of Lemmas\{~\{\[FFo], \{\[FFq].
The statements follow from definitions \(wlga), \(Wlga) of the \wtf/s and
Theorems~\[DetM], \[DetMq], \resp/.
\epf
We extend the notion of the \hint/ $I(W,w)$ and consider the \hint/ for any
\fn/ $w$ in the \fn/al space $\Fwh(z)$ of a fiber. Namely, let
$w(t,z)\in\Fwh(z)$ be a \fn/ of the form
$$
\align
P(\tell,\zn,\xin,\eta)\,\prod_{s=0}^r\,\Bigl[\,\pron\,
\pral\,{1\over(p^st_a-\xi_mz_m)\>(\xi_m t_a-p^{s+1}z_m)}\;\,\x &
\\
\nn4>
\x\prab\ {1\over(p^s\eta\)t_a-t_b)\>(t_a-p^{s+1}\eta\)t_b)}\,\Bigr]\> &
\endalign
$$
where $P$ is a Laurent \pol/. If $|z_m|=1$ for any $\mn$, the absolute values
of the parameters $\xin$ are small and the absolute value of the parameter
$\eta$ is large,
\aagood
then we define the \hint/s $I(W_\lg,w)$ by formula \(IWw).
For generic $\eta$, $\xin$, $\zn$ we define the \hint/s $I(W_\lg,w)$ by
the \anco/ \wrt/ $\eta$, $\xin$, $\zn$. Similar to Theorem~\[Ianco] one can
show that this \hint/s can be analytically continued as \hol/ univalued \fn/s
of complex \var/s $\eta$, $\xin$, $\zn$ to the region described in
Theorem~\[Ianco]. For arbitrary \fn/s $w\in\Fwh(z)$, $W\in\Fq(z)$ we define the
\hint/ $I(W,w)$ by linearity.
\Par
Let $D\Fwh(z)=\lb\)Dw\vert w\in\Fwh(z)\)\rb$.
\Lm{IWwD}
Let $0<|p|<1$. Let \(npZ)\;--\;\(assum) hold. Then the \hint/ $I(W,w)$ equals
zero for any \fn/ $w\in D\Fwh(z)$.
\endpro
\Pf.
The claim is clear if $|z_m|=1$ for any $\mn$, the absolute values of the
parameters $\xin$ are small and the absolute value of the parameter $\eta$
is large. For general $\eta$, $\xin$, $\zn$ the claim is proved by the \anco/.
\epf
\Pf of Lemma\{~\{\[IWw=0].
The first claim of the lemma follows from Lemma~\[wDw] and \[IWwD].
\par
It suffices to prove the second claim under the assumptions that $|z_m|=1$,
$\mn$, ${|\eta\)|>1}$ and the absolute values of $\xin$ are small, when the
\hint/ $I(W,w)$ is given by formula \(IWw).
\Par
Let $W\in\Qc(z)$, that is
$$
W(\tell)\,=\susi\,\Lbc\>W'(\twll)\>\Rbc_\si
$$
for a suitable \fn/ ${W'\in\Fq{[\)\koNe;\zn;\xin;\ell-1\)]}}(z)$.
Due to formula \(Phisym), we have that
$$
I(W,w)\,=\,\ell\)!\)
\int_{\){\TT^\ell}\]}\Phi(\tell)\>w(\tell)\>W'(\twll)\;\dtt
\Tag{IWw!}
$$
because the torus $\TT^\ell$ is invariant under \perm/s of the \var/s $\tell$.
\par
Since $w(0,\twll)=0$ for any $w\in\Fo(z)$, the integrand
${\Phi(t)\>w(t)\>W'(t)}$ considered as a \fn/ of $t_1$ is regular in
the disk $|t_1|\le 1$. Hence, $I(W,w)=0$ for any $w\in\Fo(z)$.
\epf
\Pf of Lemma\{~\{\[IWw=0'].
The proof is similar to the proof of Lemma~\[IWw=0]. For the proof of the first
claim Lemma~\[wDw] is to be replaced by Lemma~\[wDw']. In the proof of
the second claim the corresponding integrand is regular outside the disk
$|t_\ell|\ge 1$ decreasing as $O(t_\ell^{-2})$ at infinity.
\epf
The \hint/ defines linear \fn/als $I(W,\cdot)$ on the \fn/al space of a fiber.
Lemma~\[IWwD] means that these linear \fn/als can be considered as elements of
the top homology group $H_\ell(z)$, the dual space to the top cohomology group
of the de~Rham complex of the discrete \loc/ of the fiber.
\Pf of Theorem\{~\{\[sWz].
Recall, that in general the definition of the \hint/ depends on $z$.
In this proof we will indicate this dependence explicitly as a subscript:
$I(\cdot{,}\cdot)=I_z(\cdot{,}\cdot)$.
g\par
The section $s_W\:$ is defined by $\,{s_W\:(z)\,=\,I_z(W|_z,\cdot)}\,$ where
${W|_z\in\Fq(z)}$ denotes the restriction of the \fn/ ${W\in\Fq}$ to the fiber
over $z$. The theorem is a direct corollary of the quasiperiodicity of the
\fn/ $W$ \wrt/ each of the \var/s $\zn$:
$$
W(t,\zmn)\,=\,\xi_m^\ell\>W(t,\zn)\,,\qqq\mn\,.
$$
Namely, the periodicity of the section $s_W\:$ \wrt/ the translation
${z_m\map p\)z_m}$ means that
$$
I_{z'}(W|_{z'},w)\,=\,I_z(W|_z,(\phi_{\ell+m})|_z\>w)
\Tag{Iz}
$$
for any $w\in\Fwh(z')$. Here ${z'=(\zmn)}$ and $\phi_{\ell+m}$ is
the corresponding \conn/ of the \stype/ \loc/, see \(conns).
\par
\Wlg/ we can assume that both $w$ and $W$ are \mef/ of the parameters $\xin$
and $\eta$. So it suffices to prove \(Iz) under the assumption that the
absolute values of $\xin$ are small and the absolute value of $\eta$ is large.
Then, the \hint/ is given by formula \(IWw) and we have
$$
\align
I_{z'}(W|_{z'},w)\,=\int_{\){\TT^\ell}\]}\Phi(t,z')\>w(t)\>W(t,z')\;\dtt\,={}&
\\
\nn4>
{}=\int_{\){\TT^\ell}\]}\Phi(t,z)\>\phi_{\ell+m}(t,z)\>w(t)\>W(t,z)\;\dtt
\,={} & \,I_z(W|_z,(\phi_{\ell+m})|_z\>w)
\endalign
$$
The middle equality reflects the fact that the product ${\Phi\>W}$ is a \phf/
of \sconn/ \(conns).
\epf
\Pf of Lemmas\{~\{\[Vsing], \{\[Vsing'].
We give here a proof of Lemma~\[Vsing]. The proof of Lemma~\[Vsing'] is
similar.
\par
\Wlg/ we can assume that $\tau=\id$. Let
$$
\gather
v^\lg\>=\,\Fv\,,\qqq
b_\lg\,=\,\tpron\) q^{\lg_m(\lg_m-1)/2\,+\)\lg_m\La_m},
\\
\nn4>
\BE(t,z)\,=\,(q-\q)
\sliZ b_\lg\>w_\lg(t,z)\>E\)v^\lg\in\Fo(z)\ox(\Vox)_{\ell-1}\:\,.
\endgather
$$
Here $E$ is a generator of $\Uu$ acting in $\Vox$.
\par
By the definition of the \tenco/ $B_\id(z)$, \cf. \(Btau), and the map
$\Ib(z)$, the claim of Lemma~\[Vsing] is \eqv/ to the following statement:
$$
I(W,\BE)\,=\,0\,.
$$
Let $\,\eg(m)=(0\lc\%1_{\sss\^{$m$-th}}\lc 0)$, $\,\mn$. Since
\ifMag
$$
\align
(q & {}-\q)\>E\)v^\lg\>={}
\\
\nn4>
& {}=\,\sun\,(q^{\lg_m}\}-q^{-\lg_m})\>
(q^{\)2\La_m-\lg_m+1}\]-q^{-2\La_m+\lg_m-1})\>\tprlm q^{\)\La_l-\lg_l}
\tprml q^{-\La_l+\lg_l}\,v^{\lg-\eg(m)},
\endalign
$$
\else
$$
(q-\q)\>E\)v^\lg\>=\,\sun\,(q^{\lg_m}-q^{-\lg_m})\>
(q^{2\La_m-\lg_m+1}\!-q^{\lg_m-2\La_m-1})\>\tprllm q^{\La_l-\lg_l}
\tprml q^{\lg_l-\La_l}\cdot v^{\lg-\eg(m)},
$$
\fi
and recalling that $\eta=q^2$, $\xi_m=q^{2\La_m}$, $\mn$, we obtain
$$
\BE\,=\,-\,q^{\ell-1-\!\!\sun\!\!\La_m}\sleZ\;\Bigl(\>
\sun w_{\lg+\eg(m)}\>(1-\eta^{\lg_m+1})\>(\xi_m-\eta^{\lg_m}\xi_m\1)\>
\tprllm\!\eta^{-\lg_l}\xi_l\Bigr)\,b_\lg\> v^\lg\,.
$$
Therefore, $\BE\in\Rc(z)\ox(\Vox)_{\ell-1}\:$, see Lemma~\[wDw], and applying
Lemma~\[IWw=0] we complete the proof.
\epf
Our further strategy is as follows. First we prove Theorem~\[AsIWw] which,
together with formula \(qbeta), implies Theorem~\[asol]. Using Theorem~\[asol]
we prove that the \hpair/ is \ndeg/, \cf. Theorem~\[mu<>0]\,--\,\[mu=02].
At last, we prove Theorems~\[hmapo], \[hmapp].
\Pf of Theorem\{~\{\[AsIWw].
To simplify notations we give a proof only for the case $k=n$, so that $n_m=m$,
$\mn$. The general case is similar.
\par
Let ${w\"m\'l\in\Fo[z_m;\xi_m;l\,]}$ and
${W\"m\'l[\al]\in\Fq[\al;z_m;\xi_m;l\,]}$ be the following \fn/s:
$$
\gather
w\"m\'l(t_1\lc t_l,z_m)\,=\,\prod_{s=1}^l\,{1-\eta\over 1-\eta^s}\,
\sum_{\si\in\S^l}\,\Bigl[\,\prod_{a=1}^l\,
{t_a\over t_a-\xi_mz_m}\,\Bigr]_\si,\qqq
\Tag{wWprim}
\\
\nn8>
W\"m\'l[\al](t_1\lc t_l,z_m)\,=\,\prod_{s=1}^l\,{\tht(\eta)\over\tht(\eta^s)}\,
\sum_{\si\in\S^l}\,\LBc\,\prod_{a=1}^l\,
{\tht(\eta^{\)2a-l-1}\al\1t_a/z_m)\over\tht(\xi_m\1t_a/z_m)}\,\RBc_\si,
\endgather
$$
\cf. \(wlga), \(Wlga). We have the equalities
$$
w_\lg\,=\,w\"1\'{\lg_1}\lsym\* w\"n\'{\lg_n}\qquad\text{and}\qquad
W_\lg\,=\,W\"1\'{\lg_1}[\ka_{\lg,1}]\lsym* W\"n\'{\lg_n}[\ka_{\lg,n}]
$$
where ${\,\ka_{\lg,m}=\ka\!\prod_{1\le i<m}\!\}\eta^{-\lg_i}\xi_i\:\!
\prod_{m<i\le n}\!\}\eta^{\lg_i}\xi_i\1}$.
\vvn.1>
Therefore, we have to study the \as/s of the \hint/s $I(W_\lg,w_\mg)$.
\par
Consider the \hint/ $I(W_\lg,w_\mg)$. Due to property \(Phisym) all the terms
in formula \(Wlga) for the \fn/ $W_\lg$
\aagood
give the same contribution to the integral. So we can replace the integrand
$\Phi(t)\>w_\mg(t)\>W_\lg(t)$ by the following integrand give the same
contribution to the integral. So we can replace the integrand
$\Phi(t)\>w_\mg(t)\>W_\lg(t)$ by the following integrand
\ifAst
$$
\align
F(t)\, &{}=\,\ell\)!\,w_\mg(t)\!
\prab\;{(\eta\)t_a/t_b)\9\!\over(\eta\1t_a/t_b)\9\!}\;
\pron\,\Bigl(\,\prod_{s=1}^{\lg_m}\,{\tht(\eta)\over\tht(\eta^s)}\ \x
\\
\nn4>
{}\x\ {} & \Cph{{}={}}{\prod_{a\in\Gm_{\Rph lm}}}\ \;
{\tht(\ka_{\lg,m}\1\)t_a/z_m)\over
(p)\9\>(\xi_mt_a/z_m)\9\>(p\>\xi_mz_m/t_a)\9\!}
\;\prlm{(p\>\xi_l\1z_l/t_a)\9\over(p\>\xi_lz_l/t_a)\9}
\>\prml{(\xi_l\1t_a/z_l)\9\over(\xi_lt_a/z_l)\9}\,\Bigr)\,,
\endalign
$$
\else
$$
\align
F(t)\, &{}=\,\ell\)!\,w_\mg(t)\!
\prab\;{(\eta\)t_a/t_b)\9\!\over(\eta\1t_a/t_b)\9\!}\;
\pron\,\Bigl(\,\prod_{s=1}^{\lg_m}\,{\tht(\eta)\over\tht(\eta^s)}\ \x
\\
\nn4>
&{}\>\x\prod_{a\in\Gm_{\Rph lm}}\>{\tht(\ka_{\lg,m}\1\)t_a/z_m)\over
(p)\9\>(\xi_mt_a/z_m)\9\>(p\>\xi_mz_m/t_a)\9\!}
\;\prlm{(p\>\xi_l\1z_l/t_a)\9\over(p\>\xi_lz_l/t_a)\9}
\>\prml{(\xi_l\1t_a/z_l)\9\over(\xi_lt_a/z_l)\9}\,\Bigr)
\endalign
$$
\fi
where $\Gm_m=\lb 1+\lg^{m-1}\,\lc\lg^m\rb$, $\mn$.
\par
Assume that $\eta>1$ and $|\xi_m|<1$ for any $\mn$. If $|z_m|=1$ for all
$\mn$, then we have
$$
I(W_\lg,w_\mg)\,=\,\int_{\TT^\ell}\>F(t)\;\dtt\,.
$$
The \anco/ of $I(W_\lg,w_\mg)$ to the region $|z_1|\lsym<|z_n|$ is given by
$$
\gather
I(W_\lg,w_\mg)\,=\int_{\TT^{\lg_1}_1\lx\TT^{\lg_n}_n}\!\!\!F(t)\;\dtt
\\
\nn-3>
\Text{where}
\nn-3>
\TT^{\lg_m}_m\,=\,\lb\)(t_{1+\lg^{m-1}}\lc t_{\lg^m})\in\C^{\,\lg_m}\vert
\ |t_a|=|z_m|\,,\ \,\lg^{m-1}<a\le\lg^m\)\rb
\endgather
$$
since the integrand has no poles at the hyperplanes $t_a=p^{-s}\xi_l\1z_l$,
$s\in\Z$, for $\lg^{m-1}<a\le\lg^m$, $m>l$, has no poles at the hyperplanes
$t_a=p^s\xi_lz_l$, $s\in\Z$, for $\lg^{m-1}<a\le\lg^m$, $m<l$, and has no
poles at the hyperplanes ${t_a=p^s\eta\)t_b}$, $s\in\Z$, for $a>b$.
\par
Let $z$ tends to limit in the \azo/, $z\rto\AA$, that is
$|z_m/z_{m+1}|\to 0$ for all $\mn-1$. Consider the case $\lg=\mg$. Transform
the \hint/ $I(W_\lg,w_\lg)$ as above and replace the integrand by its
asymptotics as $z\rto\AA$. Since
$$
w_\lg(\tell)\,=\plmn\!\!\xi_l^{\lg_m}\>
\pron w\"m\'{\lg_m}(t_{\lg^{m-1}+1}\lc t_{\lg^m})\,+\,o(1)
\Tag{wAA}
$$
as $z\rto\AA$ and $t\in\TT^{\lg_1}_1\lx\TT^{\lg_n}_n$, we obtain that
\ifAst
$$
\align
I(W_\lg,w_\lg)\, &{}=\,\ell\)!\,
\pron\,\Bigl(\;\prod_{s=1}^{\lg_m}\,{\tht(\eta)\over\tht(\eta^s)}
\;\int_{\TT^{\lg_m}_m} w\"m\'{\lg_m}(t_{\lg^{m-1}+1}\lc t_{\lg_m})\,\x{}
\\
\nn6>
{}\x\ {} & \Cph{{}={}}{\prod_{a\in\Gm_{\Rph lm}}}\ \Bigl(
\bigl((\xi_mt_a/z_m)\9\>(p\>\xi_mz_m/t_a)\9\bigr)\1\!
\prod_{\tsize{b<a\atop b\in\Gm_{\Rph lm}}}
{(\eta\)t_b/t_a)\9\over(\eta\1t_b/t_a)\9\!}\;{dt_a\over t_a}\;\Bigr)\Bigr)
\,\ono
\endalign
$$
\else
$$
\align
I(W_\lg,w_\lg)\, &{}=\,\ell\)!\,
\pron\,\Bigl(\;\prod_{s=1}^{\lg_m}\,{\tht(\eta)\over\tht(\eta^s)}
\;\int_{\TT^{\lg_m}_m} w\"m\'{\lg_m}(t_{\lg^{m-1}+1}\lc t_{\lg_m})\,\x{}
\\
\nn6>
&{}\>\x\prod_{a\in\Gm_{\Rph lm}}\Bigl(
\bigl((\xi_mt_a/z_m)\9\>(p\>\xi_mz_m/t_a)\9\bigr)\1\!
\prod_{\tsize{b<a\atop b\in\Gm_{\Rph lm}}}
{(\eta\)t_b/t_a)\9\over(\eta\1t_b/t_a)\9\!}\;{dt_a\over t_a}\;\Bigr)\Bigr)
\,\ono\,,
\endalign
$$
\fi
as $z\rto\AA$. Due to \(Phisym) the integrals over $\TT^{\lg_m}_m$ are
the \hint/s $I(W\"m\'{\lg_m}[\ka_{\lg,m}],w\"m\'{\lg_m})$ up to simple factors.
Hence, we finally obtain that
$$
I(W_\lg,w_\lg)\,=\;{\ell\)!\over\lg_1!\ldots\lg_n!}\,\plmn\!\!\xi_l^{\lg_m}\,
\bigl(\)\tpron I(W\"m\'{\lg_m}[\ka_{\lg,m}],w\"m\'{\lg_m})+o(1)\bigr)\,.
$$
The \hint/ $I(W_\lg,w_\mg)$ for $\lg\ne\mg$ can be treated similarly to the
\hint/ $I(W_\lg,w_\lg)$ considered above. The final answer is
$$
\alignat2
& I(W_\lg,w_\mg)\,=\,O(1) && \for \lg\ll\mg\,,
\\
\nn2>
& I(W_\lg,w_\mg)\,=\,o(1)\,,\qqq && \text{otherwise}\,,
\endalignat
$$
which completes the proof if $\eta>1$ and $|\xi_m|<1$, $\mn$.
\Par
For general $\xin$ and $\eta$ the proof is similar. The \anco/ of
$I(W_\lg,w_\mg)$ to the region $|z_1|\lsym\ll|z_n|$ is given by
$$
\vp{\Plus_\TT}\smash{
I(W_\lg,w_\mg)\,=\int_{\TTt^{\lg_1}_1\lx\TTt^{\lg_n}_n}\kern-.8em F(t)\>\dt}
$$
where $\TTt^{\lg_m}_m$ is the respective deformation of $\TT^{\lg_m}_m$. On
every contour $\TTt^{\lg_m}_m$ the quantities $t_a/z_m$ remain bounded and
separated from zero as $z$ tends to limit, $z\rto\AA$, for all $a$ \st/
$\lg^{m-1}<a\le\lg^m$, and the rest of the proof remains the same as before.
\par
Theorem~\[AsIWw] is proved.
\epf
\Th{AsIWwt}
Let the parameters $\xin$ obey condition \(Lass). Then for any \perm/
${\tau\in\S^n}$ the \hint/ ${I(W_\lg^\tau\},w_\mg)}$ has the following \as/s
as $z$ tends to limit in the \azo/, $z\rto\AA_\tau$, so that at any moment
assumption \(assum) holds\/{\rm:}
$$
\gather
I(W_\lg^\tau\},w_\mg)\,=\;{\ell\)!\over\lg_1!\ldots\lg_n!}\,
\prod_{\tsize{1\le l<m\atop \si_l<\si_m}}\!\xi_l^{\lg_m}\!\!
\prod_{\tsize{1\le l<m\atop \si_l>\si_m}}\!\eta^{\lg_l\lg_m}\xi_l^{-\lg_m}\,
\bigl(\tpron I(W\"m\'{\lg_m}[\ka^\tau_{\lg,m}],w\"m\'{\lg_m})+o(1)\bigr)\,,
\\
\nn10>
\alignedat2
& I(W_\lg^\tau\},w_\mg)\,=\,O(1) && \for ^\tau\}\lg\ll\}{}^\tau\!\mg\,,
\\
\nn2>
& I(W_\lg^\tau\},w_\mg)\,=\,o(1)\,,\qqq && \text{otherwise}\,.
\endalignedat
\endgather
$$
Here \,${\si=\tau\1\]}$, the \fn/s $w\"m\'{\lg_m}$ and
$W\"m\'{\lg_m}[\ka^\tau_{\lg,m}]$ are defined by formulae \(wWprim),
${^\tau\}\lg=(\lg_{\tau_1}\lc\lg_{\tau_n})}$ \,and \,$\ka_{\lg,m}^\tau=
\ka\!\prlmsi\!\!\!\eta^{-\lg_l}\xi_l\:\!\prmlsi\!\eta^{\lg_l}\xi_l\1\]$.
\endpro
\nt
The proof is similar to the proof of Theorem~\[AsIWw].
\goodast
\Pf of Theorem\{~\{\[asol].
The statement follows from Theorem~\[AsIWwt] and formula \(qbeta).
\epf
\Pf of Theorem\{~\{\[mu<>0].
Since both sides of formula \(mu<>0) are \anf/s of $\xin$ and $\eta$,
it suffices to prove the formula under the assumption that ${|\eta\)|>1}$
and $|\xi_m|<1$, $\mn$.
\par
Denote by $F(z)$ the determinant $\DWwi_{\lg,\mg\in\Zln}$ and by $G(z)$ \rhs/
of formula \(mu<>0). Let $Y_\lg$ be an \adjf/ for the \ewf/ $W_\lg$ and
$\al_{\lg,m}$ be the corresponding multipliers defined by formulae \(adjust).
\par
Since for every $\lg\in\Zln$ the section $\Psi_{Y_\lg W_\lg}$ is a
$\Vl$-valued \sol/ of the \qKZe/, $F(z)$ solves the following system of \deq/s:
$$
F(\zmn)\,=\,
\Det\ell K_m(\zn)\!\tprod_{\,\,\lg\in\Zln\!\!}\al_{\lg,m}\1\>F(\zn)\,.
$$
Here $\Det\ell K_m(z)$ stands for the determinant of the operator $K_m(z)$
\(Kmz) acting in the \wt/ subspace $\Vl$. Using either formula \(Rspec) or
Theorem~\[DetM] we see that
$$
\align
& \Det\ell K_m(\zn)\!\tprod_{\,\,\lg\in\Zln\!\!}\al_{\lg,m}\1\,={}
\\
\nn2>
& \hp{\det}{}\,=\,\pros\,\Bigl(
\prlm{p\)z_m-\eta^s\xi_l\1\xi_m\1z_l\over p\)z_m-\eta^{-s}\xi_l\:\)\xi_m z_l}
\>\prml{z_l-\eta^{-s}\xi_l\:\)\xi_m z_m\over z_l-\eta^s\xi_l\1\xi_m\1 z_m}\,
\Bigr)^{\tsize{n+\ell-s-2\choose n-1}}.
\endalign
$$
Therefore, the ratio $F(z)/G(z)$ is a \p-periodic \fn/ of each of the \var/s
$\zn$:
$$
{F\over G}\)(\zmn)\,=\,{F\over G}\)(\zn)\;.
$$
Theorem~\[AsIWw] implies that the ratio $F(z)/G(z)$ tends to $1$ as $z$ tends
to limit in the \azo/, $z\rto\AA$. Hence, this ratio equals $1$
identically, which completes the proof of the determinant formula.
\Par
Let \fn/s $G_\lg$, $\lg\in\Zln$, be defined by formula \(Glg). They form a
basis in the \ehgf/ of the fiber $\Fq(z)$. Using Theorem~\[DetMq] we have that
\ifMag
$$
\gather
\Lline{\det\bigl[I(G_\lg,w_\mg)\bigr]_{\lg,\mg\in\Zln}
\,=\,\Xi\1\>(2\pii)\vpb{\ell\tsize{n+\ell-1\choose n-1}}\,
\ell\)!\vpb{\tsize{n+\ell-1\choose n-1}}\>
\eta\vpb{-n\tsize{n+\ell-1\choose n+1}}\;\x}
\\
\nn6>
\Rline{\x\;
\pron z_m^{(n-m){\tsize{n+\ell-1\choose n}}}\,
\pros\,\Bigl[\>{(\eta\1)\9^n\>(\eta^{s+1-\ell}\ka\1\}\prod\xi_m)\9\>
(p\)\eta^{s+1-\ell}\)\ka\prod\xi_m)\9\!\over
(\eta^{-s-1})\9^n\>(p)\9^{2n-1}\>\prod(\eta^{-s}\xi_m^2)\9}\ \x}
\\
\Rline{\x\plmn{1\over(p)\9\>(\eta^{-s}\xi_l\)\xi_m z_l/z_m)\9\>
(p\)\eta^{-s}\xi_l\)\xi_mz_m/z_l)\9\!}\,\Bigr]^{\tsize{n+\ell-s-2\choose n-1}}
.\!}
\endgather
$$
\else
$$
\NN8>
\align
\det\bigl[I(G_\lg,w_\mg)\bigr]_{\lg,\mg\in\Zln}
\,=\,\Xi\1\>(2\pii)\vpb{\ell\tsize{n+\ell-1\choose n-1}}\,
\ell\)!\vpb{\tsize{n+\ell-1\choose n-1}}\>
\eta\vpb{-n\tsize{n+\ell-1\choose n+1}}\>
\pron z_m^{(n-m){\tsize{n+\ell-1\choose n}}}\ \x &
\\
\x\;
\pros\,\Bigl[\>{(\eta\1)\9^n\>(\eta^{s+1-\ell}\ka\1\}\prod\xi_m)\9\>
(p\)\eta^{s+1-\ell}\)\ka\prod\xi_m)\9\!\over
(\eta^{-s-1})\9^n\>(p)\9^{2n-1}\>\prod(\eta^{-s}\xi_m^2)\9}\ \x{}
\kern8em
\\
\x\plmn{1\over(p)\9\>(\eta^{-s}\xi_l\)\xi_m z_l/z_m)\9\>
(p\)\eta^{-s}\xi_l\)\xi_mz_m/z_l)\9\!}\,\Bigr]^{\tsize{n+\ell-s-2\choose n-1}}
\!&
\endalign
$$
\fi
where $\Xi$ is the constant given in Theorem~\[DetMq]. Under the assumptions
of Theorem~\[mu<>0] we have that
${\det\bigl[I(G_\lg,w_\mg)\bigr]_{\lg,\mg\in\Zln}\ne\,0}$ which means that
the \hpair/ ${I:\Fq(z)\ox\Fo(z)\to\C}$ is \ndeg/. Theorem~\[mu<>0] is proved.
\epf
\Pf of Theorem\{~\{\[mu=0].
Since both sides of formula \(mu<>0) are \anf/s of $\xin$ and $\eta$,
it suffices to prove the formula under the assumption that ${|\eta\)|>1}$
and $|\xi_m|<1$, $\mn$.
\Par
Consider the determinant $\DWwi_{\lg,\mg\in\Zln}$ as a \fn/ of $\ka$ and denote
it by $\Dit(\ka;\ell)$. Set $\eps=\bigl(1-\ka\1\kone\bigr)$. We will show that
\ifMag
$$
\align
\kern1em
\Dit(\ka;\ell)\, &{}=\,\pros\,\Bigl(\)
{2\pii\>\eps\>\ell\over(1-\eta^{-s}\xi_1^2)}{(1-\eta)\over(1-\eta^{s+1})}
{\tht(\eta)\over\tht(\eta^{s+1})}\)\Bigr)^{\tsize{n+\ell-s-3\choose n-2}}\x{}
\Tagg{Dit}
\\
\nn4>
& {}\>\x{}\,\bigl(\Dit(\koNe;\ell-1)
\DWwi_{\%{\Rph{\lg_1=\mg_1=0}{\lg,\mg\in\Zln}}_{\ssize\lg_1=\mg_1=0}}
+\>o(1)\bigr)
\kern-1em
\endalign
$$
\else
$$
\align
\Dit(\ka;\ell)\,=\,\pros\,&\Bigl(\)
{2\pii\>\eps\>\ell\over(1-\eta^{-s}\xi_1^2)}{(1-\eta)\over(1-\eta^{s+1})}
{\tht(\eta)\over\tht(\eta^{s+1})}\)\Bigr)^{\tsize{n+\ell-s-3\choose n-2}}\x{}
\Tag{Dit}
\\
\nn4>
{}\x{} &\,\bigl(\Dit(\koNe;\ell-1)
\DWwi_{\%{\Rph{\lg_1=\mg_1=0}{\lg,\mg\in\Zln}}_{\ssize\lg_1=\mg_1=0}}
+\>o(1)\bigr)
\endalign
$$
\fi
as $\eps\to 0$, that is $\ka\to\kone$. This equality and the determinant
formula \(mu<>0) imply the determinant formula \(mu=0).
\Par
Let $\,\eg(m)=(0\lc\%1_{\sss\^{$m$-th}}\lc 0)$, $\,\mn$.
Introduce a new basis $w'_\lg$, ${\lg\in\Zln}\!$, of the \thgf/ of a fiber
by the rule: $w'_\lg=w_\lg\:$ \,for $\lg_1=0$ and
$$
\gather
\ifAst
\Cline{\aligned
\\
\nn-24>
w'_\lg\, &{}=\,w_\lg\:\>(1-\eta^{\lg_1}\])\>(1-\eta^{1-\lg_1}\xi_1^2)\;-{}
\\
\nn6>
& {}\)-\sum_{m=2}^n w_{\lg-\eg(1)+\eg(m)}\>
(1-\eta^{\lg_m+1})\>(\xi_m-\eta^{\lg_m}\xi_m\1)\>\tprllm\!\eta^{-\lg_l}\xi_l
\endaligned}
\else
\ifMag
\Lline{w'_\lg\,=\,w_\lg\:\>(1-\eta^{\lg_1}\])\>(1-\eta^{1-\lg_1}\xi_1^2)\;-
\sum_{m=2}^n w_{\lg-\eg(1)+\eg(m)}\>
(1-\eta^{\lg_m+1})\>(\xi_m-\eta^{\lg_m}\xi_m\1)\>\tprllm\!\eta^{-\lg_l}\xi_l}
\else
w'_\lg\,=\,w_\lg\:\>(1-\eta^{\lg_1}\])\>(1-\eta^{1-\lg_1}\xi_1^2)\;-
\sum_{m=2}^n w_{\lg-\eg(1)+\eg(m)}\>
(1-\eta^{\lg_m+1})\>(\xi_m-\eta^{\lg_m}\xi_m\1)\>\tprllm\!\eta^{-\lg_l}\xi_l
\fi\fi
\\
\ifAst\nn3>\fi
\nn3>
\Text{for $\lg_1>0$. We have that}
\nn2>
\ifAst\nn3>\fi
\Dit(\ka;\ell)\,=\,\pros\>\bigl((1-\eta^{s+1})\>(1-\eta^{-s}\xi_1^2)
\bigr)^{-\tsize{n+\ell-s-3\choose n-2}}\,\DWwp_{\lg,\mg\in\Zln}\,.
\endgather
$$
The main property of the new basis follows from Lemma~\[IWw=0]. Namely, we have
that $I(W_\lg,w_\mg)=O(\eps)$ as $\eps\to 0$, if either $\lg_1>0$ or $\mg_1>0$.
Therefore,
\ifAst
$$
\align
& \DWwp_{\lg,\mg\in\Zln}={}
\\
\nn4>
&\hp{\det}
{}=\,\DWwi_{\%{\Rph{\lg_1=\mg_1=0}{\lg,\mg\in\Zln}}_{\ssize\lg_1=\mg_1=0}}
\DWwp_{\%{\Rph{\lg_1>0,\mg_1>0}{\lg,\mg\in\Zln}}_{\ssize\lg_1>0,\mg_1>0}}+
\,o\bigl(\eps\vpb{\tsize{n+\ell-2\choose n-1}}\bigr)\,.
\endalign
$$
\else
\ifMag
$$
\Lline{\DWwp_{\lg,\mg\in\Zln}=\,
\DWwi_{\%{\Rph{\lg_1=\mg_1=0}{\lg,\mg\in\Zln}}_{\ssize\lg_1=\mg_1=0}}
\DWwp_{\%{\Rph{\lg_1>0,\mg_1>0}{\lg,\mg\in\Zln}}_{\ssize\lg_1>0,\mg_1>0}}+
\,o\bigl(\eps\vpb{\tsize{n+\ell-2\choose n-1}}\bigr)\,.}
$$
\else
$$
\DWwp_{\lg,\mg\in\Zln}=\,
\DWwi_{\%{\Rph{\lg_1=\mg_1=0}{\lg,\mg\in\Zln}}_{\ssize\lg_1=\mg_1=0}}
\DWwp_{\%{\Rph{\lg_1>0,\mg_1>0}{\lg,\mg\in\Zln}}_{\ssize\lg_1>0,\mg_1>0}}+
\,o\bigl(\eps\vpb{\tsize{n+\ell-2\choose n-1}}\bigr)\,.
$$
\fi\fi
\agood
If $\mg_1>0$, then by Lemmas~\[wDw] and \[wDw0]
\ifAst
$$
\align
w'_\mg & (\tell)\,={}
\\
\nn2>
& {}=\,\eps
\sum_{a=1}^\ell\,\bigl[\)w_{\mg-\eg(1)}(\twll)\)\bigr]_{(1,a)}+(1-\eps)
\sum_{a=1}^\ell\,D_a\bigl[\)w_{\mg-\eg(1)}(\twll)\)\bigr]_{(1,a)}
\endalign
$$
\else
\ifMag
$$
\Lline{w'_\mg(\tell)\,=\,\eps
\sum_{a=1}^\ell\,\bigl[\)w_{\mg-\eg(1)}(\twll)\)\bigr]_{(1,a)}
+(1-\eps)
\sum_{a=1}^\ell\,D_a\bigl[\)w_{\mg-\eg(1)}(\twll)\)\bigr]_{(1,a)}}
$$
\else
$$
w'_\mg(\tell)\,=\,
\eps\,\sum_{a=1}^\ell\,\bigl[\)w_{\mg-\eg(1)}(\twll)\)\bigr]_{(1,a)}-
(1-\eps)\>\sum_{a=1}^\ell\,D_a\bigl[\)w_{\mg-\eg(1)}(\twll)\)\bigr]_{(1,a)}
$$
\fi\fi
where $(1,a)\in\Sl$ are transpositions. Then using Lemma~\[IWwD] we see that
$I(W_\lg,w'_\mg)=\eps\>I(W_\lg,w''_{\mg-\eg(1)})$ \,where
$$
\vp{\sum}\smash{w''_{\ng}(\tell)\;=
\,\sum_{a=1}^\ell\,\bigl[\)w_{\ng}(\twll)\)\bigr]_{(1,a)}\,,\qqq\ng\in\Zll\,.}
\Tag{w''}
$$
\par
The next step is to calculate the \hint/ $I(W_\lg,w''_{\mg-\eg(1)})$ at
$\ka=\kone$. We will indicate explicitly the dependence of the \ewf/s on $\ka$.
Namely, the \ewf/ $W_\lg[\ka]$ is an element of the \ehgf/ of a fiber
${\Fq[\)\ka,\}\sun\!\lg_m](z)}$.
\par
If $\lg_1>0$, then
$$
W_\lg[\)\kone\)](\tell)\,=\;{\tht(\eta)\over\tht(\eta^{\lg_1}\})}
\,\susi\,\Lbc\>W_{\lg-\eg(1)}[\)\koNe\)](\twll)\>\Rbc_\si\,,
$$
and due to formula \(Phisym) we have that
\ifMag
$$
\gather
\Lline{I(W_\lg[\)\kone\)],w''_\mg)\,=\;\ell\;
{\tht(\eta)\over\tht(\eta^{\lg_1}\})}\ \x{}}
\\
\nn6>
\Rline{{}\x\int_{\){\TT^\ell}\]}
\Phi(\tell)\>w''_\mg(\tell)\>W_{\lg-\eg(1)}[\)\koNe\)](\twll)\;\dtt\,,}
\endgather
$$
\else
$$
I(W_\lg[\)\kone\)],w''_\mg)\,=\;\ell\;
{\tht(\eta)\over\tht(\eta^{\lg_1}\})}\,\int_{\){\TT^\ell}\]}
\Phi(\tell)\>w''_\mg(\tell)\>W_{\lg-\eg(1)}[\)\koNe\)](\twll)\;\dtt
$$
\fi
because the torus $\TT^\ell$ is invariant under \perm/s of the \var/s $\tell$.
Substitute formula \(w'') into the above integral. Similar to the proof of the
second claim of Lemma~\[IWw=0] we obtain that
$$
\gather
\ifAst
\gathered
\Lline{\int_{\){\TT^\ell}\]}\Phi(\tell)\>
\bigl[\)w_{\mg-\eg(1)}(\twll)\)\bigr]_{(1,a)}
\>W_{\lg-\eg(1)}[\)\koNe\)](\twll)\;\dtt\,={}}
\\
\Rline{{}={}\,2\pii\,\dl_{1a}\>I(W_{\lg-\eg(1)},w_{\mg-\eg(1)})\,.}
\endgathered
\else
\Cline{\aligned
\int_{\){\TT^\ell}\]}\Phi(\tell)\>
\bigl[\)w_{\mg-\eg(1)}(\twll)\)\bigr]_{(1,a)}
\>W_{\lg-\eg(1)}[\)\koNe\)] & (\twll)\;\dtt\,={}
\\
{}={}\, & 2\pii\,\dl_{1a}\>I(W_{\lg-\eg(1)},w_{\mg-\eg(1)})\,.
\endaligned}
\fi
\\
\ifMag\nn-4>\else\nn-2>\fi
\Text{Hence,}
\ifMag\ifAst\nn2>\fi\else\nn-2>\fi
I(W_\lg[\)\kone\)],w''_{\mg-\eg(1)})\,=\,2\pii\>\ell\;
{\tht(\eta)\over\tht(\eta^{\lg_1}\})}\;
I(W_{\lg-\eg(1)}[\)\koNe\)],w_{\mg-\eg(1)})
\\
\ifAst\nn2>\fi
\nn4>
\Text{and finally}
\ifMag\nn2>\else\nn-2>\fi
\ifAst
\nn-6>
\Cline{
\aligned
\DWwp_{\%{\Rph{\lg_1>0,\mg_1>0}{\lg,\mg\in\Zln}}_{\ssize\lg_1>0,\mg_1>0}}=\,
\pros\,\Bigl( & \){2\pii\>\eps\>\ell\,\tht(\eta)\over\tht(\eta^{s+1})}\)
\Bigr)^{\tsize{n+\ell-s-3\choose n-2}}\,\x{}
\\
\nn4>
{}\x{}\,\bigl( & \Dit(\koNe;\ell-1)+o(1)\bigr)
\endaligned}
\else
\Cline{
\DWwp_{\%{\Rph{\lg_1>0,\mg_1>0}{\lg,\mg\in\Zln}}_{\ssize\lg_1>0,\mg_1>0}}=\;
\pros\,\Bigl(\){2\pii\>\eps\>\ell\,\tht(\eta)\over\tht(\eta^{s+1})}\)
\Bigr)^{\tsize{n+\ell-s-3\choose n-2}}\>\bigl(\Dit(\koNe;\ell-1)+o(1)\bigr)}
\fi
\endgather
$$
as $\eps\to 0$. Formula \(Dit) is proved.
\Par
The rest of the proof is similar to the end of the proof of Theorem~\[mu<>0].
Consider the space
\ifMag\else\nl\fi
${\Fq[\)\ka\>\xi_1;z_2\lc z_n;\xi_2\lc\xi_n;\ell\,](z)}$.
It has a basis given by \fn/s $\Gti_\lg(t)$, $\lg\in\Zlm\!$, defined similarly
to formula \(Glg). Set
$$
G'_\lg(t)\,=\,\Gti_\lg(t)\,
\pral\,{\tht(\xi_1t_a/z_1)\over\tht(\xi_1\1t_a/z_1)}\;,\qqq \lg\in\Zlm.
$$
\ifMag\else\vvv.2>\fi
The \fn/s $W_\lg[\ka]$ \st/ $\lg_1=0$, $\lg\in\Zln\!$, are linear combinations
\ifMag\vvn.1>\fi
of the \fn/s $G'_\mg$, $\mg\in\Zlm\!$. Formula~\[mu=0] and Theorem~\[DetMq]
imply that
\mmgood
\ifMag
$$
\align
\det\bigl[I(G'_\lg,{} & w\'{0,\mg})\bigr]_{\lg,\mg\in\Zlm}\,={}
\\
\nn8>
&{}\!\}=\,K\>(2\pii\)\xi_1)\vpb{\ell\tsize{n+\ell-2\choose n-2}}\,
\ell\)!\vpb{\tsize{n+\ell-2\choose n-2}}\>
\eta\vpb{-n\tsize{n+\ell-2\choose n}}
\prod_{m=2}^n z_m^{(n-m){\tsize{n+\ell-2\choose n-1}}}\ \x
\\
\nn7>
& \aligned
{}\!\}\x\;
\pros\,\Bigl[\>{(\eta\1)\9^{n-1}\>(p\)\eta^{s+2-2\ell}\prod\xi_m^2)\9\>
(\eta^s\xi_1^{-2})\9\!\over(\eta^{-s-1})\9^{n-1}\>(p)\9^{2n-3}\!
\prod_{1<m\le n}\!(\eta^{-s}\xi_m^2)\9}
\ \prod_{m=2}^n\,{(\eta^s\xi_1\1\xi_m\1z_1/z_m)\9\!\over
(\eta^{-s}\xi_1\)\xi_m z_1/z_m)\9}\ \x\!{}&
\\
\nn2>
{}\x\prod_{2\le l<m\le n}{1\over(p)\9\>(\eta^{-s}\xi_l\)\xi_m z_l/z_m)\9\>
(p\)\eta^{-s}\xi_l\)\xi_mz_m/z_l)\9\!}\,\Bigr]^{\tsize{n+\ell-s-3\choose n-2}}.
& \endaligned
\endalign
$$
\vvv.2>
\else
$$
\NN8>
\align
\det\bigl[I(G'_\lg,w\'{0,\mg})\bigr]_{\lg,\mg\in\Zlm}
\,=\,K\>(2\pii\)\xi_1)\vpb{\ell\tsize{n+\ell-2\choose n-2}}\,
\ell\)!\vpb{\tsize{n+\ell-2\choose n-2}}\>
\eta\vpb{-n\tsize{n+\ell-2\choose n}}\>
\prod_{m=2}^n z_m^{(n-m){\tsize{n+\ell-2\choose n-1}}}\ \x &
\\
\x\;\pros\,\Bigl[\>{(\eta\1)\9^{n-1}\>(p\)\eta^{s+2-2\ell}\prod\xi_m^2)\9\>
(\eta^s\xi_1^{-2})\9\!\over(\eta^{-s-1})\9^{n-1}\>(p)\9^{2n-3}\!
\prod_{1<m\le n}\!(\eta^{-s}\xi_m^2)\9}
\ \prod_{m=2}^n\,{(\eta^s\xi_1\1\xi_m\1z_1/z_m)\9\!\over
(\eta^{-s}\xi_1\)\xi_m z_1/z_m)\9}\ \x {}\qquad &
\\
\x\prod_{2\le l<m\le n}{1\over(p)\9\>(\eta^{-s}\xi_l\)\xi_m z_l/z_m)\9\>
(p\)\eta^{-s}\xi_l\)\xi_mz_m/z_l)\9\!}\,\Bigr]^{\tsize{n+\ell-s-3\choose n-2}}
\! & 
\endalign
$$
\vvv.3>
\fi
where \ $\dsize K\,=\,\Bigl[\>\Rph{(p)}{(p)\9}\vpb{n(n-2)}\>
\prod_{m=1}^{n-2}\>\Bigl(\>{e^{2\pii m/(n-1)}\]-1\over
\tht(e^{2\pii m/(n-1)})}\)\Bigr)^{\!\!\raise 3pt\mbox{\ssize n-m-1}}\>
\Bigr]^{\tsize{n+\ell-2\choose n-1}}\!$.
\Par
Under the assumptions of Theorem~\[mu=0] we have that
${\det\bigl[I(G'_\lg,w\'{0,\mg})\bigr]_{\lg,\mg\in\Zlm}\ne\,0}$ which means
that the \hpair/ ${I:\Fq(z)/\Qc(z)\ox\Fo(z)/\Rc(z)\to\C}$ is \ndeg/.
Theorem~\[mu=0] is proved.
\epf
\Pf of Theorem\{~\{\[mu=02].
The proof is similar to the proof of Theorem~\[mu=0].
\epf
\Pf of Theorem\{~\{\[hmapo].
\vv.15>
Under assumptions of the theorem, for any $\ka$ in the puctured \neib/ of
$\kOne$
\vv.15>
the assumptions of Theorem~\[hmap] are valid. Therefore, the \hmap/
$I_{\tau,\tau'}(z)$ is well defined and \ndeg/ for any such $\ka$.
\vsk.2>
Introduce a matrix $X$ by the rule:
\ifMag\vv.1>\fi
$$
I_{\tau,\tau'}(z)\>\Fvqtt\,=\smiZ X_{\lg\mg}\>\Fvvt\,,\qqq\lg\in\Zln\,.
$$
\vvv-.6>
We have to show that the matrix $X$ has a finite limit as $\;\ka\to\kOne$
\ifMag\vv-.4>\else\vv.1>\fi
and $\;\lim\det X\ne 0$.
\par
Consider the \hint/ $I(W^{\tau'}_\lg\!,w^\tau_\mg)$. It is a \hof/ of $\ka$
since the corresponding integrand is a \hof/ of $\ka$ and the integration is
over a compact contour. Hence, Lemmas~\[IWw=0] and \[dimFQ] imply that
$\bigl(\ka-\kOne\bigr)\1{\>I(W^{\tau'}_\lg\!,w^\tau_\mg)}$ has a finite limit
as $\ka\to\kOne$ if $\lg_1>0$. Therefore, the entries $X_{\lg\)\mg}$ with
$\lg_1>0$ have finite limits and the entries $X_{\lg\mg}$ with $\lg_1=0$ are
regular at $\ka=\kOne$. That is the matrix $X$ has a finite limit as
$\ka\to\kOne$.
\vsk.15>
The explicit formula for the determinant $\det X$ for general $\ka$ can be
easily obtained using Theorems~\[mu<>0], \[DetM], \[DetMq], formula \(clg)
and the definition of the \hmap/ $I_{\tau,\tau'}(z)$. It follows from
the obtained expression that under the assumptions of Theorem~\[hmapo]
the limit of the determinant $\det X$ as $\ka\to\kOne$ is not equal to zero.
The theorem is proved.
\epf
\Pf of Theorem\{~\{\[hmapp].
The proof is similar to the proof of Theorem~\[hmapo].
\vvgood
\mmgood
\epf
\Pf of Theorem\{~\{\[PsivY].
\Wlg/ we can assume that the \adjm/ $Y^{\tau'}(z)$ analytically depends on
$\ka$. Hence, the same do the section
$\Psi_{\smash{v,Y^{\tau'}_{\vp|}}}^\tau\!$.
The \qKZo/s analytically depend on $\ka$ as well.
\par
For general $\ka$ the section $\Psi_{\smash{v,Y^{\tau'}_{\vp|}}}^\tau\!$ solves
the \qKZe/ with values in $\Vlt$ due to Corollary~\[qKZsol]. For special values
of $\ka$: $\;\ka=\kOne$ and $\,\ka=\kTwp\!$, \,the \qKZe/ remains valid by
the \anco/.
\vsk.15>
Theorems~\[hmap], \[hmapo] and \[hmapp] imply that under the assumptions of
each of Theorems~\[mu<>0]\,--\,\[mu=02] the sections
$\Psi_{v,Y^{\tau'}_{\vp|}}^\tau\!$, $\,v\in\Vqltt$, \;span the space of \sol/s
of the \qKZe/ over the field of \p-periodic \fn/s (quasiconstants).
Theorem~\[PsivY] is proved.
\epf

\Appendix

\ifAst\newpage\fi
\Sect[A]{Basic facts about the \tri/ \ifAst\{\nl\fi \hgf/}
Let ${\Fo=\Fo[\zn;\xin;\ell\,]}$ be the \thgf/.
\par
By construction, the \thgf/ of a fiber has the same dimension as the space of
\sym/ \pol/s in $\ell$ \var/s of degree less than $n$ in each of the \var/s,
that is
$$
\dim\Fo(z)\,=\,{n+\ell-1\choose n-1}\,.
$$
\Lm{winFo}
For any ${\lg\in\Zln}$ the \twf/ $w_\lg$ is in the \thgf/ $\Fo$.
\endpro
\Pf.
It is clear from definiton \(wlga) that the \fn/ $w_\lg(t,z)$ has the form
$$
Q(\tell,\zn)\,{\tsize\pral\}t_a}\pron\,\pral\,{1\over t_a-\xi_mz_m}\;
\prab\ \,{1\over\eta\)t_a-t_b}
$$
where $Q$ is a \pol/ which has degree less than $n+\ell-1$ in each of the
\var/s $\tell$. Furthermore, by construction the \fn/ $w_\lg$ as a \fn/ of
$\tell$ is invariant \wrt/ the action \(act) of the \symg/ $\Sl\!$, which means
that the \pol/ $Q$ is skew\sym/ \wrt/ the \var/s $\tell$. Hence, the \pol/ $Q$
is divisible by $\!\!\dsize\prab(t_a-t_b)$
\vv.1>
and the ratio is a \pol/ which is \sym/ in \var/s $\tell$ and has degree less
than $n$ in each of the \var/s $\tell$; that is the \fn/ $w_\lg$ is in the
\thgf/.
\epf
\Cr{n=1}
Let $n=1$. Then
$$
\vp\sum\smash{w\'\ell(\tell,z_1)\,=\,
\pral\,{t_a\over t_a-\xi_1z_1}\;\prab\ {t_a-t_b\over\eta\)t_a-t_b}\;.}
$$
\endpro
\Pf.
Denote by $f(t,z_1)$ the ratio ${\dsize\,w\'\ell(t,z_1)\,\Big[\pral\,
{t_a\over t_a-\xi_1z_1}\;\prab\ {t_a-t_b\over\eta\)t_a-t_b}\,\Big]\1\!}$.
Since for ${n=1}$ the \thgf/ of a fiber is \onedim/, $f(t,z_1)$ does not
depend on $t$. Let $\ts=(t_1,{\eta\1t_1}\lc\eta^{1-\ell}t_1)$. Only the term
corresponding to the identity \perm/ in \rhs/ of \(wlga) contributes to
$w(\ts\},z_1)$ and by a straightforward calculation we get $f(\ts\!,z_1)=1$.
\epf
\Lm{n=10}
Let $n=1$. Then
$$
\align
\NN6>
(1-\eta^\ell)\>& (\eta^{1-\ell}\xi_1-\xi_1\1)\>w\'\ell(\tell,z_1)\,=\,
\\
&{}\!=\,(1-\eta)\>\sum_{a=1}^\ell\,
\Bigl[\Bigl(\,{\xi_1t_1-z_1\over\xi_1\1t_1-z_1}\,
\prod_{b=2}^\ell\,{\eta\1t_1-t_b\over\eta\)t_1-t_b}\,-1\Bigr)\>
w\'{\ell-1}(\twll)\>\Bigr]\'{1,a}\,,
\\
\nn2>
(1-\eta^\ell)\>& (\xi_1-\eta^{\ell-1}\xi_1\1)\>w\'\ell(\tell,z_1)\,=\,
\\
&{}\!=\,(1-\eta)\>\sum_{a=1}^\ell\,
\Bigl[\Bigl(\,{t_1-\xi_1\1z_1\over t_1-\xi_1z_1}\,
\prod_{b=2}^\ell\,{t_1-\eta\)t_b\over t_1-\eta\1t_b}\,-1\Bigr)\>t_1\>
w\'{\ell-1}(\twll)\>\Bigr]\'{1,a}\,.
\endalign
$$
Here $(1,a)\in\Sl$ are transpositions.
\vvgood
\endpro
\Pf.
Similar to the proof of Lemma~\[winFo] one can show that \rhs/s of the both
above formulae are elements of the \thgf/. The rest of the proof is similar to
the proof of Corollary~\[n=1].
\epf
\Lm{n=11}
Let $n=1$. Then
$$
\smash{w\'\ell(\tell,z_1)\,=\,\susi\,
\Bigl[\;\pral\,{t_a\over t_a-\eta\)t_{a-1}}\;\Bigr]_\si\,.}
$$
Here $\,t_0=\eta\1\xi_1z_1\)$.
\endpro
\Pf.
\Rhs/ of the above formula is an element of the \thgf/. Comparing residues
of both sides at $t=(\xi_1z_1\lc\eta^{\ell-1}\xi_1z_1)$ completes the proof.
\epf
\Pf of Lemma\{~\{\[wlgp].
Let $\S^{\lg_1}\!\lx\S^{\lg_n}\sub\Sl$ be the subgroup of \perm/s preserving
the subsets $\Gm_m=\lb 1+\lg^{m-1}\,\lc\lg^m\rb$, $\mn$. The coset space
$\Sl/(\)\S^{\lg_1}\!\lx\S^{\lg_n})$ is in one-to-one correspondence with
the set of all \^{$\,n$-}tuples $\Gmn$ of disjoint subsets of
${\)\lb 1\lc\ell\)\rb}$ \st/ $\Gms_m$ has $\lg_m$ elements.
Namely, a permutation $\si\in\Sl$ corresponds to an \^{$\,n$-}tuple
$\si(\Gm_1)\lc\si(\Gm_n)$, and the \^{$\,n$-}tuple depends only on
the coset of the \perm/ $\si$.
\par
Perform the summation in \rhs/ of formulae \(wlga) in two steps. First take
a sum over the subgroup $\S^{\lg_1}\!\lx\S^{\lg_n}\!$. This can be done
explicitly using Corollary~\[n=1]. The remaining sum over the coset space
$\Sl/(\)S^{\lg_1}\!\lx\S^{\lg_n})$ is \eqv/ to \rhs/ of formula \(wlgp).
The lemma is proved.
\epf
\Lm{wDw0}
For any $\lg\in\Zll$ the following relations hold:
\ifAst
$$
\align
\NN8>
\sun w_{\lg+\eg(m)}\> & (1-\eta^{\lg_m+1})\>(\xi_m-\eta^{\lg_m}\xi_m\1)\>
\tprllm\!\eta^{-\lg_l}\xi_l\;=
\\
{}=\,{} & (1-\eta)\>\sum_{a=1}^\ell\,
\Bigl[\Bigl(\,\pron\,{\xi_mt_1-z_m\over\xi_m\1t_1-z_m}\,
\prod_{b=2}^\ell\,{\eta\1t_1-t_b\over\eta\)t_1-t_b}\,-1\Bigr)\>
w_\lg(\twll)\)\Bigr]_{(1,a)}\,,
\\
\ald
\nn2>
\sun w_{\lg+\eg(m)}\> & (1-\eta^{\lg_m+1})\>(\xi_m-\eta^{\lg_m}\xi_m\1)\,
z_m\tprlm\!\eta^{\lg_l}\xi_l\1\;=
\\
{}=\,{} & (1-\eta)\>
\sum_{a=1}^\ell\,\Bigl[\Bigl(\,\pron\,{t_1-\xi_m\1z_m\over t_1-\xi_mz_m}\,
\prod_{b=2}^\ell\,{t_1-\eta\)t_b\over t_1-\eta\1t_b}\,-1\Bigr)\>t_1\>
w_\lg(\twll)\)\bigr]_{(1,a)}\,,
\endalign
$$
\else
$$
\align
\NN8>
\sun w_{\lg+\eg(m)}\> & (1-\eta^{\lg_m+1})\>(\xi_m-\eta^{\lg_m}\xi_m\1)\>
\tprllm\!\eta^{-\lg_l}\xi_l\;=
\\
&{}\!=\,(1-\eta)\>\sum_{a=1}^\ell\,
\Bigl[\Bigl(\,\pron\,{\xi_mt_1-z_m\over\xi_m\1t_1-z_m}\,
\prod_{b=2}^\ell\,{\eta\1t_1-t_b\over\eta\)t_1-t_b}\,-1\Bigr)\>
w_\lg(\twll)\)\Bigr]_{(1,a)}\,,
\\
\nn2>
\sun w_{\lg+\eg(m)}\> & (1-\eta^{\lg_m+1})\>(\xi_m-\eta^{\lg_m}\xi_m\1)\,
z_m\tprlm\!\eta^{\lg_l}\xi_l\1\;=
\\
&{}\!=\,(1-\eta)\>
\sum_{a=1}^\ell\,\Bigl[\Bigl(\,\pron\,{t_1-\xi_m\1z_m\over t_1-\xi_mz_m}\,
\prod_{b=2}^\ell\,{t_1-\eta\)t_b\over t_1-\eta\1t_b}\,-1\Bigr)\>t_1\>
w_\lg(\twll)\)\bigr]_{(1,a)}
\endalign
$$
\fi
where ${(1,a)\in\Sl}$ are transpositions.
\endpro
\nt
The proof is similar to the proof of Lemma~\[wlgp] using Lemma~\[n=10] instead
of Corollary~\[n=1].
\Par
For any $\lg\in\Zln$ denote by $Q_\lg(\tell)$ the following \sym/ \pol/:
$$
Q_\lg(\tell)\,=\;{1\over\lg_1!\ldots\lg_n!}\,\susi\,
\pron\,\prod_{a\in\Gm_{\Rph lm}} t_{\si_a}^{m-1}
\Tag{Qlg}
$$
where $\Gm_m=\lb 1+\lg^{m-1}\,\lc\lg^m\rb$, $\mn$. Consider a basis in the
space $\Fo(z)$ given by \fn/s
$$
g_\lg(t,z)\,=\,Q_\lg\)(\tell)\,{\tsize\pral\}t_a}
\pron\,\pral\,{1\over t_a-\xi_mz_m}\;\prab\ {t_a-t_b\over \eta\)t_a-t_b}\;,
\ifMag\qquad\ \else\qqq\fi \lg\in\Zln\,.
$$
Define a matrix $M(z)$ by the rule:
$$
w_\lg(t,z)\,=\smiZ M_{\lg\mg}(z)\>g_\mg(t,z)\,,\qqq\lg\in\Zln.
$$
\Th{DetM}
\vvn->
$$
\det M\,=\,\pros\,\plmn\!(\eta^sz_l-\xi_l\)\xi_mz_m)
\vpb{\tsize{n+\ell-s-2\choose n-1}}.
$$
\vvgood
\endpro
\nt
The theorem is \eqv/ to Lemma 2.2 in \Cite{T}. Nevertheless, we give below
another proof of Theorem~\[DetM] which is similar to the proof of
Theorem~\[DetMq] in the elliptic case.
\Cr{F/R}
Let $\Rc(z),\,\Rc'(z)$ be the \cosub/s. Then
$$
\dim\Fo(z)/\Rc(z)\,=\,\dim\Fo(z)/\Rc'(z)\,=\,{n+\ell-2\choose n-2}
$$
provided that \ ${\xi_l\>\xi_m z_m/z_l\netr}$, \ $1\le l\le m\le n$,
\ for any $\roll$.
\endpro
\Pf.
Under assumptions of the corollary, both spaces $\Fo(z)/\Rc(z)$ and
$\Fo(z)/\Rc'(z)$ have bases induced by the set
${\lb\)w_\lg(t,z)}\vert\lg\in\Zln,\ \,{\lg_n=0\)\rb}$.
\epf
\goodast
\noPcd
\Pf of Theorem~\[DetM].
For any $\lg\in\Zln$ define a \sym/ \pol/ $P_\lg(\tell)$ by the rule:
\ifMag
$$
\align
& w_\lg(t,z)\,={}
\Tagg{Plg}
\\
& {}\}=\,P_\lg\)(\tell)\,{\tsize\pral\}t_a}\!\plmn\!\!\xi_l^{\lg_m}
\pron\,\pral\,{1\over t_a-\xi_mz_m}\;\prab\ {t_a-t_b\over \eta\)t_a-t_b}\;.
\kern-2.8em
\endalign
$$
\else
$$
w_\lg(t,z)\,=\,P_\lg\)(\tell)\,{\tsize\pral\}t_a}\!\plmn\!\!\xi_l^{\lg_m}
\pron\,\pral\,{1\over t_a-\xi_mz_m}\;\prab\ {t_a-t_b\over \eta\)t_a-t_b}\;.
\Tag{Plg}
$$
\fi
Introduce new \var/s $\xn$, $\yn$:
$$
x_m=\xi_mz_m\,,\qqq y_m=\xi_m\1z_m\,,\qquad\qquad\mn\,.
$$
Then the \pol/ $P_\lg(t)$ has the form:
$$
P_\lg(\tell)\,=\,\sum_{\Gmn}\biggl\lb\)
\prod_{\tsize{1\le m<l\le n\atop\vru1.4ex>a\in\Gms_m}}\!\!(t_a-x_l)\!
\prod_{\tsize{1\le l<m\le n\atop\vru1.4ex>a\in\Gms_m}}\!\!(t_a-y_l\:)\!
\prod_{\tsize{1\le l<m\le n\atop\vru1.4ex>a\in\Gms_l,\>b\in\Gms_{\Rph lm}}}
\!{\eta\)t_a -t_b\over t_a-t_b}\,\biggr\rb
$$
where the summation is over all \^{$\,n$-}tuples $\Gmn$ of disjoint subsets of
${\)\lb 1\lc\ell\)\rb}$ \st/ $\Gms_m$ has $\lg_m$ elements.
\par
Define a matrix $N$ by the rule:
$$
P_\lg(t)\,=\smiZ N_{\lg\mg}\>Q_\mg(t)\,,\qqq\lg\in\Zln.
$$
Then the claim of the theorem takes the form
$$
\det N\,=\,
\pros\,\plmn\!(\eta^sy_l\:-x_m)\vpb{\tsize{n+\ell-s-2\choose n-1}}.
\Tag{DetN}
$$
\par
For $\lg,\mg\in\Zln$ say that ${\lg\lle\mg}$ if
\vv-.4>
$\,\sum_{i=1}^m\lg_i\le\sum_{i=1}^m\mg_i\;$ for any $\mn-1$.
Say that $\lg\ll\mg$ if $\lg\ne\mg$ and $\lg\lle\mg$.
\vsk.1>
For any $x,y\in\Cn\!$ and $\lg\in\Zln$ set
$$
\gather
\ifMag\quad\fi
\xt\lg\,=\,(\eta^{1-\lg_1}x_1,\>\eta^{\)2-\lg_1}x_1\lc x_1,\>
\eta^{1-\lg_2}x_2\lc x_2,\,\ldots\,,\)\eta^{1-\lg_n}x_n\lc x_n)\,,
\ifMag\kern-2em\fi
\Tag{xtry}
\\
\nn4>
\ifMag\quad\fi
\yt\lg\,=\,(\eta^{\lg_1-1}y_1,\>\eta^{\lg_1-2}y_1\lc y_1,\>
\eta^{\lg_2-1}y_2\lc y_2,\,\ldots\,,\)\eta^{\lg_n-1}y_n\lc y_n)\,.
g\ifMag\kern-2em\fi
\endgather
$$
\Lm{Pxy}
$P_\lg(\xt\mg)=0$ unless $\lg\lle\mg$. $P_\lg(\yt\mg)=0$ unless $\lg\gge\mg$.
Moreover,
$$
\NN6>
\align
P_\lg(\xt\lg)\, &{} =\,\pron\,\prod_{s=0}^{\lg_m-1}\Bigl(\>
\prlm(\eta^{-s}x_m-y_l\:)\prml(\eta^{\lg_l-s}x_m-x_l)\>\Bigr)\,,
\\
P_\lg(\yt\lg)\, &{} =\,\pron\,\prod_{s=0}^{\lg_m-1}\Bigl(\>
\prlm(\eta^sy_m-\eta^{\lg_l}y_l\:)\prml(\eta^sy_m-x_l)\>\Bigr)\,.
\endalign
$$
\endpro
\nt
The proof is straightforward.
\goodbreak
\goodast
\Par
Set
$\ \dsize D(n,\ell,s)\,=\!\!\sum_{\tsize{\,r\in\Zp\!\atop 2r\le\ell-|s|-1}}^{}
\!\!{n+\ell-|s|-2r-3\choose n-2}\ifAst\vp{\sum^1}\fi$.
\Lm{detPxy}
\ifMag
$$
\NN4>
\align
\det &\bigl[P_\lg(\xt\mg)\bigr]_{\lg,\mg\in\Zln}\,=
\\
&\!{}=\,
\pros\,\plmn\!(y_l\:-\eta^{-s}x_m)\vpb{\tsize{n+\ell-s-2\choose n-1}}
\prsl\,\plmn\!(\eta^sx_m-x_l)\vpb{D(n,\ell,s)},
\\
\nn6>
\det &\bigl[P_\lg(\yt\mg)\bigr]_{\lg,\mg\in\Zln}=\,
\,\eta^{n(n-1)/2\cdot\!\!\tsize{n+\ell-1\choose n+1}}\;\x
\\
&\!{}\x\,
\pros\,\plmn\!(\eta^sy_l\:-x_m)\vpb{\tsize{n+\ell-s-2\choose n-1}}
\prsl\,\plmn\!(\eta^sy_m-y_l\:)\vpb{D(n,\ell,s)},
\endalign
$$
\else
\vvn->
$$
\NN4>
\alignat2
\det\bigl[P_\lg(\xt\mg)\bigr]_{\lg,\mg\in\Zln} &{}=\,
\pros\,\plmn\!(y_l\:-\eta^{-s}x_m)\vpb{\tsize{n+\ell-s-2\choose n-1}}
\prsl\,\plmn\!(\eta^sx_m-x_l)\vpb{D(n,\ell,s)}, &&
\\
\det\bigl[P_\lg(\yt\mg)\bigr]_{\lg,\mg\in\Zln} &{}=\,
\eta^{n(n-1)/2\cdot\!\!\tsize{n+\ell-1\choose n+1}}\>
\pros\,\plmn\!(\eta^sy_l\:-x_m)\vpb{\tsize{n+\ell-s-2\choose n-1}}\;\x &&
\\
&& \Llap{\x\,\prsl\,\plmn\!(\eta^sy_m-y_l\:)\vpb{D(n,\ell,s)},} &
\endalignat
$$
\fi
\endpro
\vsk->
\vsk0>
\Pf.
Lemma~\[Pxy] implies that
$$
\det\bigl[P_\lg(\xt\mg)\bigr]\,=\prod_{\lg\in\Zln}^{\vp.}P_\lg(\xt\lg)\qquad
\text{and}
\qquad\det\bigl[P_\lg(\yt\mg)\bigr]\,=\prod_{\lg\in\Zln}P_\lg(\yt\lg)\,.
$$
The rest of the proof consists of several applications of identity \(combi).
\epf
Consider two more determinants,
$\det\bigl[Q_\lg(\xt\mg)\bigr]_{\lg,\mg\in\Zln}$ and
$\,\det\bigl[Q_\lg(\yt\mg)\bigr]_{\lg,\mg\in\Zln}$. Since
\vvn-.5>
\ifMag
$$
\gather
\Line{
\det\bigl[P_\lg(\xt\mg)\bigr]=\det N\cdot\det\bigl[Q_\lg(\xt\mg)\bigr]
\hfil\text{and}\hfil
\det\bigl[P_\lg(\yt\mg)\bigr]=\det N\cdot\det\bigl[Q_\lg(\yt\mg)\bigr]\,,\!\!}
\\
\nn10>
\Line{\Rlap{\text{we have}}\hfil
{\det\bigl[P_\lg(\xt\mg)\bigr]\over\det\bigl[P_\lg(\yt\mg)\bigr]}\;=\;
{\det\bigl[Q_\lg(\xt\mg)\bigr]\over\det\bigl[Q_\lg(\yt\mg)\bigr]}\;,\hfil}
\endgather
$$
\else
$$
\gather
\det\bigl[P_\lg(\xt\mg)\bigr]\,=\,\det N\cdot\det\bigl[Q_\lg(\xt\mg)\bigr]
\qquad\text{and}\qquad
\det\bigl[P_\lg(\yt\mg)\bigr]\,=\,\det N\cdot\det\bigl[Q_\lg(\yt\mg)\bigr]\,,
\\
\nn10>
\Line{\Rlap{\text{we have}}\hfil
{\det\bigl[P_\lg(\xt\mg)\bigr]\over\det\bigl[P_\lg(\yt\mg)\bigr]}\;=\;
{\det\bigl[Q_\lg(\xt\mg)\bigr]\over\det\bigl[Q_\lg(\yt\mg)\bigr]}\;,\hfil}
\endgather
$$
\fi
and by standard arguments of the separation of \var/s we obtain that
$$
\align
\ifMag\kern1.7em\fi
\det\bigl[Q_\lg(\xt\mg)\bigr]_{\lg,\mg\in\Zln} &{}=\,C\>
\eta^{n(1-n)/2\cdot\!\!\tsize{n+\ell-1\choose n+1}}
\prsl\,\plmn\!(\eta^sx_m-x_l)\vpb{D(n,\ell,s)},
\ifMag\kern-1.7em\fi
\Tagg{detQxy}
\\
\nn4>
\det\bigl[Q_\lg(\yt\mg)\bigr]_{\lg,\mg\in\Zln} &{}=\,C\>
\eta^{n(n-1)/2\cdot\!\!\tsize{n+\ell-1\choose n+1}}
\prsl\,\plmn\!(\eta^sy_m-y_l\:)\vpb{D(n,\ell,s)}
\ifMag\kern-1.7em\fi
\endalign
$$
where $C$ is a nonzero constant,
\aagood
which does not depend on $\xn$, $\yn$.
Formulae \(detPxy), \(detQxy) imply formula \(DetN) up to a factor:
$$
\det N\,=\,C\1\,
\pros\,\plmn\!(\eta^sy_l\:-x_m)\vpb{\tsize{n+\ell-s-2\choose n-1}}.
\Tag{DetNC}
$$
\par
To calculate the constant $C$ we consider its dependence on $\eta$.
\Lhs/ of formula \(DetNC) is a \pol/ in $\eta$ of degree
$\dsize{n(n-1)/2\cdot\!{n+\ell-1\choose n+1}}$, the same as the double product
in \rhs/. Thus, $C$ is a \raf/ in $\eta$ with no pole at infinity. Moreover,
since $C$ does not depend on $\xn$, $\yn$, it has no zeros as a \fn/ of $\eta$.
Hence, $C$ does not depend on $\eta$ at all.
\par
Let $\eta=1$, $x_m=x_1$, $y_m=0$, $\mn$, and consider the limit $x_1\to\8$.
In this limit
$$
P_\lg(t)\,=\,(-x_1)^{\sun(n-m)\lg_m}\>\bigl(Q_\lg(t)+o(1)\bigr)\,.
$$
Therefore, $C=1$. The theorem is proved.
\epf
\ifMag\else\fixedpage\fi

\ifAst\newpage\fi
\Sect[B]{Basic facts about the \ehgf/}
Let $\om=\exp(2\pii/n)$ all over the \appx/. Fix a complex number $\al$ \st/
${\al^n\]=p}$. Let $\tht(u)=(u,p)\9(p/u,p)\9(p,p)\9$ be the Jacobi theta-\fn/.
\par
Let $A=\ka\pron z_m$. Fix a complex number $\zt$ \st/
${\zt^n\]=(-1)^{n-1}A\1\!}$. Let $\Ec[A]$ be the space of \hof/s on $\Cx\!$
\st/ $f({p\)u})=A(-u)^{-n}f(u)$. It is easy to see that $\dim\)\Ec[A]=n$,
say by Fourier series. Set
$$
\thi_l(u)\,=\,u^{l-1}\pron\tht(\zt\)\al^{l-1}\om^mu)\,,\qqq\lcn\,.
$$
\Lm*
The \fn/s $\thi_1\lc\thi_n$ form a basis in the space $\Ec[A]$.
\endpro
\Pf.
Clearly, $\thi_l\in\Ec[A]$ for any $\lcn$. Moreover,
${\thi_l(\om\)u)=\om^{l-1}\thi_l(u)}$, that is the \fn/s $\thi_1\lc\thi_n$
are eigenfunctions of the translation operator with distinct \eva/s.
Hence, they are linearly independent.
\epf
Let $\Ec_\ell[A]$ be the space of \sym/ \fn/s in \var/s $\tell$ which are
\hol/ on $\Cxl\!$ and have the property
$$
f(\tpell)\,=\,A\>(-t_a)^{-n}\>f(\tell)\,.
$$
In particular, $\Ec_1[A]=\Ec[A]$. The space $\Ec_\ell[A]$ has dimension
$\dsize{n+\ell-1\choose n-1}$. A basis in the space $\Ec_\ell[A]$ is given by
functions $\Tho_\lg(\tell)$, $\lg\in\Zln$:
$$
\Tho_\lg(\tell)\,=\;{1\over\lg_1!\ldots\lg_n!}\,\susi\,
\pron\,\prod_{a\in\Gm_{\Rph lm}}\thi_m(t_{\si_a})\,.
\Tag{Tho}
$$
Here $\Gm_m=\lb 1+\lg^{m-1}\,\lc\lg^m\rb$, $\mn$.
\Par
Let ${\Fq=\Fq[\)\ka;\zn;\xin;\ell\,]}$ be the \ehgf/.
The \ehgf/ $\Fq(z)$ of a fiber is isomorphic to the space $\Ec_\ell[A]$ by
definition. Therefore
$$
\gather
\dim\Fq(z)\,=\,{n+\ell-1\choose n-1}\,.
\\
\nn-2>
\Text{Set}
\nn-2>
G_\lg(t,z)\,=\,\Tho_\lg(\tell)\,\pron\,\pral\,{1\over\tht(\xi_m\1t_a/z_m)}\;
\prab\ {\tht(t_a/t_b)\over\tht(\eta\)t_a/t_b)}\;.
\Tag{Glg}
\endgather
$$
The \fn/s $G_\lg$, $\lg\in\Zln$ form a basis in $\Fq(z)$.
\Lm{winFq}
For any ${\lg\in\Zln}$ and $z\in\Cxn\!$ the \ewf/ $W_\lg(t,z)$ is in the \ehgf/
$\Fq(z)$ of the fiber.
\endpro
\Pf.
It is clear from definiton \(Wlga) that the \fn/ $W_\lg(t,z)$ has the form
\ifMag\else\vv-.5>\fi
$$
Q(\tell,\zn)\,\pron\,\pral\,{1\over\tht(t_a\xi_m\1/z_m)}\;
\prab\ \,{1\over\tht(\eta\)t_a/t_b)}
$$
where $Q$ is a \hof/ on $\Cxln\!$ with the properties
\ifMag
$$
\align
Q(\tpell,{} & \zn)\,={}
\\
\nn4>
{}=\,(- & t_a)^{-\ell-n}
\tprod_{b=1}^\ell t_b\]\cdot p^{1-a}\)\ka\}\tpron\!z_m\;Q(\tell,\zn)\,,
\endalign
$$
\else
$$
Q(\tpell,\zn)\,=\,(-t_a)^{-\ell-n}
\tprod_{b=1}^\ell t_b\]\cdot p^{1-a}\)\ka\}\tpron\!z_m\;Q(\tell,\zn)\,,
$$
\fi
$\aell$. Furthermore, by construction the \fn/ $W_\lg$ as a \fn/ of $\tell$ is
invariant \wrt/ the action \(ect) of the \symg/ $\Sl\!$, which means that the
\hof/ $Q$ is skew\sym/ \wrt/ the \var/s $\tell$. Hence, the ratio $\Tho$ of the
\fn/ $Q$ and the product $\!\!\dsize\prab\tht(t_a/t_b)$
\vvgood
\vv.1>
is a \hof/ on $\Cxln\!$ which is \sym/ in the \var/s $\tell$ and has
the properties
$$
\Tho (\tpell,\zn)\,=\,(-t_a)^{-n}\,\ka\}\tpron\!z_m\;\Tho(\tell,\zn)\,,
$$
$\aell$; that is the \fn/ $W_\lg$ is in the \ehgf/ of the fiber.
\epf
\Cr{n=1q}
Let $n=1$. Then
\vv-.5>
$$
W\'\ell(\tell,z_1)\,=\,\pral\,
{\tht(\ka\1t_a/z_1)\over\tht(\xi_1\1t_a/z_1)}\;
\prab\ {\tht(t_a/t_b)\over\tht(\eta\)t_a/t_b)}\;.
$$
\endpro
\nt The proof is similar to the proof of Corollary~\[n=1].
\Lm{n=11q}
Let $n=1$. Then
\ifMag
$$
\NN4>
\align
& W\'\ell(\tell,z_1)\,={}
\\
\nn4>
&{}\!=\,\pros\,{\tht(\eta^{-s}\ka\>\xi_1\1)\over
\tht(\eta^{s(s-\ell)}\ka^{\ell-s}\xi_1^{s-\ell})}\,\susi\,\LBc\;\pral\,
{\tht(\eta^{1+(a-1)(a-1-\ell)}\ka^{\ell-a+1}\xi_1^{\)a-1-\ell}t_{a-1}/t_a)\over
\tht(\eta\)t_{a-1}/t_a)}\;\RBc_\si\,.
\endalign
$$
\else
$$
W\'\ell(\tell,z_1)\,=\,\pros\,{\tht(\eta^{-s}\ka\>\xi_1\1)\over
\tht(\eta^{s(s-\ell)}\ka^{\ell-s}\xi_1^{s-\ell})}\,\susi\,\LBc\;\pral\,
{\tht(\eta^{1+(a-1)(a-1-\ell)}\ka^{\ell-a+1}\xi_1^{\)a-1-\ell}t_{a-1}/t_a)\over
\tht(\eta\)t_{a-1}/t_a)}\;\RBc_\si\,.
$$
\fi
Here $\;t_0=\eta\1\xi_1z_1$.
\endpro
\nt
The proof is similar to the proof of Lemma~\[n=11].
\makePcd
\Pf of Lemma\{~\{\[Wlgp].
The proof is similar to the proof of Lemma~\[wlgp]. The summation over the
subgroup $\S^{\lg_1}\!\lx\S^{\lg_n}\sub\Sl$ can be done explicitly using
Corollary~\[n=1q].
\epf
Let $W_\lg$, $\lg\in\Zln$, be the \ewf/s. Define a matrix $\Mq$ by the rule:
$$
W_\lg(t,z)\,=\smiZ \Mq_{\lg\mg}(z)\>G_\mg(t,z)\,,\qqq \lg\in\Zln\,.
\Tag{Mqdef}
$$
Set $\dsize d(n,m,\ell,s)\,=
\sum_{\tsize{\vp( i,j\ge 0\atop\tsize{\vp( i+j<\ell\atop i-j=s\vp(}}}
{m-1+i\choose m-1}\,{n-m-1+j\choose n-m-1}$.
\Th{DetMq}
\ifMag\else\vv->\fi
$$
\gather
\aligned
\det\Mq\}(z)\,& {}=\;\Xi\>\prsl\,\prmn\>\tht\bigl({\tsize\eta^s
\ka\1\!\!\!\prllm\!\!\xi_l\1\!\!\!\prml\!\!\xi_l\)}\bigr)^{d(n,m,\ell,s)}\ \x
\\
\nn6>
& \>{}\x\,\pron(\xi_m\1z_m)\vpb{(m-n){\tsize{n+\ell-1\choose n}}}\;
\pros\,\plmn\>\tht(\eta^s\xi_l\1\xi_m\1z_l/z_m)
\vpb{\tsize{n+\ell-s-2\choose n-1}}
\endaligned
\\
\nn2>
\Text{where}
\nn-8>
\Xi\;=\,\Bigl[\>\Rph{(p)}{(p)\9}\vpb{1-n^2}\>\prmn\>
\Bigl(\>{\tht(\om^m)\over\om^m-1}\)\Bigr)^{\!\!\raise 3pt\mbox{\ssize n-m}}\>
\Bigr]^{\tsize{n+\ell-1\choose n}}.
\endgather
$$
\endpro
\Pf.
For any $\lg\in\Zln$ define a \fn/ $J_\lg(\tell)$ by the rule:
$$
W_\lg(t,z)\,=\,J_\lg(\tell)\,\pron\,\pral\,{1\over\tht(\xi_m\1t_a/z_m)}\;
\prab\ {\tht(t_a/t_b)\over\tht(\eta\)t_a/t_b)}\;,
\Tag{Pilg}
$$
so that
$\dsize\ J_\lg(t,z)\,=\smiZ \Mq_{\lg\mg}(z)\>\Tho_\mg(t,z)$,
${\ \;\lg\in\Zln}$.\quad
\ifMag\vv.1>\fi
Introduce new \var/s $A$, $\xn$, $\yn$:
$$
A=\ka\tpron z_m\,,\qqq
x_m=\xi_mz_m\,,\qqq y_m=\xi_m\1z_m\,,\qquad\qquad\mn\,,
$$
and for any $x,y\in\Cn\!$, $\lg\in\Zln$, define $\xt\lg$, $\yt\lg$ by
formulae \(xtry).
\Lm{Jxy}
$J_\lg(\xt\mg)=0$ unless $\lg\lle\mg$.
$J_\lg(\yt\mg)=0$ unless $\lg\gge\mg$. Moreover,
\ifMag
$$
\gather
\NN6>
\Lline{J_\lg(\xt\lg)\,=\,\pron\,\prod_{s=0}^{\lg_m-1}\Bigl(\>
\tht\bigl({\tsize\eta^{s+1}A\1\!\!\!\prlm\!\!\eta^{\lg_l}y_l\:\!
\prmll\!\!\eta^{-\lg_l}x_l\)}\bigr)\;\x}
\\
\Rline{\x\,\prlm\!\}\tht(\eta^{-s}x_m/y_l\:)\!
\prml\!\}\tht(\eta^{\lg_l-s}x_m/x_l)\>\Bigr)\,,\!}
\\
\ifAst\nn6>\fi
\Lline{J_\lg(\yt\lg)\,=\plmn\!\!\eta^{\lg_l\lg_m}\,
\pron\,\prod_{s=0}^{\lg_m-1}\Bigl(\>
\tht\bigl({\tsize\eta^{-s-1}A\1\!\!\!\prllm\!\!\eta^{\lg_l}y_l\:\!
\prml\!\!\eta^{-\lg_l}x_l\)}\bigr)\;\x}
\\
\Rline{\x\,\prlm\!\}\tht(\eta^{s-\lg_l}y_m/y_l\:)\!
\prml\!\}\tht(\eta^sy_m/x_l)\>\Bigr)\,.\!}
\endgather
$$
\else
$$
\alignat2
\NN6>
J_\lg(\xt\lg)\, &{} =\,\pron\,\prod_{s=0}^{\lg_m-1}\Bigl(\>
\tht\bigl({\tsize\eta^{s+1}A\1\!\!\!\prlm\!\!\eta^{\lg_l}y_l\:\!
\prmll\!\!\eta^{-\lg_l}x_l\)}\bigr)\!
\prlm\!\}\tht(\eta^{-s}x_m/y_l\:)\!
\prml\!\}\tht(\eta^{\lg_l-s}x_m/x_l)\>\Bigr)\,, &&
\\
J_\lg(\yt\lg)\, &{} =\plmn\!\!\eta^{\lg_l\lg_m}\,
\pron\,\prod_{s=0}^{\lg_m-1}\Bigl(\>
\tht\bigl({\tsize\eta^{-s-1}A\1\!\!\!\prllm\!\!\eta^{\lg_l}y_l\:\!
\prml\!\!\eta^{-\lg_l}x_l\)}\bigr)\;\x &&
\\
&& \Llap{\x\;\prlm\!\}\tht(\eta^{s-\lg_l}y_m/y_l\:)\!
\prml\!\}\tht(\eta^sy_m/x_l)\>\Bigr)\,.} &
\endalignat
$$
\vvv->
\fi
\endpro
\nt
The proof is straightforward.
\Par
The rest of the proof is similar to the proof of Theorem~\[DetM]. Using
Lemma~\[Jxy] we calculate the determinants $\det\bigl[J_\lg(\xt\mg)\bigr]$,
$\,\det\bigl[J_\lg(\yt\mg)\bigr]$, \ \cf. Lemma~\[detPxy].
Then by the separation of \var/s we obtain the following formulae
$$
\align
\NN4>
\ifMag\quad\ \fi
\det\bigl[\Tho_\lg &(\xt\mg)\bigr]_{\lg,\mg\in\Zln}=\,K\>
\eta^{n(1-n)/2\cdot\!\!\tsize{n+\ell-1\choose n+1}}\,
\pron(-x_m)\vpb{(m-1)\tsize{n+\ell-1\choose n}}\;\x
\Tagg{detTxy}
\\
& {}\!\x\,
\pros\,\tht\bigl(\eta^{-s}A\1\!\tpron x_m\bigr)^{\!\tsize{n+s-1\choose n-1}}
\prsl\,\plmn\!\tht(\eta^sx_l/x_m)\vpb{D(n,\ell,s)},
\ifMag\kern-2em\fi
\\
\nn12>
\det\bigl[\Tho_\lg &(\yt\mg)\bigr]_{\lg,\mg\in\Zln}=\,K\>
\eta^{n(n-1)/2\cdot\!\!\tsize{n+\ell-1\choose n+1}}
\pron y\_m\vpb{(n-m)\tsize{n+\ell-1\choose n}}\;\x
\\
& {}\!\x\,
\pros\,\tht\bigl(\eta^s A\1\!\tpron y_m\bigr)^{\!\tsize{n+s-1\choose n-1}}
\prsl\,\plmn\!\tht(\eta^s y_m/y_l\:)\vpb{D(n,\ell,s)}
\ifMag\kern-2em\fi
\endalign
$$
as well as the required formula for $\det\Mq\}(z)$ with the constant $\Xi$
replaced by $K\1\!$. Here \fn/s $\Tho_\lg$ are given by formula \(Tho),
$K$ is a nonzero constant which does not depend on $\xn$, $\yn$ and
$$
D(n,\ell,s)\,=\!\!\sum_{\tsize{\,r\in\Zp\!\atop 2r\le\ell-|s|-1}}^{}
\!\!{n+\ell-|s|-2r-3\choose n-2}.
$$
\par
To calculate the constant $K$ we consider its dependence on $\eta$. Any of
formulae \(detTxy) shows that $K$ is a \hof/ in $\eta$ on $\Cx\!$ and
$K({p\)\eta})=K(\eta)$. Hence, $K$ does not depend on $\eta$ at all.
\par
Let $\om=\exp(2\pii/n)$. Take $\eta=1$, $x_m=y_m=\om^{m-1}x_1$, $\mn$.
In this case we have
$$
\Tho_\lg(\xt\mg)\,=\,Q_\lg(\xt\mg)\,
\pron\bigl(x_1^{1-m}\thi_m(x_1)\bigr)^{\lg_m}
$$
for any $\lg,\mg\in\Zln$, where the \fn/ $Q_\lg$ is defined by \(Qlg). Hence,
comparing formulae \(detQxy) and \(detTxy) we find that ${K\]=\)\Xi\1\!}$.
In calculations we use Lemma~\[thtomn] and the equality
$$
\sum_{s=1-\ell}^{\ell-1}\}D(n,\ell,s)\,=\,
{n+\ell-1\choose n}
$$
following from \(combi). The theorem is proved.
\epf
\Lm{thtomn}
Let ${\al^n\]=p}$. Then
$\ \dsize\prod_{l=1}^n\,\pron\tht(\al^{l-1}\om^mu)\,=\,
(p)\9^{\]n^2\}-1{\vrd.7ex>}}\>\tht(u^n)$.
\endpro
\nt
The proof is straightforward using the definition $\tht(u)=(u)\9(p/u)\9(p)\9$.
\goodbreak
\Pf of Lemma\{~\{\[dimFQ].
We will indicate explicitly the dependence of the \ewf/s on $\ka$. Namely,
the \ewf/ $W_\lg[\ka]$ is an element of the \ehgf/ of a fiber
${\Fq[\)\ka,\}\sun\!\lg_m](z)}$.
\par
Under the assumptions of the lemma, the \fn/s ${W_\lg[\)\kone\)]}$,
$\lg\in\Zln$, form a basis of the space $\Fq(z)={\Fq[\)\kone;\ell\,](z)}$
and the \fn/s ${W_\lg[\)\koNe\)]}$, $\lg\in\Zll$, form a basis of the space
${\Fq[\)\koNe;\ell-1\)](z)}$. Since
$$
W_\lg[\)\kone](\tell)\,=\;
{1\over\lg_2!\ldots\lg_n!}\;{\tht(\eta)\over\tht(\eta^{\lg_1})}\,
\susi\,\Lbc\>W_{\lg-\eg(1)}[\)\koNe](\twll)\>\Rbc_\si\,,
$$
for any ${\lg\in\Zln}$ \st/ ${\lg_1>0}$, the \fn/s ${W_\lg[\)\kone\)]}$,
$\,{\lg_1>0}$, ${\lg\in\Zln}$, form a basis of the \bosub/ $\Qc(z)\sub\Fq(z)$,
and the equivalence classes of the \fn/s ${W_\lg[\)\kone\)]}$, $\,{\lg_1=0}$,
${\lg\in\Zln}$, form a basis of the quotient space $\Fq(z)/\Qc(z)$.
\Par
Recall that the space ${\Fq[\)\al;l\>](z)}$ is naturally isomorphic to the
space ${\Ec_l[A]}$ of \sym/ \fn/s in \var/s $t_1\lc t_l$ which are \hol/
on $\C^{\x l}$ and have the property
$$
f(t_1\lc p\)t_a\lc t_l)\,=\,A\>(-t_a)^{-n}\>f(t_1\lc t_l)\,,
$$
where $A=\al\!\pron\!z_m$. Using this \iso/ we observe that the rank of the map
\ifMag
$$
\gather
\aligned
\Fq[\)\koNe;{} & \zn;\xin;\ell-1\)](z)\,\to{}
\\
\nn2>
& {}\!\!\to\,\Fq[\)\kone;\zn;\xin;\ell\,](z)\,,
\endaligned
\\
\nn10>
W(\telm)\,\map\susi\,\Lbc\>W(\twll)\>\Rbc_\si\,,
\endgather
$$
\else
$$
\gather
\Fq[\)\koNe;\zn;\xin;\ell-1\)](z)\,\to\,\Fq[\)\kone;\zn;\xin;\ell\,](z)\,,
\\
\nn8>
W(\telm)\,\map\susi\,\Lbc\>W(\twll)\>\Rbc_\si\,,
\endgather
$$
\fi
does not depend on $\xin$. Therefore,
$$
\dim\Qc(z)\,=\,{n+\ell-2\choose n-1}\qqq\text{and}\qqq
\dim\Fq(z)/\Qc(z)\,=\,{n+\ell-2\choose n-2}
$$
provided that ${\eta^r\!\neps}$ for any $\rtll$, $s\in\Z$.
Lemma~\[dimFQ] is proved.
\epf
\Pf of Lemma\{~\{\[dimFQ'].
The proof is similar to the proof of Lemma~\[dimFQ].
\epf
\Pf of Lemma\{~\{\[ResQ].
Let ${\dsize f(t)\,=\susi\,\Lbc\>W(\twll)\>\Rbc_\si}$
be an element of the \bosub/ $\Qc(z)$. One can see that the term of the sum
corresponding to a \perm/ $\si$ does not contribute to
$\ifMag\Res f(t)\vst{t=\xt\mg}\else\botsmash{\Res f(t)\vst{t=\xt\mg}}\fi$
unless
$$
\si_{\mg^k+1}\lsym>\si_{\mg^{k+1}}>1\,,\qqq\kon-1\,.
$$
Here $\mg^i=\sum_{j=1}^i\mg_j$. Since no \perm/
\ifMag\else\vvn-.4>\fi
satisfies above conditions, the first claim of the lemma is proved.
\Par
Consider the determinant
${\det\bigl[\Res G_\lg(t)\vst{t=\xt\mg}\)\bigr]_{\lg,\mg\in\Zln}}$ where
\fn/s $G_\lg$ are given by \(Glg). Under assumptions of the lemma the first
of formulae \(detTxy) implies that the determinant as a \fn/ of $\ka$
\vv.1>
is not zero for generic $\ka$ and it has a zero of multiplicity
$\dsize{n+\ell-2\choose n-1}$ at $\ka=\kone$.
\ifMag\else\vv-.4>\fi
Hence, for $\ka=\kone$ we have
\ifMag\else\vv->\fi
$$
\dim\>\lb\)f\in\Fq(z)\vert\Res f(t)\vst{t=\xt\mg}\!=0\,,\ \,\mg\in\Zln\)\rb\,
\le\,{n+\ell-2\choose n-1}\,.
$$
\par
Let $\ka=\kone$. We have proved already that
$$
\Qc(z)\,\sub\,
\lb\)f\in\Fq(z)\vert\Res f(t)\vst{t=\xt\mg}\!=0\,,\ \,\mg\in\Zln\)\rb\,.
$$
Since $\dsize\dim\Qc(z)={n+\ell-2\choose n-1}\,$ by Lemma~\[dimFQ],
\,the second claim of Lemma~\[ResQ] is proved.
\epf
\Pf of Lemma\{~\{\[ResQ'].
The proof is similar to the proof of Lemma~\[ResQ].
\epf

\ifAst\newpage\fi
\Sect[S]{The \Spair/s of the \hgf/s of fibers}
In this \appx/ we define pairings of the \hgf/s of fibers which provides
a geometric interpretation for the coefficients $c^\tau_\lg$, \cf. \(clg),
used in the definition of the \tenco/ on the \ehgf/ of a fiber, see formula
\(c=1/N).
\par
In this \appx/ we always suppose that assumptions \(npZ)\;--\;\(assum) hold.
Set
$$
\Omq(\tell)\,=\,{\tsize\pral t_a\1}
\pron\,\pral\,{\tht(\xi_m\1t_a/z_m)\over\tht(\xi_mt_a/z_m)}\;
\prab\ {\tht(\eta\)t_a/t_b)\over\eta\)\tht(\eta\1t_a/t_b)}\;.
\Tag{Omq}
$$
Let $\xt\lg$ and $\yt\lg$, $\lg\in\Zln\!$, be the points defined by \(xtyt).
Recall that we define the multiple residue $\Res f(t)\vst{t=\ts}\!$ by
formula \(Rests).
\par
Consider the \ehgf/s of a fiber ${\Fq(z)=\Fq[\)\ka;\ell\,](z)}$ and
\ifMag\nl\fi
${\Fqt(z)=\Fq[\)\ka\1;\ell\,](z)}$. For any \fn/s $W\in\Fq(z)$ and
$\Wti\in\Fqt(z)$ set
$$
\Sq(W,\Wti)\;=\smiZ
\Res\bigl(\Omq(t)\>W(t)\>\Wti(t)\bigr)\vst{t=\xt\mg}\,.
\Tag{SWWx}
$$
\Lm{SWWy}
For any \fn/s $W\in\Fq(z)$ and $\Wti\in\Fqt(z)$ we have that
$$
\Sq(W,\Wti)\,=\,(-1)^\ell\smiZ
\Res\bigl(\Omq(t)\>W(t)\>\Wti(t)\bigr)\vst{t=\yt\mg}\,.
$$
\endpro
\nt
The statement follows from Lemma~\[Resi].
\Par
The pairing ${\Sq:\Fq(z)\ox\Fqt(z)\>\to\>\C}$ is called the \em{\Spair/}
of the \ehgf/s of a fiber.
\Par
We will indicate explicitly the dependence of the \ewf/s on $\ka$. Namely,
the \ewf/ $W^\tau_\lg[\ka]$ is an element of the \ehgf/ of a fiber
${\Fq[\)\ka;\ell\,](z)}$.
\Th{SWW}
Let $\om\in\S^n$ be the permutation of the maximal length.
Let \(npZ)\;--\;\(assum) hold. Then for any \perm/ $\tau\in\S^n$ and any
$\lg,\mg\in\Zln$ we have that
$$
\Sq\bigl(W^\tau_\lg[\ka]\),W^{\tau\om}_\mg[\ka\1]\)\bigr)\,=\,
\dl_{\lg\mg}\>N^\tau_\lg\,.
$$
Here $\dl_{\lg\mg}$ is the Kronecker symbol and
$N^\tau_\lg=N_{{\vp(}^\tau\!\lg}\:(\xi_{\tau_1}\lc\xi_{\tau_n})$ where
$^\tau\}\lg=(\lg_{\tau_1}\lc\lg_{\tau_n})$ and
\ifMag\else\vv->\fi
\ifAst
\aagood
$$
\align
& \aligned
N_\lg(\xin)\, &{}=\,\pron\,\prod_{s=1}^{\lg_m}\,\Bigl(\>
{\tht(\eta)\over\tho\>\tht(\eta^s)\>\tht(\eta^{1-s}\xi_m^2)}\ \x{}
\\
\nn8>
&{}\>\x\>\tht\bigl(\eta^s\ka\1\!\!\!\tprlm\!\!\eta^{\)\lg_l}\xi_l\1\!\!\!
\tprmll\!\!\eta^{-\lg_l}\xi_l\:\bigr)\>\tht\bigl(\eta^s\ka\!
\tprllm\!\!\eta^{-\lg_l}\xi_l\:\!\tprml\!\!\eta^{\)\lg_l}\xi_l\1\bigr)
\}\Bigr)\,,\!\!\!
\endaligned
\\
\nn6>
& \;\tho\,=\;{d\over du}\tht(u)\vst{u=1}\>=\,-\)(p)^3_\8\,.
\endalign
$$
\else
$$
\gather
\aligned
N_\lg(\xin)\,=\,\pron\,\prod_{s=1}^{\lg_m}\,\Bigl(\> &
{\tht(\eta)\over\tho\>\tht(\eta^s)\>\tht(\eta^{1-s}\xi_m^2)}\ \x{}
\\
\nn8>
{}\x{} & \>\tht\bigl(\eta^s\ka\1\!\!\!\tprlm\!\!\eta^{\)\lg_l}\xi_l\1\!\!\!
\tprmll\!\!\eta^{-\lg_l}\xi_l\:\bigr)\>\tht\bigl(\eta^s\ka\!
\tprllm\!\!\eta^{-\lg_l}\xi_l\:\!\tprml\!\!\eta^{\)\lg_l}\xi_l\1\bigr)
\}\Bigr)\,,\!\!\!
\endaligned
\\
\ifMag\nn4>\else\nn-2>\fi
\Lline{\tho\,=\;{d\over du}\tht(u)\vst{u=1}\>=\,-\)(p)^3_\8\,.}
\endgather
$$
\fi
\endpro
\ifAst
\else
\vsk->
\vsk0>
\fi
\goodbreak
\Pf.
For $\lg,\mg\in\Zln$ say that $\,{\lg\lle\mg}\,$ if
$\;\sum_{i=1}^m\lg_i\le\sum_{i=1}^m\mg_i\;$ for any $\mn-1$. Formulae \(Wlgp)
and \(Wtau) for the \ewf/s imply that
\ifAst
$$
\alignat2
& \Res\bigl(\Omq(t)\>W^\tau_\lg[\ka](t)\>W^{\tau\om}_\mg[\ka\1](t)\bigr)
\vst{t=\xt\ng}\> &&{}=\,0\,,
\\
\nn8>
& \Res\bigl(\Omq(t)\>W^\tau_\mg[\ka](t)\>W^{\tau\om}_\lg[\ka\1](t)\bigr)
\vst{t=\yt\ng}\> &&{}=\,0\,,
\\
\nn4>
\Text{unless $\lg\lle\ng\lle\mg$, and}
\nn4>
& \Res\bigl(\Omq(t)\>W^\tau_\lg[\ka](t)\>W^{\tau\om}_\lg[\ka\1](t)\bigr)
\vst{t=\xt\lg}\>\Rlap{{}=\,N_\lg^\tau\,.} &&
\endalignat
$$
\else
$$
\alignat2
& \Res\bigl(\Omq(t)\>W^\tau_\lg[\ka](t)\>W^{\tau\om}_\mg[\ka\1](t)\bigr)
\vst{t=\xt\ng}\> &&{}=\,0\,,
\\
\nn6>
& \Res\bigl(\Omq(t)\>W^\tau_\mg[\ka](t)\>W^{\tau\om}_\lg[\ka\1](t)\bigr)
\vst{t=\yt\ng}\> &&{}=\,0\,,
\\
\Text{unless $\lg\lle\ng\lle\mg$, and}
& \Res\bigl(\Omq(t)\>W^\tau_\lg[\ka](t)\>W^{\tau\om}_\lg[\ka\1](t)\bigr)
\vst{t=\xt\lg}\>\Rlap{{}=\,N_\lg^\tau\,.} &&
\endalignat
$$
\fi
Therefore, by formula \(SWWx) we have that
$\Sq\bigl(W^\tau_\lg[\ka]\),W^{\tau\om}_\mg[\ka\1]\)\bigr)=0$, unless
$\lg\lle\mg$ and
$$
\Sq\bigl(W^\tau_\lg[\ka]\),W^{\tau\om}_\lg[\ka\1]\)\bigr)\,=\,N^\tau_\lg\,.
$$
Similarly, by formula \(SWWy) we have that
$\Sq\bigl(W^\tau_\lg[\ka]\),W^{\tau\om}_\mg[\ka\1]\)\bigr)=0$, unless
$\lg\gge\mg$. The theorem is proved.
\epf
\Rem
The coefficients $c_\lg\:(\xin)$ defined by formula \(clg) are inverse to
the coefficients $N_\lg(\xin)$ defined in Theorem~\[SWW] up to a common factor.
Namely,
$$
c_\lg\:(\xin)\,=\;{\eta^{\)\ell(1-\ell)/2}\>\ka^\ell\over
(p)^{3\ell}_\8\,N_\lg(\xin)}\;\pron\xi_m^\ell\,.
\Tag{c=1/N}
$$
\enddemo
The \Spair/ of the \thgf/s of a fiber is defined similarly to the \Spair/ of
the \ehgf/s of a fiber. Assumptions \(npZ)\;--\;\(assum) can be replaced by
weaker assumptions
$$
\NN4>
\gather
\eta^r\!\ne 1\,,\qqq\rell\,,
\Tag{ne1}
\\
\xi_m^2\ne\,\eta^r\,,\qqq \mn\,,\qquad\rll\,,
\\
\nn3>
\xi_l^{\pm1}\xi_m^{\pm1}z_l/z_m\,\ne\,\eta^r\,,\qqq l,\mn\,,\quad l\ne m\,,
\qquad\rll\,.
\\
\ald
\ifMag\nn4>\fi
\Text{Let}
\nn-10>
\Om(\tell)\,=\,{\tsize\pral t_a^{-2}}
\pron\,\pral\,{t_a-\xi_mz_m\over\xi_mt_a-z_m}\;
\prab\ {\eta\)t_a-t_b\over t_a-\eta\)t_b}\;.
\Tag{Omtri}
\endgather
$$
For any \fn/s $w,\wti$ in the \thgf/ of a fiber $\Fo(z)$ set
$$
S(w,\wti)\;=\smiZ
\Res\bigl(\Om(t)\>w(t)\>\wti(t)\bigr)\vst{t=\xt\mg}\,.
$$
\Lm{Swwy}
For any \fn/s $w,\wti\in\Fo(z)$ we have that
$$
S(w,\wti)\,=\,(-1)^\ell\smiZ
\Res\bigl(\Om(t)\>w(t)\>\wti(t)\bigr)\vst{t=\yt\mg}\,.
$$
\endpro
\goodbreak
\Pf.
Due to Lemma~\[wbasis] it suffices to prove the statement if both \fn/s
$w,\wti$ are \twf/s. Moreover, since $S(w_\lg,w_\mg)$ depends analytically
on parameters $\eta$, $\xin$, $\zn$, it is enough to prove the lemma under
the assumptions: ${|\eta\)|>1}$ and $|z_m|=1$, $|\xi_m|<1$, $\mn$.
\par
\ifMag\vsk.5>\fi
Consider the integral
${\dsize\int_{\){\TT^\ell}\]}\Om(t)\>w(t)\>\wti(t)\;d^\ell t}$
\vv.2>
\;where $\,{\TT^\ell=\lb\)t\in\Cl\vert\ }{|t_1|=1}\llc{|t_\ell|=1\)\rb}$.
\alb
Similarly to Theorems~\[xsum] and~\[ysum] one can show that under the above
assumptions
$$
\gather
\int_{\){\TT^\ell}\]}\Om(t)\>w(t)\>\wti(t)\;d^\ell t\;=\,(2\pii)^\ell\>\ell\)!
\smiZ\Res\bigl(\Om(t)\>w(t)\>\wti(t)\bigr)\vst{t=\xt\mg}
\\
\nn-3>
\Text{and}
\nn-3>
\int_{\){\TT^\ell}\]}\Om(t)\>w(t)\>\wti(t)\;d^\ell t\;=\,(-2\pii)^\ell\>\ell\)!
\smiZ\Res\bigl(\Om(t)\>w(t)\>\wti(t)\bigr)\vst{t=\yt\mg}\,.
\endgather
$$
The lemma is proved.
\epf
\goodbreak
\Th{Sww}
Let $\om\in\S^n$ be the permutation of the maximal length. Let assumptions
\(ne1) hold. Then for any \perm/ $\tau\in\S^n$ and any $\lg,\mg\in\Zln$ we
have that
$$
S(w^\tau_\lg,w^{\tau\om}_\mg)\,=\,\dl_{\lg\mg}\,\pron\,\prod_{s=1}^{\lg_m}\,
{1-\eta\over z_m\>(1-\eta^s)\>(\xi_m^2-\eta^{s-1})}
$$
where $\dl_{\lg\mg}$ is the Kronecker symbol.
\endpro
\nt
The proof is similar to the proof of Theorem~\[SWW].
\Par
The pairing ${S:\Fo(z)\ox\Fo(z)\>\to\>\C}$ is called the \em{\Spair/}
of the \thgf/s of a fiber.
\Par
The \tri/ \Spair/ $S$ can be considered as a degeneration either of the
elliptic \Spair/ $\Sq$ or the \hpair/
\ifAst $I:\Fq(z)\ox\Fo(z)\alb\to\>\C$. \else
${I:\Fq(z)\ox\Fo(z)\>\to\>\C}$. \fi
Namely, let $p\to 0$ and after that $\ka\to\8$. Then in this limit we have that
$$
(u)\9\to\> 1-u\,,\qqq \tht(u)\>\to\> 1-u\,,\qqq (p)\9\to\> 1\,,
$$
the \ewf/s tend to the respective \twf/s:
\ifMag
$$
\gather
W_\lg[\ka](t)\,\to\,w_\lg(t)\,(-1)^\ell\tpral t_a\1\pron z_m^{\lg_m}\!\!
\plmn\!\!\xi_l^{\lg_m}\,,
\\
\nn6>
 \ka^{-\ell}\,W_\lg[\ka\1](t)\,\to\,w_\lg(t)\!\plmn\!\!\xi_m^{\lg_l}\,,
\endgather
$$
\else
$$
W_\lg[\ka](t)\,\to\,w_\lg(t)\,(-1)^\ell\tpral t_a\1\pron z_m^{\lg_m}\!\!
\plmn\!\!\xi_l^{\lg_m}\,,\qqq
\ka^{-\ell}\,W_\lg[\ka\1](t)\,\to\,w_\lg(t)\!\plmn\!\!\xi_m^{\lg_l}\,,
$$
\fi
the \fn/ $\Omq$ defined in \(Omq) and the short \phf/ $\Phi$ defined in \(Phis)
turn into the \fn/ $\Om$ given by \(Omtri):
$$
\Omq(t)\,\to\,\Om(t)\>\tpral t_a\>\pron\xi_m^{-\ell}\,,\qqq
\Phi(t)\,\to\,\Om(t)\;\eta^{\ell(\ell-1)/2}\tpral t_a^2\>\pron\xi_m^{-\ell}\,,
$$
and both the elliptic \Spair/ and the \hpair/ become the \tri/ \Spair/:
$$
\gather
\ka^{-\ell}\,\Sq\bigl(W_\lg[\ka]\),W_\mg[\ka\1]\)\bigr)\,\to\,S(w_\lg,w_\mg)\,
(-1)^\ell\pron z_m^{\lg_m}\!\!\plmn\!\!\xi_m^{\mg_l-\lg_l}\,,
\\
\nn6>
I\bigl(W_\lg[\ka]\),w_\mg\bigr)\,\to\,S(w_\lg,w_\mg)\,(-1)^\ell\>
\eta^{\ell(\ell-1)/2}\pron z_m^{\lg_m}\!\!\plmn\!\!\xi_m^{-\lg_l}\,.
\endgather
$$
\par
In the rest of the \appx/ we formulate and prove Lemma~\[Resi] which implies
Lemma~\[SWWy].
\par
Let $u_1\lc u_k$ be nonzero complex numbers. Consider a \fn/ $F(\tell)$ which
is quasiperiodic \wrt/ each of the \var/s $\tell$:
$$
\gather
F(\tpell)\,=\,p\1\)F(\tell)\,,
\\
\Text{and which has the form}
\nn-4>
f(\tell)\,\prod_{j=1}^k\,\pral\,{1\over\tht(t_k/u_j)}\;
\prod_{a=1}^\ell\,\prod_{\tsize{b=1\atop b\ne a}}^\ell
\ {\tht(t_a/t_b)\over\tht(\eta\)t_a/t_b)}
\Tag{Fform}
\endgather
$$
where $f(\tell)$ is a \sym/ \hof/ on ${\Cxl\!}$. We call the points
$u_1\lc u_k$ the ``root'' singularities of the \fn/ $F$.
\par
For any $\lg\in\Zp^k$, $\sum_{j=1}^k\lg_m=\ell$, introduce the points
$\uxt\lg,\,\uyt\lg\in\Cxl\!$ by the rule:
$$
\gather
\uxt\lg\,=\,(\eta^{1-\lg_1}u_1,\>\eta^{\)2-\lg_1}u_1\lc u_1,
\>\eta^{1-\lg_2} u_2\lc u_2,\,\ldots\,,\)\eta^{1-\lg_k}u_k\lc u_k)\,,
\\
\nn4>
\uyt\lg\,=\,(\eta^{\lg_1-1}u_1,\>\eta^{\lg_1-2}u_1\lc u_1,
\>\eta^{\lg_2-1} u_2\lc u_2,\,\ldots\,,\)\eta^{\lg_k-1}u_k\lc u_k)\,,
\endgather
$$
\cf. \(xtyt). For any $\vp{\sum^1} i=1\lc k\,$ let
\vv-.5>
$$
\gather
\Zps_i\)=\,\lb\,\lg\in\Zp^k\vert\ \lg_{i+1}\lsym=\lg_k=0\,,
\ \ \tsum_{j=1}^k\lg_j=\ell\,\rb
\\
\nn-4>
\Text{and}
\nn-6>
\Zms_i\)=\,\lb\,\lg\in\Zp^k\vert\ \lg_1\lsym=\lg_i=0\,,
\ \ \tsum_{j=1}^k\lg_j=\ell\,\rb\,.
\endgather
$$
\vvv-.5>
In particular, $\>\Zps_0=\Zms_k=\Empty\,\vp\sum$ and
$\>\Zps_k=\Zms_0=\lb\>\lg\in\Zp^k\vert\sum_{j=1}^k\lg_j=\ell\>\rb$.
\Lm{Resi}
For any $i=1\lc k$ we have that
$$
\NN3>
\gather
\sum_{\mg\in\Rph Z{\Zps_i}}\,\)\Res F(t)\vst{t=\uxt\mg}\,=\,(-1)^\ell\,
\sum_{\mg\in\Rph Z{\Zms_i}}\,\)\Res F(t)\vst{t=\uyt\mg}\,,
\\
\nn3>
\Text{provided that the singularity hyperplanes}
\nn4>
\aligned
& t_a=u_j\,,
\\
& t_a=\eta\)t_b\,,
\endaligned
\qqq
\gathered
a=1\lc\ell\,,\quad j=1\lc k\,,
\\
a,b=1\lc\ell\,,\quad a\ne b\,,
\endgathered
\endgather
$$
have no multiple intersections at the points $\,\uxt\mg$, $\mg\in\Zps_i\!$,
\,and at the points $\,\uyt\mg$, $\mg\in\Zms_i\!$.
\endpro
\Pf.
We prove the lemma by induction \wrt/ $\ell$. For $\ell=1$ the statement is
standard, which provides the base of induction.
\par
Let $\ell>1$. To avoid complicated notations we first give the idea of
the proof in the general case and then explain technical details for
the example $\ell=3$, $k=3$, $i=1$.
\par
The \fn/ $F$ considered as a \fn/ of $\twll$ for a fixed $t_1$ has a form
similar to \(Fform) with the ``root'' singularities at the points
${\eta\)t_1,\,\eta\1 t_1,\,{}}u_1\lc u_k$. Using the induction assumption we
apply the lemma in this case taking in \lhs/ of the formula the sum of residues
corresponding to the ``root'' singularities $u_1\lc u_i,{\,\eta\1 t_1}$ and
taking in \rhs/ of the formula the sum of residues corresponding to the
``root'' singularities $u_{i+1}\lc u_k,{\,\eta\)t_1}$. To complete the proof we
apply the lemma to the \fn/ of $t_1$ given by the sum of residues in \lhs/ of
the formula obtained at the previous step.
\Par
Consider the example $\ell=3$, $k=3$, $i=1$. At the first step we get
\ifMag
$$
\NN8>
\gather
\Lline{
\aligned\quad
\Res\bigl(\)\Res F(t)\vst{t_3=u_1}\bigr)\vst{t_2=\eta\1u_1}\! & {}+\>
\Res\bigl(\)\Res F(t)\vst{t_3=u_1}\bigr)\vst{t_2=\eta\1t_1}+{}
\\
& {}+\>
\Res\bigl(\)\Res F(t)\vst{t_3=\eta\1t_1}\bigr)\vst{t_2=\eta^{-2}t_1}={}
\endaligned}
\\
\ald
\nn4>
\Rline{
\alignedat2
{}=\,\Res\bigl(\) &\Rlap{\Res F(t)\vst{t_3=u_2}\bigr)\vst{t_2=\eta\)u_2}\!+\>
\Res\bigl(\)\Res F(t)\vst{t_3=u_3}\bigr)\vst{t_2=\eta\)u_3}+{}}
\\
{}+{}\> &\!\Res\bigl( &&\!\Res F(t)\vst{t_3=u_3}\bigr)\vst{t_2=u_2}\!+\>
\Res\bigl(\)\Res F(t)\vst{t_3=u_2}\bigr)\vst{t_2=\eta\)t_1}+{}
\\
&& \Llap{{}+{}\>}&\!
\Res\bigl(\)\Res F(t)\vst{t_3=u_3}\bigr)\vst{t_2=\eta\)t_1}\!+\>
\Res\bigl(\)\Res F(t)\vst{t_3=\eta\)t_1}\bigr)\vst{t_2=\eta^2\]t_1}.\quad
\endalignedat}
\endgather
$$
\else
$$
\Lline{
\aligned
\Res\bigl(\)\Res F(t)\vst{t_3=u_1}\bigr)\vst{t_2=\eta\1u_1}\!+\>
\Res\bigl(\)\Res F(t)\vst{t_3=u_1}\bigr)\vst{t_2=\eta\1t_1}\!+\>
\Res\bigl(\)\Res F(t)\vst{t_3=\eta\1t_1}\bigr)\vst{t_2=\eta^{-2}t_1}={} &
\\
\nn10>
{}=\,\Res\bigl(\)\Res F(t)\vst{t_3=u_2}\bigr)\vst{t_2=\eta\)u_2}\!+\>
\Res\bigl(\)\Res F(t)\vst{t_3=u_3}\bigr)\vst{t_2=\eta\)u_3}\!+\>
\Res\bigl(\)\Res F(t)\vst{t_3=u_3}\bigr)\vst{t_2=u_2}+{} &
\\
\nn10>
\Llap{{}+\,{}}
\Rph{\Res F(t)\vst{t_3=u_2}\bigr)\vst{t_2=\eta\)u_2}\!+\>
\Res\bigl(\)\Res F(t)\vst{t_3=u_3}\bigr)\vst{t_2=\eta\)u_3}\!+\>
\Res\bigl(\)\Res F(t)\vst{t_3=u_2}\bigr)\vst{t_2=u_3}+{}}
{\Res\bigl(\)\Res F(t)\vst{t_3=u_2}\bigr)\vst{t_2=\eta\)t_1}\!+\>
\Res\bigl(\)\Res F(t)\vst{t_3=u_3}\bigr)\vst{t_2=\eta\)t_1}\!+\>
\Res\bigl(\)\Res F(t)\vst{t_3=\eta\)t_1}\bigr)\vst{t_2=\eta^2\]t_1}.}&
\endaligned}
$$
\fi
Denote by $G(t_1)$ the sum of residues in \lhs/ of the above formula, or
\eqv/ly, the sum of residues in \rhs/ of the above formula. The \fn/ $G$ can
have poles only at the following points:
$p^s\eta^{-r}u_1,\,p^s\eta^ru_2,\,p^s\eta^ru_3$, $r=0,1,2$, $s\in\Z$,
because at any other point at least one of the defining expressions
for the \fn/ $G$ has no pole. Hence, we have that
$$
-\sum_{s=0}^2\>\Res G(t_1)\vst{t_1=\eta^{-s}\]u_1}=\,
\sum_{s=0}^2\>\bigl(\>\Res G(t_1)\vst{t_1=\eta^s\]u_2}\!+\>
\Res G(t_1)\vst{t_1=\eta^s\]u_3}\)\bigr)\,.
$$
Substituting \resp/ the left (right) hand side definition of the \fn/ $G$ into
the left (right) hand side of the above formula, we obtain that
\ifMag
$$
\NN8>
\alds
\align
{}-\,\Tres{t_3=u_1}{t_2=\eta\1u_1}{t_1=\eta^{-2}u_1}\,-\,
\Tres{t_3=u_1}{t_2=\eta\1t_1}{t_1=\eta\1u_1}\,-\,
\Tres{t_3=\eta\1t_1}{t_2=\eta^{-2}t_1}{t_1=u_1}\ =\qquad
\hp{\Tres{t_3=u_3}{t_2=\eta\)u_3}{t_1=\eta^{\)2}u_3}\,+{}} &
\\
\ifAst\nn2>\fi
\nn2>
{}=\ \Tres{t_3=u_2}{t_2=\eta\)u_2}{t_1=\eta^{\)2}u_2}\,+\,
\Tres{t_3=u_2}{t_2=\eta\)u_2}{t_1=u_3}\,+\,
\Tres{t_3=u_3}{t_2=\eta\)u_3}{t_1=u_2}\,+\,
\Tres{t_3=u_3}{t_2=\eta\)u_3}{t_1=\eta^{\)2}u_3}\,+\quad{} &
\\
\ifAst\nn4>\fi
{}+\,\Tres{t_3=u_3}{t_2=u_2}{t_1=\eta\)u_2}\,+\,
\Tres{t_3=u_3}{t_2=u_2}{t_1=\eta\)u_3}\,+\,
\Tres{t_3=u_2}{t_2=\eta\)t_1}{t_1=\eta\)u_2}\,+\,
\Tres{t_3=u_2}{t_2=\eta\)t_1}{t_1=u_3}\,+{} &
\\
\ifAst\nn4>\fi
{}+\,\Tres{t_3=u_3}{t_2=\eta\)t_1}{t_1=u_2}\,+\,
\Tres{t_3=u_3}{t_2=\eta\)t_1}{t_1=\eta\)u_3}\,+\,
\Tres{t_3=\eta\)t_1}{t_2=\eta^{\)2}t_1}{t_1=u_2}\,+\,
\Tres{t_3=\eta\)t_1}{t_2=\eta^{\)2}t_1}{t_1=u_3} & \,,
\endalign
$$
\else
$$
\alds
\align
{}-\,\Tres{t_3=u_1}{t_2=\eta\1u_1}{t_1=\eta^{-2}u_1}\,-\,
\Tres{t_3=u_1}{t_2=\eta\1t_1}{t_1=\eta\1u_1}\,-\,
\Tres{t_3=\eta\1t_1}{t_2=\eta^{-2}t_1}{t_1=u_1}\ =
\ \Tres{t_3=u_2}{t_2=\eta\)u_2}{t_1=\eta^{\)2}u_2}\,+\,
\Tres{t_3=u_2}{t_2=\eta\)u_2}{t_1=u_3}\,+{}\quad &
\\
\nn10>
{}+\,\Tres{t_3=u_3}{t_2=\eta\)u_3}{t_1=u_2}\,+\,
\Tres{t_3=u_3}{t_2=\eta\)u_3}{t_1=\eta^{\)2}u_3}\,+\,
\Tres{t_3=u_3}{t_2=u_2}{t_1=\eta\)u_2}\,+\,
\Tres{t_3=u_3}{t_2=u_2}{t_1=\eta\)u_3}\,+\,
\Tres{t_3=u_2}{t_2=\eta\)t_1}{t_1=\eta\)u_2}\,+{}\, &
\\
\nn10>
{}+\,\Tres{t_3=u_2}{t_2=\eta\)t_1}{t_1=u_3}\,+\,
\Tres{t_3=u_3}{t_2=\eta\)t_1}{t_1=u_2}\,+\,
\Tres{t_3=u_3}{t_2=\eta\)t_1}{t_1=\eta\)u_3}\,+\,
\Tres{t_3=\eta\)t_1}{t_2=\eta^{\)2}t_1}{t_1=u_2}\,+\,
\Tres{t_3=\eta\)t_1}{t_2=\eta^{\)2}t_1}{t_1=u_3}&
\endalign
$$
\fi
\ifAst\vsk.2>\nt\fi
where we use the notation
$\dsize\ \Tres{t_3=c}{t_2=b}{t_1=a}\,=\,\Res\bigl(\)\Res\bigl(\)\Res F(t)
\vst{t_3=c}\bigr)\vst{t_2=b}\bigr)\vst{t_1=a}$.
\ifAst\Par\else\vsk.25>\fi
To complete the proof we have to transform the multiple residues to the form
\(Rests). This is straightforward under the assumptions of the lemma because
at each step we have to calculate residues only at simple poles and the \fn/
$F$ is \sym/. The transformation can be done term by term, for instance,
\ifAst
\fixedpage\nt
$$
\align
\Res\bigl(\)\Res\bigl(\) & \Res F(t)\vst{t_3=u_2}\bigr)
\vst{t_2=\eta\)t_1}\bigr)\vst{t_1=u_3}={}
\\
\ald
\nn8>
{}=\,\Res\bigl(\) & \Res\bigl(\)\Res F(t)\vst{t_3=u_2}\bigr)
\vst{t_1=u_3}\bigr)\vst{t_2=\eta\)u_3}={}
\\
\ald
\nn8>
{}=\,{} & \Res\bigl(\)\Res\bigl(\)\Res F(t)\vst{t_3=u_3}\bigr)
\vst{t_2=\eta\)u_3}\bigr)\vst{t_1=u_2}.
\endalign
$$
\else
$$
\align
\Res\bigl(\)\Res\bigl(\)\Res F(t)\vst{t_3=u_2}\bigr)
\vst{t_2=\eta\)t_1}\bigr)\vst{t_1=u_3}=\,
\Res\bigl(\)\Res\bigl(\)\Res F(t)\vst{t_3=u_2}\bigr)
\vst{t_1=u_3}\bigr)\vst{t_2=\eta\)u_3}={} &
\\
\nn8>
{}=\,\Res\bigl(\)\Res\bigl(\)\Res F(t)\vst{t_3=u_3}\bigr)
\vst{t_2=\eta\)u_3}\bigr)\vst{t_1=u_2} & .
\endalign
$$
\fi
The lemma is proved.
\epf

\ifAst\newpage\fi
\Sect[C]{The \{\qSelberg/ integral}
In this \appx/ we give proofs of formulae \(qbeta) and \(qbetaR). The last
formula is \eqv/ to the formula for the \qSelberg/ integral, see
\Cite{AK, Theorem 3.2}.
\par
Denote
by ${\Fabc}$ the integrand of the integral \(qbeta):
$$
\Fabc\,=\,\prkl\,{\tht(c\)t_k)\over t_k\>(at_k)\9\>(b/t_k)\9\!}\;
\pjkl{(t_j/t_k)\9\over(x\)t_j/t_k)\9\!}\;.
\Tag{Fabc}
$$
\Pf of formula \{\(qbeta).
Consider the integral in \lhs/ of \(qbeta) as a \fn/ of $c$ and denote it
by $S(c)$. Let $f(c)$ be the ratio of $S(c)$ and \rhs/ of formula \(qbeta).
\par
The \fn/ $S(c)$ satisfies a \deq/
$$
S(p\)c)\,=\,S(c)\,\pros\,{1-x^sa/c\over1-x^sbc}\;,
\Tag{c-eq}
$$
\cf. Corollary~\[Sdeq]. \Rhs/ of formula \(qbeta) solves the same \deq/ \wrt/
$c$. Therefore, $f(c)$ is a \p-periodic \fn/: ${f(p\)c)=f(c)}$.
\par
$S(c)$ is a \hof/ of $c$ on $\Cx$, since the integrand ${\Fabc}$ is a \hof/
of $c$ on $\Cx$ and the integration contour is compact. So, the \fn/
$f(c)$ is regular in the annulus ${|p\)a|<|c|<|b|\1}$ of width greater than
$|p|$. Hence, $f(c)$ is a \hof/ of $c$ on $\Cx\!$, and therefore, $f(c)$ is
a constant \fn/.
\par
We will show that $f=1$ by induction \wrt/ $\ell$. We will indicate
the dependence of $\ell$ explicitly, that is $S(c)=S_\ell(c)$ and $f=f_\ell$.
The trivial case $\ell=0$ provides the base of the induction.
\par
Formulae \(Wlga) and \(n=1q) after a suitable change of notation give the
following identity:
$$
\prkl\>\tht(c\)t_k)\;=\,\prosl\,{\tht(x)\over\tht(x^s)}\,
\susi\,\Bigl(\,\prkl\>\tht(x^{\ell-2k+1}c\)t_{\si_k})
\prjk\,{\tht(xt_{\si_k}/t_{\si_j})\over\tht(t_{\si_k}/t_{\si_j})}\,\Bigr)\,.
$$
Replace the product ${\prkl\tht(c\)t_k)}$ in the integrand $F(\tell)$ by \rhs/
of the identity. Since the rest of the integrand is a \sym/ \fn/ in $\tell$
and the integration contour is invariant \wrt/ \perm/s of the \var/s $\tell$,
we can keep in the sum only the term corresponding to the identity \perm/ and
then multiply the result of the integration by ${\ell\)!}$.
Therefore, we have that
$$
S_\ell(c)\,=\,\ell\)!\,\prosl\,{\tht(x)\over\tht(x^s)}\;
\int_{\)\TT^\ell\]}\;\prkl\,
{\tht(x^{\ell-2k+1}c\)t_k)\over t_k\>(at_k)\9\>(b/t_k)\9\!}\;
\prjk{(1-t_j/t_k)\>(p\)x\1t_j/t_k)\9\over(x\)t_j/t_k)\9}\;\dt\,.
$$
\par
Let ${c=p\)x^{1-\ell}b\1}$. Then the integrand of the above integral is regular
in the punctured disk $0<|t_1|\le 1$ and has a simple pole at $t_1=0$.
Perfoming the integration \wrt/ $t_1$ we obtain
$$
S_\ell(p\)x^{1-\ell}b\1)\,=\,2\pii\>\ell\>(p)\9\;{\tht(x)\over\tht(x^\ell)}\;
S_{\ell-1}(p\)x^{-\ell}b\1)\,.
$$
\Rhs/ of formula \(qbeta) satisfies the same recurrence relation \wrt/ $\ell$.
Hence, $f_\ell=f_{\ell-1}$, which completes the proof.
\epf
\Lm{X1Xl}
Let $\ {\dsize X_k\;=\int_{\TT^\ell}t_1\ldots t_k\>\Fabc\dt}$, \ $\koll$.
\vv.5>\ifMag\nl\fi
The following recurrence relation holds\/{\rm:}
$$
X_k\,=\,X_{k-1}\;{k\>(1-x^{\ell-k+1})\>(p-x^{k-1}bc)\over
(\ell-k+1)\>(1-x^k)\>(p\)x^{\ell-k}a-c)}\;,\qqq\kell\,.
$$
\endpro
\Pf.
Consider the integrals
$$
\Xti_k\;=\int_{\TT^\ell}(1-at_1)\>t_2\lc t_k\>
\prod_{j=2}^\ell\,{xt_1-t_j\over t_1-t_j}\,\,\Fabc\dt\,,
\Tag{Xtia}
$$
$\kell$. Notice that the integrands are regular on $\TT^\ell$ since ${\Fabc}$
vanishes at all diagonals $t_i=t_j$.
\par
Replacing $t_1$ by $t_1/p$ in the integrand and using the explicit formula
for $\Fabc$ we obtain that
$$
\Xti_k\;=\!\int_{\,\TT_p\x\TT^{\ell-1}\!}\!\!p\1\)c\>(b-t_1)\>t_2\lc t_k\>
\prod_{j=2}^\ell\,{t_1-xt_j\over t_1-t_j}\,\,\Fabc\dt
$$
where ${\TT_p=\lb\)t_1\in\C\vert\ |t_1|=p\>\rb}$. The integrand considered as
a \fn/ of $t_1$ is regular in the annulus $p\le|t_1|\le 1$. Therefore, we can
replace the integration contour ${\TT_p\x\TT^{\ell-1}}$ in the above integral
by $\TT^\ell$ without changing the integral:
$$
\Xti_k\;=\int_{\TT^\ell}p\1\)c\>(b-t_1)\>t_2\lc t_k\>
\prod_{j=2}^\ell\,{t_1-xt_j\over t_1-t_j}\,\,\Fabc\dt\,.
\Tag{Xtib}
$$
Since the integration contour $\TT^\ell$ is invariant \wrt/ \perm/s of the
\var/s $\tell$, we can symmetrize the integrands in the formulae for
$X_0\lc X_\ell$ and $\Xti_1\lc\Xti_\ell$. Then formula \(Xtia) and
the first two identities \(xident) imply that
$$
(1-x)\>\Xti_k\,=\;{1-x^{\ell-k+1}\over \ell-k+1}\,X_{k-1}\,-\;
{x^{\ell-k}-x^\ell\over k}\,a\>X_k\,,\qqq \kell\,,
$$
\goodast\nt
and formula \(Xtib) and the last two identities \(xident) imply that
$$
(1-x)\>\Xti_k\,=\;{x^{k-1}-x^\ell\over \ell-k+1}\,p\1b\)c\>X_{k-1}\,-\;
{1-x^k\over k}\,p\1c\>X_k\,,\qqq \kell\,,
$$
because the \fn/ $\Fabc$ is \sym/ in the \var/s $\tell$. The rest of the proof
is obvious. The lemma is proved.
\epf
\Cr{Sdeq}
The \deq/ \(c-eq) holds.
\endpro
\Pf.
The statement is clear since $S(c)=X_0$ and ${S(c/p)=(-c/p)^\ell\>X_\ell}$.
\epf
\Lm{xident}
The following identities hold\/{\rm:}
$$
\align
k\>(1-x)\>\susi\>\Bigl(\tsik\>\prod_{j=2}^\ell\,
{xt_{\si_1}-t_{\si_j}\over t_{\si_1}-t_{\si_j}}\,\Bigr)\, &{}=\;
(x^{\ell-k}-x^\ell)\>\susi\tsik\,,
\\
(\ell-k)\>(1-x)\>\susi\>\Bigl(\tsik\>\prod_{j=1}^{\ell-1}\,
{xt_{\si_1}-t_{\si_j}\over t_{\si_1}-t_{\si_j}}\,\Bigr)\,
&{}=\;(1-x^{\ell-k})\>\susi\tsik\,,
\\
\nn3>
k\>(1-x)\>\susi\>\Bigl(\tsik\>\prod_{j=2}^\ell\,
{t_{\si_1}-xt_{\si_j}\over t_{\si_1}-t_{\si_j}}\,\Bigr)\, &{}=\;
(1-x^k)\>\susi\tsik\,,
\\
(\ell-k)\>(1-x)\>\susi\>\Bigl(\tsik\>\prod_{j=1}^{\ell-1}\,
{t_{\si_1}-xt_{\si_j}\over t_{\si_1}-t_{\si_j}}\,\Bigr)\,
&{}=\;(x^k-x^\ell)\>\susi\tsik\,.
\endalign
$$
\endpro
\Pf.
\Lhs/s of the above formulae are homogeneous \sym/ \pol/s in the \var/s $\tell$
of the homogeneous degree $k$ and of degree one in each of the \var/s $\tell$.
Hence, they are proportional to $\susi^{\vp.}\tsik$.
\par
Restrict the polynomials in question to the line $t_j=x^{j-1} t_1$,
$j=1\lc\ell$, and use the following identity
\vvn-.8>
$$
\sum_{0\le r_1\lsym< r_k\le\ell}\!\!\!x^{\)r_1\lsym+r_k}=\,
\prod_{s=0}^{k-1}\,{x^s-x^\ell\over 1-x^{s+1}}\;.
\Tag{xfact}
$$
Then the calculation of the proportionality coefficients is straighforward.
Identity \(xfact) can be proved by induction \wrt/ $\ell$.
\epf
Let
\vvn-1.5>\nl\strut
$$
S(\tell)\,=\,\prkl\,{(p\)t_k)\9\over(\al t_k)\9}
\prjk\){(1-t_k/t_j)\>(p\)x\1t_k/t_j)\9\over(x\)t_k/t_j)\9}\;.
$$
\Th{qbetaR}
Let $|u|<\min\>(1,|x^{\ell-1}|)$. Then
$$
\align
\sum_{k=1}^\ell\,\sum_{r_k=0}^\8 u^{\>\sum_{i=1}^\ell(\ell-i+1)\)r_i}
x^{-\!\sum_{i=1}^\ell(i-1)(\ell-i+1)\)r_i}
S(p^{r_1},p^{r_1+r_2}x & {}\lc p^{r_1\lsym+r_\ell}x^{\ell-1})\,={}
\\
\nn-4>
{}={}\] & \)\,\pros\;
{(x)\9\>(x^s\al u)\9\>(p)\9\over(x^{s+1})\9\>(x^s\al)\9(x^{-s}u)\9\!}
\endalign
$$
provided the parameters $\al$ and $x$ are \st/ all the terms of the sum are
regular.
\endpro
\Pf.
The sum in \lhs/ of the formula is absolutely convergent and, hence, defines
an \anf/ of the parameters $\al, u, x$. Therefore, it suffices to prove the
formula under the assumption $|\al|<1$, $|x|<1$.
\par
Under this assumption formula \(qbetaR) follows from formula \(qbeta) and
Lemma~\[qbetaS], since the sum in formula \(qbetaR) is proportional term by
term to the sum in formula \(qbetaS) if we identify ${\al=ab}$,
${u=p\)b\1c\1}\!$ and the proportionality coefficient equals
$\,{\pros\,\bigl((p)\9^2/\tht(x^{-s}u)\bigr)}$.
\epf
\Lm{qbetaS}
Let $\,|a|<1$, $|b|<1$, $|x|<1$ and $\,|x^{\ell-1}b\)c|>|p|$. Then
\ifMag
$$
\gather
\Lline{\int_{\TT^\ell}\Fabc\dt\,=\,(2\pii)^\ell\>\ell\)!\,\x{}}
\\
\nn-6>
\Rline{\x\,
\sum_{k=1}^\ell\,\sum_{r_k=0}^\8 \Res\bigl(\)\Res\bigl(\,\ldots\,
\Res\Fabc\vst{t_\ell=p^{r_\ell}xt_{\ell-1}}\!\ldots\,\bigr)
\vst{t_2=p^{r_2}xt_1}\)\bigr)\vst{t_1=p^{r_1}b}\,.\!}
\endgather
$$
\else
$$
\align
\int_{\TT^\ell} & \Fabc\dt\,={}
\\
\nn-6>
& \!\}{}=\,(2\pii)^\ell\>\ell\)!\;
\sum_{k=1}^\ell\,\sum_{r_k=0}^\8 \Res\bigl(\)\Res\bigl(\,\ldots\,
\Res\Fabc\vst{t_\ell=p^{r_\ell}xt_{\ell-1}}\!\ldots\,\bigr)
\vst{t_2=p^{r_2}xt_1}\)\bigr)\vst{t_1=p^{r_1}b}\,.
\endalign
$$
\fi
\endpro
\Pf.
We begin the proof with the next identity which follows from formulae
\(n=1q), \(n=11q) after a suitable change of notations:
\ifAst
$$
\align
\prkl\,{\tht(c\)t_k)\over\tht(t_k/b)}\;\,=\;\, &
\pros\,{\tht(p\)x^sbc)\over\tht(x^{s(s-\ell)}(bc)^{s-\ell})}\ \x{}
\\
\nn6>
{}\x\susi\,\Bigl(\, & \prkl
{\tht(x^{1+(k-1)(k-1-\ell)}(bc)^{k-1-\ell}t_{\si_{k-1}}/t_{\si_k})\over
\tht(xt_{\si_{k-1}}/t_{\si_k})}
\prjk\,{\tht(xt_{\si_j}/t_{\si_k})\over\tht(t_{\si_j}/t_{\si_k})}\,\Bigr)
\endalign
$$
\else
$$
\align
\prkl\,{\tht(c\)t_k)\over\tht(t_k/b)}\;\, &{}=\;\,
\pros\,{\tht(p\)x^sbc)\over\tht(x^{s(s-\ell)}(bc)^{s-\ell})}\ \x{}
\\
\nn6>
& {}\>\x\susi\,\Bigl(\,\prkl
{\tht(x^{1+(k-1)(k-1-\ell)}(bc)^{k-1-\ell}t_{\si_{k-1}}/t_{\si_k})\over
\tht(xt_{\si_{k-1}}/t_{\si_k})}
\prjk\,{\tht(xt_{\si_j}/t_{\si_k})\over\tht(t_{\si_j}/t_{\si_k})}\,\Bigr)
\endalign
$$
\fi
where $\si_0=0$. Here and below we set $t_0=x\1b$.
Using the identity we obtain that
\ifMag
$$
\NN6>
\gather
\int_{\TT^\ell}\Fabc\dt\;=\;\ell\)!\>\int_{\TT^\ell}\Fti(\tell)\dt\,,
\Tag{FFti}
\\
\Lline{\Fti(\tell)\,=\,(-b)^{-\ell}\>(p)\9^\ell\,
\pros\,{\tht(p\)x^sbc)\over\tht(x^{s(s-\ell)}(bc)^{s-\ell})}\ \x}
\\
\Rline{{}\x\prkl
{\tht(x^{1+(k-1)(k-1-\ell)}(bc)^{k-1-\ell}t_{k-1}/t_k)\>(p\>t_k/b)\9
\over\tht(xt_{k-1}/t_k)\>(at_k)\9\!}
\prjk{(1-t_k/t_j)\>(p\)x\1t_k/t_j)\9\over(x\)t_k/t_j)\9}\;.\!}
\endgather
$$
\else
$$
\NN6>
\gather
\int_{\TT^\ell}\Fabc\dt\;=\;\ell\)!\>\int_{\TT^\ell}\Fti(\tell)\dt\,,
\Tag{FFti}
\\
\aligned
\Fti(\tell)\, &{}=\,(-b)^{-\ell}\>(p)\9^\ell\,
\pros\,{\tht(p\)x^sbc)\over\tht(x^{s(s-\ell)}(bc)^{s-\ell})}\ \x
\\
& {}\>\x\prkl
{\tht(x^{1+(k-1)(k-1-\ell)}(bc)^{k-1-\ell}t_{k-1}/t_k)\>(p\>t_k/b)\9\over
\tht(xt_{k-1}/t_k)\>(at_k)\9\!}
\prjk{(1-t_k/t_j)\>(p\)x\1t_k/t_j)\9\over(x\)t_k/t_j)\9}\;.
\endaligned
\endgather
$$
\fi
This step is similar to the transformation of the integral in the proof
of formula \(qbeta).
\par
Suppose that the \var/s $\telm$ are fixed and $|t_k|=1$, $\kell-1$. Then poles
of the \fn/ $\Fti(\tell)$ in the punctured disk $0<|t_\ell|<1$ form the
following set: ${\lb\)p^rb\vert r\in\Zp\)\rb}$. Hence,
\ifMag
$$
\align
\int_{\TT^\ell}\Fti(\tell)\dt\; &{}=\int_{\TT^{\ell-1}\x\TT_{p^s\]x\eps}}\!\!\!
\Fti(\tell)\dt\;+{}
\\
\nn3>
&\){}+\,2\pii\;\sum_{r=0}^{s-1}\ \;\int_{\TT^{\ell-1}}\!
\Res\Fti(\tell)\vst{t_\ell=p^r\}xt_{\ell-1}}d^{\ell-1}t
\endalign
$$
\else
$$
\int_{\TT^\ell}\Fti(\tell)\dt\;=\int_{\TT^{\ell-1}\x\TT_{p^s\]x\eps}}\!\!\!
\Fti(\tell)\dt\;+\, 2\pii\;\sum_{r=0}^{s-1}\ \;\int_{\TT^{\ell-1}}\!
\Res\Fti(\tell)\vst{t_\ell=p^r\}xt_{\ell-1}}d^{\ell-1}t\,,
$$
\fi
where ${\TT_y=\lb\)t_\ell\in\C\vert\ |t_\ell|=y\>\rb}$, $\;\eps$ is an
arbitrary positive number between $|p|$ and $1$, and the residues are
calculated \wrt/ the \var/ $t_\ell$, all other \var/s being fixed.
Due to the \fn/al relation
$$
{\Fti(\telm,p\>t_\ell)\over\Fti(\tell)}\;=\,x^{1-\ell}\>(bc)\1\;
{1-at_\ell\over1-p\>b\1t_\ell}\;\prod_{k=1}^{\ell-1}\,
{(t_k-p\>t_\ell)\>(t_k-xt_\ell)\over(t_k-t_\ell)\>(t_k-p\)x\1t_\ell)}
$$
we see that $\Fti(\tell)=O\bigl((x^{\ell-1}bc)^{-s}\bigr)\;$ uniformly at
$\;\TT^{\ell-1}\x\TT_{p^sx\eps}$ \ and
$$
\gather
\Res\Fti(\tell)\vst{t_\ell=p^r\}xt_{\ell-1}}=\,
O\bigl(p^s(x^{\ell-1}bc)^{-s}\bigr)
\\
\Text{as $s\to\8$. Therefore,}
\nn6>
\int_{\TT^\ell}\Fti(\tell)\dt\;=\,2\pii\;\sum_{r=0}^\8\ \;
\int_{\TT^{\ell-1}}\!
\Res\Fti(\tell)\vst{t_\ell=p^r\}xt_{\ell-1}}d^{\ell-1}t\,.
\endgather
$$
\par
Similarly to the previous consideration we transform the integrals in \rhs/ of
the above formula and obtain that
\ifMag
$$
\align
\int_{\TT^\ell} & \,\Fti(\tell)\dt\;={}
\\
\nn-6>
& {}\!\}=\,(2\pii)^2\;\sum_{r_\ell=0}^\8\,
\sum_{r_{\ell-1}=0}^\8\ \;\int_{\TT^{\ell-2}}\!\Res\bigl(\)\Res\Fti(\tell)
\vst{t_\ell=p^{r_\ell}xt_{\ell-1}}\)\bigr)
\vst{t_{\ell-1}=p^{r_{\ell-1}}xt_{\ell-2}}\)d^{\)\ell-2}t\,.
\endalign
$$
\else
$$
\int_{\TT^\ell}\,\Fti(\tell)\dt\;=\,(2\pii)^2\;\sum_{r_\ell=0}^\8\,
\sum_{r_{\ell-1}=0}^\8\ \;\int_{\TT^{\ell-2}}\!\Res\bigl(\)\Res\Fti(\tell)
\vst{t_\ell=p^{r_\ell}xt_{\ell-1}}\)\bigr)
\vst{t_{\ell-1}=p^{r_{\ell-1}}xt_{\ell-2}}\)d^{\)\ell-2}t\,.
$$
\fi
It is clear that the order of summation is irrelevant.
Repeating the procedure $\ell$ times we get the following formula
$$
\align
\int_{\TT^\ell} & \,\Fti(\tell)\dt\,={}
\\
\nn-6>
& \!\}{}=\,(2\pii)^\ell\;
\sum_{k=1}^\ell\,\sum_{r_k=0}^\8 \Res\bigl(\)\Res\bigl(\,\ldots\,
\Res\Fti(\tell)\vst{t_\ell=p^{r_\ell}xt_{\ell-1}}\!\ldots\,\bigr)
\vst{t_2=p^{r_2}xt_1}\)\bigr)\vst{t_1=p^{r_1}b}
\endalign
$$
which is \eqv/ to formula \(qbetaS) because of relation \(FFti) and
the equality
$$
\align
\Res\bigl(\)\Res\bigl(\,\ldots\,{} & \Res\Fabc
\vst{t_\ell=p^{r_\ell}xt_{\ell-1}}\!\ldots\,\bigr)
\vst{t_2=p^{r_2}xt_1}\)\bigr)\vst{t_1=p^{r_1}b}\;=
\\
\nn6>
{}=\,{} & \Res\bigl(\)\Res\bigl(\,\ldots\,\Res\Fti(\tell)
\vst{t_\ell=p^{r_\ell}xt_{\ell-1}}\!\ldots\,\bigr)
\vst{t_2=p^{r_2}xt_1}\)\bigr)\vst{t_1=p^{r_1}b}\,.
\endalign
$$
The lemma is proved.
\epf

\ifAst\newpage\fi
\Sect[D]{The multidimensional Askey-Roy formula and
\ifAst\{\nl\fi Askey's \ifAst\else\ifMag\{\nl\fi\fi conjecture}
In this \appx/ we give proofs of formula \(ARl) and Askey's conjecture
\Cite{As, Conjecture 8}, see formula \(Ascj).
\Pf of formula \{\(ARl).
Let $n=2$ and $\;\ka=p\1\eta^{\ell-1}\xi_1\1\xi_2\1\!$. Assume that parameters
$\xi_1,\xi_2,z_1,z_2$ and $\eta$ are generic. Consider formula \(mu=02).
It takes the form
$$
I(W_0,w)\,=\;\pros\,{(\eta\1)\9\>(p\)\eta^{s+2-\ell}\xi_1^2\xi_2^2)\9
\>(\eta^s\xi_2^{-2})\9\>(\eta^s\xi_1\1\xi_2\1z_1/z_2)\9\over
(\eta^{-s-1})\9\>(p)\9\>(\eta^{-s}\xi_1^2)\9\>(\eta^s\xi_1\xi_2z_1/z_2)\9}\;.
\Tag{IW0}
$$
Here and after we use the simplified notations: $w_0\:=w\'{\ell,0}$,
$W_k=W\'{\ell-k,k}$, $k=0\lc\ell$, \,and $I$ denotes the \hint/ \(IWwt) as
usual.
\par
In the case in question the quotient space $\Fq(z)/\Qc'(z)$ is \onedim/ and
is spanned by the equivalence class of the \fn/ $W_0$, see Lemma~\[dimFQ'].
Let $y\in\Cxl$ be the following point:
$$
y\,=\,(\eta^{\ell-1}\xi_1\1z_1,\eta^{\ell-2}\xi_1\1z_1\lc\xi_1\1z_1)\,.
$$
Lemma~\[ResQ'] means that for any element $W(t)$ of the \ehgf/ $\Fq(z)$ of
a fiber the \fn/ ${W(t)\>W_0(y)-W_0(t)\>W(y)}$ is an element of the \bosub/
$\Qc'(z)$. Therefore, by Lemma~\[IWw=0'] for any $W\in\Fq(z)$ we have that
$$
I(W,w_0\:)\,=\,I(W_0,w_0\:)\;{W(y)\over W_0(y)}\;.
\Tag{IWW0}
$$
\par
Let $\zt$ be a nonzero complex number. The next function is an element of
the space $\Fq(z)$:
\ifMag\vv-.8>\fi
$$
W(\tell)\,=\,\pral\,
{\tht(p\)\eta^{1-\ell}\xi_1\xi_2\zt\1t_a/z_1)\>\tht(\zt\)t_a/z_2)\over
\tht(\xi_1\1t_a/z_1)\>\tht(\xi_2\1t_a/z_2)}
\prab\ {\tht(t_a/t_b)\over\tht(\eta\)t_a/t_b)}\;.
\Tag{Wal}
$$
In particular, for $\zt=\xi_2\1$ we get the \fn/ $W_0$.
\goodbm
\Par
Let $|z_1|=|z_2|=1$, ${|\)\xi_1|<1}$, ${|\)p\)\xi_2|<1}$ and ${|\eta\)|>1}$.
Under these assumptions for any \fn/ $W\in\Fq(z)$ the \hint/ $I(W,w_0\:)$
is given by formula \(IWw):
$$
\vp\int
{\smash I(W,w_0\:)\,=\int_{\){\TT^\ell}\]}\Phi(t)\>w_0\:(t)\>W(t)\;\dtt\,.}
$$
Calculating the \hint/ $I(W,w_0\:)$
\aagood
for a \fn/ $W$ of the form \(Wal) via
formulae \(IW0) and \(IWW0) we obtain formula \(ARl) for generic values
of parameters $a,b,c,\al,\bt,x$ up to a change of notations:
$$
\ifMag
a=\xi_1/z_1\,,\quad b=\xi_2/z_2\,,\quad c=p\)z_2/\zt\,,\quad
\al=\xi_1z_1\,,\quad \bt=p\>\xi_2z_2\,,\quad x=\eta\1\,.
\else
a=\xi_1/z_1\,,\qquad b=\xi_2/z_2\,,\qquad c=p\)z_2/\zt\,,\qquad
\al=\xi_1z_1\,,\qquad \bt=p\>\xi_2z_2\,,\qquad x=\eta\1\,.
\fi
\Tag{abcx}
$$
Formula \(ARl) extends to arbitrary values of parameters $a,b,c,\al,\bt,x$ by
the \anco/.
\epf
Formula \(ARl) admits the following modifications.
\goodbreak
\Lm{ARl'}
Let $|a|<1$, $|b|<1$, $|\al|<1$, $|\bt|<1$, $|x|<1$. Then
$$
\align
\int_{\)\TT^\ell\]}\,\prkl\,{\tht(p\)t_k/c)\>\tht(x^{\ell-1}abc\)t_k)\over
t_k\>(at_k)\9\>(b\)t_k)\9\>(\al/t_k)\9\>(\bt/t_k)\9\!}\;\prjk
\ {\tht(t_j/t_k)\9\over(x\)t_j/t_k)\9\>(p\)x\)t_k/t_j)\9\!}\;\dt\,=\qqq &
\\
\nn4>
=\;(2\pii)^\ell\,(p)\9^{\ell(\ell-3)/2}\;\pros\,
{(p\)x)\9\>(x^{\ell+s-1}ab\al\bt)\9\>\tht(x^s ac)\>\tht(x^s bc)\over
(p\)x^{s+1})\9\>(x^s a\al)\9\>(x^s a\bt)\9\>(x^s b\al)\9\>(x^s b\bt)\9\!}&\;.
\endalign
$$
\endpro
\Lm{ARl''}
Let $|a|<1$, $|b|<1$, $|\al|<1$, $|\bt|<1$, $|x|<1$. Then
$$
\align
\int_{\)\TT^\ell\]}\,\prkl\,{\tht(p\)t_k/c)\>\tht(x^{\ell-1}abc\)t_k)\over
t_k\>(at_k)\9\>(b\)t_k)\9\>(\al/t_k)\9\>(\bt/t_k)\9\!}\;\prjk
\ {\tht(t_k/t_j)\9\over(x\)t_k/t_j)\9\>(p\)x\)t_j/t_k)\9\!}\;\dt\,=\qqq &
\\
=\;(2\pii)^\ell\,(p)\9^{\ell(\ell-3)/2}\;\pros\,
{(p\)x)\9\>(x^{\ell+s-1}ab\al\bt)\9\>\tht(x^s ac)\>\tht(x^s bc)\over
(p\)x^{s+1})\9\>(x^s a\al)\9\>(x^s a\bt)\9\>(x^s b\al)\9\>(x^s b\bt)\9\!}&\;.
\endalign
$$
\endpro
\noPcd
\Pf of Lemma\{~\{\[ARl'].
Denote by $f(\tell)$ the integrand in formula \(ARl). Then the integrand in
the first formula of the lemma equals
$$
f(\tell)\>(p)\9^{\ell(\ell-1)/2}\!\!\prjk\;{t_j-x\)t_k\over t_j-t_k}\;.
$$
Since the integration contour $\TT^\ell$ is invariant \wrt/ \perm/s of the
\var/s $\tell$, the first formula follows from formula \(ARl) and the next
identity:
$$
\susi\;\prjk\;{t_{\si_j}-x\)t_{\si_k}\over t_{\si_j}-t_{\si_k}}\ =\;
\prosl\,{1-x^s\over 1-x}\;.
$$
The identity is \eqv/ to Corollary~\[n=1] up to a change of notations.
\epf
\nt
The proof of Lemma~\[ARl''] is similar.
\goodbm
\goodast
\Par
Introduce points $v_s\in\C$, $s\in\Z$, by the rule:
$$
\gather
v_s=p^s b\qquad\for\ s\ge 0\,,\qqq v_s=p^{-s-1}a\qquad\for\ s<0\,.
\Tag{points}
\\
\nn8>
\Text{Set}
\nn-16>
A(\uell;x)\,=\,\prkl\,
{u_k\>(p\)u_k/a)\9\>(p\)u_k/b)\9\over(\al u_k)\9\>(\bt u_k)\9}
\prjk\){(p\)x\1u_k/u_j)\9\over(p\)x\)u_k/u_j)\9}\;.
\endgather
$$
\Th{Ascj}
Let $m$ be a nonnegative integer. Then
\ifMag
\vv-.2>
$$
\gather
\Lline{\sum_{k=1}^\ell\,\sum_{r_k\in\Z}
\sgn(r_k)\>A(v_{r_1}\lc v_{r_\ell};p^m)\prkl v_{r_k}^{2m(\ell-k)}\,={}}
\\
\nn4>
\Rline{{}=\,
p^{m^2{\tsize{\ell\choose 3}}-{\tsize{m\choose 2}{\ell\choose 2}}}\,
\pros\,{(p^{m+1})\9\>(p^{m(\ell+s-1)}ab\al\bt)\9\>(-ab)^{ms}\>b\>\tht(a/b)
\over(p^{m(s+1)+1})\9\>(p^{ms}a\al)\9\>(p^{ms}a\bt)\9\>(p^{ms}b\al)\9\>
(p^{ms}b\bt)\9\!}}
\endgather
$$
\else
$$
\align
\sum_{k=1}^\ell\,\sum_{r_k\in\Z} & \sgn(r_k)\>A(v_{r_1}\lc v_{r_\ell};p^m)
\prkl v_{r_k}^{2m(\ell-k)}\,={}
\\
\nn4>
&\,{}=\,p^{m^2{\tsize{\ell\choose 3}}-{\tsize{m\choose 2}{\ell\choose 2}}}\,
\pros\,{(p^{m+1})\9\>(p^{m(\ell+s-1)}ab\al\bt)\9\>(-ab)^{ms}\>b\>\tht(a/b)
\over(p^{m(s+1)+1})\9\>(p^{ms}a\al)\9\>(p^{ms}a\bt)\9\>(p^{ms}b\al)\9\>
(p^{ms}b\bt)\9\!}
\endalign
$$
\fi
provided the parameters $a,b,\al,\bt$ are \st/ all the terms of the sum are
regular. Here $\sgn(r)=1\,$ for $r\ge 0$, and $\sgn(r)=-1\,$ for $r<0$.
\endpro
\Rem
The above formula was conjectured by Askey \Cite{As, Conjecture 8}. There is
the following correspondence of notations:
\ifAst
$$
\Line{p=q\,,\hfil a=-c\,,\hfil b=d\,,\hfil \al=-q^x\}/c\,,\hfil \bt=q^y\}/d\,,
\hfil m=k\,,\hfil \ell=n}
$$
\else
$$
p=q\,,\qquad a=-c\,,\qquad b=d\,,\qquad \al=-q^x\}/c\,,\qquad \bt=q^y\}/d\,,
\qquad m=k\,,\qquad \ell=n
$$
\fi
where notations in \lhs/ are from this paper and notations in \rhs/s are from
\Cite{As}. The above formula differs from the conjectured formula in \Cite{As}
by the factor
$p^{m^2{\tsize{\ell\choose 3}}-{\tsize{m\choose 2}{\ell\choose 2}}}\!$.
\enddemo
\Rem
After this paper was written we found out that formula \(Ascj) was proved in
\Cite{E}.
\enddemo
\Pf of Theorem\{~\{\[Ascj].
The sum in \lhs/ of the formula is absolutely convergent since
\ifMag\vv-.7>\fi
$$
A(v_{r_1}\lc v_{r_\ell};p^m)\prkl v_{r_k}^{2m(\ell-k)}\,=\,
O(\)p^{\)|r_1|\lsym+|r_\ell|})
$$
as $|r_1|\lsym+|r_\ell|\to\8$. Hence, the sum defines an \anf/ of the
parameters $a,b,\al,\bt$. Therefore, it suffices to prove the formula under
the assumptions $|a|<1$, $|b|<1$, $|\al|<1$, $|\bt<1|$.
\par
Consider the integral in \lhs/ of formula \(ARl'') for $x=p^m$, and denote
the integrand by $f(\tell)$. The poles of the integrand are located at
the hyperplanes
\ifMag\vv.1>\fi
$$
t_j=p^{-s}a\1\,,\qquad t_j=p^{-s}b\1\,,\qquad t_j=p^s\]\al\,,\qquad
t_j=p^s\]\bt\,,
$$
\ifMag\vvv.3>\fi
$j=1\lc\ell$, $s\in\Zp$. Due to Lemma~\[ARl''] the formula \(Ascj) is \eqv/ to
\aagood
the next formula
\ifMag
$$
\gather
\Lline{\qqq\,\ \int_{\TT^\ell} f(\tell)\dt\,={}}
\Tag{resfv}
\\
\nn-6>
\Rline{{}=\,(-2\pii)^\ell\;
\sum_{k=1}^\ell\,\sum_{r_k\in\Z} \Res\bigl(\)\Res\bigl(\,\ldots\,
\Res f(\tell)\vst{t_\ell=v_{r_\ell}\1}\!\ldots\,\bigr)
\vst{t_2=v_{r_2}\1}\)\bigr)\vst{t_1=v_{r_1}\1}\,,}
\endgather
$$
\else
$$
\align
\int_{\TT^\ell} & f(\tell)\dt\,={}
\Tagg{resfv}
\\
\nn-6>
& \!\}{}=\,(-2\pii)^\ell\;
\sum_{k=1}^\ell\,\sum_{r_k\in\Z} \Res\bigl(\)\Res\bigl(\,\ldots\,
\Res f(\tell)\vst{t_\ell=v_{r_\ell}\1}\ldots\,\bigr)
\vst{t_2=v_{r_2}\1}\)\bigr)\vst{t_1=v_{r_1}\1}\,,
\endalign
$$
\fi
since the sum of residues in \rhs/ coincides with the sum in \lhs/ of formula
\(Ascj). Formula \(resfv) can be proved by standard arguments, \cf. definition
\(points) of points $v_s$, $s\in\Z$. The theorem is proved.
\epf
Theorem~\[Ascj] admits the following generalization. Set
$$
\At(\uell;x)\,=\,\prkl\,
{u_k\>(p\)u_k/a)\9\>(p\)u_k/b)\9\over(\al u_k)\9\>(\bt u_k)\9}
\prjk\){(1-u_k/u_j)\>(p\)x\1u_k/u_j)\9\over(x\)u_k/u_j)\9}\;.
$$
\Th{Ascj'}
Let $|\)p\)x^{\ell-1}|<1$. Then
\ifMag
$$
\NN6>
\gather
\Lline{\sum_{j=0}^\ell\,\sum_{k=1}^\ell\,\sum_{r_k=0}^\8\,
(-1)\vpb{\>j}\>x^{\)\sum_{i=1}^\ell(\ell-i)(\ell-i+1)\)r_i\>-\>
(\ell-j-1)(\ell-j)(1\)+\>2\]\sum_{i=1}^j\]r_i)/2}
\ \prod_{s=0}^{\ell-j-1}\,{\tht(x^{\)j+s}a/b)\over\tht(x^{\)j-s}a/b)}\ \x}
\\
\Rline{{}\x\,\At(\)p^{r_1}a,p^{r_1+r_2}x\)a\lc{}p^{r_1\lsym+r_j}x^{\)j-1}a,
p^{r_{j+1}}b\lc p^{r_{j+1}\lsym+r_\ell}x^{\ell-j-1}b;x)\,={}\ \ }
\\
\Rline{{}=\,\pros\,{(x)\9\>(x^{\ell+s-1}ab\al\bt)\9\>b\>\tht(x^s a/b)\over
(x^{s+1})\9\>(x^s a\al)\9\>(x^s a\bt)\9\>(x^s b\al)\9\>(x^s b\bt)\9\!}}
\endgather
$$
\else
$$
\NN6>
\align
\sum_{j=0}^\ell\,\sum_{k=1}^\ell\,\sum_{r_k=0}^\8\,(-1)\vpb{\>j}\>
x^{\)\sum_{i=1}^\ell(\ell-i)(\ell-i+1)\)r_i\>-\>
(\ell-j-1)(\ell-j)(1\)+\>2\]\sum_{i=1}^j\]r_i)/2}
\ \prod_{s=0}^{\ell-j-1}\,{\tht(x^{\)j+s}a/b)\over\tht(x^{\)j-s}a/b)}\ \x\qqq &
\\
{}\x\,\At(\)p^{r_1}a,p^{r_1+r_2}x\)a\lc p^{r_1\lsym+r_j}x^{\)j-1}a,
p^{r_{j+1}}b\lc p^{r_{j+1}\lsym+r_\ell}x^{\ell-j-1}b;x)\,={}\quad &
\\
{}=\,\pros\,{(x)\9\>(x^{\ell+s-1}ab\al\bt)\9\>b\>\tht(x^s a/b)\over
(x^{s+1})\9\>(x^s a\al)\9\>(x^s a\bt)\9\>(x^s b\al)\9\>(x^s b\bt)\9\!}&
\endalign
$$
\fi
provided the parameters $a,b,\al,\bt,x$ are \st/ all the terms of the sum are
regular.
\endpro
\Pf.
The proof is similar to the proof of the previous Theorem. The terms of the sum
behave as ${O\bigl((\)p\)x^{\ell-1})\vpb{\)r_1\lsym+r_\ell}\bigr)}$ for
$r_1\lsym+r_\ell$ going to infinity. Hence, the sum is absolutely convergent
and defines an \anf/ of the parameters $a,b,\al,\bt,x$. Therefore, it suffices
to prove the formula under the assumptions $|a|<1$, $|b|<1$, $|\al|<1$,
$|\bt<1|$, $|x|<1$.
\par
Consider the integral in \lhs/ of formula \(ARl) and denote the integrand
by $f(\tell)$. The poles of the integrand are located at the hyperplanes
$$
t_j=p^{-s}a\1\,,\qquad t_j=p^{-s}b\1\,,\qquad t_j=p^s\]\al\,,\qquad
t_j=p^s\]\bt\,,
$$
$j=1\lc\ell$, $s\in\Zp$, and at the hyperplanes
$$
t_j=p^sx\)t_k\,,\qqq j,k=1\lc\ell,\quad j\ne k\,,\qquad s\in\Zp\,.
$$
Moreover, the integrand vanishes at the hyperplanes
$$
t_j=p^st_k\,,\qqq j,k=1\lc\ell,\quad j\ne k\,,\qquad s\in\Z\,,
$$
and it is a \sym/ \fn/ of the integration \var/s $\tell$. Using these
properties of the integrand by standard arguments we obtain the following
formula
\ifMag
$$
\gather
\Lline{\kern4em \int_{\TT^\ell} f(\tell)\dt\,=\,(2\pii)^\ell\>\ell\)!\,\x{}}
\Tag{ysum2}
\\
\ifAst
\Rline{
\aligned
{}\x\,\sum_{j=0}^\ell\;\sum_{k=1}^\ell\,\sum_{r_k=0}^\8
\Res\bigl(\, & \ldots\,\Res\bigl(\)\Res\bigl(\,\ldots\,
\Res f(\tell)\vst{t_\ell=p^{-r_\ell}x\1t_{\ell-1}}
\\
\nn-1>
& \ldots\,\bigr)\vst{t_{j+1}=p^{-r_{j+1}}b\1}\)\bigr)
\vst{t_j=p^{-r_j}x\1t_{j-1}}\!\ldots\,\bigr)\vst{t_1=p^{-r_1}a\1}\)\bigr)\,,\!
\endaligned}
\else
\Rline{
\aligned
{}\x\,\sum_{j=0}^\ell\;\sum_{k=1}^\ell\,\sum_{r_k=0}^\8
\Res\bigl(\,\ldots\,\Res\bigl(\)\Res\bigl(\, & \ldots\,
\Res f(\tell)\vst{t_\ell=p^{-r_\ell}x\1t_{\ell-1}}
\\
\nn-1>
& \ldots\,\bigr)\vst{t_{j+1}=p^{-r_{j+1}}b\1}\)\bigr)
\vst{t_j=p^{-r_j}x\1t_{j-1}}\!\ldots\,\bigr)\vst{t_1=p^{-r_1}a\1}\)\bigr)\,,\!
\endaligned}
\fi
\endgather
$$
\else
$$
\alignat2
\int_{\TT^\ell} f &(\tell)\dt\,=\,(2\pii)^\ell\>\ell\)!\,\x{} &&
\Tag{ysum2}
\\
\nn-6>
&{}\x\,\sum_{j=0}^\ell\;\sum_{k=1}^\ell\,\sum_{r_k=0}^\8
\Res\bigl(\,\ldots\,\Res\bigl(\)\Res\bigl( && \,\ldots\,
\Res f(\tell)\vst{t_\ell=p^{-r_\ell}x\1t_{\ell-1}}
\\
\nn-1>
&&& \ \;\ldots\,\bigr)\vst{t_{j+1}=p^{-r_{j+1}}b\1}\)\bigr)
\vst{t_j=p^{-r_j}x\1t_{j-1}}\!\ldots\,\bigr)\vst{t_1=p^{-r_1}a\1}\)\bigr)\,,
\kern-3em
\endalignat
$$
\fi
\ifMag\vvv.4>\fi
Formula \(ysum2) is a particular case of formula \(ysum). We will give
a detailed proof of a more general formula in Appendix~\SNo{E}.
\ifMag\vsk-.5>\goodbm\vsk.5>\fi
\par
Due to formula \(ARl) the formula \(Ascj') is \eqv/ to formula \(ysum2) since
the sum of residues in \rhs/ of formula \(ysum2) coisides with the sum in \lhs/
of formula \(Ascj').
\epf

\ifAst\newpage\fi
\Sect[E]{The Jackson integrals via the \hgeom/ \ifAst\{\nl\fi integrals}
In this \appx/ we present two theorems which connect the \hint/s described
in this paper with \sym/ \Atype/ Jackson integrals.
\par
Let $\sb=(s_1\lc s_\ell)$ be a vector with integer components. For any
$\lg\in\Zln\!$ define the points $\xt(\lg,\sb),\,\yt(\lg,\sb)\in\Cxl\!$
as follows:
\ifAst
$$
\NN4>
\gather
{\align
\xt(\lg,\sb)\,=\,
(\)p^{s_1\lsym+s_{\lg_1}}\eta^{1-\lg_1}\xi_1z_1,
\)p^{s_2\lsym+s_{\lg_1}}\eta^{\)2-\lg_1}\xi_1z_1\lc \)p^{s_{\lg_1}}\xi_1z_1 &,
\\
p^{s_{\lg_1\]+1}\lsym+s_{\lg_1\]+\lg_2}}\eta^{1-\lg_2}\xi_2z_2\lc
\)p^{s_{\lg_1\]+\lg_2}}\xi_2z_2 &,\,\ldots\,,
\\
p^{s_{\ell-\lg_n\]+1}\lsym+s_\ell}\eta^{1-\lg_n}\xi_nz_n \,\,\)&\!\!\]\lc
\)p^{s_\ell}\xi_nz_n)\,,
\endalign}
\\
\nn2>
{\align
\yt(\lg,\sb)\,=\,(\)p^{s_1\lsym+s_{\lg_1}}\eta^{\lg_1-1}\xi_1\1z_1,\)
p^{s_2\lsym+s_{\lg_1}}\eta^{\lg_1-2}\xi_1\1z_1\lc\)p^{s_{\lg_1}}\xi_1\1\]z_1 &,
\\
p^{s_{\lg_1\]+1}\lsym+s_{\lg_1\]+\lg_2}}\eta^{\lg_2-1}\xi_2\1z_2\lc
\)p^{s_{\lg_1\]+\lg_2}}\xi_2\1\]z_2 &,\,\ldots\,,
\\
p^{s_{\ell-\lg_n\]+1}\lsym+s_\ell}\eta^{\lg_n-1}\xi_n\1z_n \quad\>&\kern-1em\}
\lc\)p^{s_\ell}\xi_n\1z_n)\,,
\endalign}
\endgather
$$
\else
\ifMag
$$
\NN4>
\gather
{\align
\xt(\lg,\sb)\,=\,
(\)p^{s_1\lsym+s_{\lg_1}}\eta^{1-\lg_1}\xi_1z_1,
\)p^{s_2\lsym+s_{\lg_1}}\eta^{\)2-\lg_1}\xi_1z_1\lc \)p^{s_{\lg_1}}\xi_1z_1 &,
\\
p^{s_{\lg_1\]+1}\lsym+s_{\lg_1\]+\lg_2}}\eta^{1-\lg_2}\xi_2z_2\lc
\)p^{s_{\lg_1\]+\lg_2}}\xi_2z_2 &,\,\ldots\,,
\\
p^{s_{\ell-\lg_n\]+1}\lsym+s_\ell}\eta^{1-\lg_n}\xi_nz_n &\lc
\)p^{s_\ell}\xi_nz_n)\,,
\endalign}
\\
\nn2>
{\align
\yt(\lg,\sb)\,=\,(\)p^{s_1\lsym+s_{\lg_1}}\eta^{\lg_1-1}\xi_1\1z_1,\)
p^{s_2\lsym+s_{\lg_1}}\eta^{\lg_1-2}\xi_1\1z_1\lc\)p^{s_{\lg_1}}\xi_1\1\]z_1 &,
\\
p^{s_{\lg_1\]+1}\lsym+s_{\lg_1\]+\lg_2}}\eta^{\lg_2-1}\xi_2\1z_2\lc
\)p^{s_{\lg_1\]+\lg_2}}\xi_2\1\]z_2 &,\,\ldots\,,
\\
p^{s_{\ell-\lg_n\]+1}\lsym+s_\ell}\eta^{\lg_n-1}\xi_n\1z_n & \lc
\)p^{s_\ell}\xi_n\1z_n)\,,
\endalign}
\endgather
$$
\else
$$
\NN3>
\gather
{\align
\xt(\lg,\sb)\,=\,
({}\) & p^{s_1\lsym+s_{\lg_1}}\eta^{1-\lg_1}\xi_1z_1,
\)p^{s_2\lsym+s_{\lg_1}}\eta^{\)2-\lg_1}\xi_1z_1\lc \)p^{s_{\lg_1}}\xi_1z_1,
\\
& p^{s_{\lg_1\]+1}\lsym+s_{\lg_1\]+\lg_2}}\eta^{1-\lg_2}\xi_2z_2\lc
\)p^{s_{\lg_1\]+\lg_2}}\xi_2z_2,\,\ldots\,,
\)p^{s_{\ell-\lg_n\]+1}\lsym+s_\ell}\eta^{1-\lg_n}\xi_nz_n\lc
\)p^{s_\ell}\xi_nz_n)\,,
\endalign}
\\
\nn3>
{\align
\yt(\lg,\sb)\,=\,({}\) & p^{s_1\lsym+s_{\lg_1}}\eta^{\lg_1-1}\xi_1\1z_1,
\)p^{s_2\lsym+s_{\lg_1}}\eta^{\lg_1-2}\xi_1\1z_1\lc\)p^{s_{\lg_1}}\xi_1\1\]z_1,
\\
& p^{s_{\lg_1\]+1}\lsym+s_{\lg_1\]+\lg_2}}\eta^{\lg_2-1}\xi_2\1z_2\lc
\)p^{s_{\lg_1\]+\lg_2}}\xi_2\1\]z_2,\,\ldots\,,
\)p^{s_{\ell-\lg_n\]+1}\lsym+s_\ell}\eta^{\lg_n-1}\xi_n\1z_n\lc
\)p^{s_\ell}\xi_n\1z_n)\,,
\endalign}
\endgather
$$
\fi\fi
\cf. \(xtyt). In particular, for $\sb=(0\lc 0)$ we have $\xt(\lg,\sb)=\xt\lg$
and $\yt(\lg,\sb)=\alb\yt\lg$.
\Par
Recall, that the short \phf/ $\Phi(t,z)$ is given by formula \(Phis) and we
define the multiple residue by formula \(Rests). Set
$\,{\Pht(t,z)\>=\>t_1\1\!\ldots t_\ell\1\)\Phi(t,z)}$.
\Th{xsum}
Let ${|\)p\)\ka\zy|<\min\>(1,|\eta\)|^{1-\ell})}$. Let \(npZ)\;--\;\(assum)
hold. Then for any \fn/s $w\in\Fo(z)$ and $W\in\Fq(z)$ we have that
$$
I(W,w)\,=\,(2\pii)^\ell\>\ell\)!\smiZ\!\sum_{\;\sb\in\Zp^\ell\!\!}
\Res\bigl(\Pht(t)\>w(t)\>W(t)\bigr)\vst{t=\xt(\mg,\sb)}\,.
\vv-1.2>
$$
\endpro
\Th{ysum}
Let ${|\)\ka\yz|>\max\>(1,|\eta\)|^{\ell-1})}$. Let \(npZ)\;--\;\(assum)
hold. Then for any \fn/s $w\in\Fo(z)$ and $W\in\Fq(z)$ we have that
$$
I(W,w)\,=\,(-2\pii)^\ell\>\ell\)!\smiZ\!\sum_{\;\sb\in\Zn^\ell\!\!}
\Res\bigl(\Pht(t)\>w(t)\>W(t)\bigr)\vst{t=\yt(\mg,\sb)}\,.
\vv->
$$
\endpro
\Pf of Theorem~\[xsum].
We prove the theorem assuming that the \fn/s $w$ and $W$ are the \tri/ and
\ewf/, \resp/. For arbitrary \fn/s $w$ and $W$ the statement holds by
linearity.
\vsk.1>
Let ${\zt=\max\>(1,|\eta\)|^{\ell-1})\>|\)p\)\ka\zy|}$, that is ${0<\zt<1}$
\vvgood
\ifMag\vv.15>\fi
under the assumptions of the theorem. Using the \fn/al relation
\ifMag
\mmgood
$$
\gather
\Lline{\ifAst\quad\else\ \)\fi\qqq
{\Phi(\tpell)\>W(\tpell)\over\Phi(\tell)\>W(\tell)}\;={}}
\Tag{PWPW}
\\
\nn6>
\Rline{{}=\,\ka\pron{\xi_mt_a-z_m\over t_a-\xi_mz_m}\,
\prod_{a<b\le\ell}{t_a-\eta\)t_b\over\eta\)t_a-t_b}\;
\prod_{1\le b<a}{t_a-p\)\eta\)t_b\over p\)\eta\)t_a-t_b}\;,}
\endgather
$$
\else
$$
\Rline{{\Phi(\tpell)\>W(\tpell)\over\Phi(\tell)\>W(\tell)}\;=\,
\ka\pron{\xi_mt_a-z_m\over t_a-\xi_mz_m}\,
\prod_{a<b\le\ell}{t_a-\eta\)t_b\over\eta\)t_a-t_b}\;
\prod_{1\le b<a}{t_a-p\)\eta\)t_b\over p\)\eta\)t_a-t_b}\;,}
\Tag{PWPW}
$$
\fi
$\aell$, we estimate the residues: ${\Res\bigl(\Pht(t)\>w(t)\>W(t)\bigr)
\vst{t=\xt(\mg,\sb)}\}=\>O(\zt^{s_1\lsym+s_\ell})}$, \,and obtain that the sum
is absolutely convergent. Therefore, the sum defines an \anf/ of the parameters
$\eta$, $\xin$, $\zn$, and it suffices to prove the formula under
the assumptions that ${|\eta\)|>1}$ and $|z_m|=1$, $|\xi_m|<1$, $\mn$.
\par
Under these assumptions the \hint/ $I(W,w)$ is given by formula \(IWw)
and the claim of the theorem is \eqv/ to Lemma~\[xsumf] for $f=1$.
Theorem~\[xsum] is proved.
\epf
\Lm{xsumf}
Let $\,{|\eta\)|>1}$ and $\,{|\xi_lz_l|<|\xi_m\1z_m|}$ \,for all $l,\mn$. Let
\vv.1>
\ifMag\nl\fi
$\,{|\)p\)\ka\zy|<1}$. Let \(npZ)\;--\;\(assum) hold. Let $\,\al$ be
a real number \st/
\ifMag\vv.1>\nl\fi
$\,{\max\>(|\xi_1z_1|\)}\lc|\xi_nz_n|)<\al
<{\min\>(|\xi_1\1z_1|\)}\lc|\xi_n\1z_n|)$. Then,
$$
\int_{\){\TT^\ell_\al}\]}\Phi(t)\>w(t)\>W(t)\>f(t)\;\dtt\,=\,
(2\pii)^\ell\>\ell\)!\smiZ\!\sum_{\;\sb\in\Zp^\ell\!\!}
\Res\bigl(\Pht(t)\>w(t)\>W(t)\>f(t)\bigr)\vst{t=\xt(\mg,\sb)}
$$
for any \sym/ \fn/ $f(t)$ regular inside the torus
$\,\TT^\ell_\al={\lb\)t\in\Cl\vert\ |t_1|=\al}\llc{|t_\ell|=\al\)\rb}$.
\endpro
\Pf.
We prove the lemma by induction \wrt/ $\ell$, the number of integrations.
\par
Notice that the integral does not change if we replace the integration contour
in the integral by $\TT_{\al_1}\!\lx\TT_{\al_\ell}$, where $\al_1\lc\al_\ell$
are \pd/ real numbers close to $\al$.
\par
Similar to the proof of Theorem~\[xsum] we see that the sum of residues is
absolutely convergent. Hence, the sum defines a \hof/ of the parameters
$\eta$, $\xin$, $\zn$, and the same does the integral in \lhs/ of formula
\(xsumf). Therefore, \wlg/ we can assume that all the numbers $|p^s\al_m|$,
$|p^s\xi_mz_m|$, $\mn$, $s\in\Zp$, are \pd/.
\par
Replace the integration contour in the integral by
${\TT_{\eps\al_1}\!\lx\TT_{\eps\al_\ell}}$ for real $\eps$, so that at $\eps=1$
we have the initial integral. If $\eps$ is decreasing starting from $1$,
the integral is not changing until the integration contour touches one of the
hyperplanes where the integrand has a pole. And every time when the integration
\aagood
contour crosses the singularity hyperplane we pick up the integral of dimension
$\ell-1$ of the corresponding residue. Notice that during the described
deformation the integration contour can touch only the singularity hyperplanes
of the form $t_a=p^s\xi_mz_m$, $s\in\Zp$. Finally, for any positive integer
$r$ we get
\ifMag
$$
\gather
\Lline{\int_{\){\TT^\ell_\al}\]}\Phi(t)\>w(t)\>W(t)\>f(t)\;\dtt\;=
\int_{\){\TT^\ell_{p^r\}\al}}\]}\Phi(t)\>w(t)\>W(t)\>f(t)\;\dtt\,+{}}
\\
\nn2>
\Rline{{}+\,2\pii\,\sum\nolimits_{\tsize\>(r)}\int\Res\bigl(t_a\1\)
\Phi(t)\>w(t)\>W(t)\>f(t)\bigr)\vst{t_a=p^s\xi_mz_m}(dt/t)^{\ell-1}\,,}
\endgather
$$
\else
$$
\align
\int_{\){\TT^\ell_\al}\]}\Phi(t)\>w(t)\>W(t)\>f(t)\;\dtt\;=
\int_{\){\TT^\ell_{p^r\}\al}}\]} & \Phi(t)\>w(t)\>W(t)\>f(t)\;\dtt\,+{}
\\
&\!{}+\,2\pii\,\sum\nolimits\'{\tsize r}\int\Res\bigl(t_a\1\)
\Phi(t)\>w(t)\>W(t)\>f(t)\bigr)\vst{t_a=p^s\xi_mz_m}(dt/t)^{\ell-1}\,,
\endalign
$$
\fi
where each term of the sum $\sum\nolimits\'r$ corresponds to a passing of the
integration contour through a singularity hyperplane when $\eps$ goes from $1$
to $p^r$.
\par
Using relation \(PWPW) it is easy to show that
$$
\smash{\int_{\){\TT^\ell_{p^r\}\al}}\]}\Phi(t)\>w(t)\>W(t)\>f(t)\;\dtt\,=\;
O\bigl(\)|\)p\)\ka\zy|\vpb{\)\ell\)r}\)\bigr)}\vp{\int_.}
$$
as $r\to\8$. Hence, the integral disappears as $r\to\8$.
\par
All integrals of smaller dimension in the sum $\sum\nolimits\'r$ are of the
same form as the initial integral and can be replaced by the sums of residues
according to formula \(xsumf) because of the induction assumption. Since
${\Pht(t)\>w(t)\>W(t)\>f(t)}$ is a \sym/ \fn/ of $\tell$, it is straightforward
to show that the resulting sum of residues can be transformed to the sum in
\rhs/ of formula \(xsumf). The lemma is proved.
\vvgood
\epf
\ifAst\else\ifMag\vsk-.8>\fi\fi
The proof of Theorem~\[ysum] is similar to the proof of Theorem~\[xsum].
\Par
The sums of residues in formulae \(xsum) and \(ysum) coincide with \sym/
\Atype/ Jackson integrals, see for example \Cite{AK}. This means that
the \traf/s between \asol/s of the \qKZe/ coincide with the connection matrices
for \sym/ \Atype/ Jackson integrals. The Gauss decomposition of the connection
matrices studied in \Cite{AK} is related to Lemmas~\[Jxy] and \[Xixy] in this
paper.

\ifAst\newpage\fi
\Sect[X]{One useful identity}
Particular specializations of the following identity:
$$
\sum_{a=0}^l\>{j+a\choose j}{j+k+a\choose k}{l+m-a\choose m}\,=\,
{j+k\choose k}{j+k+l+m+1\choose j+k+m+1}\,,
\Tag{combi}
$$
are often used in the paper. The identity can be proved by induction \wrt/
$l$ and $m$.

\ifAst\newpage\fi
\myRefs
\widest{WW}

\ref\Key A
\by \Aomoto/
\paper $q$-analogue of de~Rham cohomology associated with Jackson integrals, I
\jour Proceedings of Japan Acad.{} \vol 66 {\rm Ser\&A} \yr 1990
\pages 161--164
\moreref \paper II \jour Proceedings of Japan Acad.{} \vol 66 {\rm Ser\&A}
\yr 1990 \pages 240--244
\endref

\ref\Key AK
\by \Aomoto/ and Y\]\&Kato
\paper Gauss decomposition of connection matrices for \sym/ \Atype/ Jackson
integrals \jour \SMNS/ \vol 1 \yr 1995 \issue 4 \pages 623--666
\endref

\ref\Key As
\by R\&Askey
\paper Some basic \hgeom/ extensions of integrals of Selberg and Andrews
\jour SIAM J.\ Math.\ Anal.{} \vol 11 \yr 1980 \pages 938--951
\endref

\ref\Key CP
\by V\]\&Chari and A\&Pressley
\book A guide to quantum groups
\yr 1994 \publ \CUP/
%% \pages 651
\endref

\ref\Key D1
\by \Dri/
\paper Quantum groups
\inbook in Proceedings of the ICM, Berkley, 1986
\ifAst\else\ifMag\adjustnext\nl\fi\fi
\ed A\&M\&Gleason \yr 1987 \pages 798--820 \publ \AMSa/
\endref

\ref\Key D2
\by \Dri/
\paper Quasi-Hopf algebras
\jour \LMJ/ \vol 1 \yr 1990 \pages 1419--1457
\endref

\ref\Key DF
\by J\&Ding and \Fre/
\paper Isomorphism of two realizations of \qaff/
$U_q\bigl(\widehat{\frak{gl}}(n)\bigr)$
\jour \CMP/ \vol 156 \yr 1993 \pages 277--300
\endref

\ref\Key E
\by R\&J\&Evans
\paper Multidimensional beta and gamma integrals
\jour Contemp.\ Math.{} \vol 166 \yr 1994 \pages 341--357
\endref

\ref\Key F
\by \Feld/
\paper Conformal field theory and integrable systems associated to \ecurve/s
\inbook in Proceedings of the ICM, Z\"urich, 1994
\publ \Birk/ \yr 1994 \pages 1247--1255
%% \jour Preprint {\sl hep-th/9407154} \yr 1994
\endref

\ref
\by\refin \Feld/
\paper Elliptic quantum groups
\inbook in Proceedings of the ICMP, Paris 1994 \ed D\&Iagolnitzer %% XI-th
\yr 1995 \publ Intern.\ Press \publaddr Boston \pages 211--218
%% \jour Preprint {\sl hep-th/9412207} \yr 1994
\endref

\ref\Key FR
\by \Fre/ and \Reshy/
\paper Quantum affine algebras and holonomic \dif/ \eq/s
\jour \CMP/ \vol 146 \yr 1992 \pages 1--60
\endref

\ref\Key FTV1
\by \Feld/, \VT/ and \Varch/
\paper Solutions of elliptic \qKZBe/s and \Ba/ I
\jour \AMS/ Transl.,\ Ser\&\)2 \vol 180 \yr 1997 \pages 45--75
%% \jour Preprint ENSLAPP\~~L\~~598/96 \yr 1996 \pages 1--22
\endref

\ref\Key FTV2
\by \Feld/, \VT/ and \Varch/
\paper Monodromy of \sol/s of the elliptic quantum \KZvB/ \deq/s
\jour Preprint \yr 1997 \pages 1--20
\endref

\ref\Key FV
\by \Feld/ and \Varch/
\paper On \rep/s of the elliptic \qg/ $E_{\tau,\eta}(sl_2)$
\jour \CMP/ \vol 181 \yr 1996 \issue 3 \pages 741--761
\endref

\ref\Key GR
\by G\&Gasper and M\&Rahman
\book Basic \hgeom/ series \bookinfo Encycl.\ Math.\ Appl.{}
\yr 1990 \publ \CUP/
\endref

\ref\Key IK
\by A\&Izergin and \Kor/
\paper The quantum scattering method approach to correlation \fn/s
\jour \CMP/ \vol 94 \yr 1984 \pages 67--92
\endref

\ref\Key J
\by M\&Jimbo
\paper Quantum \Rm/ for the generalized Toda system
\jour \CMP/ \vol 102 \yr 1986 \pages 537--547
\endref

\ref\Key JM
\by M\&Jimbo and T\]\&Miwa
\paper Algebraic analysis of solvable lattice models
\jour CBMS Regional Conf.\ Series in Math.{} \vol 85 \yr 1995
\endref

\ref\Key K
\by T\]\&Kohno
\paper Monodromy \rep/s of braid groups and \YB/s
\jour Ann.\ Inst.\ Fourier \vol 37 \yr 1987 \pages 139--160
\endref

\ref
\by\refin T\]\&Kohno
\paper Linear \rep/s of braid groups and classical \YB/s
\jour Contemp.\ Math.{} \vol 78 \yr 1988 \pages 339--363
\endref

\goodbreak

\ref\Key KL
\by \Kazh/ and \Lusz/
\paper Affine Lie algebras and \qg/s
\jour Intern.\ Math.\ Research Notices \vol 2 \yr 1991 \pages 21--29
\endref

\ref
\by\refin \Kazh/ and \Lusz/
\paper Tensor structures arising from affine Lie algebras, I
\jour J.\ \AMS/ \vol 6 \yr 1993 \pages 905--947
\moreref \paper II
\jour J.\ \AMS/ \vol 6 \yr 1993 \pages 949--1011
\endref

\ref\Key KS
\by \Kazh/ and Ya\&S\&\&Soibelman
\paper Representation theory of \qaff/s
\jour \SMNS/ \vol 3 \yr 1995 \issue 3 \pages 537--595
\endref

\ref\Key L
\by F\]\&Loeser
\paper Arrangements d'hyperplans et sommes de Gauss
\jour Ann.\ Scient.\ \'Ecole Normale Super.{} \vol {\rm 4-e serie, t\&24}
\yr 1991 \pages 379-400
\endref

\ref\Key Lu
\by S\&Lukyanov
\paper Free field \rep/ for massive integrable models
\jour \CMP/ \vol 167 \yr 1995 \pages 183--226
\endref

\ref\Key RS
\by \Reshy/ and M\&A\&Semenov-Tian-Shansky
\paper Central extensions of the quantum current groups
\jour \LMP/ \vol 19 \yr 1990 \pages 133--142
\endref

\ref\Key S
\by F\]\&A\&Smirnov
\book Form factors in completely integrable models of quantum field theory
\bookinfo Advanced Series in Math.\ Phys., vol\&14
\yr 1992 \publ \WSa/
\endref

\ref\Key SV1
\by \SchV/ and \Varn/
\paper Hypergeometric \sol/s of \KZv/ \eq/s
\jour \LMP/ \vol 20 \yr 1990 \pages 279--283
\endref

\goodbm

\ref\Key SV2
\by \SchV/ and \Varn/
\paper Arrangements of hyperplanes and Lie algebras homology
\jour Invent.\ Math.{} \vol 106 \yr 1991 \pages 139--194
\endref

\ref
\by\refin \SchV/ and \Varn/
\paper Quantum groups and homology of local systems
\inbook in Algebraic Geometry and Analytic Geometry, Proceedings of the ICM
Satellite Conference, Tokyo, 1990
\yr 1991 \publ \Spria/ \pages 182--191
\endref

\ref\Key T
\by \VoT/
\paper Irreducible monodromy matrices for the \Rm/ of the {\sl XXZ}-model
and lattice local quantum Hamiltonians
\jour \TMP/ \vol 63 \yr 1985 \pages 440--454
\endref

\ifAst
\ref\Key TV1
\by \VoT/ and \Varn/
\paper Jackson integral \rep/s of \sol/s of the quantized \KZv/ \eq/
\jour \LpMJ/ \vol 6 \issue 2 \yr 1995 \pages 275--313
\endref
\else
\ref\Key TV1
\by \VoT/ and \Varn/
\paper Jackson integral \rep/s of \sol/s of the quantized \KZv/ \eq/
\jour \LpMJ/ \vol 6 \issue 2 \yr 1995 \pages 275--313
\endref
\fi

\ref\Key TV2
\by \VT/ and \Varch/
\paper Asymptotic \sol/ to the quantized \KZv/ \eq/ and \Bv/s
\jour \AMS/ Transl.,\ Ser\&\)2 \vol 174 \yr 1996 \pages 235--273
\endref

\ref\Key TV3
\by \VT/ and \Varch/
\paper Geometry of $q$-hypergeometric \fn/s as a bridge between Yangians and
\qaff/s \jour \Inv/ \yr 1997 \vol 128 \issue 3 \pages 501--588
\endref

\ref\Key V1
\by \Varn/
\paper The Euler beta-\fn/, the Vandermonde determinant, Legendre's equation,
and critical values of linear \fn/s on a configuration of hyperplanes, I
\jour Math.\ USSR, Izvestia \vol 35 \yr 1990 \pages 543--571
\moreref \paper II
\jour Math.\ USSR Izvestia \vol 36 \yr 1991 \pages 155-168
\endref

\ref
\by\refin \Varn/
\paper Determinant formula for Selberg type integrals
\jour Funct.\ Anal.\ Appl.{} \vol 4 \yr 1991 \pages 65--66
\endref

\ref\Key V2
\by \Varch/
\book Multidimensional \hgeom/ \fn/s and \rep/ theory of Lie algebras and
\qg/s \bookinfo Advanced Series in Math.\ Phys., vol\&21
\yr 1995 \publ \WSa/
%% \pages 371
\endref

\ref\Key V3
\by \Varch/
\paper Quantized \KZv/ \eq/s, quantum \YB/, and \deq/s for $q$-\hgeom/ \fn/s
\ifAst\adjustnext\nl\fi
\jour \CMP/ \vol 162 \yr 1994 \pages 499--528
\endref

\ref\Key V4
\by \Varch/
\paper Asymptotic \sol/s to the \KZv/ \eq/ and crystal base
\jour \CMP/ \vol 171 \yr 1995 \pages 99-137
\endref

\endRefs

\ifAst
\vskip 0pt plus 1filll
\rightline{
\vtop{\hbox{\smc \VT/}
\hbox{\SPb/ Branch of}
\hbox{Steklov Mathematical Institute}
\hbox{Fontanka 27, \SPb/ \,191011}
\hbox{\smc Russia}
\hbox{\sl E-mail\/{\rm:} \homemail/}
\vsk.25>
\hbox{\it and}
\vsk.25>
\hbox{Department of Mathematics}
\hbox{Faculty of Science, Osaka University}
\hbox{Toyonaka, Osaka \,560}
\hbox{\smc Japan}
\vsk1.5>
\hbox{\smc \Varch/}
\hbox{Department of Mathematics}
\hbox{University of North Carolina}
\hbox{\ChH/, NC 27599}
\hbox{USA}
\hbox{\sl E-mail\/{\rm:} \avemail/}}}
\else
\fixedpage
\nopagenumber
\ifx\Pcde\empty
\ifMag
\contents
\ContM
\endco
\else
\contents
\ContS
\endco
\fi
\fixedpage
\nopagenumber
\else
\ifMag
\raggedbottom
\contents
\ContMP
\endco
\else
\let\Appncd\Appencd\def\Appencd{\fixedpage\Appncd}
\contents
\ContSP
\endco
\fi
\nopagenumber
\ifMag\vskip0pt plus 1filll \let\abstxt\abstext
\def\abstext{\abstxt\vsk-1.2>\lline{\vp1}}\else
\vsk2>
\fi
\fi
\abstext
\fi

\bye